\documentclass[11pt]{book}
\usepackage{graphicx}
\usepackage{amssymb}
\usepackage{amsmath}
\usepackage{slashbox}
\usepackage[latin1]{inputenc}

\usepackage[colorlinks=true, pdfstartview=FitV, linkcolor=blue, 
            citecolor=blue, urlcolor=blue]{hyperref}

\newfont{\bbb}{msbm10 scaled 500}

\newfont{\bb}{msbm10 scaled 1100}

\newcommand{\beq}{\begin{eqnarray}}
\newcommand{\eeq}{\end{eqnarray}}
\newcommand{\beqs}{\begin{eqnarray*}}
\newcommand{\eeqs}{\end{eqnarray*}}

\title{ \bf{\Huge{Divergences.\\Scale invariant Divergences.\\
 Applications to linear inverse problems.\\
- N.M.F. - \\
- Blind deconvolution -}}  \\[2cm]
 Version 1.0\\[1cm]  \bf{Construction of the divergences.\\Problems related to their use.}}

\author{Henri Lantéri}

\begin{document}

\frontmatter
\maketitle

\tableofcontents
\pagenumbering{roman} 
\listoffigures
\listoftables

\pagenumbering{arabic}

\setcounter{table}{0}  \setcounter{equation}{0}  \setcounter{figure}{0} \setcounter{chapter}{1} \setcounter{section}{0} 
\chapter{Résumé.}
Ce livre traite des fonctions permettant d'exprimer la différence entre deux champs de données ou "Divergence", pour des applications à des problèmes linéaires inverses.
La plupart des divergences trouvées dans la littérature sont utilisées en théorie de l'information pour quantifier la différence entre deux densités de probabilité, c'est-à-dire des quantités positives dont les sommes sont égales à 1. Dans ce contexte, elles prennent des formes simplifiées qui les rendent inadaptées aux problèmes considérés ici, dans lesquels les champs de données ont des valeurs positives quelconques.\\
De manière systématique, les divergences classiques sont réexaminés et on en donne des formes qui peuvent être utilisées pour les problèmes inverses que nous considérons.
Cela nous amène à préciser les méthodes de construction des écarts et à proposer un certain nombre de généralisations.
La résolution des problèmes inverses implique systématiquement la minimisation d'une divergence entre les mesures physiques et un modèle dépendant de paramètres inconnus. Dans le cadre de la reconstruction d'images, le modèle est généralement linéaire et c'est un problème de minimisation sous la contrainte de la non-négativité et (éventuellement) de la somme constante des paramètres inconnus.\\
Afin de prendre en compte la contrainte de somme de façon simple, nous introduisons la classe des divergences invariantes par changement d'échelle sur les paramètres du modèle (divergences affines invariantes) et nous montrons des propriétés intéressantes de ces divergences.\\
Une extension de ces divergences permet d'obtenir l'invariance par changement d'échelle par rapport aux deux arguments intervenant dans les divergences ; ceci autorise l'utilisation de ces divergences dans la régularisation des problèmes inverses par contrainte de lissage.
Des méthodes algorithmiques de minimisation des divergences sont développées, sous contraintes de non-négativité et de somme des composantes de la solution. Les méthodes présentées sont basées sur les conditions de Karush-Kuhn-Tucker qui doivent être satisfaites à l'optimum. La régularisation au sens de Tikhonov est prise en compte dans ces méthodes.
Le chapitre 11 associé à l'annexe 9 traite des applications de la MMF, tandis que le chapitre 12 est consacré à la déconvolution à l'aveugle.
 Dans ces deux chapitres, l'accent est mis sur l'intérêt des divergences invariantes.

\chapter{Abstract.}
This book deals with functions allowing to express the dissimilarity (discrepancy) between two data fields or "divergence functions" with the aim of applications to linear inverse problems.\\
Most of the divergences found in the litterature are used in the field of information theory to quantify the difference between two probability density functions, that is between positive data whose sum is equal to one. In such context, they take a simplified form that is not adapted to the problems considered here, in which the data fields are non-negative but with a sum not necessarily equal to one.\\
In a systematic way, we reconsider the classical divergences and we give their forms adapted to inverse problems.\\
To this end, we will recall the methods allowing to build such divergences, and propose some generalizations.\\
The resolution of inverse problems implies systematically the minimisation of a divergence between physical measurements and a model depending of the unknown parameters.\\
In the context image reconstruction, the model is generally linear and the constraints that must be taken into account are the non-negativity  as well as (if necessary) the sum constraint of the unknown parameters.\\
To take into account in a simple way the sum constraint, we introduce the class of scale invariant or affine invariant divergences. Such divergences remains unchanged when the model parameters are multiplied by a constant positive factor. We show the general properties of the invariance factors, and we give some interesting characteristics of such divergences.\\
An extension of such divergences allows to obtain the property of invariance with respect to both the arguments of the divergences; this characteristic can be used to introduce the smoothness regularization of inverse problems, that is a regularisation in the sense of Tikhonov.\\
We then develop in a last step, minimisation methods of the divergences subject to non-negativity and sum constraints on the solution components. These methods are founded on the Karush-Kuhn-Tucker conditions that must be fulfilled at the optimum. The Tikhonov regularization is considered in these methods.\\
Chapter 11 associated with Appendix 9 deal with the application to the NMF, while Chapter 12 is dedicated to the Blind Deconvolution problem.\\
In these two chapters, the interest of the scale invariant divergences is highlighted.

\chapter{Preface - Outline of the document.}
Solving inverse problems generally involves minimizing, with respect to unknown parameters, a gap function or discrepancy or divergence between measurements and a model describing the physical phenomenon under consideration. The divergences found in the literature are for the most part dedicated to problems related to information theory and are unsuitable for solving inverse problems. In order to have divergences adapted to these problems, we therefore re-examine the classical divergences, specify their construction methods, give some generalisations and present the invariant forms of these divergences. Some algorithmic methods of minimization under constraints are proposed.\\
Chapter 1 deals with some general considerations on inverse problems and on deviation functions or divergences. In particular, the problems concerning the influence of the inevitable noise on the measurements are briefly discussed. The types of divergences to be considered are specified.\\
Since these discrepancies are essentially based on the properties of the differentiable convex functions, some basic properties of the convex functions are recalled in Chapter 2 and the \textit{``standard convex functions"} are defined, as opposed to the \textit{``simple convex functions"}.\\
We then indicate the rules of construction of the Csiszar, Bregman and Jensen divergences based on such functions, then we analyze the convexity of the resulting ``\textit{separable}" divergences and we give some relations between the different types of divergences.\\
In Chapter 3, we introduce the ``scale-invariant divergences", discuss in detail the methods of obtaining the invariance factors and indicate some remarkable properties of the resulting invariant divergences.\\
Chapter 4 is devoted to divergences based on classical entropy functions and some extensions of those divergences.\\
Chapter 5 deals with the ``Alpha and Beta divergences" and we analyzes how they are constructed.\\
The forms of these divergences invariant by scale change are then given; this allows us to obtain the well known ``Gamma divergence" and some divergences similar to this one.\\
In Chapter 6, a certain number of classical divergences are grouped together and extensions of these divergences are proposed in case the data fields considered are not probability densities.\\
Chapter 7 deals with divergences based on inequalities between weighted and unweighted generalized means, and extensions are given to the logarithmic forms of these divergences.\\
Chapter 8 deals with divergences between one of the two data fields and weighted mixtures of the two fields. Generalizations of these divergences are also discussed.\\
In Chapter 9, we study the application of scale invariant divergences to the problem of regularization of inverse problems by smoothness constraint and we show the interest of invariant divergences by change of scale on the 2 arguments for certain forms of regularization of this type.\\
Finally, in Chapter 10, we discuss algorithmic methods. After recalling the S.G.M. method, we develop extensions of this method in order to take into account simultaneously constraints of non negativity and sum of the unknown parameters. We develop in particular the case of scale invariant divergences  and we show the interest of their specific properties to take into account the sum constraint.\\
Chapter 11 is devoted to Non-Negative Matrix Factoring (NMF), and more specifically to the introduction of sum constraints on unknown matrix columns.\\
 Concerning this particular constraint, in the case where the divergence used does not have an invariance property, we proceed by changing variables.\\
 On the other hand, if the considered divergence is invariant by change of scale, we show that the specific properties of these divergences lead to new particularly interesting algorithms.\\
Annex 9 deals with the problem of regularization in MMF and provides some clarifications and corrections to the classical techniques.\\
Chapter 12 deals with Blind Deconvolution.\\
It is shown in this chapter that the use of scale invariant divergences is of great interest insofar as the commutativity of the convolution product makes it possible to fully exploit the properties of such divergences when taking into account the sum constraints on the solution images.\\
As a comparison, the use of non-invariant divergences is also considered.\\
In these two chapters, the interest of the scale invariant divergences is clearly highlighted insofar as they make it possible to take into account the sum constraints, in particular with regard to the Blind Deconvolution.

\chapter{chapter 1- Introduction}  \label{chptr::chapitre1}
\section{Some Preliminary Considerations.}

The goal of this book is to analyze the functions allowing to express the gap between 2 data fields and some methods to minimize these functions in order to apply them to inverse (linear) problems\cite{bertero1998}.\\
Two aspects will thus be examined successively: on the one hand the inverse problems aspect on which we will only make a few brief reminders, and on the other hand the ``\textit{gap functions}" aspect, i.e. the ``\textit{divergences functions}" on which we will extend much more extensively and which constitutes a large part of this work.

 \subsection{Brief reminders on inverse problems and the use of divergence functions.}
 
 In the general case, we have experimental measurements of a physical quantity dependent on a temporal, spatial or other variable, i.e. $y\left(t\right)$ ("$t$" not necessarily being time) and a model $m\left(t,x\right)$ to describe this quantity; this model depends on a certain number of parameters $x=\left(x_{1},x_{2},...x_{n}\right)$ whose values are unknown.\\
  We'll move to the special case where the model is linear with respect to the parameters.\\ 
	
It is further assumed that the measurements and the model are sampled with a step $\Delta t$ so that $t_{i}=i\; \Delta t$; thus a set of measurements $y_{i}=y\left(t_{i}\right)$ is available and the corresponding values of the model are noted $m_{i}\left(x\right)=m\left(t_{i},x\right)$.\\  

To fix the ideas, we consider the problem of image restoration (reconstruction), and it is more precisely the problem of image deconvolution that will be considered throughout this book.\\
 In this case, the sampling is implicitly linked to the pixels of the observed (convolved) image, the measurements are the intensities in the pixels of this image, the (linear) model corresponds to the convolution of the unknown original image with an impulse response that we will assume to be known; the unknown parameters are the intensities of the pixels of the original image.\\   
In the context of inverse problems, the objective is to determine the optimal values of the unknown parameters of the model, i.e. those which will make it possible to bring the physical model as close as possible to the experimental measurements; one is thus led to compare the measurements with the model.\\ 
In order to do this, we need a function to quantify the difference between the two quantities (the measurements and the model).\\

\textbf{In terms of inverse problems, the ``discrepancy (gap) function'' as we just defined expresses ``\textit{the data attachment}" or ``\textit{the data consistency}".}\\

 Obviously, in this exercise, the measurements are fixed quantities and the values of the model parameters are varied until the gap between the measurements and the model is minimal.\\
It is an inverse problem with all the difficulties that can be inherent in this type of problem \cite{bertero1998} \cite{Idier01a}.\\
Assuming the model based on physical considerations is known, two steps occur in sequence:
\begin{itemize}
\item the choice of a gap function
\item the choice of a method allowing to minimize such function with respect to parameters of the physical model (the unknowns).
\end{itemize}

\subsection{The gap function $F(p,q)$ and its arguments ``$p$" and ``$q$".}

An important point should be stressed immediately: except in some special cases, the variables ``$p$" and ``$q$" considered in this book are implicitly non-negative quantities.
\subsubsection{* The data ``$p$".}
 
 The first of the arguments in the deviation function concerns measurements. In order to adopt classical notations in this area, they are referred to as ``$p_{i}$'', (these are the ``$y_{i}$" from the previous section).\\
 With rare exceptions, these measurements will be corrupted by errors that are attributed to noise whose nature is assumed to be known, which is not always true.\\
 A first question immediately comes to mind: due to the presence of noise, do the measurements have the properties of the physical quantity being studied ?\\
 Specifically, if the physical quantity under consideration is non-negative, ``are the measurements non-negative ?".\\
If not, then the gap function must be able to take into account the negative values of the measurements. If the discrepancy function can only take into account non-negative quantities, then the measurements must be modified either by setting the negative values to zero or by shifting all the measurements to become non-negative, but in doing so, the statistical properties of the noise are changed.\\
  For example, in astronomical imaging, corrections related to ``dead" sensor pixels, or corrections related to differences in the gain of the elements of a CCD detector array change the nature of the noise present in the raw measurements.\\
  In general, the range of definition of the deviation function should include the range of values of the measurements.\\
 Clearly, the problem of non-negative measurement is a crucial point.\\
 
\subsubsection{* The physical model ``$q$".} 

The second argument which is used in the discrepancy function is the model, let us designate it by ``$q_{i}$'' (they are the $m_{i}(x)$); these values depend on the unknown parameters, which we note ``$x_{l}$''.\\
 The question is, what does the physics say about the possible values of the ``$q_{i}$'' and of the ``$x_{l}$''?\\
 We can assume that the model is such that if the values of the ``$x_{l}$'' are chosen within a physically acceptable range, then the values of the ``$q_{i}$'' will be acceptable, that is, in particular, they will be non-negative.\\
  In any event, that means not proposing any values for the unknowns ``$x_{l}$'', it will be essential to introduce some constraints on these values.\\
It is quite obvious, moreover, that a choice of  the parameters ``$x_{l}$''  that satisfies the constraints, must lead to a value ``$q_{i}$'' of the model belonging to the domain of definition of the discrepancy function.\\
For example, in image reconstruction, the unknowns ``$x_{l}$'' must be non-negative, and the transformation operation (convolution) must lead to non-negative images.\\

\subsection{Discrepancy (gap) functions - Divergences.}
	
	Most of the ``discrepancy functions'' or ``divergences" that appear in the literature are either Csiszär's divergences \cite{csiszar1967} \cite{csiszar1974}, Ali-Silvey\cite{ali1966}, either Bregman  divergences \cite{bregman1967}, or Jensen  divergences \cite{jensen1906} and, more generally , when they do not strictly belong to one of these categories (this is the case when considering the generalizations of these divergences), they are based on the use of strictly convex differentiable basic functions, and rely on the properties of such functions.\\
	This classification into different categories is purely formal, and it can be observed that a divergence may belong  simultaneously to more than one of these categories.\\
Furthermore, non-differentiable divergences such as variational distance will not be considered in this book.\\
The major difficulty with the divergences found in the literature lies in the fact that many works have been developed within the framework of the information theory \cite{{basseville1989},{basseville1996},{pardo2005},{taneja2001},{salicru1994}}. In this case, the data fields ``\textit{p}" and ``\textit{q}" are probability densities, which implies that the quantities considered are both non-negative and of sum equal to 1.\\
 Thus, many of the divergences used in this context are constructed specifically for such applications, or have simplifications related to this particular case.\\
On the other hand, for the applications considered in this book, two important points must be taken into account: on the one hand, the data fields involved are almost never probability densities, so that simplifications concerning the sum of the arguments are not necessary, on the other hand, it is not simply a question of quantifying a difference between the two data fields, but of minimizing a difference between ``\textit{p}" and ``\textit{q}" with respect to the unknown parameters in one of the two fields considered: the model ``\textit{q}".\\
 In addition, situations should be considered where the measurements ``\textit{p}" may be punctually negative due to noise, mainly in the case of Gaussian additive noise; if this is the case, it is important that the divergences used and the basic convex functions on which they are based, be defined for the negative values of the variable.\\
  If this is not the case, the alternative is to pre-process the data to make it non-negative.\\
These remarks imply to analyze precisely the divergences found in the literature and to make sure that they are well adapted to the considered problem.

\subsubsection{* Csiszär's divergences.}

The above remarks apply in particular to Csiszär divergences, the use of which requires some precautions. Indeed, their mode of construction implies the use of a strictly convex basis function $f\left(x\right)$, but moreover, for applications related to information theory (and more precisely if we are dealing with probability densities), we impose on the basis function to have the property $f\left(1\right)=0$ (this will correspond to particular convex functions that we will designate by ``\textit{simple convex functions}").\\
On the contrary, in our problems, the fields ``$p$" and ``$q$" have almost never spontaneously equal sums, (this will be one of the constraints of the problem), no more than sums equal to 1.\\
Then, in order to construct Csiszär divergences usable for our applications, the basic functions, in addition to the preceding properties, will have to be such that $f'\left(1\right)=0$ (this will correspond to convex functions which we will designate by ``\textit{standard convex functions}"); we will come back to this point in the following sections.

\subsubsection{* Bregman's and Jensen's divergences.} 

Unlike Csiszär's divergences, Bregman and Jensen's divergences, which are what one might call ``\textit{convexity measures}", imply only the use of a strictly convex base function $f\left(x\right)$ . Except in the particular case where simplifications have been introduced after construction of the divergence, their use does not pose a problem for our applications, as long as the particular divergence considered is convex with respect to the true unknowns of the problem. Indeed, as we will see in the following sections, the fact that a Bregman or Jensen divergence is constructed on a strictly convex base function does not necessarily imply its convexity.

\subsubsection{General properties of divergences.} 

Subject to differentiability, the divergences considered shall have the following properties:\\
\begin{itemize}
\item They'll have to be positive.\\
\item They'll have to be zero when $p_{i}=q_{i}$ $\forall i$.\\
\end {itemize}

Those two conditions are enough if you're comparing two different probability densities.\\ 
However, for applications to inverse problems that involve minimizing deviations from the model parameters ``$q$", additional properties are required:\\
 \begin{itemize}
\item Their gradient with respect to ``$q$" must be zero for $p_{i}=q_{i}$ $\forall i$.\\
\item They'll have to be convex to the true unknowns "$x$". (usually going through the convexity in relation to ``$q$").\\
\end {itemize}
Generally speaking:
\begin{itemize}
\item They are not necessarily symmetrical; if the need for this property arises for a particular problem, one can always build a ``\textit{ad hoc}" divergence.\\
\item They do not necessarily meet the triangular inequality, so they are not necessarily distances (but this is not prohibited).\\
\end{itemize}

The consequence of these remarks is that a few checks are required before using a divergence, and this explains in part the analyses that are carried out in each case.

\subsection{Choice of a divergence - Algorithmic aspects.}
 
From a physical point of view, when considering the data attachment term, the choice of the gap function may be dictated by the statistical properties of the measurement noise, as is the case for maximum likelihood methods \cite{Idier01a}.\\
When the measurements are simulated, this point does not present any difficulty because the noise statistics are perfectly known. On the other hand, in the case of real measurements and in particular when a pre-treatment has been carried out, it is necessary to estimate the probability density of the residual noise before use.\\
In the Bayesian framework \cite{Idier01a}, given the lack of information on the real nature of the noise, it is generally proposed to model the probability density of the noise using a Gaussian law.\\
This is justified by the fact that it is the ``\textit{least compromising}" choice given the lack of information on the exact nature of the measurement noise; one is thus led to minimize a divergence which is the root mean square deviation.\\
However, in some cases such as low photon count imaging, the use of a binomial distribution or Poisson's law seems more appropriate, leading to the minimization of a Kullback-Leibler divergence.   
Similarly, in \cite {fevotte2011}, it is proposed to use the ``$Beta$" divergence; a special case of this divergence is the Itakura-Saito divergence which seems to be better suited to the problem of musical signal processing.\\
Generally speaking, depending on the problem at hand, one could consider using any one of the various divergences discussed in this book to express the data attachment term.\\ 
As far as the minimization problem is concerned, the choice of method depends on the properties of the divergence used, and again the question is, which divergence to select?\\
From a strict optimization point of view, the simplest answer is to choose a strictly convex and differentiable discrepancy function, so that gradient-type descent methods can be used; however, for the chosen function, there must be a minimum corresponding to $p_{i}=q_{i}\ \forall i$, obtained by nulling the gradient of the discrepancy function (in the simplest case, i.e. without constraints). This is not always the case, as we will see in the following.\\
In addition, during the minimization process, constraints on the solution properties must be taken into account.\\
In this framework, properties reflecting strict constraints, such as the non-negativity of the components of the solution, or a constraint on the sum of these components, will be easily taken into account by conventional Lagrangian techniques of constrained minimization.\\
However, in the case of ill-posed problems  \cite{bertero1998}  \cite{Idier01a}, solutions obtained by estimating the maximum likelihood under constraints, using iterative methods (usually), reveal instabilities as the number of iterations increases; an acceptable solution can only be obtained by empirically limiting the number of iterations. This implicitly regularizes the problem.\\
The classical alternative to solve this instability problem is to perform an explicit regularization.\\
 The characteristic of these regularization methods consists in searching for a solution to make a tradeoff between the fidelity to the measured data and a fidelity to an information ``\textit{a priori}" \cite{titterington1985}.\\
To this end, one seeks to minimize a composite criterion made by introducing, next to the data attachment term , a penalty term allowing to reinforce certain suitable properties of the solution and which summarize our knowledge ``\textit{a priori}"; the relative importance of the two terms of the composite criterion thus made up being adjusted by a ``\textit{regularization factor}".\\
This ``\textit{energetic}" point of view is the one we'll adopt; the data attachment term and the penalty term will be expressed by two divergences of the type studied in this book; the penalty term expressing a ``\textit{gap}" between the current solution and a ``\textit{default}" solution reflecting the desirable  property(ies) of the solution.\\
Such a regularisation can be interpreted in a Bayesian framework cited by \cite{Idier01a} in which one must model a ``\textit{a priori}" law of probability of the solution, which should make it possible to take into account the desirable and known properties of the solution.\\
Applying Bayes' theorem yields the ``\textit{a posteriori}" distribution.\\
The estimation of the maximum of the ``\textit{a posteriori}"  law is equivalent to finding the minimum of the composite criterion.\\

\subsection{A few general considerations about the divergences.}

Basically, a divergence is used to express the difference between 2 data fields: Field 1 ($C1$) and Field 2 ($C2$), and the divergence is written $D\left(C1\|C2\right)$. With the notations indicated previously, the basic fields (if one can say so) are ``$p$" and ``$q$", but many authors have introduced a 3rd field which is the weighted sum of the basic fields, i.e. $\alpha p+\left(1-\alpha\right)q\ \ 0\leq \alpha \leq 1$; thus, 3 different fields are involved, and one can easily imagine the variants of the divergences related to the assignment of these 3 fields on the 2 fields appearing in the divergence.\\
Add to that the fact that the divergences are generally not symmetrical, which adds to the number of possible divergences\\
In any case, all of these divergences will always express a gap between the two basic fields ``$p$" and ``$q$."\\
Finally, to further add to the variety of possible divergences, we will rely on the fact that if a $D\left(p\|q\right)$ divergence is expressed as a difference of 2 positive terms $D\left(p\|q\right)=A\left(p,q\right)-B\left(p,q\right)$, a generalization can be introduced by applying to each of the terms $A\left(p,q\right)$ and $B\left(p,q\right)$ an increasing function (e.g. the generalized logarithm or the logarithm) without changing the (positive) sign of the divergence.\\
From these few remarks, one can see the wide variety of divergences that can appear.\\
This being said, the important problem when minimizing a divergence remains that of the convexity of the divergence considered in relation to the true unknowns of the problem.

\section{Application to the regularization.}

The resolution of the inverse problem as we have presented it, implies the minimization (under constraints) with respect to ``$x$", of a divergence between the measurements and the corresponding physical model $D_{1}\left(y\|m\left(x\right)\right)$.\
This divergence reflects the ``attachment to the data."\\
However, because of the ill-posed character of the inverse problems, one is led to carry out a regularization of the problem \cite{bertero1998,Idier01a,titterington1985} by introducing, next to the data attachment term, a ``\textit{penalty}" or ``\textit{regularization}"  term which makes it possible to give particular properties to the solution. In the most classical case, where (for example) a certain ``\textit{smoothness}" is imposed on the solution, this term is expressed as a discrepancy (divergence) $D_{2}\left(x\|x_{d}\right)$ between the current solution and a ``default solution $x_{d}$" having the required properties (in this case a smoothed version of the solution); on the other hand, it is reasonable to think that this default solution must also fulfill the constraints imposed on the solution.\\
We are thus led to minimize with respect to ``$x$", under constraints, an expression of the form:
\begin{equation}
	J\left(x\right)=D_{1}\left(y\|m\left(x\right)\right)+\gamma D_{2}\left(x\|x_{d}\right)
\end{equation}

In this expression, ``$\gamma$" is the positive regularization factor.\\
The constraints considered in this book are the non-negativity constraint and the sum constraint of the components of the solution:
\begin{equation}
	x_{i}\geq 0\ \ \forall i\ \ \ \ \ \ ;\ \ \ \ \ \ \sum_{i}x_{i}= C
\end{equation}

The terms $D_{1}$ and $D_{2}$ are divergences that will be analyzed in this work.\\
After having analyzed a certain number of classical divergences, and in order to simply take into account the sum constraint, we will introduce the scale invariant divergences  and we will show the advantages of the latter to satisfy this constraint.
In particular, we will show that scale invariant divergences with respect to both arguments are particularly adapted to the problem of regularization by smoothness constraint in the sense of Tikhonov \cite{tikhonov1974methods}, and, more specifically, when the default solution ``$x_{d}$" depends explicitly on the variable ``$x$".

\setcounter{table}{0}  \setcounter{equation}{0}  \setcounter{figure}{0} \setcounter{chapter}{2} \setcounter{section}{0} 
\chapter{chapter 2 -\\Convex functions and divergences.}  \label{chptr::chapitre2}

In this chapter, we recall some properties of the convex functions \cite{bertsekas1995, boyd2004, hiriart2012, rockafellar2015} that will be useful in constructing the divergences.\\
We then look at the constructive modes of different types of classical divergences. We study their convexity and give some relationships that link them.\\

\section{Convex functions - Some properties.}
\subsection{Generalities.}
From a strictly convex function $f\left(x\right)$, we define the ``mirror function" (Basseville \cite{basseville1989, basseville1996}) or ``dual function" or ``*conjugated function" (Osterreicher \cite{osterreicher2002}) by:
\begin{equation}
	\breve{f}\left(x\right)=x f\left(\frac{1}{x}\right)
\end{equation}

This function will help build the Csiszär dual divergences.\\
The properties of the dual function will be:\\
* if $f\left(x\right)$ is convex, then $\breve{f}\left(x\right)$ is convex .\\
* if $f\left(1\right)=0$, then we have $\breve{f}\left(1\right)=0$.\\
A convex function with this property will be referred to as a \textit{``\textbf{simple convex function.}"}\\
* if $f\left(1\right)=$0, we have: $\breve{f^{'}}\left(1\right)=-f^{'}\left(1\right)$, which leaves the possibility of having $\breve{f^{'}}\left(1\right)=-f^{'}\left(1\right)=0$ in the case of ``self-mirror" functions.\\
From functions $f\left(x\right)$ and $\breve{f}\left(x\right)$, we define the function:
\begin{equation}
	\hat{f}\left(x\right)=\frac{f\left(x\right)+\breve{f}\left(x\right)}{2}
	\label{eq.fchap}
\end{equation}

This function will allow us to construct symmetrical Csiszär divergences.\\
For this function, we have the following properties:\\
* si $f\left(1\right)=0$, on a $\hat{f}\left(1\right)=0$ \\
* si $f\left(1\right)=0$, on a $\hat{f}^{'}\left(1\right)=0$ \\

This last property will be fundamental to our purpose.\\
This leads us, starting from a function $f\left(x\right)$ with no other property than the strict convexity, to define the function $f_{c}\left(x\right)$ which we can designate by \textit{``\textbf{standard convex function}"}:
\begin{equation}
	f_{c}\left(x\right)=f\left(x\right)-f\left(1\right)-\left(x-1\right)f'\left(1\right)
	\label{eq.fc}
\end{equation}

This function will be strictly convex, with $f_{c}\left(1\right)=0$ and $f'_{c}\left(1\right)=0$, so it will be positive or zero.\\
If we divide the right-hand side of (\ref{eq.fc}) by $f''\left(1\right)$, then we get $f_{c}''\left(1\right)=$1; the latter property is not necessary for our purposes.\\

Contrary to the \textit{``simple convex functions"}, the \textit{``standard convex functions"} will allow us to construct Csiszär divergences that can be used with data fields that are not probability densities, the so-called by Zhang \cite{zhang2004} ``measure invariant divergences" or ``measures extended to allow denormalized densities".\\

\textbf{\*Remarks}
\begin{itemize}
\item a function of the type $\hat{f}\left(x\right)$ is a \textit{``standard convex function"}.
\item the reciprocal is not true.
\item for the ``auto-mirror" functions, we have $f\left(x\right)=\breve{f}\left(x\right)=\hat{f}\left(x\right)=f_{c}\left(x\right)$.
\end{itemize}
\subsection{Some properties of convex functions.}

For a convex $f\left(x\right)$ function, provided that we remain in its definition domain, that its derivative is defined and continuous on the same domain, we have the following properties:
\begin{equation}
	\left(q-p\right)f'\left(p\right)\leq f\left(q\right)-f\left(p\right)
	\label{eq.fcv1}
\end{equation}
\begin{equation}
	f\left(q\right)-f\left(p\right)\leq \left(q-p\right)f'\left(q\right)
	\label{eq.fcv2}
\end{equation}

\begin{figure}[ht!]
\centering
\includegraphics[width=\linewidth]{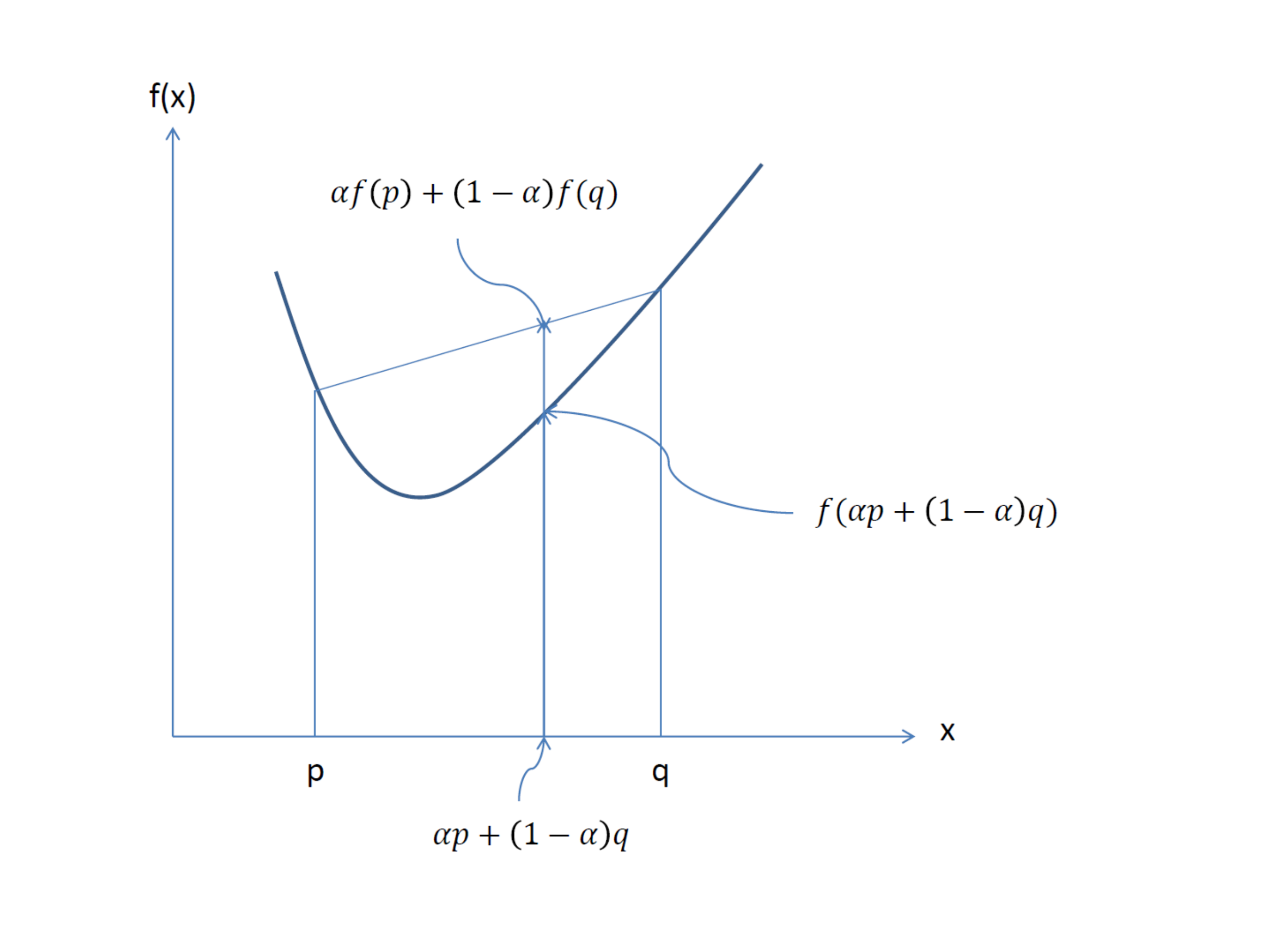}
\caption{Jensen divergence}
\label{fig:jensen}
\end{figure}

We also have Jensen's inequality \cite{jensen1906}:  for $\alpha_{i}>0\ \forall i$ and $\sum_{i}\alpha_{i}=1$.  
\begin{equation}
	\sum_{i}\alpha_{i}f\left(x_{i}\right)\geq f\left(\sum_{i}\alpha_{i}x_{i}\right)
	\label{eq.jen1}
\end{equation}
This is writen with 2 points (see figure (\ref{fig:jensen}):
\begin{equation}
	\alpha f\left(p\right)+\left(1-\alpha\right)f\left(q\right)\geq f\left[\alpha p+\left(1-\alpha\right)q\right]
	\label{eq.jen2}
\end{equation}

\begin{figure}[ht!]
\centering
\includegraphics[height=\linewidth, angle=-90]{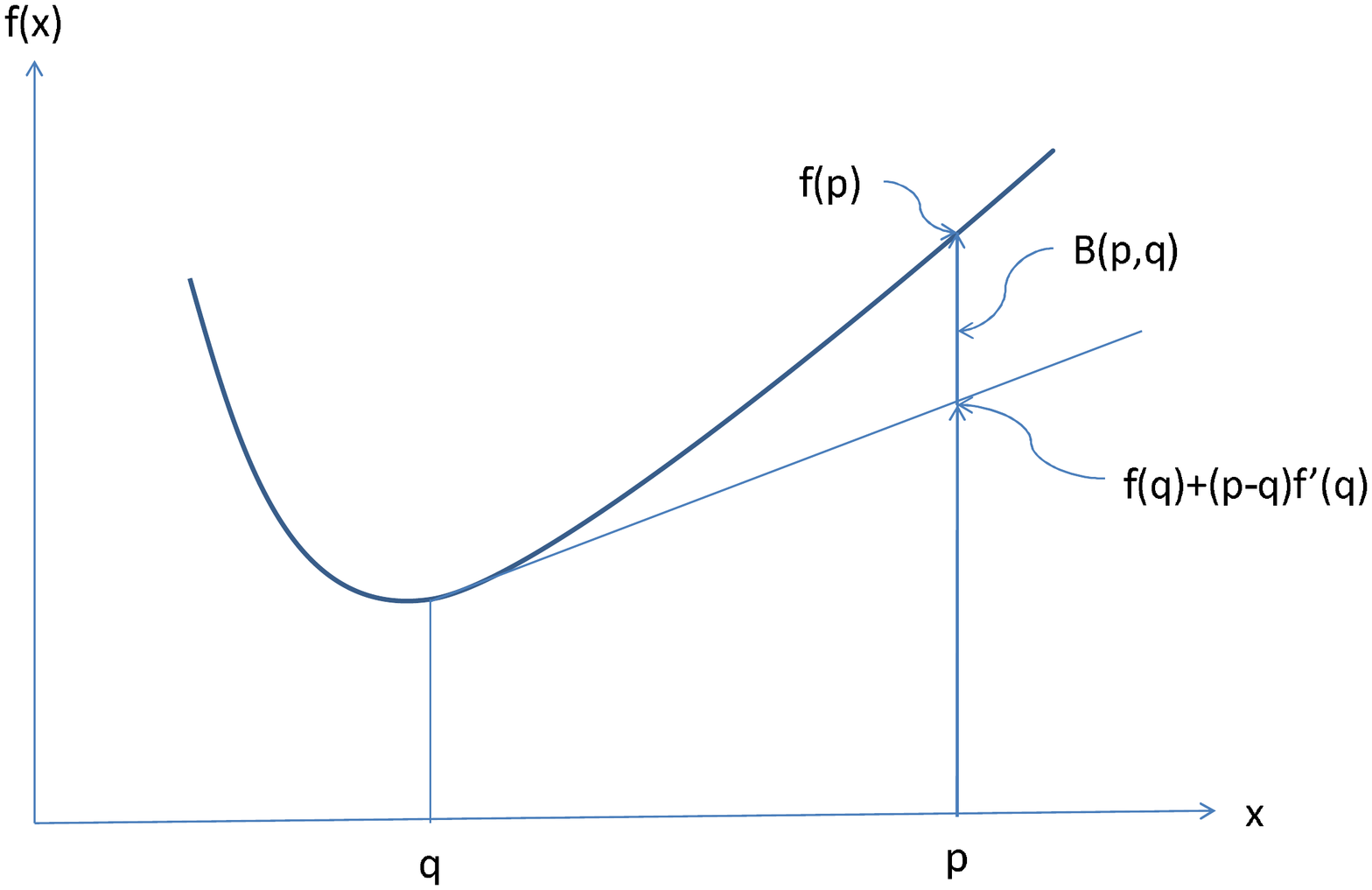}
\caption{Bregman divergence}
\label{fig:bregman}
\end{figure}

In the basic case, $\alpha=$1/2, it comes:
\begin{equation}
	\frac{f\left(p\right)+f\left(q\right)}{2}\geq f\left(\frac{p+q}{2}\right)
	\label{eq.jen3}
\end{equation}
On the other hand, from (\ref{eq.fcv2}), we can write:
\begin{equation}
	f\left(p\right)-f\left(q\right)-\left(p-q\right)f'\left(q\right)\geq 0
	\label{eq.fcv2bis}
\end{equation}

This inequality will be at the origin of Bregman's divergences, \cite{bregman1967}, see figure(\ref{fig:bregman}).\\
Similarly, from (\ref{eq.fcv1}), we have:
\begin{equation}
	f\left(q\right)-f\left(p\right)-\left(q-p\right)f'\left(p\right)\geq 0
	\label{eq.fcv1bis}
\end{equation}

Which will allow us to build the dual Bregman divergence of the previous one.
Finally, by summing up (\ref{eq.fcv1bis}) and (\ref{eq.fcv2bis}), we have:
\begin{equation}
	\left(q-p\right)\left[f'\left(q\right)-f'\left(p\right)\right]\geq 0
	\label{eq.fcv3}
\end{equation}

This inequality will intervene in the Burbea-Rao divergences \cite{burbea1982}.\\
Note that Burbea and Rao define a more general expression that is written:
\begin{equation}
	\left(q-p\right)\left[\frac{\phi\left(q\right)}{q}-\frac{\phi\left(p\right)}{p}\right]\geq 0
\end{equation}

where $\frac{\phi\left(x\right)}{x}$ is an increasing function, which coincides with (\ref{eq.fcv3}) if $f\left(x\right)$ is convex.

\textbf{More generally, it can be considered that any inequality can give rise to a divergence.}\\

The question of the definition domain is essential, because, as we shall see, this is what often restricts the field of use of the divergences to non-negative variables (for example).\\
In the general case that will be ours, the variables ``$p$" and ``$q$" are discretized (with the same step size); they are considered as vectors with components ``$p_{i}$" and ``$q_{i}$".\\
A divergence between ``$p$" and ``$q$" will be written in the general case as the separable form:
\begin{equation}
	D\left(p\|q\right)=\sum_{i}d\left(p_{i}\|q_{i}\right)
	\label{eq.D}
\end{equation}
And, from time to time, in the non-separable form, which is more delicate to analyse:
\begin{equation}
	D\left(p\|q\right)=h\left[\sum_{i}d\left(p_{i}\|q_{i}\right)\right]
\end{equation}
If we limit ourselves (for that time) to the separable form (\ref{eq.D}), the convexity of this divergence passes through the convexity of one of the terms of the sum: $d\left(p_{i}\|q_{i}\right)$.\\
We rely on the following definition:\\

\textbf{Definition}: The term $d\left(p_{i}\|q_{i}\right)$ is jointly convex with respect to ``$p_{i}$" and ``$q_{i}$" if its Hessian matrix $H$ is positive defined ; it is convex with respect to ``$p_{i}$" if $h_{11}$ is positive and convex with respect to ``$q_{i}$" if $h_{22}$ is positive.

\section{Csiszär's divergences - f (or I) divergences.}

These divergences were introduced by Csiszär \cite{csiszar1967} and simultaneously by Ali and Silvey \cite{ali1966}.\\
Lets consider a strictly convex function $f\left(x\right)$, a Csiszär divergence between two data fields "$p$" and "$q$" constructed on the function $f\left(x\right)$ is written:
\begin{equation}
	C_{f}\left(p\|q\right)=\sum_{i}c_{f}\left(p_{i}\|q_{i}\right)=\sum_{i}q_{i}f\left(\frac{p_{i}}{q_{i}}\right)
	\label{eq.Cf}
\end{equation}

For applications in information theory, the fields ``\textit{p}" and ``\textit{q}" are probability densities, they are positive quantities with equal sums (to 1 moreover); in this context, we simply impose on the basic function, the property: $f\left(1\right)=0$ \cite{pardo2005} \cite{osterreicher2002}.\\
As previously stated, a convex function possessing only this property will be referred to as: ``\textit{simple convex function}".\\
Divergences built on such a function are not usable outside of this context without some precautions.\\
It should also be noted that if in a Csiszâr divergence, one explicitly introduces the fact that the variables are specifically summed to 1, the resulting simplified divergence is no longer a Csiszâr divergence in the strict sense, since there is no convex function to construct it according to the relation (\ref{eq.Cf}).\\
A symmetrical Csiszär divergence can always be obtained in Jeffreys' sense \cite{jeffreys1946}  applying the constructive process (\ref{eq.Cf}) using the base function $\hat{f}\left(x\right)$ (\ref{eq.fchap}).\\
In our problems, as noted above, the ``\textit{standard convex function}" will allow us to construct Csiszär's divergences that can be used with data fields that are not probability densities, what Zhang \cite{zhang2004} calls ``measure invariant divergences" or ``measures extended to allow denormalized densities".\\
To show the necessity of using standard convex functions, let's consider Csiszär's divergences built on a simple convex function $f\left(x\right)$ on the one hand, and on the standard convex function $f_{c}\left(x\right)$ associated with it on the other hand.
In the first case, the gradient with respect to ``$q_{j}$'' is written:
\begin{equation}
	\frac{\partial C_{f}\left(p\|q\right)}{\partial q_{j}}=f\left(\frac{p_{j}}{q_{j}}\right)-\frac{p_{j}}{q_{j}}f^{'}\left(\frac{p_{j}}{q_{j}}\right)
\end{equation}
If we want this gradient to be zero for $p_{j}=q_{j}\ \forall j$, it requires that $f\left(1\right)=f^{'}\left(1\right)$.\\
Given the property $f\left(1\right)=0$, that implies that one must have: $f^{'}\left(1\right)=0$.\\
This observation obviously leads to the standard convex function deduced from $f\left(x\right)$:
\begin{equation}
	f_{c}\left(x\right)=	f\left(x\right)-f\left(1\right)-\left(x-1\right)f^{'}\left(1\right)
\end{equation}
The gradient with respect to ``$q_{j}$'' of the Csiszär's divergence built on $f_{c}\left(x\right)$ will be written:
\begin{equation}
	\frac{\partial C_{f_{c}}\left(p\|q\right)}{\partial q_{j}}=f\left(\frac{p_{j}}{q_{j}}\right)-\frac{p_{j}}{q_{j}}f^{'}\left(\frac{p_{j}}{q_{j}}\right)-f\left(1\right)+f^{'}\left(1\right)
\end{equation}
It will be spontaneously zero for $p_{j}=q_{j}\ \forall j$.\\
In the next section, we give two examples to illustrate this point.\\
Let's note that in the context of inverse problems, in general, if we want to use simplified divergences, we must take into account the simplifications that have been introduced, for example by changing variables, that is, by introducing normalized variables.\\

\subsection{A few examples to illustrate these difficulties.}
\subsubsection{Exemple 1.}
We consider the standard convex function:
\begin{equation}
	f_{c}\left(x\right)=\left(x-1_{}\right)^{2}
\end{equation}
The corresponding Csiszär divergence is known as Neyman's Chi2:
\begin{equation}
	\chi_{N}^{2}\left(p\|q\right)=\sum_{i}\frac{\left(p_{i}-q_{i}\right)^{2}}{q_{i}}
\end{equation}
Its gradient with respect to ``$q$" will be written:
\begin{equation}
	\frac{\partial \chi_{N}^{2}\left(p\|q\right)}{\partial q_{j}}=1-\frac{p_{j}^{2}}{q_{j}^{2}}
\end{equation}
The components of the gradient can cancel for $p_{j}=q_{j}$.\\
Let's now move to the simplified form in which, we introduced: $\sum_{i}p_{i}=\sum_{i}q_{i}$; the divergence becomes:
\begin{equation}
		\chi_{N}^{2}\left(p\|q\right)=\sum_{i}\frac{p_{i}^{2}}{q_{i}}-p_{i}
\end{equation}
This divergence is constructed in the sense of Csiszär on the simple convex function $f\left(x\right)$, associated with $f_{c}\left(x\right)$ which is written:
\begin{equation}
	f\left(x\right)=x^{2}-x
\end{equation}
The corresponding gradient will be::
\begin{equation}	
	\frac{\partial \chi_{N}^{2}\left(p\|q\right)}{\partial q_{j}}=-\frac{p_{j}^{2}}{q_{j}^{2}}
	\label{eq.gradki2}
\end{equation}
Now, the gradient doesn't cancel out for $p_{j}=q_{j}$, it doesn't even cancel out at all unless $q_{j}\rightarrow \infty$, hence the problem appear.\\
We can associate with the same standard convex function another simple convex function that's written:
\begin{equation}
	f\left(x\right)=x^{2}-1	
\end{equation}
The corresponding Csiszär divergence is given by: 
\begin{equation}
		\chi_{N}^{2}\left(p\|q\right)=\sum_{i}\frac{p_{i}^{2}}{q_{i}}-q_{i}
\end{equation}
The corresponding gradient will be:
\begin{equation}	
	\frac{\partial \chi_{N}^{2}\left(p\|q\right)}{\partial q_{j}}=-\frac{p_{j}^{2}}{q_{j}^{2}}-1
\end{equation}
As in the previous example, it never cancels.\\
Let's go even further and set the sum of the variables to 1: $\sum_{i}p_{i}=\sum_{i}q_{i}=1$; then, the divergence becomes:
\begin{equation}
		\chi_{N}^{2}\left(p\|q\right)=\left[\sum_{i}\frac{p_{i}^{2}}{q_{i}}\right]-1
\end{equation}
This function taken as such is no longer a Csiszär divergence, there is no basic convex function to obtain it according to the relation (\ref{eq.Cf}), its gradient with respect to ``$q$" is given by (\ref{eq.gradki2}).

\subsubsection{Example 2.}
Lets consider the standard convex function:
\begin{equation}
	f_{c}\left(x\right)=x\log x+x-1
\end{equation}
The Csiszär divergence built on this function is the Kullback-Leibler divergence \cite{kullback1951}.\begin{equation}
	KL\left(p\|q\right)=\sum_{i}p_{i}\log\frac{p_{i}}{q_{i}}+q_{i}-p_{i}
	\label{eq.KL}
\end{equation}
Its gradient with respect to ``$q$" will be:
\begin{equation}
	\frac{\partial 	KL\left(p\|q\right)}{\partial q_{j}}=1-\frac{p_{j}}{q_{j}}
\end{equation}
It will be zero for $p_{j}=q_{j}$.\\
On the other hand, if we introduce in the expression (\ref{eq.KL}), the simplification $\sum_{i}p_{i}=\sum_{i}q_{i}$, we obtain ``Kullback information" \cite{basseville1996}:
\begin{equation}
	IKL\left(p\|q\right)=\sum_{i}p_{i}\log\frac{p_{i}}{q_{i}}
\end{equation}
This expression is a Csiszär divergence built on the simple convex function:
\begin{equation}
	f\left(x\right)=x\log x
\end{equation}
Its gradient with respect to ``$q$" will be:
\begin{equation}
	\frac{\partial 	IKL\left(p\|q\right)}{\partial q_{j}}=-\frac{p_{j}}{q_{j}}
\end{equation}
He'll never be equal to zero.\\
Note that here, the extreme simplification $\sum_{i}p_{i}=\sum_{i}q_{i}=1$ will not add anything.

\subsection{Consequences of these examples.}

In our problem, it is a matter of making one of the two fields (the one that represents the model) evolve through the variations of the model parameters (the true unknowns), until the model is as close as possible to the measurements, in the sense of the divergence considered.\\

	The divergences used being convex with respect to the unknown parameters, the classical optimization methods always imply to look for the set of unknown parameters that corresponds to "zero gradient", (even in constrained problems, it provides part of the solution).\\
	
	 It is quite obvious from the previous remarks and examples that if the divergence to be minimized, although convex, does not have a finite minimum, or if the minimum obtained by this method does not have a suitable physical meaning, it is inappropriate for our problem.\\
	
Therefore, in our problems, if the divergence used is a Csiszär divergence, it is imperative that we consider only those that are constructed on the basis of a standard convex function.\\
Furthermore, the divergences, regardless of how they are constructed, must not have undergone any simplification related to a particular application.\\ 
 
\subsection{Convexity of Csiszär's divergences.}
Taking into account (\ref{eq.D}) and (\ref{eq.Cf}), we must calculate the Hessian of a term of the form:
\begin{equation}
c_{f}\left(p\|q\right)=q f\left(\frac{p}{q}\right)
\label{eq.termebase}
\end{equation}
We classically consider that "$p$" and "$q$" are positive, and we have:
\begin{equation}
	\frac{\partial^{2}c_{f}\left(p\|q\right)}{\partial p^{2}}=\frac{1}{q}f"\left(\frac{p}{q}\right)>0
	\label{eq.d2p}
\end{equation}
\begin{equation}
	\frac{\partial^{2}c_{f}\left(p\|q\right)}{\partial q^{2}}=\frac{p^{2}}{q^{3}}f"\left(\frac{p}{q}\right)>0
	\label{eq.d2q}
\end{equation}
\begin{equation}
	\frac{\partial^{2}c_{f}\left(p\|q\right)}{\partial q \partial p}=\frac{\partial^{2}c_{f}\left(p\|q\right)}{\partial p \partial q}=\frac{p}{q^{2}}f"\left(\frac{p}{q}\right)>0
\end{equation}
Separate convexity is clear from (\ref{eq.d2p}) and (\ref{eq.d2q}).\\
The joint convexity is analyzed by considering the Hessian determinant that is written:
\begin{equation}
\frac{p^{2}}{q^{4}}\left[f"\left(\frac{p}{q}\right)\right]^{2}-\left[\frac{p}{q^{2}}f"\left(\frac{p}{q}\right)\right]^{2}=0 	
\end{equation}

One of the eigenvalues is zero, the other is equal to the ``Trace" of the Hessian matrix, therefore positive; the expression (\ref{eq.termebase}) is therefore jointly convex.\\
The Csiszär divergence is therefore jointly convex as a sum of jointly convex terms.

\section{Bregman's divergences .}
These divergences are typically  convexity measurements.\\
They are based on a property of convex functions; they therefore imply the use of a basis convex function.\\
From this function, nothing is required other than strict convexity.\\
The property used to construct these divergences is expressed as:
\begin{itemize}
\item a convex curve is always located above any tangent to that curve; the Bregman divergence \cite{bregman1967},\cite{censor1997} is the difference between the curve and the tangent (taken in that order).
\item we can also say that for a strictly convex $f\left(x\right)$ function, the Bregman divergence is the difference between the function and its first order Taylor's development.
\end{itemize}
\begin{equation}	B_{f}\left(p\|q\right)=\sum_{i}b_{f}\left(p_{i}\|q_{i}\right)=\sum_{i}\left[f\left(p_{i}\right)-f\left(q_{i}\right)-\left(p_{i}-q_{i}\right)f^{'}\left(q_{i}\right)\right]
	\label{eq.Bf}
\end{equation}

Of course, since it is a basic property of convexity, a simple convex function and the corresponding standard convex function, which even have a second derivative, lead to the same Bregman's divergence, in other words, two functions that differ from each other by a linear function lead to the same Bregman's divergence.\\
Therefore, as one would expect, whether the divergence is constructed on the simple convex function or on the associated standard convex function, the gradients with respect to ``$q$" lead to the same expression being written:
\begin{equation}
	\frac{\partial B_{f}\left(p\|q\right)}{\partial q_{j}}=\frac{\partial B_{f_{c}}\left(p\|q\right)}{\partial q_{j}}=\left(q_{j}-p_{j}\right)f^{''}\left(q_{j}\right)=\left(q_{j}-p_{j}\right)f_{c}^{''}\left(q_{j}\right)
\end{equation}
That expression is obviously zero for $p_{j}=q_{j}\ \forall j$.\\ 
The convexity of Bregman's divergences is studied using the Hessian for one of the terms of the sum.\\

\textbf{Proposal:} If $f\left(x\right)$ is convex, $B_{f}\left(p\|q\right)$ is always convex with respect to the first ``$p$" argument, but may be non convex with respect to the second argument ``$q$".\\

\textbf{Demonstration:} From one of the sum terms appearing in (\ref{eq.Bf}), we calculate the elements of the Hessian.\\
\begin{equation}
	\frac{\partial^{2}b_{f}\left(p\|q\right)}{\partial p^{2}}=f^{''}\left(p\right)\geq 0
	\label{eq.d2bp}
\end{equation}
\begin{equation}
	\frac{\partial^{2}b_{f}\left(p\|q\right)}{\partial q^{2}}=f^{''}\left(q\right)+\left(q-p\right)f^{'''}\left(q\right)
	\label{eq.d2bq}
\end{equation}
\begin{equation}
	\frac{\partial^{2}b_{f}\left(p\|q\right)}{\partial p\partial q}=\frac{\partial^{2}b_{f}\left(p\|q\right)}{\partial q\partial p} =-f^{''}\left(q\right)
\end{equation}

Thus, (\ref{eq.d2bp}) implies convexity with respect to ``$p$''.\\
The sign of (\ref{eq.d2bq}) depends of course on $f\left(x\right)$, so the convexity with respect to ``$q$" depends on $f\left(x\right)$.\\
For joint convexity, we express the determinant of the Hessian matrix:
\begin{equation}	Det_{H}\left(p,q\right)=-\left[f^{''}\left(q\right)\right]^{2}+f^{''}\left(p\right)f^{''}\left(q\right)-\left(p-q\right)f^{''}\left(p\right)f^{'''}\left(q\right)
\end{equation}
Dividing by $f^{''}\left(p\right)\left[f^{''}\left(q\right)\right]^{2}> 0$, we have:
\begin{equation}
\frac{Det_{H}\left(p,q\right)}{f^{''}\left(p\right)\left[f^{''}\left(q\right)\right]^{2}}=-\frac{1}{f^{''}\left(p\right)}+\frac{1}{f^{''}\left(q\right)}-\left(p-q\right)\left[-\frac{1}{f^{''}\left(q\right)}\right]^{'}	
\end{equation}
This quantity is a Bregman's divergence built on the function $\left[-\frac{1}{f^{''}\left(x\right)}\right]$, it will be positive if $\left[-\frac{1}{f^{''}\left(x\right)}\right]$ is convex, i.e. if $\left[\frac{1}{f^{''}\left(x\right)}\right]$ is concave.\\
 If so, the Hessian determinant will be positive, the Hessian will be positive defined  and the corresponding Bregman divergence will be ``jointly convex".\\

\textbf{Property:} For a strictly convex $f\left(x\right)$ function, Bregman's divergence is jointly convex with respect to ``$p$" and ``$q$" if and only if $\left[1/f^{''}\left(x\right)\right]$ is concave.\\
 We'll get a similar property for Jensen's divergences.\\
From the expression (\ref{eq.fcv1}), we obtain the adjoint (dual) form of the Bregman divergence:
\begin{equation}	B_{f}\left(q\|p\right)=\sum_{i}b_{f}\left(q_{i}\|p_{i}\right)=\sum_{i}\left[f\left(q_{i}\right)-f\left(p_{i}\right)-\left(q_{i}-p_{i}\right)f^{'}\left(p_{i}\right)\right]
\end{equation}
Following the same reasoning as for $B_{f}\left(p\|q\right)$ we can easily show that $B_{f}\left(q\|p\right)$ is jointly convex to ``$p$" and ``$q$" if and only if $\left[1/f^{''}\left(x\right)\right]$ is concave.\\
Similarly, from the expression (\ref{eq.fcv3}), we can deduce a form of divergence proposed by Burbea and Rao \cite{burbea1982} which is expressed as follows:
\begin{equation}	
A_{f}\left(p,q\right)=\sum_{i}a_{f}\left(p_{i},q_{i}\right)=\sum_{i}\left[\left(p_{i}-q_{i}\right)\left(f^{'}\left(p_{i}\right)-f^{'}\left(q_{i}\right)\right)\right]
\end{equation}
This symmetrical divergence is related to Bregman's $B_{f}\left(p\|q\right)$ and $B_{f}\left(q\|p\right)$ divergences by the relationship:
\begin{equation}
	A_{f}\left(p,q\right)=B_{f}\left(p\|q\right)+B_{f}\left(q\|p\right)
\end{equation}

By relying on the convexity of Bregman's divergences, $A_{f}\left(p,q\right)$ is jointly convex with respect to ``$p$" and ``$q$" if and only if $\left[1/f^{''}\left(x\right)\right]$ is concave.\\ 
It should be noted, however, that $A_{f}\left(p,q\right)$ is not a Bregman divergence in that there is no convex function to construct it directly.

\subsection{Example.}
We consider the functions $f_{1}\left(x\right)=x^{2}-1$ , $f_{2}\left(x\right)=x^{2}-x$ and $f_{c}\left(x\right)=\left(x-1\right)^{2}$; these convex functions defined for any $x$, differ from each other by a linear function; only $f_{c}\left(x\right)$ is a standard convex function, but all lead to the same Bregman's divergence which is the mean square deviation (which is otherwise symmetric and respects the triangular inequality; it is therefore a distance).

\subsection{Some possible variations.}
Considering the Bregman divergence, based on a strictly convex function $f\left(x\right)$, between a convex combination of the variables ``$\alpha p+\left(1-\alpha\right)q$" and ``$p$" or ``$q$", and taking into account the fact that the order of the arguments can be reversed, we can consider the divergence:
\begin{equation}
	B^{\alpha}_{f}\left(p\|q\right)=B_{f}\left[\alpha p+\left(1-\alpha\right)q\|q\right]
\end{equation}
which is written in more detail:
\begin{equation}
B^{\alpha}_{f}\left(p\|q\right)=\sum_{i}f\left[\alpha p_{i}+\left(1-\alpha\right)q_{i}\right]-f\left(q_{i}\right)-\alpha \left(p_{i}-q_{i}\right)\nabla f\left(q_{i}\right)
\end{equation}
But we can also look at the divergences:\\
* $B_{f}\left[\alpha p+\left(1-\alpha\right)q\|p\right]$\\
* $B_{f}\left[p\|\alpha p+\left(1-\alpha\right)q\right]$\\
* $B_{f}\left[q\|\alpha p+\left(1-\alpha\right)q\right]$\\

That being said, one may ask if it's of much interest...

\section{Jensen's divergences.}
These divergences are applications of Jensen's inequality \cite{jensen1906} (\ref{eq.jen1}), (\ref{eq.jen2}) or (\ref{eq.jen3}), with simple or standard convex functions.\\
 Since it is, as for Bregman's divergences, a convexity measure, two convex functions having  equal  second derivatives (i.e. differing from each other by a linear function) will lead to the same Jensen's divergence.\\
For a strictly convex $f\left(x\right)$ base function and for $0\leq \alpha\leq 1$, we define:
\begin{equation}
	j_{f}^{\alpha}\left(p_{i}\|q_{i}\right)=\alpha f\left(p_{i}\right)+\left(1-\alpha\right)f\left(q_{i}\right)-f\left[\alpha p_{i}+\left(1-\alpha\right)q_{i}\right]
\end{equation}
The corresponding Jensen's divergence is then:
\begin{equation}
J_{f}^{\alpha}\left(p\|q\right)=\sum_{i}j_{f}^{\alpha}\left(p_{i}\|q_{i}\right)
	\label{eq.jenalfa}	
\end{equation}
Therefore, as expected, whether the divergence is built on the simple convex function or on the associated standard convex function, the gradients with respect to ``$q$" lead to the same expression which is written:
\begin{equation}
	\frac{\partial J^{\alpha}_{f}\left(p\|q\right)}{\partial q_{j}}=\frac{\partial J^{\alpha}_{f_{c}}\left(p\|q\right)}{\partial q_{j}}=\left(1-\alpha\right)\left[f^{'}\left(q_{j}\right)-f^{'}\left(\alpha p_{j}+\left(1-\alpha\right)q_{j}\right)\right]
\end{equation}
That expression is obviously zero for $p_{j}=q_{j}\ \forall j$.\\
So the Jensen's 1/2 divergence will be expressed as follows:
\begin{equation}
J_{f}^{\frac{1}{2}}\left(p\|q\right)=\sum_{i}\left\{\frac{1}{2} f\left(p_{i}\right)+\frac{1}{2}f\left(q_{i}\right)-f\left[\frac{ p_{i}+q_{i}}{2}\right]\right\}
\label{eq.J1/2}	
\end{equation}
This last divergence is symmetrical.\\
Its convexity is subject to the following rule:\\

\textbf{Rule:} For a strictly convex $f\left(x\right)$ function, Jensen's ``$\alpha$" divergence  is jointly convex with respect to ``$p$" and ``$q$" if and only if $\left[1/f''\left(x\right)\right]$ is concave. (see. Burbea-Rao \cite{burbea1982}).\\

\textbf{Démonstration}: From the expression (\ref{eq.jenalfa}), we calculate the elements of the Hessian matrix of one of the terms of the sum:
\begin{equation}
	\frac{\partial^{2}j_{f}^{\alpha}\left(p,q\right)}{\partial p^{2}}=\alpha f^{''}\left(p\right)-\alpha^{2} f^{''}\left[\alpha p+\left(1-\alpha\right)q\right]
\end{equation}
\begin{equation}
	\frac{\partial^{2}j_{f}^{\alpha}\left(p,q\right)}{\partial q^{2}}=\left(1-\alpha\right)f^{''}\left(q\right)-\left(1-\alpha\right)^{2}f^{''}\left[\alpha p+\left(1-\alpha\right)q\right]
\end{equation}
\begin{equation}
	\frac{\partial^{2}j_{f}^{\alpha}\left(p,q\right)}{\partial p\partial q}=\frac{\partial^{2}j_{f}^{\alpha}\left(p,q\right)}{\partial q\partial p}=-\alpha\left(1-\alpha\right)f^{''}\left[\alpha p+\left(1-\alpha\right)q\right]
\end{equation}
The Hessian will be positive defined  if its determinant is positive, i.e. if:
\begin{equation}
	\begin{split}
\left\{f^{''}\left(p\right)-\alpha f^{''}\left[\alpha p+\left(1-\alpha\right)q\right]\right\}	\left\{f^{''}\left(q\right)-\left(1-\alpha\right) f^{''}\left[\alpha p+\left(1-\alpha\right)q\right]\right\}  \\ 
	-\alpha\left(1-\alpha\right)\left\{f^{''}\left[\alpha p+\left(1-\alpha\right)q\right]\right\}^{2}>0
	\end{split}
\end{equation}
Which gives immediately:
\begin{equation}
	f^{''}\left(p\right)f^{''}\left(q\right)-\left[\left(1-\alpha\right)f^{''}\left(p\right)+\alpha f^{''}\left(q\right)\right]f^{''}\left[\alpha p+\left(1-\alpha\right)q\right]>0
\end{equation}
And finally:
\begin{equation}
	\frac{\alpha}{f^{''}\left(p\right)}+\frac{1-\alpha}{f^{''}\left(q\right)}-\frac{1}{f^{''}\left[\alpha p+\left(1-\alpha\right)q\right]}<0
\end{equation}
Which expresses the concavity of $\left[1/f''\left(x\right)\right]$.\\
This result is identical to what we found with the Bregman's divergences.

\section{Relationship between these divergences.}
\subsection{Between Csiszär's and Bregman's divergences.}
In Censor and Zenios \cite{censor1997}, we find the relationship:
\begin{equation}
	C_{f}\left(p\|q\right)=\sum_{i}q_{i}\;f\left(\frac{p_{i}}{q_{i}}\right)=\sum_{i}q_{i}\;b_{f}\left(\frac{p_{i}}{q_{i}}\|1\right)
	\label{eq.relCB}
\end{equation}
Or, similarly:
\begin{equation}
	f\left(\frac{p_{i}}{q_{i}}\right)=b_{f}\left(\frac{p_{i}}{q_{i}}\|1\right)
	\label{eq.relCBS}
\end{equation}
We can show this easily, by expressing $b_{f}\left(\frac{p_{i}}{q_{i}}\|1\right)$ and considering that $f\left(1\right)=f^{'}\left(1\right)=0$, which is the case for the \textit{``standard convex functions"}, on the other hand, the relations (\ref{eq.relCB}) and (\ref{eq.relCBS}) are not valid for \textit{``simple convex functions"}.

\subsection{Between Jensen's and Bregman's divergences.}
* - In the direction: $Bregman\rightarrow Jensen$.\\
A first relation is given by Basseville \cite{basseville1989}, \cite{basseville1996}:
\begin{equation}
	J^{\alpha}_{f}\left(p\|q\right)=\alpha B_{f}\left[p\|\alpha p+\left(1-\alpha\right)q\right]+\left(1-\alpha\right)B_{f}\left[q\|\alpha p+\left(1-\alpha\right)q\right]
\end{equation}
But we can also establish another relationship that can be expressed as:
\begin{equation}
	J^{\alpha}_{f}\left(p\|q\right)=\alpha B_{f}\left(p\|q\right)-B_{f}\left[\alpha p+\left(1-\alpha\right)q\|q\right]
	\label{eq.jalpha}
\end{equation}
* - In the direction: $Jensen\rightarrow Bregman$.\\
Basseville \cite{basseville1989} gives the following relationship:
\begin{equation}
	B_{f}\left(p\|q\right)=\lim _{\alpha_\rightarrow 0}\;	J^{\alpha}_{f}\left(p\|q\right)
\end{equation}
But we can also establish the relationship:
\begin{equation}	B_{f}\left(p\|q\right)=\sum^{\infty}_{n=1}\frac{1}{\alpha^{n}}J^{\alpha}_{f}\left[\alpha^{n-1}p+\left(1-\alpha^{n-1}\right)q\|q\right]
\label{eq.devBf}
\end{equation}
To obtain this relationship, we start from expression (\ref{eq.jalpha}) we write:
\begin{equation}
	 B_{f}\left(p\|q\right)=\frac{1}{\alpha}J^{\alpha}_{f}\left(p\|q\right)+\frac{1}{\alpha}B_{f}\left[\alpha p+\left(1-\alpha\right)q\|q\right]
	 \label{eq.jalphabis}	
\end{equation}
Then replace $B_{f}\left[\alpha p+\left(1-\alpha\right)q\|q\right]\equiv B_{f}\left[z\|q\right]$ with its expression deduced from (\ref{eq.jalphabis}), that is:
\begin{equation}
B_{f}\left(z\|q\right)=\frac{1}{\alpha}J^{\alpha}_{f}\left(z\|q\right)+\frac{1}{\alpha}B_{f}\left[\alpha z+\left(1-\alpha\right)q\|q\right]	
\end{equation}
hence, by replacing ``$z$" with its expression:
\begin{equation}
B_{f}\left(\alpha p+(1-\alpha)q\|q\right)=\frac{1}{\alpha}J^{\alpha}_{f}\left(\alpha p+(1-\alpha)q\|q\right)+\frac{1}{\alpha}B_{f}\left[\alpha^{2} p+\left(1-\alpha^{2}\right)q\|q\right]	
\end{equation}
This expression is reported in (\ref{eq.jalphabis}) which gives:
\begin{equation}
B_{f}\left(p\|q\right)=\frac{1}{\alpha}J^{\alpha}_{f}\left(p\|q\right)+\frac{1}{\alpha^{2}}J^{\alpha}_{f}\left(\alpha p+(1-\alpha)q\|q\right)+\frac{1}{\alpha^{2}}B_{f}\left[\alpha^{2} p+\left(1-\alpha^{2}\right)q\|q\right]	
\end{equation}
and so on and so forth, recursively.\\

\subsection{Between Jensen's and Csiszär's divergences.}
Based on the previous results, we can establish the relationship that yields Csiszär's divergences when we know Jensen's divergences.\\
From (\ref{eq.devBf}), we have:
\begin{equation}	b_{f}\left(\frac{p_{i}}{q_{i}}\|1\right)=\sum^{\infty}_{n=1}\frac{1}{\alpha^{n}}J^{\alpha}_{f}\left[\alpha^{n-1}\frac{p_{i}}{q_{i}}+\left(1-\alpha^{n-1}\right).1\|1\right]
\end{equation}
Or also:
\begin{equation}	b_{f}\left(\frac{p_{i}}{q_{i}}\|1\right)=\sum^{\infty}_{n=1}\frac{1}{\alpha^{n}}J^{\alpha}_{f}\left[\frac{\alpha^{n-1}p_{i}+\left(1-\alpha^{n-1}\right)q_{i}}{q_{i}}\|1\right]
\end{equation}
Then, with (\ref{eq.relCB}):
\begin{equation}	C_{f}\left(p\|q\right)=\sum_{i}q_{i}\sum^{\infty}_{n=1}\frac{1}{\alpha^{n}}J^{\alpha}_{f}\left[\frac{\alpha^{n-1}p_{i}+\left(1-\alpha^{n-1}\right)q_{i}}{q_{i}}\|1\right]
\end{equation}

\setcounter{table}{0}  \setcounter{equation}{0}  \setcounter{figure}{0} \setcounter{chapter}{3} \setcounter{section}{0} 
\chapter{chapter 3 -\\Scale change invariance}  \label{chptr::chapitre3}

In this chapter, we are interested in the divergences invariant by change of scale and we indicate their construction mode; we also specify for this type of divergence, some useful properties to build minimization algorithms under sum and non-negativity constraints.
We will consider in this chapter that the variables ``$p$" and ``$q$" involved in the expression of divergences are non-negative. 
\section{Introduction.}
In the context of linear inverse problems, one is generally led to minimize with respect to the true unknown ``$x$", a divergence between measures $y\equiv p$ and a linear model $Hx\equiv q$; this minimization is frequently associated with a non-negativity constraint $x_{i}\geq 0\ \forall i$ and a sum constraint of the type $\sum_{i}x_{i}=Cte$.\\
In order to simply take into account this last constraint, we develop a class of divergences which are invariant by changing the scale with respect to ``$q$", i.e. such that:
\begin{equation}
	DI\left(p\|q\right)=DI\left(p\|aq\right)
	\label{eq.divinv}
\end{equation}
where ``$DI$" denotes an invariant divergence and ``$a$" is a positive scalar.\\
The underlying idea being that during the iterative process of minimization, after each iteration, we will be able to renormalize ``$x$" and respect the sum constraint, without changing the value of the divergence to be minimized; indeed, according to (\ref{eq.divinv}) we will obviously have:
\begin{equation}
DI\left(y\|Hx\right)=DI\left(y\|a\left(Hx\right)\right)=D\left(y\|H\left(ax\right)\right)	
\end{equation}
This being the case, we will show that the divergences invariant by change of scale possess a property that will be particularly interesting when the non-negativity constraint is associated with the sum constraint in a minimization problem such as the one mentioned above.\\
Some examples of applications of invariant divergences have been proposed in \cite{lanteri2013}, \cite{lanteri2015}, \cite{lanteri2015x} .

\section{Invariance factor $K$ - General properties.}

The general idea is that starting with a divergence of any kind, we're going to transform it into a divergence of $DI\left(p\|q\right)$ (related to $D\left(p\|q\right)$), which is invariant with respect to the variable ``$q$".\\
More precisely, in order to make a divergence scale invariant  with respect to ``$q$", we look for the expression of a \textbf{positive scalar factor} $K\left(p,q\right)$ such that the divergence $D\left(p\|Kq\right)=DI\left(p\|q\right)$ remains unchanged when the components of ``$q$" are multiplied by a positive scalar.\\

\textbf{The solution of this problem is not unique, indeed all expressions of $K\left(p,q\right)$ having the following properties are possible solutions of this problem:\\}
 
\textbf{1 - $K\left(p,q\right)$ must be scalar and positive.},\\
   
 \textbf{2 - In order to obtain a divergence which is invariant by scale change with respect to ``$q$", the vector $\left[K\left(p,q\right).\ q\right]$ must be invariant when multiplying ``$q$" by a constant.}\\

Finally, in general terms:\\

\textbf{* If $p\rightarrow q$, we must have $K\left(p,q\right)\rightarrow1$.}\\
 consequently:\\
 
 \textbf{ * If $p\ \rightarrow\ q$, $D\left(p\|Kq\right)\ \rightarrow\ D\left(p\|q\right)\ \rightarrow\ 0$.}\\

We will note that the constraints imposed on the $K\left(p,q\right)$ factor do not necessarily relate it to any given divergence.\\

\textbf{These observations lead to the following very important remark: an invariance factor having the properties listed above will make any divergence invariant with respect to ``$q$".}\\

\subsection{Calculation of the nominal invariance factor.} 
In order to obtain an expression of $K\left(p,q\right)$, the method defined in \cite{eguchi2010} consists of calculating the invariance factor: \\
\begin{equation}
	K_{0}\left(p,q\right)=\arg\min_{K>0}D\left(p\|Kq\right)
	\label{eq.defK0}
\end{equation}
In this mode of calculation, the invariance factor $K_{0}\left(p,q\right)$, which we will designate by
``Nominal invariance factor" is specifically associated with the divergence $D\left(p\|q\right)$ and therefore implicitly associated with a basic convex function and a constructive mode of divergence.\\
So we have to solve with respect to $K$ (positive), the equation:
\begin{equation}
	\frac{\partial D\left(p\|Kq\right)}{\partial K}=0
	\label{eq.resolK0}
\end{equation}

If we put the resulting $K$ expression, $K_{0}\left(p,q\right)$, into the divergence $D\left(p\|q\right)$, we get the invariant divergence $DI\left(p\|q\right)=D\left(p\|K_{0}q\right)$.\\
In this calculation, $K$ is considered a scalar quantity.\\

\textbf{However, it should be noted that (\ref{eq.resolK0}) does not necessarily have an explicit solution.}\\

\textbf{* Property.}\\

Taking into account the definition (\ref{eq.defK0}), the resulting invariant divergence $D\left(p\|K_{0}q\right)$ is less than or equal to any other invariant divergence derived from $D\left(p\|q\right)$ using an invariance factor $K_{1}\left(p,q\right)$ different from $K_{0}\left(p,q\right)$.
\begin{equation}
 D\left(p\|K_{0}q\right)\leq D\left(p\|K_{1}q\right)	
\end{equation}
With equality if $p_{i}=q_{i}\ \forall i$.\\
Furthermore, when $q\rightarrow p$, we have:\\
 $K_{0}\left(p,q\right)\ \rightarrow K_{1}\left(p,q\right)\rightarrow 1$\\
therefore\\
 $D\left(p\|K_{0}q\right)\ \rightarrow D\left(p\|K_{1}q\right)\ \rightarrow \ D\left(p\|q\right)\rightarrow 0$.\\
 
 \textbf{Examples allowing to show this more precisely are given in Annex 7.}
\section{A few comments regarding the invariance factor $K$.}
For divergences that remain invariant by change of scale, whatever their form, we note:
\begin{equation}
D\left(p\|Kq\right)=\sum_{i}d\left(p_{i}\|Kq_{i}\right)	
\end{equation}

\subsection{Fundamental property.}
For divergences invariant by scale change on ``$q$", we will first establish a fundamental property which is writtens:
\begin{equation}
	\sum_{j}q_{j}\frac{\partial D\left(p\|Kq\right)}{\partial q_{j}}=0
	\label{eq.Pfond}
\end{equation}

In a first step, we show that this relation is verified when the invariance factor is a nominal factor, i.e. when it is calculated explicitly by resolution of (\ref{eq.resolK0}); this invariance factor $K_{0}\left(p,q\right)$ is thus directly associated to the divergence $D\left(p\|q\right)$ considered.\\

\textbf{Démonstration.}\\

Knowing that $K$ is a function of ``$p$" and ``$q$", the gradient of $D\left(p\|Kq\right)$ with respect to ``$q$" is written:
\begin{equation}
	\frac{\partial D\left(p\|Kq\right)}{\partial q_{j}}=\sum_{i}\frac{\partial d\left(p_{i}\|Kq_{i}\right)}{\partial \left(Kq_{i}\right)}\frac{\partial \left(Kq_{i}\right)}{\partial q_{j}}
\end{equation}
But, we have:
\begin{equation}
\frac{\partial \left(Kq_{i}\right)}{\partial q_{j}}=\frac{\partial \left(Kq_{i}\right)}{\partial K}\frac{\partial K}{\partial q_{j}}+\frac{\partial \left(Kq_{i}\right)}{\partial q_{i}}\frac{\partial q_{i}}{\partial q_{j}}=q_{i}\frac{\partial K}{\partial q_{j}}+K \delta_{ij}
\label{eq.dKqi}	
\end{equation}
which leads to:
\begin{equation}
\frac{\partial D\left(p\|Kq\right)}{\partial q_{j}}=\sum_{i}\frac{\partial d\left(p_{i}\|Kq_{i}\right)}{\partial \left(Kq_{i}\right)}\frac{\partial K}{\partial q_{j}}q_{i}+K	\frac{\partial d\left(p_{j}\|Kq_{j}\right)}{\partial \left(Kq_{j}\right)}
\label{eq.Pfond0}
\end{equation}
From which it can be deduced:
\begin{equation}
\sum_{j}q_{j}\frac{\partial D\left(p\|Kq\right)}{\partial q_{j}}=\left[\sum_{i}\frac{\partial d\left(p_{i}\|Kq_{i}\right)}{\partial \left(Kq_{i}\right)}q_{i}\right]\left[\sum_{j}q_{j}\frac{\partial K}{\partial q_{j}}+K\right]
\label{eq.Pfondbis}
\end{equation}
With:
\begin{equation}
	q_{i}=\frac{\partial \left(Kq_{i}\right)}{\partial K}
\end{equation}
one can also write (\ref{eq.Pfondbis}) in the form:
\begin{equation}
\sum_{j}q_{j}\frac{\partial D\left(p\|Kq\right)}{\partial q_{j}}=\left[\sum_{i}\frac{\partial d\left(p_{i}\|Kq_{i}\right)}{\partial K}\right]\left[\sum_{j}q_{j}\frac{\partial K}{\partial q_{j}}+K\right]
\label{eq.bilan2}	
\end{equation}
The fundamental relationship (\ref{eq.Pfond}) is thus made explicit in the form:
\begin{equation}
\left[\sum_{i}\frac{\partial d\left(p_{i}\|Kq_{i}\right)}{\partial K}\right]\left[\sum_{j}q_{j}\frac{\partial K}{\partial q_{j}}+K\right]=0
\label{eq.prod}	
\end{equation}
Or also:
\begin{equation}
	\left[\sum_{i}\frac{\partial d\left(p_{i}\|Kq_{i}\right)}{\partial \left(Kq_{i}\right)}q_{i}\right]\left[K+\sum_{j}q_{j}\frac{\partial K}{\partial q_{j}}\right]=0
	\label{eq.Base0}
\end{equation}

So it's verified if one of the terms in the product (\ref{eq.prod}) is zero.\\
We examine each of these two terms in turn.

\subsubsection {First term of the product.}
The nominal invariance factor $K_{0}$ is calculated by solving with respect to $K$ the equation:
\begin{equation}
	\frac{\partial D\left(p\|Kq\right)}{\partial K}=\sum_{i}\frac{\partial d\left(p_{i}\|Kq_{i}\right)}{\partial K}=0
\end{equation}
that is:
\begin{equation}
	\sum_{i}q_{i}\frac{\partial d\left(p_{i}\|Kq_{i}\right)}{\partial \left(Kq_{i}\right)}=0
	\label{eq.resolK}
\end{equation}
For such an invariance factor $K_{0}\left(p,q\right)$, specifically associated with the divergence $D\left(p\|q\right)$  considered, the first term of the product (\ref{eq.prod}) is null; we can therefore conclude that:\\

\textbf{The relationship (\ref{eq.Pfond}) is verified if $K\left(p,q,\right)=K_{0}\left(p,q,\right)$.\\}

\textbf{This property is of fundamental importance as we will see in the chapter dealing with the algorithmic developments of minimization under non-negativity and sum constraints.}\\

We first examine the expressions of the term (\ref{eq.resolK}) corresponding to the different forms of classical divergences.

\subsubsection{Relationship with the different construction modes of divergences (Csiszär, Bregman, Jensen).}

The factor $K\left(p,q\right)$, when calculated by resolving (\ref{eq.resolK0}) (\ref{eq.resolK}), is specifically related to a given $D\left(p\|q\right)$ divergence, so it is associated with both a basic convex function and the constructive mode of the divergence (Csiszär, Bregman or Jensen); generalizations outside this framework will not be considered here. \\
We will therefore establish the relations that link the $K$ factor and the convex functions that allows us to construct the divergences. To do so, we will give the particular expressions of (\ref{eq.resolK}) for each of the three constructive modes of divergences.
\subsubsection{1 - Csiszär's divergences.}
The function $f\left(x\right)$ considered here is a standard convex function; the Csiszär divergence $C_{f}\left(p\|q\right)$ is constructed according to the relationship:
\begin{equation}
	C_{f}\left(p\|q\right)=\sum_{i}q_{i}f\left(\frac{p_{i}}{q_{i}}\right)
	\label{eq.dcsisz}
\end{equation}
The factor $K\left(p,q\right)$ making this divergence invariant with respect to ``$q$" is associated with the function ``$f$" and Csiszär's constructive mode. The invariant divergence corresponding to (\ref{eq.dcsisz}) will be written:
\begin{equation}
	C_{f}\left(p\|K q\right)=\sum_{i}c_{f}\left(p_{i}\|K q_{i}\right)=\sum_{i}K q_{i}\;f\left(\frac{p_{i}}{K q_{i}}\right)
\end{equation}
After a few simple calculations, we have:
\begin{equation}
	\frac{\partial c_{f}\left(p_{i}\|Kq_{i}\right)}{\partial \left(Kq_{i}\right)}=f\left(\frac{p_{i}}{Kq_{i}}\right)-\frac{p_{i}}{Kq_{i}}f'\left(\frac{p_{i}}{Kq_{i}}\right)
\end{equation}
And the equation (\ref{eq.resolK}) is written:
\begin{equation}
\sum_{i}q_{i}\left[f\left(\frac{p_{i}}{Kq_{i}}\right)-\frac{p_{i}}{Kq_{i}}f'\left(\frac{p_{i}}{Kq_{i}}\right)\right]=0
\label{eq.BaseC}
\end{equation}

\subsubsection{2 - Bregman's divergences.}
The invariant divergence associated with this constructive mode is by definition written:
\begin{equation}
	B_{f}\left(p\|Kq\right)=\sum_{i}b_{f}\left(p_{i}\|Kq_{i}\right)
\end{equation}
that is:
\begin{equation}	B_{f}\left(p\|Kq\right)=\sum_{i}\left[f\left(p_{i}\right)-f\left(Kq_{i}\right)-\left(p_{i}-Kq_{i}\right)f'\left(Kq_{i}\right)\right]
\end{equation}
we then obtain:
\begin{equation}
\frac{\partial b_{f}\left(p_{i}\|Kq_{i}\right)}{\partial \left(Kq_{i}\right)}=\left(Kq_{i}-p_{i}\right)f''\left(Kq_{i}\right)
\end{equation}
And the relationship (\ref{eq.resolK}) is written:
\begin{equation}
\sum_{i}q_{i}\left(Kq_{i}-p_{i}\right)f''\left(Kq_{i}\right)=0
\label{eq.BaseB}
\end{equation}

\subsubsection{3 - Jensen's divergences.}
By definition, we have:
\begin{equation}
	J_{f}\left(p\|Kq\right)=\sum_{i}j_{f}\left(p_{i}\|Kq_{i}\right)
\end{equation}
that is:
\begin{equation}
	J_{f}\left(p\|Kq\right)=\sum_{i}\left\{\alpha f\left(p_{i}\right)+\left(1-\alpha\right)f\left(Kq_{i}\right)-f\left[\alpha p_{i}+\left(1-\alpha\right)Kq_{i}\right]\right\}
\end{equation}
Therefore:
\begin{align}
\frac{\partial j_{f}\left(p_{i}\|Kq_{i}\right)}{\partial \left(Kq_{i}\right)} =
\left(1-\alpha\right)\left[f'\left(Kq_{i}\right)-f'\left(\alpha p_{i}+\left(1-\alpha\right)Kq_{i}\right)\right]
\end{align}
And the relationship (\ref{eq.resolK}) is written:
\begin{equation}
\sum_{i}q_{i}\left[f'\left(Kq_{i}\right)-f'\left(\alpha p_{i}+\left(1-\alpha\right)Kq_{i}\right)\right]=0
\label{eq.BaseJ}
\end{equation}

\subsubsection{4 - Overview}
If we make an analysis of the relationships corresponding to the 3 types of divergences considered, we can observe that the relation (\ref{eq.resolK}) will be satisfied for Csiszär divergences, if:
\begin{equation}
\sum_{i}q_{i}\left[f\left(\frac{p_{i}}{Kq_{i}}\right)-\frac{p_{i}}{Kq_{i}}f'\left(\frac{p_{i}}{Kq_{i}}\right)\right]=0
\label{eq.relC} 	
\end{equation}
Similarly, this equation will be satisfied for Bregman's divergences if:
\begin{equation}
\sum_{i}q_{i}\left(Kq_{i}-p_{i}\right)f''\left(Kq_{i}\right)=0
\label{eq.relB} 	
\end{equation}
Finally, this relationship will be satisfied for Jensen's divergences, if:
\begin{equation}
\sum_{i}q_{i}\left[f'\left(Kq_{i}\right)-f'\left(\alpha p_{i}+\left(1-\alpha\right)Kq_{i}\right)\right]=0
\label{eq.relJ} 	
\end{equation}
It is quite obvious that the equations (\ref{eq.relC}), (\ref{eq.relB}) and (\ref{eq.relJ}) are nothing else but the translation of the equation (\ref{eq.resolK}) corresponding to the 3 types of divergences.\\

The resolution of (\ref{eq.relC}), (\ref{eq.relB}) or (\ref{eq.relJ}) depending on the type of divergence and more generally of (\ref{eq.resolK}), when possible, allows us to obtain an expression of the nominal invariance factor $K_{0}\left(p,q\right)$ and thus to satisfy (\ref{eq.Pfond}).\\

\textbf{The question is then: is the relation (\ref{eq.Pfond}) still true if we use an invariance factor different from $K_{0}$?}
\subsection{Extension of the fundamental property to "non-nominal" invariance factors.}
The fundamental property (\ref{eq.Pfond}) has been expressed in the form:
\begin{equation}
	\left[\sum_{i}\frac{\partial d\left(p_{i}\|Kq_{i}\right)}{\partial \left(Kq_{i}\right)}q_{i}\right]\left[K+\sum_{j}q_{j}\frac{\partial K}{\partial q_{j}}\right]=0
	\label{eq.Base0bis}	
\end{equation}
Or also:
\begin{equation}
	\left[\sum_{i}\frac{\partial d\left(p_{i}\|Kq_{i}\right)}{\partial K}\right]\left[K+\sum_{j}q_{j}\frac{\partial K}{\partial q_{j}}\right]=0
	\label{eq.prodbis}
\end{equation}
The effect of the first product term was discussed in the previous section.\\

\textbf{The effect of the second term of the product is examined here.}\\

The relations (\ref{eq.Base0bis}) (\ref{eq.prodbis}) are always satisfied if the factor $K\left(p,q\right)$ is a solution to the differential equation:
\begin{equation}
K+\sum_{j}q_{j}\frac{\partial K}{\partial q_{j}}=0
\label{eq.DiffK}	
\end{equation}

This differential equation doesn't depend on the divergence under consideration.\\

\textbf{Its resolution allows to extend the fundamental property (\ref{eq.Pfond}) to invariance factors different from $K_{0}$, i.e. non-nominal.}

\subsection{Some precisions on the differential equation (\ref{eq.DiffK}).}
Clearly, only the "$q$" dependency is exhibited in this equation. That means that the set of solutions of this equation contains the $K\left(p,q\right)$ expressions which satisfy (\ref{eq.Pfond})(\ref{eq.prod})(\ref{eq.Base0}).\\
 In order for a solution  $K\left(p,q\right)$ of the differential equation (\ref {eq.DiffK}) to be an acceptable invariance factor, it must also possess other properties that have already been mentioned:\\
 
  \textbf{* $K\left(p,q\right)$ must be a positive scalar},\\
	
  on the other hand, and, this is a crucial point already indicated:\\
   
 \textbf{* In order to obtain an invariant divergence by scale-change on ``$q$", the vector ``$\left[K\left(p,q\right).\ q\right]$" must be invariant when multiplying ``$q$" by a constant.}\\
 
Finally, in a general way:\\

 \textbf{* If $p\rightarrow q$, then we must have $K\left(p,q\right)\rightarrow1$.}\\

However, the general solution of (\ref{eq.DiffK}) mentioned in \cite{Polyanin:2001ul} (p.94), is written $\forall j$:
\begin{equation}
	K\left(p,q\right)=\frac{1}{q_{j}}\Phi\left(p,\frac{q_{1}}{q_{j}},\frac{q_{2}}{q_{j}},...,\frac{q_{j-1}}{q_{j}},1,\frac{q_{j+1}}{q_{j}},...,\frac{q_{n}}{q_{j}}\right)
	\label{eq.exprK}
\end{equation}
Where $\Phi$ is any function.\\
This means that any function $K\left(p,q\right)$ of the form (\ref{eq.exprK}) will satisfy (\ref{eq.prod})(\ref{eq.Base0}), i.e. (\ref{eq.Pfond}).\\
We will observe that with expressions of $K\left(p,q\right)$ of this form, the vector $\left[K\left(p,q\right).\ q\right]$ is invariant when multiplying ``$q$" by a constant. However, among the solutions of the form (\ref{eq.exprK}), only the positive scalar expressions of $K\left(p,q\right)$ can play the role of invariance factors.\\
	
Finally, one will be able to check that all the expressions of $K_{0}\left(p,q\right)$ calculated (when possible), by explicit resolution of (\ref{eq.resolK}) that is to say, according to the case, of (\ref{eq.relC}), (\ref{eq.relB}) where (\ref{eq.relJ}) will be of the form (\ref{eq.exprK}).\\

These observations induce the following remarkable property:\\

\textbf{An expression of  $K\left(p,q\right)$, positive scalar, as long as it is a solution of the differential equation (\ref{eq.DiffK}), (i.e. as long as it is of the form (\ref{eq.exprK})), will make invariant any divergence (because the vector $\left[K\left(p,q\right).\ q\right]$ is invariant with respect to ``$q$"), and the relation (\ref{eq.Pfond}) will be satisfied.}\\

\textbf{We can therefore say, in order to globalize these observations concerning the property (\ref{eq.Pfond}), that two cases can arise:}\\

\textbf{- either the invariance factor is computed by explicit resolution of (\ref{eq.resolK}), the invariance factor is then referred to as $K_{0}\left(p,q\right)$, this will be the ``nominal" invariance factor for the considered divergence.}\\

\textbf{- or the factor $K\left(p,q\right)$ is solution of (\ref{eq.DiffK}) and has the specific properties of the invariance factors, without being solution of (\ref{eq.resolK}), i.e. without any relation with the starting divergence, then, one always obtains an invariant divergence as shown on the examples of the following paragraph.}\\

\textbf{- the nominal invariance factors belong to the set of solutions of the differential equation (\ref{eq.DiffK}).}\\

\textbf{Exemple 1: Kullback-Leibler divergence.}\\

We can consider that this divergence is constructed in the sense of Csiszär on the standard convex function:
\begin{equation}
	f_{c}\left(x\right)=x\log x+1-x
\end{equation}
It is written:
\begin{equation}
	KL\left(p\|q\right)=\sum_{i}\left[p_{i}\log\frac{p_{i}}{q_{i}}+q_{i}-p_{i}\right]
\end{equation}
 The calculation of the nominal invariance factor leads to the explicit solution: 
\begin{equation}
K_{0}\left(p,q\right)=\frac{\sum_{j}p_{j}}{\sum_{j}q_{j}}
\label{eq.Kpart}	
\end{equation}
This factor will be the solution of the equation (\ref{eq.relC}), but this divergence can also be obtained in the Bregman sense by relying on the same convex function (it is the common point between the 2 types of divergences), it will thus be made invariant by the same factor $K_{0}$ which is the solution of the equation (\ref{eq.relB}).\ \
The expression of $K_{0}\left(p,q\right)$ (\ref{eq.Kpart}) can be put in the form (\ref{eq.exprK}); it is thus the solution of the differential equation (\ref{eq.DiffK}) and thus makes it possible to make any divergence invariant.\\
This can be verified, for example, on the Euclidean distance:
\begin{equation}
	EQM\left(p\|q\right)=\sum_{i}\left(p_{i}-q_{i}\right)^{2}
\end{equation}
With the invariance factor (\ref{eq.Kpart}) \textbf{\textsl{``which is not the nominal invariance factor for this divergence"}}, we obtain:
\begin{equation}
	EQM\left(p\|K_{0}q\right)=\sum_{i}\left(p_{i}-\frac{\sum_{j}p_{j}}{\sum_{j}q_{j}}q_{i}\right)^{2}
\end{equation}
This divergence is invariant under scale change on ``$q$ "; it can also be written as follows:
\begin{equation}	
EQM\left(p\|K_{0}q\right)=\left(\sum_{j}p_{j}\right)^{2}\sum_{i}\left(\frac{p_{i}}{\sum_{j}p_{j}}-\frac{q_{i}}{\sum_{j}q_{j}}\right)^{2}
\end{equation}
or also:
\begin{equation}	
EQM\left(p\|K_{0}q\right)=\left(\sum_{j}p_{j}\right)^{2}\sum_{i}\left(\bar{p_{i}}-\bar{q_{i}}\right)^{2}
\end{equation}
When using this very particular form of the invariance factor $K$ given by (\ref{eq.Kpart}), the resulting invariant divergence is analogous (except for one factor that depends only on ``$p_{j}$") to the initial divergence, where ``$p_{i}$" has been replaced with $\bar{p}_{i}=\frac{p_{i}}{\sum_{j}p_{j}}$ and ``$q_{i}$" has been replaced with $\bar{q}_{i}=\frac{q_{i}}{\sum_{j}q_{j}}$.\\
We can check that property on any discrepancies that come up later.\\

\textbf{Exemple 2: Mean square deviation (E.Q.M.).}\\

We consider here, the Mean square deviation (Euclidean distance):
\begin{equation}
	EQM\left(p\|q\right)=\sum_{i}\left(p_{i}-q_{i}\right)^{2}
\end{equation} 
This can be considered as being a Bregman divergence based on the standard convex function:
\begin{equation}
	f_{c}\left(x\right)=\left(x-1\right)^{2}
\end{equation}
 It is made invariant to a scale change on ``$q$" by calculating the nominal factor that is written: 
\begin{equation}
	K_{0}\left(p,q\right)=\frac{\sum_{i}p_{i}q_{i}}{\sum_{i}q^{2}_{i}}
	\label{eq.KEQM}
\end{equation}
But the E.Q.M. is also a Jensen (1/2) divergence based on the same convex function.\\
With this expression of $K\left(p,q\right)$, the equation (\ref{eq.relB}) will be satisfied because the Euclidean distance is a Bregman divergence, but simultaneously the equation (\ref{eq.relJ}) will be satisfied because it is also a Jensen divergence, and of course, this expression of $K\left(p,q\right)$ is a solution of the differential equation (\ref{eq.DiffK}).\\
If we now use the expression of $K_{0}\left(p,q\right)$ given in (\ref{eq.KEQM}) as the invariance factor in the Kullback-Leibler divergence, we can see that, although this expression is not the nominal invariance factor for this divergence, we obtain an invariant form by scale change on ``$q$".

\subsection{General form of the invariance factor $K\left(p,q\right)$.}
Given these observations, and all the constraints on the invariance factor $K\left(p,q\right)$, a general expression acceptable for $K$ can be written in the form:
\begin{equation}	
K\left(p,q\right)=\left(\frac{\sum_{i}p^{\alpha}_{i}q^{\beta}_{i}}{\sum_{i}p^{\delta}_{i}q^{\gamma}_{i}}\right)^{\mu}\ \ \ \ ;\ \ \ \alpha+\beta=\delta+\gamma \ \ \ \ ; \ \ \ \ \mu\left(\gamma-\beta\right)=1
\label{eq.Kgene}
\end{equation}
Indeed, this expression which can be put in the form (\ref{eq.exprK}) represents a family of solutions to the differential equation (\ref{eq.DiffK}).\\
Note that the constraints expressed in (\ref{eq.Kgene}) reflect the following properties:\\
\begin{equation}
	\alpha+\beta=\delta+\gamma\ \ \ \ \Leftrightarrow\ \ \ K\left(p,q\right)\ is\ a\ scalar
\end{equation}
\begin{equation}
\mu\left(\gamma-\beta\right)=1\ \ \ \Leftrightarrow\ \ \ \left[K\left(p,q\right).\;q\right]\ \ \ is\ invariant\ with\ respect\ to\ q	
\end{equation}
Taking into account the constraints mentioned in (\ref{eq.Kgene}), this general expression of $K$ can also be written with a smaller number of parameters, in the form: 
\begin{equation}	
K\left(p,q\right)=\left(\frac{\sum_{i}p^{\alpha}_{i}q^{\beta}_{i}}{\sum_{i}p^{\alpha-\frac{1}{\mu}}_{i}q^{\beta+\frac{1}{\mu}}_{i}}\right)^{\mu}
\label{eq.Kgene2}
\end{equation}
But, another expression of $K\left(p,q\right)$ is more explicit; indeed, taking into account the relations between the parameters, it can be written, after some calculations, in the form:
\begin{equation}
K\left(p,q\right)=\left[\sum_{i}\left(\frac{p_{i}}{q_{i}}\right)^{\alpha-\delta}\frac{p^{\delta}_{i}q^{\gamma}_{i}}{\sum_{j}p^{\delta}_{j}q^{\gamma}_{j}}\right]^{\frac{1}{\alpha-\delta}}
\label{eq.Kgene4}
\end{equation}\\
These expressions are of the form:
\begin{equation}
K\left(p,q\right)=\left[\sum_{i}\left(\frac{p_{i}}{q_{i}}\right)^{t}\frac{p^{\delta}_{i}q^{\gamma}_{i}}{\sum_{j}p^{\delta}_{j}q^{\gamma}_{j}}\right]^{\frac{1}{t}}
\label{eq.Kgene5}
\end{equation}\\
These are, for quantities of the form $\left(p_{i}/q_{i}\right)$, a weighted generalized mean of the order ``$t$" with the exponent $t=\alpha-\delta$ (see.Appendix 2), and weighting coefficients $w_{i}$ such that $\sum_{i}w_{i}=1$ which are given by: 
\begin{equation}
	w_{i}=\frac{p^{\delta}_{i}q^{\gamma}_{i}}{\sum_{j}p^{\delta}_{j}q^{\gamma}_{j}}
	\label{eq.wpond}
\end{equation}
In the special case of the dual Kullback -Leibler divergence, the nominal invariance factor will be expressed as follows:
\begin{equation}
	K_{0}=\exp\left[\sum_{i}w_{i}\log\left(\frac{p_{i}}{q_{i}}\right)\right]
	\label{eq.KKLDuale}
\end{equation}
It is a weighted generalized mean of $p_{i}/q_{i}$, based on the function $\psi\left(x\right)=\log x$, with weighting factors $\omega_{i}=\frac{p_{i}}{\sum_{j}p_{j}}$.
\begin{equation}
	K_{0}=\psi^{-1}\left[\sum_{i}w_{i}\psi\left(\frac{p_{i}}{q_{i}}\right)\right]
\end{equation}

It can be verified after some calculations that expressions of the invariance factor given by (\ref{eq.Kgene5}) and (\ref{eq.KKLDuale}), associated with weighting coefficients of the form (\ref{eq.wpond}), satisfy the differential equation (\ref{eq.DiffK}).\\

\subsection{Remarks.}
Based on this general form, we can go even further in the discussion; by returning to the expression $K\left(p,q\right)=\frac{\sum_{j}p_{j}q_{j}}{\sum_{j}q^{2}_{j}}$ calculated explicitly for K.L.'s divergence, we can imagine an expression for the invariance factor written $K\left(p,q\right)=\left[\frac{\sum_{j}p^{2}_{j}}{\sum_{j}q^{2}_{j}}\right]^{\frac{1}{2}}$ which is of the form (\ref{eq.Kgene}), but which does not correspond "\textit{a priori}" to any divergence. \\
This expression of $K\left(p,q\right)$ is a solution of the differential equation, so it makes invariant any divergence.\\
Similarly, from the expression of the invariance factor corresponding to the mean square deviation: $K\left(p,q\right)=\frac{\sum_{j}p_{j}q_{j}}{\sum_{j}q^{2}_{j}}$, we can imagine an expression of the invariance factor $K\left(p,q\right)=\frac{\sum_{j}p_{j}q^{2}_{j}}{\sum_{j}q^{3}_{j}}$ which is of the form (\ref{eq.Kgene}), but which does not correspond ``\textit{a priori}" to any divergence; this expression of $K\left(p,q\right)$ is the solution of the differential equation (\ref{eq.DiffK}); this factor makes invariant any divergence.\\

\section{Some properties of the Gradient of an invariant divergence.}
We recall the relation (\ref{eq.Pfond0}) giving the expression of the gradient with respect to the variable ``$q$" of an invariant divergence $D\left(p\|Kq\right)$ for an invariance factor $K\left(p,q\right)$:
\begin{equation}
\frac{\partial D\left(p\|Kq\right)}{\partial q_{j}}=\sum_{i}\frac{\partial d\left(p_{i}\|Kq_{i}\right)}{\partial \left(Kq_{i}\right)}\frac{\partial K}{\partial q_{j}}q_{i}+K	\frac{\partial d\left(p_{j}\|Kq_{j}\right)}{\partial \left(Kq_{j}\right)}
\label{eq.GradDpKQ1}
\end{equation}
This expression can also be written:
\begin{equation}
\frac{\partial D\left(p\|Kq\right)}{\partial q_{j}}=\frac{\partial K}{\partial q_{j}}\sum_{i}q_{i}\frac{\partial d\left(p_{i}\|Kq_{i}\right)}{\partial \left(Kq_{i}\right)}+K	\frac{\partial d\left(p_{j}\|Kq_{j}\right)}{\partial \left(Kq_{j}\right)}
\label{eq.GradDpKQ2}
\end{equation}
Or also:
\begin{equation}
\frac{\partial D\left(p\|Kq\right)}{\partial q_{j}}=\frac{\partial K}{\partial q_{j}}\sum_{i}\frac{\partial d\left(p_{i}\|Kq_{i}\right)}{\partial K}+K	\frac{\partial d\left(p_{j}\|Kq_{j}\right)}{\partial \left(Kq_{j}\right)}
\label{eq.GradDpKQ3}
\end{equation}
For a given divergence, suppose that the expression of the invariance factor can be calculated explicitly, namely $K_{0}\left(p,q\right)$ this expression; it is of course of the general form (\ref{eq.Kgene2}).\\
Let us now consider another expression of the invariance factor $K_{1}\left(p,q\right)$ respecting the general form (\ref{eq.Kgene2}), but which is not in correspondence with the divergence considered.\\
 
The question is: what happens to the gradient expressions (\ref{eq.GradDpKQ1}), (\ref{eq.GradDpKQ2}) or (\ref{eq.GradDpKQ3}) depending on whether we use $K_{0}\left(p,q\right)$ or $K_{1}\left(p,q\right)$? \\
In the first case, $K_{0}\left(p,q\right)$ is the solution to the differential equation (\ref{eq.DiffK}), then the first term of the second member of the equation (\ref{eq.GradDpKQ3}) is zero and it simply remains:
\begin{equation}
\frac{\partial D\left(p\|Kq\right)}{\partial q_{j}}=K_{0}	\left[\frac{\partial d\left(p_{j}\|Kq_{j}\right)}{\partial \left(Kq_{j}\right)}\right]_{K=K_{0}}
\label{eq.GradDpKQ3S}
\end{equation}
If $K=K_{1}$, the first term of the second member of the equation (\ref{eq.GradDpKQ3}) is no longer zero, and we have:
\begin{equation}
\frac{\partial D\left(p\|Kq\right)}{\partial q_{j}}=\frac{\partial K_{1}}{\partial q_{j}}\left[\sum_{i}\frac{\partial d\left(p_{i}\|Kq_{i}\right)}{\partial K}\right]_{K=K_{1}}+K_{1}\left[\frac{\partial d\left(p_{j}\|Kq_{j}\right)}{\partial\left(Kq_{j}\right)}\right]_{K=K_{1}}
\label{eq.GradDpKQ3C}
\end{equation}
However, in both cases, one still has the fundamental property:
\begin{equation}
	\sum_{j}q_{j}\frac{\partial D\left(p\|Kq\right)}{\partial q_{j}}=\sum_{j}q_{j}\frac{\partial DI\left(p\|q\right)}{\partial q_{j}}=0
\end{equation}

\textbf{Some examples to show this in more detail are given in Appendix 4.}

\section{A very special case of the invariance factor.}
We saw that a particular expression of the invariance factor could be:
\begin{equation}
	K\left(p,q\right)=\frac{\sum_{j}p_{j}}{\sum_{j}q_{j}}
\end{equation}
This expression is the explicit result of calculating the invariance factor for a Kullback-Leibler divergence.\\
Using this factor, any divergence can be made invariant.\\

\textbf{Forms of the divergences related to this invariance factor.}\\

By introducing this invariance factor with respect to ``$q$", in any divergence, one makes it appear in the resulting invariant divergence, 
 normalized variables noted $\bar{p_{i}}=\frac{p_{i}}{\sum_{j}p_{j}}$ and $\bar{q_{i}}=\frac{q_{i}}{\sum_{j}q_{j}}$, so that we systematically obtain invariant divergences corresponding to $\bar{p}$ and $\bar{q}$ variables of sum equal to $1$.\\ 
 
Except for a multiplicative factor (which depends only on $\sum_{j}p_{j}$), the resulting invariant divergences have the same expression as the initial divergences, with the normalized variables simply replacing the initial variables.
This operation being carried out, some simplifications can appear.\\
 One thus obtains invariant divergences similar to the simplified divergences applicable to densities of probabilities, provided that one introduces explicitly in the simplified divergences, normalized variables.\\
 
\textbf{Moreover, if we disregard the multiplicative factor, these divergences are invariant not only with respect to the ``$q$" variable, but also with respect to the ``$p$" variable, as can be seen from the examples given in Appendix 6.}\\

\setcounter{table}{0}  \setcounter{equation}{0}  \setcounter{figure}{0} \setcounter{chapter}{4} \setcounter{section}{0} 
\chapter{chapter 4 -\\ Divergences and Entropies}  \label{chptr::chapitre4}

This chapter presents the divergences that can be related to the various forms of entropy found in the literature. In most cases they will be constructed in the sense of Csiszär on the basis of a convex function; in this constructive mode, we will be led to distinguish the case of \textit{``standard convex functions"}. The relation with the entropies will of course be rather related to the \textit{``simple convex functions"}.\\
 Finally, the extensions of the divergences by using the function ``Generalized Logarithm" (see Appendix 1), will lead us to deviate from the constructive mode of Csiszär, but will make it possible to make the connection with the entropies of Sharma-Mittal \cite{sharma1975} and Renyi \cite{renyi1961} \cite{renyi1955}.\\
  In a final section, we'll discuss Jensen's differences based on Entropies.\\
General references for this chapter will be given in: \cite{arndt2001},\cite{basseville1989},\cite{basseville1996},\cite{pardo2005} et \cite{taneja2001}.

\section{Shannon Entropy related divergences.}
In this section, the divergences linked to Shannon's entropy \cite{shannon1948} are analyzed.
\subsection{Direct form.}
\begin{figure}[h!]
\centering
\includegraphics[width=0.7\linewidth]{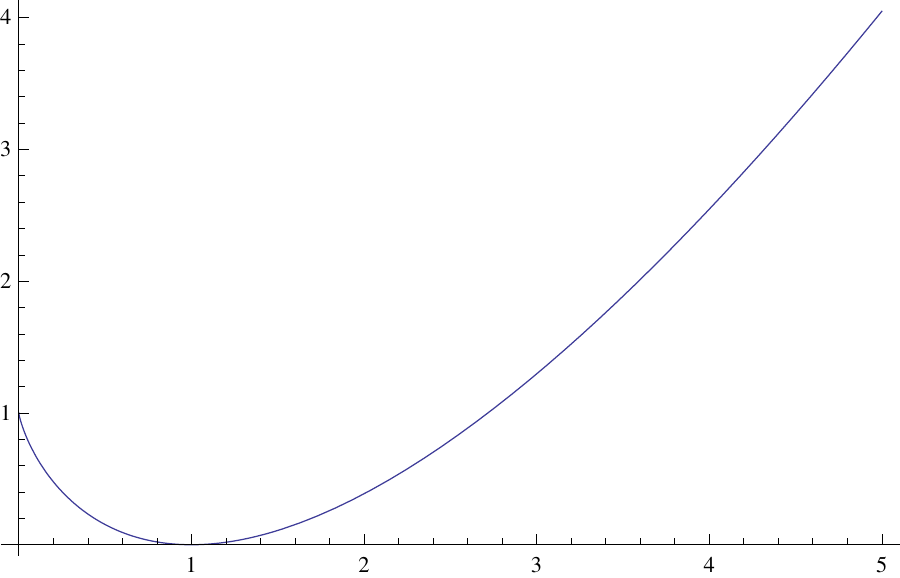}
\caption{Function $f_{c}\left(x\right)=x\log x+1-x$}
\label{fig:SHfc}
\end{figure}

The standard convex function used in this case, see.figure (\ref{fig:SHfc}), is written:
\begin{equation}
	f_{c}\left(x\right)=x\log x+1-x
	\label{eq.fcKL}
\end{equation}

\begin{figure}[h!]
\centering
\includegraphics[width=0.7\linewidth]{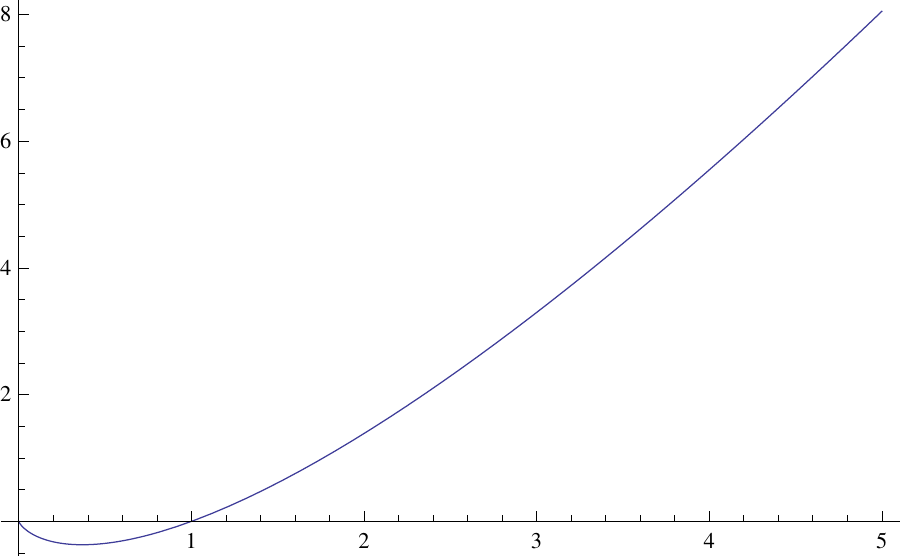}
\caption{Function $f\left(x\right)=x\log x$}
\label{fig:SHf}
\end{figure}

The corresponding Csiszär's divergence is the Kullback-Leibler  divergence \cite{kullback1951}:
\begin{equation}
	KL\left(p\|q\right)=\sum_{i}p_{i}\log\frac{p_{i}}{q_{i}}+q_{i}-p_{i}
	\label{eq.DKL}
\end{equation}
The gradient with respect to ``$q$" is written:
\begin{equation}
	\frac{\partial KL\left(p\|q\right)}{\partial q_{j}}=-\frac{p_{j}}{q_{j}}+1
\end{equation}
It will be zero for $p_{j}=q_{j}\ \forall j$.\\
If, at this point, we consider that we have $\sum_{i}p_{i}=\sum_{i}q_{i}$, the divergence (\ref{eq.DKL}) simplifies and is written:
\begin{equation}
	IKL\left(p\|q\right)=\sum_{i}p_{i}\log\frac{p_{i}}{q_{i}}
	\label{eq.IKL}
\end{equation}
That's Kullback's information \cite{basseville1989}.\\
The gradient with respect to ``$q$" is written:
\begin{equation}
	\frac{\partial IKL\left(p\|q\right)}{\partial q_{j}}=-\frac{p_{j}}{q_{j}}
\end{equation}
We can see that this gradient will never be equal to zero, therefore $ IKL\left(p\|q\right)$ is not usable in our problem without special precautions; indeed, to use such a divergence, it will be necessary to explicitly introduce the fact that $\sum_{i}p_{i}=\sum_{i}q_{i}$.\\
This difficulty will appear every time one wants to use simplified divergences, that is to say, those built on ``simple" convex functions.\\
An additional simplification such as $\sum_{i}p_{i}=\sum_{i}q_{i}=1$ will not bring any additional simplification, but if we don't make this assumption, the divergence (\ref{eq.IKL}) will probably not be positive.\\
The form (\ref{eq.IKL}) is mainly used in works dealing with probability densities; this form is deduced in the Csiszär sense from the simple convex function (see figure \ref{fig:SHf}):
\begin{equation}
	f\left(x\right)=x\log x
\end{equation}

\subsubsection{Invariance by change of scale.}
The invariance factor $K_{0}\left(p,q\right)$ corresponding to the Kullback-Leibler divergence (\ref{eq.DKL}) can be calculated according to the method indicated in Chapter 3; it is expressed in explicit form:
\begin{equation}
	K_{0}=\frac{\sum_{j}p_{j}}{\sum_{j}q_{j}}
\end{equation}
It is a special case already mentioned in Chapter 3 that leads to the invariant divergence:
\begin{equation}
	KLI\left(p\|q\right)=\sum_{j}p_{j}\sum_{i}\bar{p}_{i}\log\frac{\bar{p}_{i}}{\bar{q}_{i}}
	\label{eq.KLI}
\end{equation}
with $\bar{p}_{j}=\frac{p_{j}}{\sum_{l}p_{l}}$ et $\bar{q}_{j}=\frac{q_{j}}{\sum_{l}q_{l}}$.\\
We can observe that we get a divergence equivalent to (\ref{eq.IKL}), but here we explicitly have $\sum_{i}\bar{p_{i}}=\sum_{i}\bar{q_{i}}\left(=1\right)$.\\
On the other hand, disregarding the multiplicative factor $\sum_{j}p_{j}$, we can observe that the divergence obtained is invariant not only with respect to ``$q$" but also with respect ``$p$".\\
Its gradient with respect to ``$q$" is written:
\begin{equation}
\frac{\partial KLI\left(p\|q\right)}{\partial q_{l}}=\frac{1}{\sum_{j}q_{j}}\left(1-\frac{\bar{p}_{l}}{\bar{q}_{l}}\right)
\label{eq.gradKLI}	
\end{equation}

\subsection{Dual form.}
\begin{figure}[h!]
\centering
\includegraphics[width=0.7\linewidth]{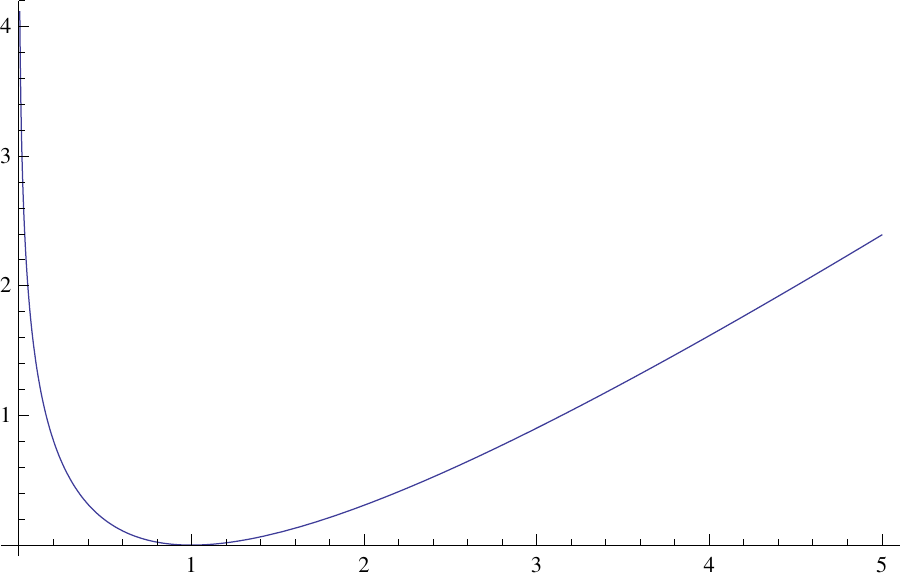}
\caption{Function $\breve{f}_{c}\left(x\right)=\log\frac{1}{x}+x-1$}
\label{fig:SHfct}
\end{figure}

To obtain the dual form, we are using the mirror function of (\ref{eq.fcKL}) (see figure \ref{fig:SHfct}), which is written:
\begin{equation}
	\breve{f}_{c}\left(x\right)=\log\frac{1}{x}+x-1
\end{equation}

\begin{figure}[h!]
\centering
\includegraphics[width=0.7\linewidth]{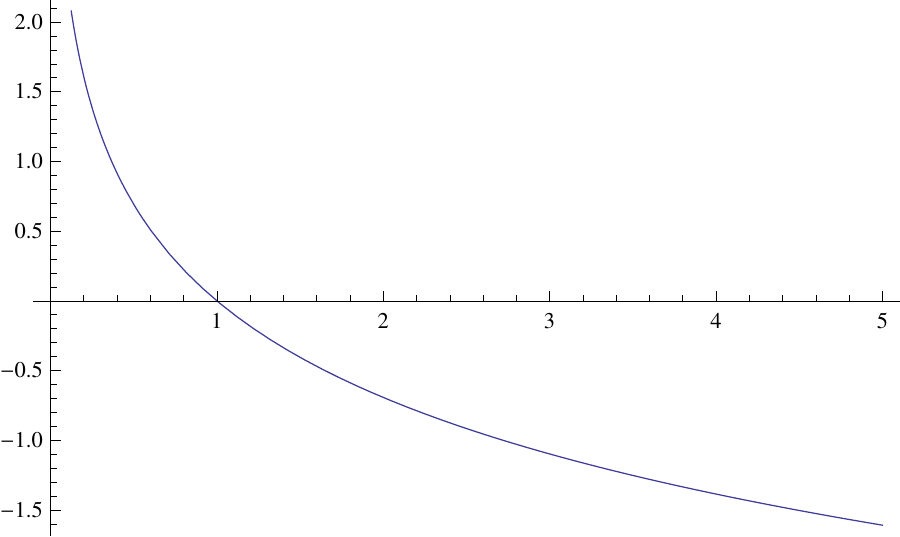}
\caption{Function $\breve{f}\left(x\right)=\log\frac{1}{x}$}
\label{fig:SHft}
\end{figure}

So, we obtain:
\begin{equation}
	KL\left(q\|p\right)=\sum_{i}q_{i}\log\frac{q_{i}}{p_{i}}+p_{i}-q_{i}
\label{eq.KLduale}
\end{equation}
The gradient with respect to ``$q$" is written:
\begin{equation}
	\frac{\partial KL\left(q\|p\right)}{\partial q_{j}}=\log\frac{q_{j}}{p_{j}}
\end{equation}
It will be zero for $p_{j}=q_{j}\ \forall j$.\\
Again, assuming that $\sum_{i}p_{i}=\sum_{i}q_{i}$, we have the simplified form:
\begin{equation}
	IKL\left(q\|p\right)=\sum_{i}q_{i}\log\frac{q_{i}}{p_{i}}
	\label{eq.IKLduale}
\end{equation}
The gradient with respect to ``$q$" is written:
\begin{equation}
	\frac{\partial IKL\left(q\|p\right)}{\partial q_{j}}=\log\frac{q_{j}}{p_{j}}+1
\end{equation}
This gradient will not be zero for $p_{j}=q_{j}$, therefore $IKL\left(q\|p\right)$ is not usable in our problem without special precautions.\\
An additional specification such as $\sum_{i}p_{i}=\sum_{i}q_{i}=1$ will not provide any additional simplification, but if we don't make that assumption, the divergence (\ref{eq.IKLduale}) is unlikely to be positive.\\
The divergence (\ref{eq.IKLduale}) is constructed in the Csiszär sense on the simple convex function (see figure (\ref{fig:SHft}):
\begin{equation}
	\breve{f}\left(x\right)=\log\frac{1}{x}
\end{equation}

\subsubsection{Invariant form of dual divergence (\ref{eq.KLduale}).}
The derivation of the nominal invariance factor for this divergence leads to the expression:
\begin{equation}
	K_{0}=\exp \sum_{i}\frac{q_{i}}{\sum_{j}q_{j}}\log \frac{p_{i}}{q_{i}}
\end{equation}
It is a weighted generalized average of terms of the form ($p_{i}/q_{i}$) that is written:
\begin{equation}
	K_{0}=\psi^{-1}\left[\sum_{i}w_{i}\psi\left(\frac{p_{i}}{q_{i}}\right)\right]
\end{equation}
With $\psi\left(x\right)=\ln \left(x\right)$ and weighting factors $w_{i}=\frac{q_{i}}{\sum_{j}q_{j}}$.\\
By introducing this invariance factor in the divergence (\ref{eq.KLduale}), we obtain after some simple calculations, the invariant form that can be written:
\begin{equation}
	KLI\left(q\|p\right)=\sum_{i}\left[p_{i}-q_{i}\exp \sum_{j}\frac{q_{j}}{\sum_{l}q_{l}}\log \frac{p_{j}}{q_{j}}\right]=\sum_{i}\left[p_{i}-K_{0}q_{i}\right]
\end{equation}
It is easily verified that this divergence is not modified if ``$q$" is multiplied by a positive constant.

\subsection{Symmetrical form.}
\begin{figure}[h!]
\centering
\includegraphics[width=0.7\linewidth]{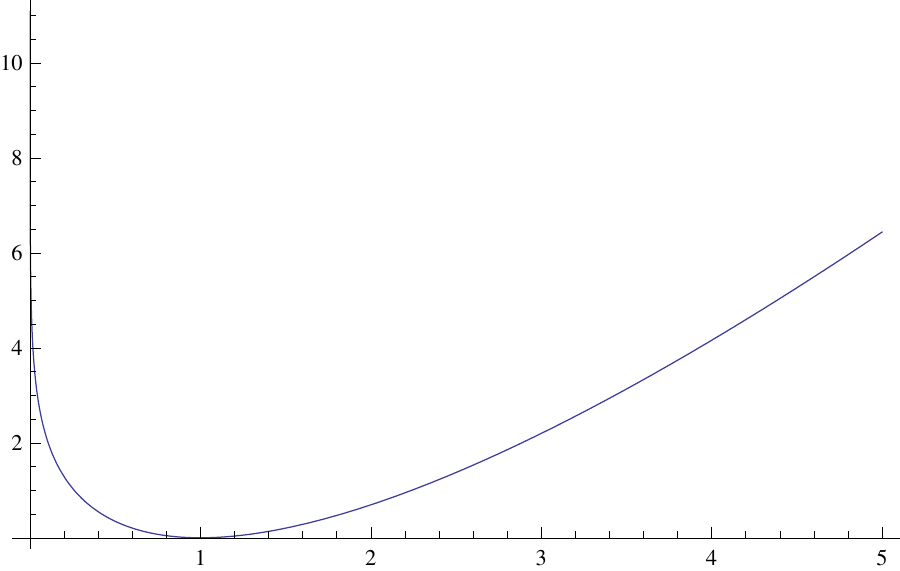}
\caption{Function $\hat{f}_{c}\left(x\right)=x \log x-\log x$}
\label{fig:SHsym}
\end{figure}

It is obtained in the sense of Jeffreys \cite{jeffreys1946} by using Csiszär's constructive mode, based on the standard convex function (see figure \ref{fig:SHsym}):
\begin{equation}	\hat{f}_{c}\left(x\right)=\frac{1}{2}\left[f_{c}\left(x\right)+\breve{f}_{c}\left(x\right)\right]=\frac{1}{2}\left[f\left(x\right)+\breve{f}\left(x\right)\right]
\end{equation}

So, we have:
\begin{equation}
	KL\left(p,q\right)=\sum_{i}\left(p_{i}-q_{i}\right)\log\frac{p_{i}}{q_{i}}
\end{equation}
Whose gradient with respect to ``$q$" is written: 
\begin{equation}
	\frac{\partial KL\left(p,q\right)}{\partial q_{j}}=\frac{q_{j}-p_{j}}{q_{j}}-\log\frac{p_{j}}{q_{j}}
\end{equation}
This expression is of course equal to zero for $p_{j}=q_{j}\ \forall j$.

\section{Havrda-Charvat Entropy related divergences.} 
This section deals with Csiszär's divergences related to the entropy of Havrda-Charvat \cite{havrda1967} and more specifically based on the standard convex function:
\begin{equation}
	f_{c}\left(x\right)=\frac{1}{\alpha\left(\alpha-1\right)}\left[x^{\alpha}-\alpha x+\left(1-\alpha\right)\right]
	\label{eq.fcHC}
\end{equation}
To these divergences we will associate those built on the corresponding simple convex functions:
\begin{equation}
	f_{1}\left(x\right)=\frac{1}{\alpha\left(\alpha-1\right)}\left[x^{\alpha}- x\right]
\end{equation}
and
\begin{equation}
		f_{2}\left(x\right)=\frac{1}{\alpha\left(\alpha-1\right)}\left[x^{\alpha}- 1\right]
\end{equation}
We will see in the \textbf{chapter 5} that the Csiszär divergences built on these convex functions, in particular on the function (\ref{eq.fcHC}), belong to the class of ``Alpha divergences" of Amari \cite{amari2009}.\\
The divergence associated with $f_{c}\left(x\right)$ will be written:
\begin{equation}	
HC_{\alpha}\left(p\|q\right)=\frac{1}{\alpha\left(\alpha-1\right)}\left\{\sum_{i}p_{i}^{\alpha}q_{i}^{1-\alpha}-\sum_{i}\left[\alpha p_{i}+\left(1-\alpha\right)q_{i}\right]\right\}
\label{eq.HCa}
\end{equation}
It's clearly a difference between a generalized geometric mean and a generalized arithmetic mean.\\
The gradient with respect to ``$q$" is given by:
\begin{equation}
\frac{\partial HC_{\alpha}\left(p\|q\right)}{\partial q_{j}}=\frac{1}{\alpha}\left(1-p_{j}^{\alpha}q_{j}^{-\alpha}\right)	
\end{equation}
It will be zero for $p_{j}=q_{j}\ \forall j$.\\
The divergences related to $f_{1}\left(x\right)$ and $f_{2}\left(x\right)$ will be written respectively:
\begin{equation}	
A_{1}\left(p\|q\right)=\frac{1}{\alpha\left(\alpha-1\right)}\left[\sum_{i}p_{i}^{\alpha}q_{i}^{1-\alpha}-p_{i}\right]
\label{eq.A1}
\end{equation}
and
\begin{equation}	A_{2}\left(p\|q\right)=\frac{1}{\alpha\left(\alpha-1\right)}\left[\sum_{i}p_{i}^{\alpha}q_{i}^{1-\alpha}-q_{i}\right]
\label{eq.A2}
\end{equation}
We can notice that if in the expression (\ref{eq.HCa}), one makes $\sum_{i}p_{i}=\sum_{i}q_{i}$, one finds (\ref{eq.A1}) or (\ref{eq.A2}).\\
The corresponding gradients with respect to ``$q$" will be written respectively:
\begin{equation}
	\frac{\partial A_{1}\left(p\|q\right)}{\partial q_{j}}=\frac{1}{\alpha}\left(p_{j}^{\alpha}q_{j}^{-\alpha}\right)
	\label{eq.ga1}
\end{equation}
that will never be zero, and:
\begin{equation}
\frac{\partial A_{2}\left(p\|q\right)}{\partial q_{j}}=\frac{1}{\alpha\left(\alpha-1\right)}\left[\left(1-\alpha\right)p_{j}^{\alpha}q_{j}^{-\alpha}-1\right]	
\end{equation}
which may be zero, but for $p_{j}\neq q_{j}$.\\
Furthermore, if we make the additional assumption that $\sum_{i}p_{i}=\sum_{i}q_{i}=1$, we get 3 identical divergences that are written:
 \begin{equation}	 HCS_{\alpha}\left(p\|q\right)=\frac{1}{\alpha\left(\alpha-1\right)}\left[\sum_{i}p_{i}^{\alpha}q_{i}^{1-\alpha}-1\right]
 \label{eq.HCSa}
\end{equation}
whose the gradient relative to ``$q$" given by (\ref{eq.ga1}).
This last divergence (\ref{eq.HCSa}) is the Havrda-Charvat divergence \cite{havrda1967} quoted by Arndt \cite{arndt2001} and Basseville \cite{basseville1996}.\\
The dual divergences are constructed by using as convex functions, the mirror functions of the previous ones; they will be indicated in \textbf{chapter 5}.

\section{Sharma-Mittal Entropy related divergences.}

The various divergences highlighted in this section are no longer Csiszär's divergences, in fact, they cannot be obtained with this constructive method, but they can be considered as a form of generalisation of the divergences obtained in the previous section. They are related to the entropy of Sharma-Mittal \cite{sharma1975}.\\
This generalization can be summarized by the following simple rule.\\

\textbf{Rule:} When a divergence is expressed as the difference between two positive terms, it can be generalized by applying to each of the terms an increasing function (which will not change the sign of the divergence obtained); the increasing function used is often the ``Generalized Logarithm" function (see Appendix 1), which makes it possible to change, by action on a single parameter, from the linear function (which leaves the initial divergence unchanged) to the logarithmic function itself.\\
Of course, at the end of this operation, we obtain a new divergence whose convexity properties are not guaranteed, even if the initial divergence was convex.\\
If we apply this rule on the divergence $HC_{\alpha}\left(p\|q\right)$ (\ref{eq.HCa}), using the generalized logarithm with the exponent $1-d=\frac{s-1}{\alpha-1}$, we obtain:
\begin{equation}	SM_{\alpha,s}\left(p\|q\right)=\frac{1}{\alpha\left(s-1\right)}\left\{\left[\sum_{i}p_{i}^{\alpha}q_{i}^{1-\alpha}\right]^{\frac{s-1}{\alpha-1}}-\left[\sum_{i}\alpha p_{i}+\left(1-\alpha\right)q_{i}\right]^{\frac{s-1}{\alpha-1}}\right\}
\label{eq.SMas}
\end{equation}
Similarly, if we apply this rule on the very simplified divergence (\ref{eq.HCSa}), we get:
\begin{equation}	SMS_{\alpha,s}\left(p\|q\right)=\frac{1}{\alpha\left(s-1\right)}\left\{\left[\sum_{i}p_{i}^{\alpha}q_{i}^{1-\alpha}\right]^{\frac{s-1}{\alpha-1}}-1\right\}
\label{eq.SMSa}
\end{equation}
It's this last divergence that is commonly referred to as the Sharma-Mittal divergence \cite{sharma1975} \cite{arndt2001}.\\
Calculating the gradient of $SM_{\alpha,s}\left(p\|q\right)$ (\ref{eq.SMas}) with respect to ``$q$" gives:
\begin{equation}
	\frac{\partial SM_{\alpha,s}\left(p\|q\right)}{\partial q_{j}}=\frac{1}{\alpha}\left\{\left[\sum_{i}\alpha p_{i}+\left(1-\alpha\right)q_{i}\right]^{\frac{s-\alpha}{\alpha-1}}-\left[\sum_{i}p_{i}^{\alpha}q_{i}^{1-\alpha}\right]^{\frac{s-\alpha}{\alpha-1}}p_{j}^{\alpha}q_{j}^{-\alpha}\right\}
	\label{eq.GSMa}
\end{equation}
He'll cancel for $p_{i}=q_{i}\ \forall i$.\\
On the other hand, the calculation of the gradient of $SMS_{\alpha}\left(p\|q\right)$ (\ref{eq.SMSa}) with respect to ``$q$" yields:
\begin{equation}
	\frac{\partial SMS_{\alpha,s}\left(p\|q\right)}{\partial q_{j}}=-\frac{1}{\alpha}\left[\sum_{i}p_{i}^{\alpha}q_{i}^{1-\alpha}\right]^{\frac{s-\alpha}{\alpha-1}}p_{j}^{\alpha}q_{j}^{-\alpha}
	\label{eq.GSMsa}
\end{equation}
It will never cancel.

\subsubsection{Scale invariance with respect to ``$q$".}
The divergence given by the relation (\ref{eq.SMas}) is made invariant using the invariance factor $K^{*}(p,q)=\frac{\sum_{j}p_{j}}{\sum_{j}q_{j}}$ which is not the nominal invariance factor; it is written after simplifications:
\begin{equation}	SM_{\alpha,s}I\left(p\|q\right)=\frac{1}{\alpha\left(s-1\right)}\left\{\left[\sum_{i}\bar{p}_{i}^{\alpha}\bar{q}_{i}^{1-\alpha}\right]^{\frac{s-1}{\alpha-1}}-\left[\sum_{i}\alpha \bar{p}_{i}+\left(1-\alpha\right)\bar{q}_{i}\right]^{\frac{s-1}{\alpha-1}}\right\}
\label{eq.SMasI}
\end{equation}
Its gradient with respect to ``$q$" is written:
\begin{equation}
\frac{\partial SM_{\alpha,s}I\left(p\|q\right)}{\partial q_{l}}=\frac{1}{\alpha\sum_{j}q_{j}}\left[\sum_{i}\bar{p}_{i}^{\alpha}\bar{q}_{i}^{1-\alpha}\right]^{\frac{s-\alpha}{\alpha-1}}\left(\sum_{i}\bar{p}^{\alpha}_{i}\bar{q}^{1-\alpha}_{i}-\bar{p}^{\alpha}_{l}\bar{q}^{\alpha}_{l}\right)
\label{eq.gradSMasI}	
\end{equation}

\section{Renyi Entropy related divergences.}
These divergences related to Renyi's entropy \cite{renyi1955} \cite{renyi1961} \cite{arndt2001} \cite{bercher2008}, correspond to the limit $d\rightarrow 1$ in the ``generalized Logarithm", that is $s\rightarrow 1$ in the divergences of the previous section.\\
If we perform this operation on the expressions (\ref{eq.SMas}) and (\ref{eq.SMSa}), we get respectively:
\begin{equation}	
R_{\alpha}\left(p\|q\right)=\frac{1}{\alpha\left(\alpha-1\right)}\left\{\log\sum_{i}p_{i}^{\alpha}q_{i}^{1-\alpha}-\log\sum_{i}\alpha p_{i}+\left(1-\alpha\right)q_{i}\right\}
\label{eq.Ra}
\end{equation}
The form (\ref{eq.Ra}) can be seen as an extension of Renyi's divergence to data fields whose sum is not equal to 1.\\
In the case of probability densities, i.e. with $\sum_{i}p_{i}=\sum_{i}q_{i}=1$, it comes:
\begin{equation}	RS_{\alpha}\left(p\|q\right)=\frac{1}{\alpha\left(\alpha-1\right)}\left\{\log\sum_{i}p_{i}^{\alpha}q_{i}^{1-\alpha}\right\}
\label{eq.RSa}
\end{equation}
The expression (\ref{eq.RSa}) is Renyi's divergence in its classical form related to Renyi's entropy and probability densities  \cite{arndt2001}.\\
The gradients with respect to ``$q$" can be deduced from the gradient expressions (\ref{eq.GSMa}) and (\ref{eq.GSMsa}) by making $s=1$; they are written respectively:
\begin{equation}
	\frac{\partial R_{\alpha}\left(p\|q\right)}{\partial q_{j}}=\frac{1}{\alpha}\left\{\left[\sum_{i}\alpha p_{i}+\left(1-\alpha\right)q_{i}\right]^{-1}-\left[\sum_{i}p_{i}^{\alpha}q_{i}^{1-\alpha}\right]^{-1} p_{j}^{\alpha}q_{j}^{-\alpha}\right\}
\end{equation}
and
\begin{equation}
	\frac{\partial RS_{\alpha}\left(p\|q\right)}{\partial q_{j}}=-\frac{1}{\alpha}\left[\sum_{i}p_{i}^{\alpha}q_{i}^{1-\alpha}\right]^{-1}p_{j}^{\alpha}q_{j}^{-\alpha}
\end{equation}
We can observe that $\frac{\partial R_{\alpha}\left(p\|q\right)}{\partial q_{j}}$ will be equal to zero if $p_{i}=q_{i}\ \forall i$, whereas $\frac{\partial RS_{\alpha}\left(p\|q\right)}{\partial q_{j}}$ could never be zero.

\subsubsection{Scale invariance with respect to ``$q$".}
For this divergence given by the relation (\ref{eq.Ra}), the invariance factor cannot be calculated explicitly, so we use as invariance factor the expression $K^{*}(p,q)=\frac{\sum_{j}p_{j}}{\sum_{j}q_{j}}$.\\
In these conditions, the invariant divergence is written after simplifications:
\begin{equation}	
RI_{\alpha}\left(p\|q\right)=\frac{\log\sum_{j}p_{j}}{\alpha\left(\alpha-1\right)}\left[\log\sum_{i}\bar{p}_{i}^{\alpha}\bar{q}_{i}^{1-\alpha}\right]
\label{eq.RaI}
\end{equation}
In this expression, the multiplicative factor $\log\sum_{j}p_{j}$ can be omitted, and the gradient with respect to ``$q$" is given by:
\begin{equation}
	\frac{\partial RI_{\alpha}\left(p\|q\right)}{\partial q_{l}}=\frac{1}{\alpha\sum_{j}q_{j}}\left[1-\frac{\bar{p}^{\;\alpha}_{l}\bar{q}^{\;-\alpha}_{l}}{\sum_{j}\bar{p}^{\alpha}_{j}\bar{q}^{1-\alpha}_{j}}\right]
	\label{eq.gradRaI}
\end{equation}
Of course, it will be equal to zero if $p_{i}=q_{i}\ \ \forall i$.
 
\section{Arimoto entropy related divergences.}

\subsection{Direct form.}
\begin{figure}[h!]
\centering
\includegraphics[width=0.7\linewidth]{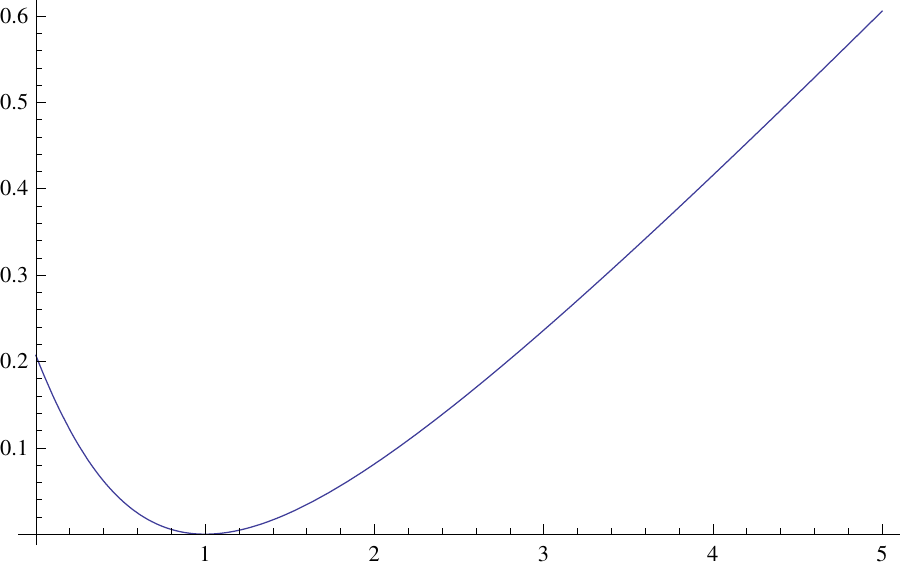}
\caption{Function $f_{c\delta}\left(x\right)=\frac{1}{\delta-1}\left[\left(\frac{1+x^{\delta}}{2}\right)^{\frac{1}{\delta}}-\left(\frac{1+x}{2}\right)\right], \delta=2$}
\label{fig:ARfc}
\end{figure}

\begin{figure}[ht!]
\centering
\includegraphics[width=0.7\linewidth]{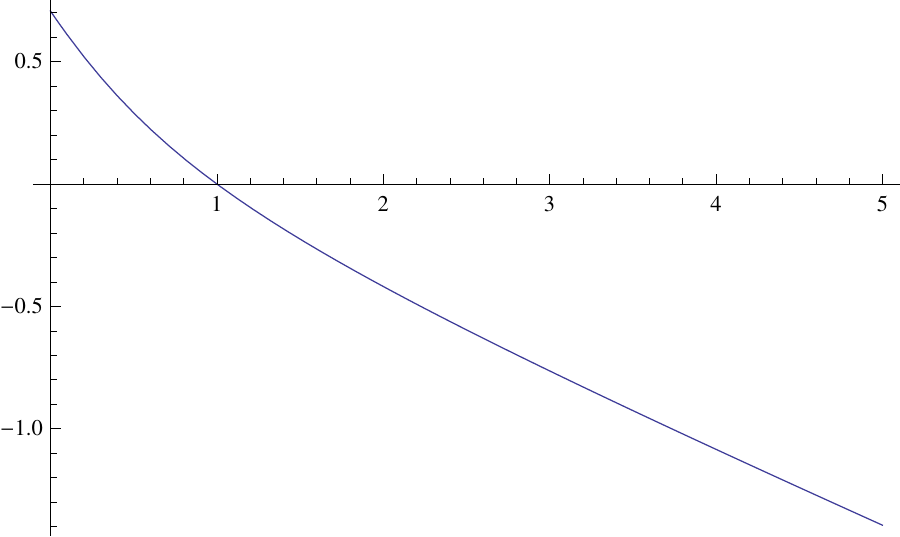}
\caption{Function $f_{1,\delta}\left(x\right)=\frac{1}{\delta-1}\left[\left(\frac{1+x^{\delta}}{2}\right)^{\frac{1}{\delta}}-x\right], \delta=2$}
\label{fig:ARf1}
\end{figure}

\begin{figure}[ht!]
\centering
\includegraphics[width=0.7\linewidth]{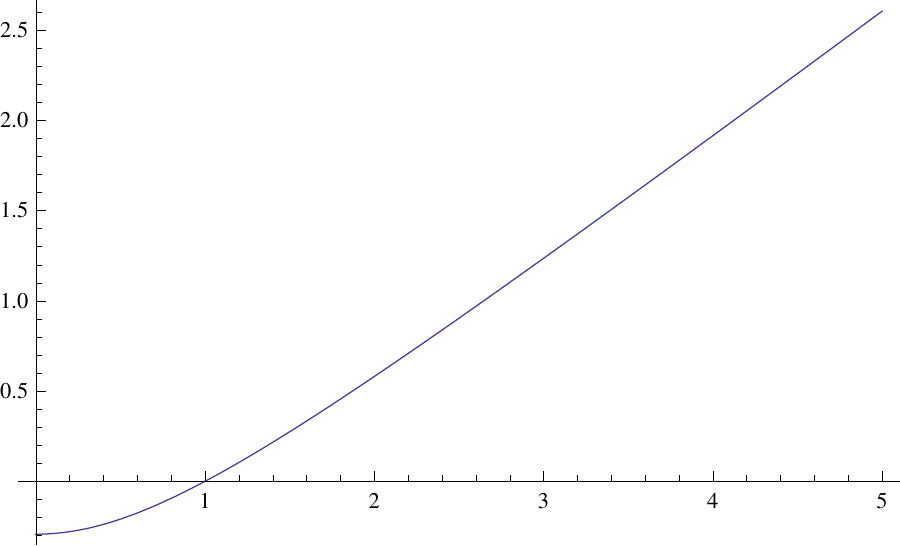}
\caption{Function $f_{2,\delta}\left(x\right)=\frac{1}{\delta-1}\left[\left(\frac{1+x^{\delta}}{2}\right)^{\frac{1}{\delta}}-1\right], \delta=2$}
\label{fig:ARf2}
\end{figure}

These divergences developed by Osterreicher \cite{osterreicher1996} \cite{osterreicher2013} rely on the use of generalized averages as in the Arimoto Entropy \cite{arimoto1971}. In their initial form, they are constructed in the Csiszär sense on the basis of the standard convex function shown in figure (\ref{fig:ARfc}) for $\delta=2$:
\begin{equation}	f_{c\delta}\left(x\right)=\frac{1}{\delta-1}\left[\left(\frac{1+x^{\delta}}{2}\right)^{\frac{1}{\delta}}-\left(\frac{1+x}{2}\right)\right]\ \ \ \delta\neq1
\end{equation}
This leads to the divergence:
\begin{equation}	AR\left(p\|q\right)=\frac{1}{\delta-1}\left[\sum_{i}\left(\frac{p_{i}^{\delta}+q_{i}^{\delta}}{2}\right)^{\frac{1}{\delta}}-\sum_{i}\left(\frac{p_{i}+q_{i}}{2}\right)\right]
\label{eq.divarimoto}
\end{equation}
This divergence is symmetrical, so it can be denoted as $AR\left(p,q\right)$.\\
The second term is clearly the unweighted arithmetic mean of the two data fields, whereas the first term corresponds, according to the value of ``$\delta$", to the different unweighted means between the two fields; indeed:\\
* if $\delta=2$, the first term is the square root mean.\\
* if $\delta=$1, the first term is the harmonic mean.\\
* if $\delta\rightarrow 0$, the first term is the geometric mean. \\
This approach allows us to find the divergences between averages, which will be developed in \textbf{chapter 7}.\\
The gradient with respect to ``$q$" will be written:
\begin{equation}
	\frac{\partial AR\left(p\|q\right)}{\partial q_{j}}=\frac{1}{2\left(\delta-1\right)}\left[\left(\frac{p_{j}^{\delta}+q_{j}^{\delta}}{2}\right)^{\frac{1-\delta}{\delta}}q_{j}^{\delta-1}-1\right]
\end{equation}
This gradient is equal to zero if $p_{i}=q_{i}\ \forall i$.\\
Since the basic convex function is a standard convex function, we can try to show the associated simple convex functions; after a few thinking about them, we can show 2 of these functions (as always), which are represented on figures (\ref{fig:ARf1}) and (\ref{fig:ARf2}) for $\delta=2$::
\begin{equation}	f_{1,\delta}\left(x\right)=\frac{1}{\delta-1}\left[\left(\frac{1+x^{\delta}}{2}\right)^{\frac{1}{\delta}}-x\right]
\end{equation}

and
\begin{equation}	f_{2,\delta}\left(x\right)=\frac{1}{\delta-1}\left[\left(\frac{1+x^{\delta}}{2}\right)^{\frac{1}{\delta}}-1\right]
\end{equation}

They will respectively lead to the divergences:
\begin{equation}
ARS1\left(p\|q\right)=\frac{1}{\delta-1}\left[\sum_{i}\left(\frac{p_{i}^{\delta}+q_{i}^{\delta}}{2}\right)^{\frac{1}{\delta}}-\sum_{i}p_{i}\right]
\end{equation}
and
\begin{equation}	
ARS2\left(p\|q\right)=\frac{1}{\delta-1}\left[\sum_{i}\left(\frac{p_{i}^{\delta}+q_{i}^{\delta}}{2}\right)^{\frac{1}{\delta}}-\sum_{i}q_{i}\right]
\end{equation}
If in these divergences as well as in $AR\left(p\|q\right)$ one introduces the simplification $\sum_{i}p_{i}=\sum_{i}q_{i}=1$, one obtains the divergence of Arimoto \cite{arndt2001}:
\begin{equation}	
ARS\left(p\|q\right)=\frac{1}{\delta-1}\left[\sum_{i}\left(\frac{p_{i}^{\delta}+q_{i}^{\delta}}{2}\right)^{\frac{1}{\delta}}-1\right]
\end{equation}
Of course, such a divergence can only be used to compare data fields whose sum is explicitly equal to $1$.\\
The gradients with respect to ``$q$" of the $ARS1\left(p\|q\right)$ and $ARS\left(p\|q\right)$ divergences will never cancel, while the gradient of $ARS2\left(p\|q\right)$ may become zero, but not for $p_{i}=q_{i}\ \forall i$.

\subsection{Related dual divergences.}
They are built on the dual convex functions of those used in the previous section.\\
We can note that $\breve{f}_{1,\delta}\left(x\right)=f_{2,\delta}\left(x\right)$ and that $\breve{f}_{2,\delta}\left(x\right)=f_{1,\delta}\left(x\right)$.

\subsection{Symmetrical divergence.}
Taking into account the remark in the previous section, it is constructed, in the sense of Jeffreys, on the basis of $\hat{f}_{\delta}\left(x\right)=\frac{\breve{f}_{1,\delta}\left(x\right)+f_{1,\delta}\left(x\right)}{2}=\frac{\breve{f}_{2,\delta}\left(x\right)+f_{2,\delta}\left(x\right)}{2}=f_{c,\delta}\left(x\right)$, and of course leads to the divergence $AR\left(p,q\right)$.

\subsection{Weighted versions of these divergences.}
In a first simple variant of these divergences, weighted versions of these divergences can be introduced by replacing the first term of the divergence $AR\left(p,q\right)$ (\ref{eq.divarimoto}) by a weighted generalized mean and second term by a weighted arithmetic mean.\\
This is equivalent to constructing a Csiszâr's divergence by using the basic standard convex function:
\begin{equation}	f_{c,\delta,\alpha}\left(x\right)=\frac{1}{\left(1-\alpha\right)\left(\delta-1\right)}\left\{\left[\alpha x^{\delta}+\left(1-\alpha\right)\right]^{\frac{1}{\delta}}-\left[\alpha x+\left(1-\alpha\right)\right]\right\}
\end{equation}
With $0\leq\alpha\leq1$, this leads to the divergence:
\begin{equation}	AR_{\delta,\alpha}\left(p\|q\right)=\frac{1}{\left(1-\alpha\right)\left(\delta-1\right)}\left\{\sum_{i}\left[\alpha p_{i}^{\delta}+\left(1-\alpha\right)q_{i}^{\delta}\right]^{\frac{1}{\delta}}-\sum_{i}\left[\alpha p_{i}+\left(1-\alpha\right)q_{i}\right]\right\}
\end{equation}
Of course, by varying the ``$\delta$" values as indicated above, we can review in the first term, the various weighted means.\\
Note that the symmetry property of the unweighted version has disappeared; it only exists for $\alpha=1/2$.\\
The gradient with respect to ``$q$" will be written:
\begin{equation}
\frac{\partial AR_{\delta,\alpha}\left(p\|q\right)}{\partial q_{j}}=\frac{1}{\delta-1}\left\{\left[\alpha p_{j}^{\delta}+\left(1-\alpha\right)q_{j}^{\delta}\right]^{\frac{1-\delta}{\delta}}q_{j}^{\delta-1}-1\right\}	
\end{equation}
It will be zero if $p_{j}=q_{j}\ \forall j$.\\
One can observe that a simplified version of this divergence corresponding to $\sum_{i}p_{i}=\sum_{i}q_{i}=1$ exists, but the gradient of this simplified form will not be equal to zero if $p_{j}=q_{j}\ \forall j$. 

\subsection{First type of extension.}
In order to retrieve any divergences based on differences between the weighted means developed in \textbf{chapter 7} (the unweighted case will be obtained immediately with $\alpha=\frac{1}{2}$), we can consider, as proposed in \cite{lanteri2011}, the standard convex function:
\begin{equation}	f_{c,\delta,\gamma,\alpha}\left(x\right)=\frac{1}{\left(1-\alpha\right)\left(\delta-1\right)}\left\{\left[\alpha x^{\delta}+\left(1-\alpha\right)\right]^{\frac{1}{\delta}}-\left[\alpha x^{\gamma}+\left(1-\alpha\right)\right]^{\frac{1}{\gamma}}\right\}
\end{equation}
With such a function, the resulting Csiszär divergence is written:
\begin{align}	AR_{\delta,\gamma,\alpha}\left(p\|q\right)=&\frac{1}{\left(1-\alpha\right)\left(\delta-1\right)}\left\{\sum_{i}\left[\alpha p_{i}^{\delta}+\left(1-\alpha\right)q_{i}^{\delta}\right]^{\frac{1}{\delta}}\right. \nonumber \\  & \left.  -\sum_{i}\left[\alpha p_{i}^{\gamma}+\left(1-\alpha\right)q_{i}^{\gamma}\right]^{\frac{1}{\gamma}}\right\}
\label{eq.AR}
\end{align}
The divergences obtained for different values of the parameters
``$\gamma$" and ``$\delta$" are summarized in the following table.\\

 \begin{table}[h]
\begin{center}
\begin{tabular}{|c|c|c|c|c|}
\hline
\backslashbox{$\delta$}{$\gamma$}
 & $-1$ & $0$ & $1$ & $2$ \\  \hline
$-1$ & $0$ & GM-HM & AM-HM & SM-HM  \\  \hline
$0$ & HM-GM &$0$  & AM-GM & SM-GM  \\  \hline
$1$ &  &  & Jensen-Shannon &   \\  \hline
$2$ & SM-HM & SM-GM & SM-AM & $0$  \\  \hline
\end{tabular}
\caption{Divergences resulting from various values of ``$\gamma$" and ``$\delta$". }
\end{center}
\end{table}

The divergences not listed in this table are either negative or non-convex; moreover the Jensen-Shannon divergence implies the computation of limits $\delta\rightarrow 1$, $\gamma\rightarrow 1$.\\
There is no simplified version of the divergence (\ref{eq.AR}); its gradient with respect to ``$q$" is written:
\begin{align}
\frac{\partial AR_{\delta,\gamma,\alpha}\left(p\|q\right)}{\partial q_{j}}=&\frac{1}{\delta-1}\left\{\left[\alpha p_{j}^{\delta}+\left(1-\alpha\right)q_{j}^{\delta}\right]^{\frac{1-\delta}{\delta}}q_{j}^{\delta-1}\right. \nonumber \\  & \left. -\left[\alpha p_{j}^{\gamma}+\left(1-\alpha\right)q_{j}^{\gamma}\right]^{\frac{1-\gamma}{\gamma}}q_{j}^{\gamma-1}\right\}	
\end{align}
Such a gradient will be zero pour $p_{i}=q_{i}\ \forall i$

\subsection{Extension in the sense of the Generalized Logarithm.}
Since the divergence (\ref{eq.AR}) shows a difference of 2 positive terms, an increasing function can be applied to each of the two terms, for example, the generalized logarithm, with the exponent ``$1-d$" noted $1-d=\frac{s-1}{\delta-1}$; this gives the divergence:
\begin{align}	
ARSM_{\delta,\gamma,\alpha,s}\left(p\|q\right)=&\frac{1}{\left(1-\alpha\right)\left(s-1\right)}\left\{\left[\sum_{i}\left[\alpha p_{i}^{\delta}+\left(1-\alpha\right)q_{i}^{\delta}\right]^{\frac{1}{\delta}}\right]^{\frac{s-1}{\delta-1}} 
\right. \nonumber \\  & \left. -\left[\sum_{i}\left[\alpha p_{i}^{\gamma}+\left(1-\alpha\right)q_{i}^{\gamma}\right]^{\frac{1}{\gamma}}\right]^{\frac{s-1}{\delta-1}}
\right\}
\end{align}
From this expression, the operation $s\rightarrow\delta$ ($d\rightarrow0$) leads to the divergence of the previous section, while the transition $s\rightarrow1$ ($d\rightarrow1$) leads to the following logarithmic form:
\begin{align}
	LARSM_{\delta,\gamma,\alpha}\left(p\|q\right)=&\frac{1}{\left(1-\alpha\right)\left(\delta-1\right)}\left\{\log\sum_{i}\left[\alpha p_{i}^{\delta}+\left(1-\alpha\right)q_{i}^{\delta}\right]^{\frac{1}{\delta}} 
\right. \nonumber \\  & \left. -\log\sum_{i}\left[\alpha p_{i}^{\gamma}+\left(1-\alpha\right)q_{i}^{\gamma}\right]^{\frac{1}{\gamma}}
\right\}
\end{align}
The specific cases corresponding to the various values of ``$\delta$'' and ``$\gamma$'' will be found in \textbf{chapter 7}.\\
The gradients with respect to ``$q$" of these two divergences are written respectively:
\begin{align}	
\frac{\partial ARSM_{\delta,\gamma,\alpha,s}\left(p\|q\right)}{\partial q_{j}}=&\frac{1}{\delta-1}\left\{\left[\sum_{i}\left[\alpha p_{i}^{\delta}+\left(1-\alpha\right)q_{i}^{\delta}\right]^{\frac{1}{\delta}}\right]^{\frac{s-\delta}{\delta-1}}\left[\alpha p_{j}^{\delta}+\left(1-\alpha\right)q_{j}^{\delta}\right]^{\frac{1-\delta}{\delta}}q_{j}^{\delta-1} 
\right. \nonumber \\  & \left. -\left[\sum_{i}\left[\alpha p_{i}^{\gamma}+\left(1-\alpha\right)q_{i}^{\gamma}\right]^{\frac{1}{\gamma}}\right]^{\frac{s-\delta}{\delta-1}}\left[\alpha p_{j}^{\gamma}+\left(1-\alpha\right)q_{j}^{\gamma}\right]^{\frac{1-\gamma}{\gamma}}q_{j}^{\gamma-1}
\right\}
\label{eq.GARSM}
\end{align}
and
\begin{align}	
\frac{\partial LARSM_{\delta,\gamma,\alpha,s}\left(p\|q\right)}{\partial q_{j}}=&\frac{1}{\delta-1}\left\{\left[\sum_{i}\left[\alpha p_{i}^{\delta}+\left(1-\alpha\right)q_{i}^{\delta}\right]^{\frac{1}{\delta}}\right]^{-1}\left[\alpha p_{j}^{\delta}+\left(1-\alpha\right)q_{j}^{\delta}\right]^{\frac{1-\delta}{\delta}}q_{j}^{\delta-1} 
\right. \nonumber \\  & \left. -\left[\sum_{i}\left[\alpha p_{i}^{\gamma}+\left(1-\alpha\right)q_{i}^{\gamma}\right]^{\frac{1}{\gamma}}\right]^{-1}\left[\alpha p_{j}^{\gamma}+\left(1-\alpha\right)q_{j}^{\gamma}\right]^{\frac{1-\gamma}{\gamma}}q_{j}^{\gamma-1}
\right\}
\end{align}
This last expression can be obtained from (\ref{eq.GARSM}) by simply doing $s=1$.\\
The components of these gradients will be zero for $p_{j}=q_{j}\ \forall j$.

\section{Jensen divergences based on Entropies.}
We use as basic convex functions, the opposite of the Shannon, Havrda-Charvat and Renyi Entropies.

\subsection{Shannon's entropy.}
On the basis of the simple convex  function which is the opposite of Shannon's entropy, that is:
\begin{equation}
	f\left(x\right)=x \log x
\end{equation}
we treat the case of Jensen's divergence weighted by $0<\beta<$1, the classic case $\beta=1/2$ is immediately inferred.\\
We have:
\begin{align}
	J_{S}\left(p\|q\right)=&\beta\sum_{i}p_{i}\log p_{i}+\left(1-\beta\right)\sum_{i}q_{i}\log q_{i}
	 \nonumber \\  & -\sum_{i}\left[\beta p_{i}+\left(1-\beta\right)q_{i}\right]\log\left[\beta p_{i}+\left(1-\beta\right)q_{i}\right]
	 \label{eq.JS}
\end{align}
The gradient with respect to ``$q$" is expressed as:
\begin{equation}
	\frac{\partial J_{S}\left(p\|q\right)}{\partial q_{j}}=\left(1-\beta\right)\left\{\log q_{j}-\log\left[\beta p_{j}+\left(1-\beta\right)q_{j}\right]\right\}
\end{equation}
Of course these components will be zero for $p_{j}=q_{j}\ \forall j$.

\subsubsection{Scale invariance with respect to ``$q$"}
For this divergence given by the relation (\ref{eq.JS}), the invariance factor can not be calculated explicitly, so we use as invariance factor the expression $K^{*}(p,q)=\frac{\sum_{j}p_{j}}{\sum_{j}q_{j}}$.\\
In these conditions, the invariant divergence is written after simplifications:
\begin{align}
	J_{S}I\left(p\|q\right)=&\beta\sum_{i}\bar{p}_{i}\log \bar{p}_{i}+\left(1-\beta\right)\sum_{i}\bar{q}_{i}\log \bar{q}_{i}
	 \nonumber \\  & -\sum_{i}\left[\beta \bar{p}_{i}+\left(1-\beta\right)\bar{q}_{i}\right]\log\left[\beta \bar{p}_{i}+\left(1-\beta\right)\bar{q}_{i}\right]
	 \label{eq.JSI}
\end{align}
The gradient with respect to ``$q$" is expressed as:
\begin{equation}
\frac{\partial 	J_{S}I\left(p\|q\right)}{\partial q_{l}}=\frac{1-\beta}{\sum_{j}q_{j}}\left[\log\frac{\bar{q}_{l}}{\beta \bar{p}_{l}+\left(1-\beta\right)\bar{q}_{l}}-\sum_{i}\bar{q}_{i}\log\frac{\bar{q}_{i}}{\beta \bar{p}_{i}+\left(1-\beta\right)\bar{q}_{i}}\right]
 \label{eq.gradJSI}	
\end{equation}

\subsection{Havrda-Charvat entropy.}
Here, the convex base function is the opposite of the Havrda-Charvat entropy; we treat the case of the Jensen divergence weighted by $0<\beta<1$, the classic case $\beta=1/2$ is deduced immediately.\\
We have:
\begin{equation}	J_{HC}\left(p\|q\right)=\frac{1}{\alpha\left(\alpha-1\right)}\left\{\beta\sum_{i}p_{i}^{\alpha}+\left(1-\beta\right)\sum_{i}q_{i}^{\alpha}-\sum_{i}\left[\beta p_{i}+\left(1-\beta\right)q_{i}\right]^{\alpha}\right\}
\label{eq.JHC}
\end{equation}
The gradient with respect to ``$q$" is written:
\begin{equation}
	\frac{\partial J_{HC}\left(p\|q\right)}{\partial q_{j}}=\frac{1-\beta}{\alpha-1}\left\{q_{j}^{\alpha-1}-\left[\beta p_{j}+\left(1-\beta\right)q_{j}\right]^{\alpha-1}\right\}
\end{equation}

\subsubsection{Scale invariance with respect to ``$q$".}
For this divergence given by the relation (\ref{eq.JHC}), the invariance factor can not be calculated explicitly, so we use as invariance factor the expression $K^{*}(p,q)=\frac{\sum_{j}p_{j}}{\sum_{j}q_{j}}$.\\
Under these conditions, the invariant divergence is written after simplifications:
\begin{equation}	J_{HC}I\left(p\|q\right)=\frac{1}{\alpha\left(\alpha-1\right)}\left\{\beta\sum_{i}\bar{p}_{i}^{\alpha}+\left(1-\beta\right)\sum_{i}\bar{q}_{i}^{\alpha}-\sum_{i}\left[\beta \bar{p}_{i}+\left(1-\beta\right)\bar{q}_{i}\right]^{\alpha}\right\}
\label{eq.JHCI}
\end{equation}
The gradient with respect to ``$q$" is expressed as:
\begin{align}
\frac{\partial 	J_{HC}I\left(p\|q\right)}{\partial q_{l}}=&\frac{1-\beta}{\left(\alpha-1\right)\sum_{j}q_{j}}\left\{\bar{q}^{\alpha-1}_{l}-\sum_{i}\bar{q}^{\alpha}_{i}-\left[\beta \bar{p}_{l}+\left(1-\beta\right)\bar{q}_{l}\right]^{\alpha-1}\right. \nonumber \\  & \left. 
+\sum_{i}\bar{q}_{i}\left[\beta \bar{p}_{i}+\left(1-\beta\right)\bar{q}_{i}\right]^{\alpha-1}\right\}
\label{eq.gradJHCI}	
\end{align}
If, in the expression of the divergence (\ref{eq.JHC}), we replace the arithmetic mean by a generalized geometric mean, which will not change the sign of the expression, we do what O. Michel \cite{michel1994} has proposed for Renyi's divergence and on which we will come back later; the divergence thus constructed will be written:
\begin{equation}	J2_{HC}\left(p\|q\right)=\frac{1}{\alpha\left(\alpha-1\right)}\left\{\beta\sum_{i}p_{i}^{\alpha}+\left(1-\beta\right)\sum_{i}q_{i}^{\alpha}-\sum_{i}\left[ p_{i}^{\beta}q_{i}^{1-\beta}\right]^{\alpha}\right\}
\label{eq.J2HC}
\end{equation}
The gradient with respect to ``$q$" will be:
\begin{equation}
	\frac{\partial J2_{HC}\left(p\|q\right)}{\partial q_{j}}=\frac{1-\beta}{\alpha-1}\left\{q_{j}^{\alpha-1}-\left[ p_{j}^{\beta}q_{j}^{1-\beta}\right]^{\alpha-1}\left(\frac{p_{j}}{q_{j}}\right)^{\beta}\right\}
\end{equation}

\subsubsection{Scale invariance with respect to ``$q$".}
For this extension proposed by O.Michel \cite{michel1994} given by (\ref{eq.J2HC}), using the invariance factor $K^{*}\left(p,q\right)$ previously defined, we obtain after simplifications the following invariant divergence:
\begin{equation}	J2_{HC}I\left(p\|q\right)=\frac{1}{\alpha\left(\alpha-1\right)}\left\{\beta\sum_{i}\bar{p}_{i}^{\alpha}+\left(1-\beta\right)\sum_{i}\bar{q}_{i}^{\alpha}-\sum_{i}\left[ \bar{p}^{\beta}_{i}\bar{q}^{1-\beta}_{i}\right]^{\alpha}\right\}
\label{eq.J2HCI}
\end{equation}
The gradient with respect to ``$q$" will be:
\begin{equation}
\frac{\partial 	J2_{HC}I\left(p\|q\right)}{\partial q_{l}}=\frac{1-\beta}{\left(\alpha-1\right)\sum_{j}q_{j}}\left\{\bar{q}^{\alpha-1}_{l}-\sum_{i}\bar{q}^{\alpha}_{i}-\left[ \bar{p}^{\beta}_{l}\bar{q}^{1-\beta}_{l}\right]^{\alpha}\bar{p}^{\beta}_{l}\bar{q}^{-\beta}_{l}+\sum_{i}\left[\bar{p}^{\beta}_{i}\bar{q}^{1-\beta}_{i}\right]^{\alpha+1}\right\}
\label{eq.gradJ2HCI}	
\end{equation}
This being the case, we can also continue to follow the reasoning of O. Michel \cite{michel1994} in his work on Renyi's divergence, and introduce a Jensen-Havrda Charvat divergence between a generalized arithmetic mean of the 2 fields and one or the other of the fields, then, one can always invert the order of the arguments insofar as the divergences thus formed are not necessarily symmetrical.

\subsection{Renyi's entropy.}
Here, the convex base function is the opposite of Renyi's entropy, these leads to Arndt's $^{1}R_{\alpha}^{1}$ \cite{arndt2001}. We treat the case of Jensen's divergence weighted by $0<\beta<1$, the classic case $\beta=1/2$ is inferred immediately.\\
The corresponding divergence is written as:
\begin{align}
	J_{R}^{\beta}\left(p\|q\right)=&\frac{1}{\alpha\left(\alpha-1\right)}\left\{\beta \log\sum_{i}p_{i}^{\alpha}+\left(1-\beta\right)\log\sum_{i}q_{i}^{\alpha}\right. \nonumber \\  & \left. -\log\sum_{i}\left[\beta p_{i}+\left(1-\beta\right)q_{i}\right]^{\alpha}\right\}
\label{eq.JRab}
\end{align}
The gradient with respect to ``$q$" will be:
\begin{equation}
	\frac{\partial J_{R}^{\beta}\left(p\|q\right)}{\partial q_{j}}=\frac{1-\beta}{\alpha-1}\left\{\frac{q_{j}^{\alpha-1}}{\sum_{i}q_{i}^{\alpha}}-\frac{\left[\beta p_{j}+\left(1-\beta\right)q_{j}\right]^{\alpha-1}}{\sum_{i}\left[\beta p_{i}+\left(1-\beta\right)q_{i}\right]^{\alpha}}\right\}	
\end{equation}

\subsubsection{Scale invariance with respect to ``$q$".}
For such a divergence given by the relation (\ref{eq.JRab}), we use the $K^{*}(p,q)$ invariance factor already mentioned, which leads to the invariant divergence:
\begin{align}
	J_{R}^{\beta}I\left(p\|q\right)=&\frac{1}{\alpha\left(\alpha-1\right)}\left\{\beta \log\sum_{i}\bar{p}_{i}^{\alpha}+\left(1-\beta\right)\log\sum_{i}\bar{q}_{i}^{\alpha}\right. \nonumber \\  & \left. -\log\sum_{i}\left[\beta \bar{p}_{i}+\left(1-\beta\right)\bar{q}_{i}\right]^{\alpha}\right\}
\label{eq.JRabI}
\end{align}
The gradient with respect to ``$q$" is written as:
\begin{align}
\frac{\partial 	J_{R}^{\beta}I\left(p\|q\right)}{\partial q_{l}}=\frac{1-\beta}{\left(\alpha-1\right)\sum_{j}q_{j}}&\left\{\frac{\bar{q}^{\alpha-1}_{l}}
{\sum_{j}q^{\alpha}_{j}}-1-\frac{\left[\beta \bar{p}_{l}+\left(1-\beta\right)\bar{q}_{l}\right]^{\alpha-1}}{\sum_{i}\left[\beta \bar{p}_{i}+\left(1-\beta\right)\bar{q}_{i}\right]^{\alpha}}\right. \nonumber \\  & \left.
+\frac{\sum_{i}\bar{q}_{i}\left[\beta \bar{p}_{i}+\left(1-\beta\right)\bar{q}_{i}\right]^{\alpha-1}}{\sum_{i}\left[\beta \bar{p}_{i}+\left(1-\beta\right)\bar{q}_{i}\right]^{\alpha}}\right\}
\label{eq.gradJRabI}	
\end{align}

If in the expression (\ref{eq.JRab}) the generalized arithmetic mean is replaced by a generalized geometric mean, as proposed by O. Michel in \cite{michel1994}, the divergence is written:
\begin{align}
	J2_{R}^{\beta}\left(p\|q\right)=&\frac{1}{\alpha\left(\alpha-1\right)}\left\{\beta \log\sum_{i}p_{i}^{\alpha}+\left(1-\beta\right)\log\sum_{i}q_{i}^{\alpha}\right. \nonumber \\  & \left. -\log\sum_{i}\left( p_{i}^{\beta}q_{i}^{1-\beta}\right)^{\alpha}\right\}
	\label{eq.J2R}
\end{align}
Whose gradient is:
\begin{equation}
	\frac{\partial J2_{R}^{\beta}\left(p\|q\right)}{\partial q_{j}}=\frac{1-\beta}{\alpha-1}\left\{\frac{q_{j}^{\alpha-1}}{\sum_{i}q_{i}^{\alpha}}-\frac{\left[ p_{j}^{\beta}q_{j}^{1-\beta}\right]^{\alpha-1}}{\sum_{i}\left[ p_{i}^{\beta}q_{i}^{1-\beta}\right]^{\alpha}}\frac{p_{j}^{\beta}}{q_{j}^{\beta}}\right\}
\label{eq.gradJ2R}		
\end{equation}

\subsubsection{Scale invariance with respect to ``$q$".}
From this modification proposed by O.Michel \cite{michel1994} given by (\ref{eq.J2R}); the corresponding scale invariant divergence is deduced by introducing the invariance factor $K^{*}(p,q)$:
\begin{align}
	J2_{R}^{\beta}I\left(p\|q\right)=&\frac{1}{\alpha\left(\alpha-1\right)}\left\{\beta \log\sum_{i}\bar{p}_{i}^{\alpha}+\left(1-\beta\right)\log\sum_{i}\bar{q}_{i}^{\alpha}\right. \nonumber \\  & \left. -\log\sum_{i}\left( \bar{p}_{i}^{\beta}\bar{q}_{i}^{1-\beta}\right)^{\alpha}\right\}
	\label{eq.J2RI}
\end{align}
The gradient with respect to ``$q$" is expressed as:
\begin{equation}
\frac{\partial 	J2_{R}^{\beta}I\left(p\|q\right)}{\partial q_{l}}=\frac{1-\beta}{\left(\alpha-1\right)\sum_{j}q_{j}}\left[\frac{\bar{q}^{\alpha-1}_{l}}{\sum_{j}\bar{q}^{\alpha}_{j}}-\frac{\left(\bar{p}^{\beta}_{l}\bar{q}^{1-\beta}_{l}\right)^{\alpha-1}\bar{p}^{\beta}_{l}\bar{q}^{-\beta}_{l}}{\sum_{j}\left(\bar{p}^{\beta}_{j}\bar{q}^{1-\beta}_{j}\right)^{\alpha}}\right]
\label{eq.gradJ2RI}	
\end{equation}
In the same article \cite{michel1994}, another variation of this type of divergence has been suggested. To fully understand how to obtain it, we'll proceed in successive steps.\\
 First we write the dual divergence of $J2_{R}^{\beta}\left(p\|q\right)$ which gives:
\begin{align}
	J2_{R}^{\beta}\left(q\|p\right)=&\frac{1}{\alpha\left(\alpha-1\right)}\left\{\beta \log\sum_{i}q_{i}^{\alpha}+\left(1-\beta\right)\log\sum_{i}p_{i}^{\alpha}\right. \nonumber \\  & \left. -\log\sum_{i}\left( q_{i}^{\beta}p_{i}^{1-\beta}\right)^{\alpha}\right\}
\end{align}
Then we construct a divergence from the previous general form, but by replacing $q$ by a linear combination of ``$p$'' and ``$q$'' with a weighting factor ``$\gamma$'', i.e. $\gamma p+\left(1-\gamma\right)q$, which gives:
\begin{align}
	J3_{R}^{\beta,\gamma}\left(p\|q\right)=\frac{1}{\alpha\left(\alpha-1\right)}&\left\{\beta \log\sum_{i}\left[\gamma p_{i}+\left(1-\gamma\right)q_{i}\right]^{\alpha}+\left(1-\beta\right)\log\sum_{i}p_{i}^{\alpha} \right. \nonumber \\ & \left.
	-\log\sum_{i}\left( \left[\gamma p_{i}+\left(1-\gamma\right)q_{i}\right]^{\beta}p_{i}^{1-\beta}\right)^{\alpha}\right\}
\end{align}
Of course, this little exercise can go on for a very long time and may not be of much interest...\\
We can still calculate the gradient relative to ``$q$'':
\begin{align}
	\frac{\partial J3_{R}^{\beta,\gamma}\left(p\|q\right)}{\partial q_{j}}&=\beta\frac{1-\gamma}{\alpha-1}\left\{\frac{\left[\gamma p_{j}+\left(1-\gamma\right) q_{j}\right]^{\alpha-1}}{\left[\sum_{i}\gamma p_{i}+\left(1-\gamma\right) q_{i}\right]^{\alpha}} 
	\right. \nonumber \\ & \left.	
	-\frac{\left( \left[\gamma p_{j}+\left(1-\gamma\right)q_{j}\right]^{\beta}p_{j}^{1-\beta}\right)^{\alpha}}{\left(\sum_{i} \left[\gamma p_{i}+\left(1-\gamma\right)q_{i}\right]^{\beta}p_{i}^{1-\beta}\right)^{\alpha}}\frac{1}{\gamma p_{j}+\left(1-\gamma\right)q_{j}}\right\}
\end{align}
\setcounter{table}{0}  \setcounter{equation}{0}  \setcounter{figure}{0} \setcounter{chapter}{5} \setcounter{section}{0} 
\chapter{chapter 5 - Alpha, Beta and Gamma divergences}  \label{chptr::chapitre5}
\section{Introduction.}
\begin{figure}[h!]
\centering
\includegraphics[width=0.7\linewidth]{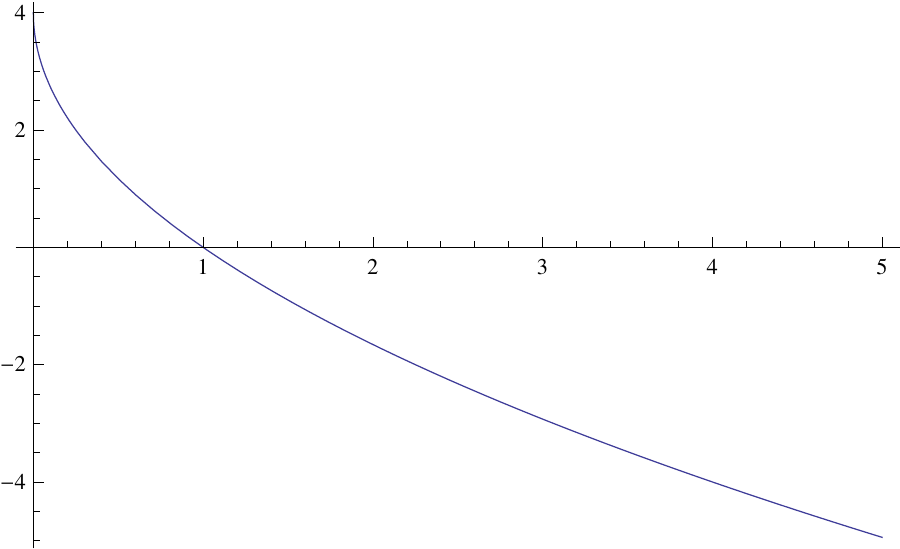}
\caption{Function $f_{1}\left(x\right)=\frac{1}{\lambda\left(\lambda-1\right)}\left(x^{\lambda}-1\right), \lambda=0.5$}
\label{fig:ABGf1}
\end{figure}

\begin{figure}[h!]
\centering
\includegraphics[width=0.7\linewidth]{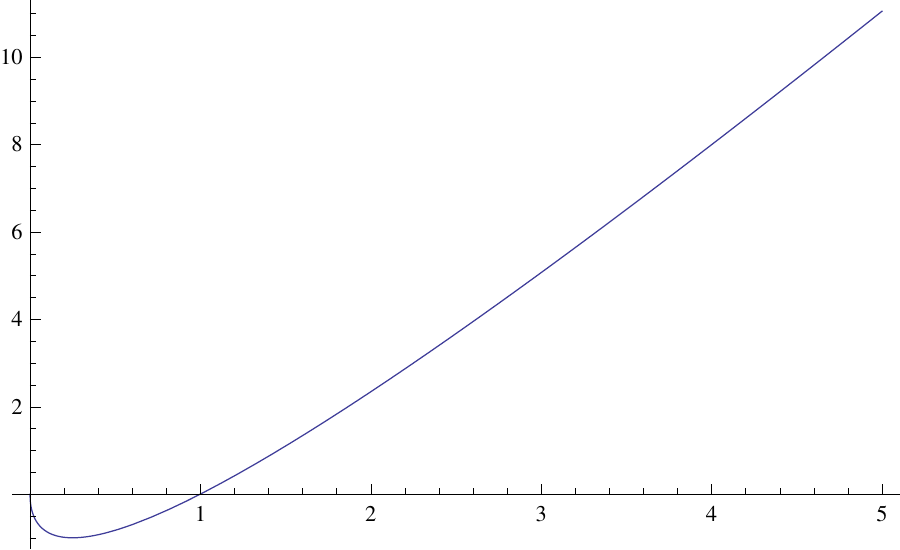}
\caption{Function $f_{2}\left(x\right)=\frac{1}{\lambda\left(\lambda-1\right)}\left(x^{\lambda}-x\right), \lambda=0.5$}
\label{fig:ABGf2}
\end{figure}

In this chapter, we present the information corresponding to the divergences constructed on the basis of functions of the type:
\begin{equation}
	f\left(x\right)=\frac{x^{\lambda}}{\lambda-1}
\end{equation}
The variable ``$x$" is assumed to be positive and we use the parameter ``$\lambda$" which hasn't been used yet and can become $\alpha,\beta$ or $\gamma$ as needed.\\
This function is convex, for $\lambda>0$; to make it convex regardless of ``$\lambda$", we have to write:
\begin{equation}
	f\left(x\right)=\frac{x^{\lambda}}{\lambda\left(\lambda-1\right)}
\end{equation}
To get moreover $f\left(1\right)=0$, we have to turn it into:
\begin{equation}
	f_{1}\left(x\right)=\frac{1}{\lambda\left(\lambda-1\right)}\left(x^{\lambda}-1\right)
	\label{eq.ABGf1}
\end{equation}
This function is shown on the figure (\ref{fig:ABGf1}),
\\
or else:
\begin{equation}
	f_{2}\left(x\right)=\frac{1}{\lambda\left(\lambda-1\right)}\left(x^{\lambda}-x\right)
	\label{eq.ABGf2}
\end{equation}
which is represented on the figure: (\ref{fig:ABGf2}).

\begin{figure}[h!]
\centering
\includegraphics[width=0.7\linewidth]{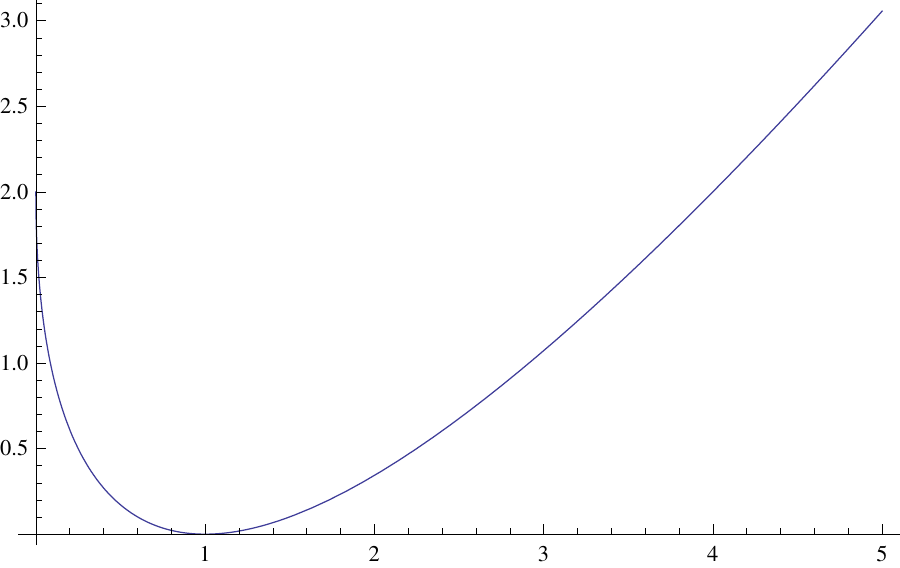}
\caption{Function $f_{c}\left(x\right)=\frac{1}{\lambda\left(\lambda-1\right)}\left[x^{\lambda}-\lambda x-\left(1-\lambda\right)\right], \lambda=0.5$}
\label{fig:ABGfc}
\end{figure}

The functions $f_{1}\left(x\right)$ and $f_{2}\left(x\right)$ are \textit{``simple convex functions"}.\\
If we want to have, in addition, a null derivative in $x=1$, that is to say if we want to obtain a \textit{``standard convex function"}, we obtain in both cases the same function $f_{c}\left(x\right)$ which is represented on the figure (\ref{fig:ABGfc}) and which is written: 

\begin{equation}
	f_{c}\left(x\right)=\frac{1}{\lambda\left(\lambda-1\right)}\left[x^{\lambda}-\lambda x-\left(1-\lambda\right)\right]
	\label{eq.ABGfc}
\end{equation}

In all these figures, $\lambda=0.5$.\\

An important point to keep in mind is that these functions are convex whatever $``\lambda"$ is due to the $``1/\lambda"$ factor that's been introduced, which may not always be found in the divergences proposed in the literature, but it will allow the $\lambda\rightarrow 0$ limit to be computed when necessary.\\
The question now is, what divergences can be built using such functions?\\
We will consider two constructive modes that have already been dealt with in the literature: Csiszär's divergences and Bregman's divergences.\\ 
In addition, we will also consider Jensen's divergences.

\section{Csiszär's divergences - Alpha divergences.}
They have been largely developed in the works of Amari \cite{amari2009} and Cichocki \cite{cichocki2010}; they are based on the functions $f_{1}\left(x\right)$, $f_{2}\left(x\right)$ and $f_{c}\left(x\right)$.
\subsection{On $f_{1}\left(x\right)$.}
\begin{figure}[h!]
\centering
\includegraphics[width=0.7\linewidth]{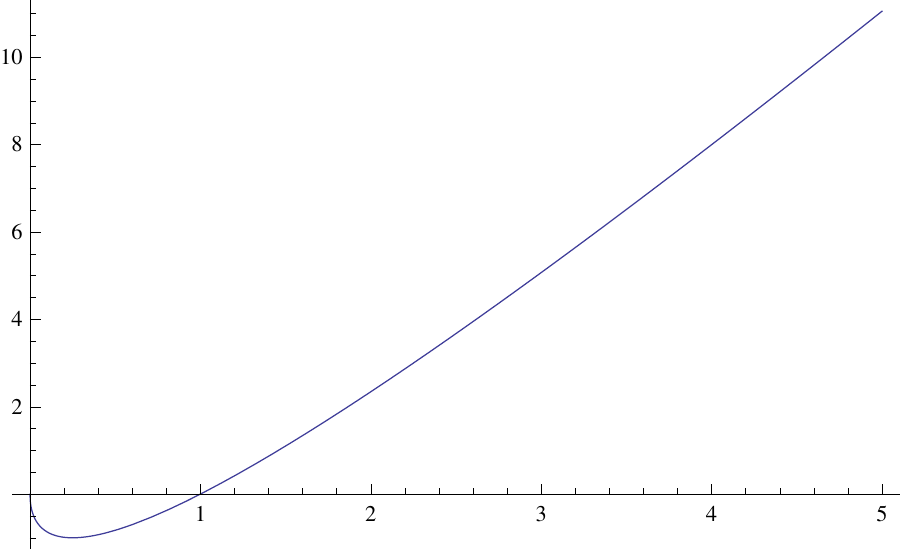}
\caption{Function $\breve{f}_{1}\left(x\right)=\frac{1}{\lambda\left(\lambda-1\right)}\left(x^{1-\lambda}-x\right), \lambda=0.5$}
\label{fig:ABGft1}
\end{figure}

In this case we obtain the divergence:
\begin{equation}
	A_{1}\left(p\|q\right)=\frac{1}{\lambda\left(\lambda-1\right)}\sum_{i}\left(p_{i}^{\lambda}q_{i}^{1-\lambda}-q_{i}\right)
\end{equation}
If we calculate the gradient with the intention to minimize with respect to ``$q$", we get:
\begin{equation}
	\frac{\partial A_{1}\left(p\|q\right)}{\partial q_{j}}=\frac{1}{\lambda\left(\lambda-1\right)}\left[\left(1-\lambda\right)\left(\frac{p_{j}}{q_{j}}\right)^{\lambda}-1\right]
\end{equation}
This seems suitable for building a descent algorithm, but the problem will be as always, that the gradient doesn't cancel out for $p_{i}=q_{i}\ \forall i$.\\
Moreover, we can notice that if we operate on probability densities ($\sum_{i}p_{i}=\sum_{i}q_{i}=1$), the divergence simplifies, but in this case the corresponding gradient will never be zero.\\
The dual divergence is built on the mirror function represented in figure (\ref{fig:ABGft1})
\begin{equation}
	\breve{f}_{1}\left(x\right)=\frac{1}{\lambda\left(\lambda-1\right)}\left(x^{1-\lambda}-x\right)
	\label{eq.ABGft1}
\end{equation}
It is expressed as follows:
\begin{equation}		A_{1}\left(q\|p\right)=\frac{1}{\lambda\left(\lambda-1\right)}\sum_{i}\left(q_{i}^{\lambda}p_{i}^{1-\lambda}-p_{i}\right)
\end{equation}
On this divergence, we can see that the gradient with respect to ``$q$" never cancels out, hence a problem if we try to minimize.
\subsection{On $f_{2}\left(x\right)$.}
\begin{figure}[h!]
\centering
\includegraphics[width=0.7\linewidth]{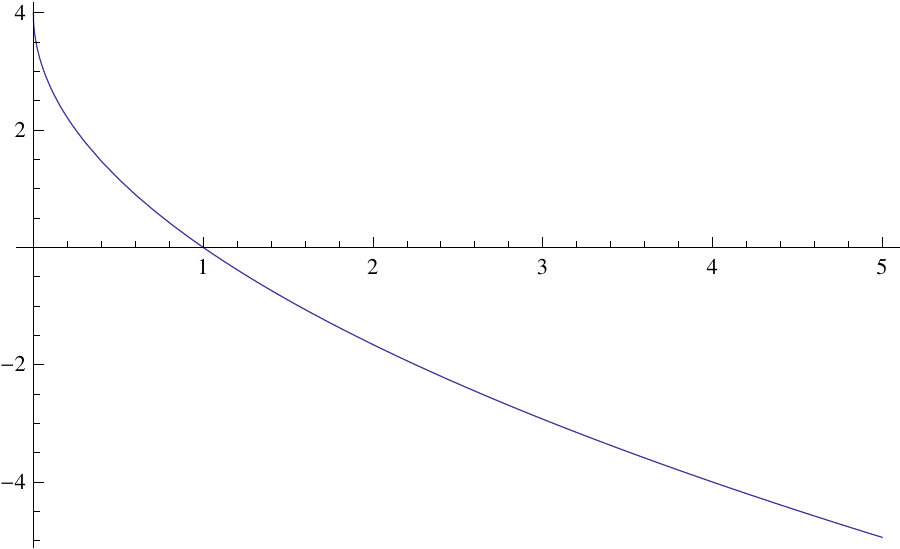}
\caption{Function $\breve{f}_{2}\left(x\right)=\frac{1}{\lambda\left(\lambda-1\right)}\left(x^{1-\lambda}-1\right), \lambda=0.5$}
\label{fig:ABGft2}
\end{figure}

In this case we obtain the divergence:
\begin{equation}
	A_{2}\left(p\|q\right)=\frac{1}{\lambda\left(\lambda-1\right)}\sum_{i}\left(p_{i}^{\lambda}q_{i}^{1-\lambda}- p_{i}\right)
\end{equation}
The gradient will be:
\begin{equation}
\frac{\partial A_{2}\left(p\|q\right)}{\partial q_{j}}=-\frac{1}{\lambda}\left(\frac{p_{j}}{q_{j}}\right)^{\lambda}	
\end{equation}
We can see that the gradient with respect to ``$q$" never cancels out.\\
If we operate on probability densities, the divergence becomes simpler, but the problem is again, the gradient will not cancel out.\\
The dual divergence will be built on the mirror function represented in figure (\ref{fig:ABGft2}):
\begin{equation}
	\breve{f}_{2}\left(x\right)=\frac{1}{\lambda\left(\lambda-1\right)}\left(x^{1-\lambda}-1\right)
\end{equation}
It is expressed as follows:
\begin{equation}
	A_{2}\left(q\|p\right)=\frac{1}{\lambda\left(\lambda-1\right)}\sum_{i}\left(q_{i}^{\lambda}p_{i}^{1-\lambda}- q_{i}\right)
\end{equation}
This divergence is made simpler if we operate on probability densities.\\
The gradient with respect to ``$q$" is:
\begin{equation}
	\frac{\partial A_{2}\left(q\|p\right)}{\partial q_{j}}=\frac{1}{\lambda\left(\lambda-1\right)}\left[\lambda\left(\frac{p_{j}}{q_{j}}\right)^{1-\lambda}-1\right]
\end{equation}
Of course, the problem is still that the gradient with respect to ``$q$" will not be zero for $p_{i}=q_{i}$.

\subsection{On $f_{c}\left(x\right)$.}
\begin{figure}[h!]
\centering
\includegraphics[width=0.7\linewidth]{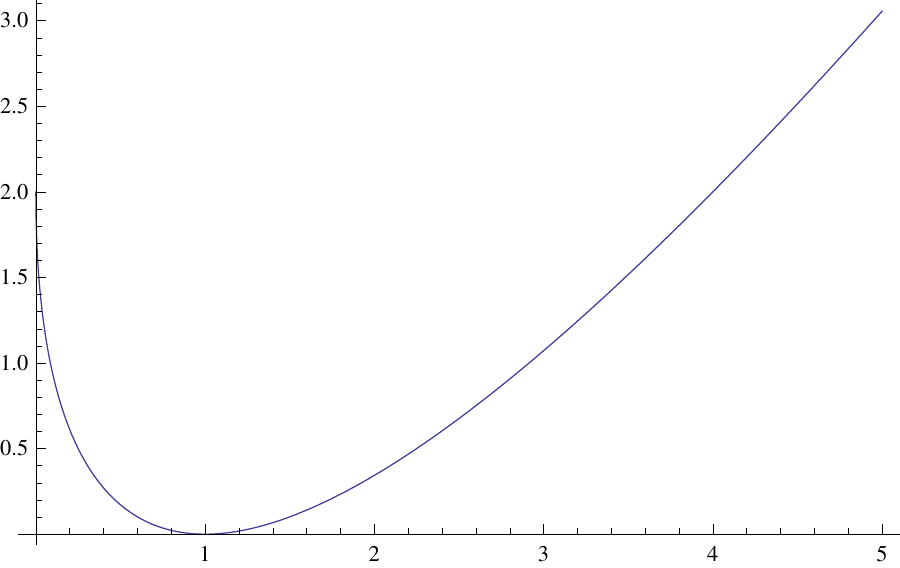}
\caption{Function $\breve{f}_{c}\left(x\right)=\frac{1}{\lambda\left(\lambda-1\right)}\left[x^{1-\lambda}-\lambda-\left(1-\lambda\right)x\right], \lambda=0.5$}
\label{fig:ABGftc}
\end{figure}

In this case we obtain the divergence:
\begin{equation}
	A\left(p\|q\right)=\frac{1}{\lambda\left(\lambda-1\right)}\sum_{i}\left(p_{i}^{\lambda}q_{i}^{1-\lambda}-\lambda p_{i}-\left(1-\lambda\right)q_{i}\right)
	\label{eq.AC}
\end{equation}
This divergence \cite{havrda1967} can also be thought of as a arithmetic-geometric divergence \cite{taneja2001}, it is also the ``Alpha divergence" of Amari \cite{cichocki2009}($\lambda\leftrightarrow \alpha$).\\
Its gradient with respect to ``$q$" is  written as follows:
\begin{equation}
		\frac{\partial A\left(p\|q\right)}{\partial q_{j}}=\frac{1}{\lambda}\left[1-\left(\frac{p_{j}}{q_{j}}\right)^{\lambda}\right]
\label{eq.gradAC}
\end{equation}
It will be zero if $p_{i}=q_{i}\ \forall i$.\\
The dual divergence will be constructed on the mirror function shown in figure (\ref{fig:ABGftc}):
\begin{equation}	
\breve{f}_{c}\left(x\right)=\frac{1}{\lambda\left(\lambda-1\right)}\left[x^{1-\lambda}-\lambda-\left(1-\lambda\right)x\right]
\end{equation}
It will be written:
\begin{equation}
	A\left(q\|p\right)=\frac{1}{\lambda\left(\lambda-1\right)}\sum_{i}\left(q_{i}^{\lambda}p_{i}^{1-\lambda}-\lambda q_{i}-\left(1-\lambda\right)p_{i}\right)
	\label{eq.ACduale}
\end{equation}
Its gradient with respect to ``$q$" will be:
\begin{equation}
		\frac{\partial A\left(q\|p\right)}{\partial q_{j}}=\frac{1}{\lambda-1}\left[\left(\frac{p_{j}}{q_{j}}\right)^{1-\lambda}-1\right]
\end{equation}
It will be zero if $p_{i}=q_{i}\ \forall i$.

\subsection{Some classical particular cases.}
If $\lambda\rightarrow 1$ we have $A\left(p\|q\right)=KL\left(p\|q\right)$.\\
If $\lambda\rightarrow 0$ we have $A\left(p\|q\right)=KL\left(q\|p\right)$.\\
If $\lambda=1/2$ we have $A\left(p\|q\right)=H \left(p\|q\right)$, that is $\sum_{i}\left(\sqrt{p_{i}}-\sqrt{q_{i}}\right)^{2}$ (Hellinger's divergence \cite{hellinger1909}).\\
If $\lambda=2$ we have $A\left(p\|q\right)=\chi^{2}_{N}\left(p\|q\right)$, that is $\frac{1}{2}\sum_{i}\frac{\left(p_{i}-q_{i}\right)^{2}}{q_{i}}$ (Neyman's Chi2).\\
If $\lambda=-1$ we have $A\left(p\|q\right)=\chi^{2}_{P}\left(p\|q\right)$, that is $\frac{1}{2}\sum_{i}\frac{\left(p_{i}-q_{i}\right)^{2}}{p_{i}}$ (Pearson's Chi2).\\

\section{Bregman's divergences - Beta divergences.}
Since the Bregman divergences are convexity measures, it is obvious that the functions
$f_{1}\left(x\right)$, $f_{2}\left(x\right)$ or $f_{c}\left(x\right)$ which differ from each other only by a linear function (i.e. which have the same second derivative) lead to the same divergence which can be written as follows:
\begin{equation}
	B\left(p\|q\right)=\frac{1}{\lambda\left(\lambda-1\right)}\sum_{i}\left[p_{i}^{\lambda}-\lambda p_{i}q_{i}^{\lambda-1}-\left(1-\lambda\right)q_{i}^{\lambda}\right]
	\label{eq.BC}
\end{equation}
With ($\lambda \leftrightarrow \beta$) we obtain the ``Beta" divergence of Mihoko and Eguchi \cite{mihoko2002} or Uchida and Shioya \cite{uchida2005}; we also find the divergence of BHHJ \cite{basu1998} by deleting in (\ref{eq.BC}) the multiplicative factor ``$1/\lambda$".\\ 
The gradient with respect to ``$q$" will be:
\begin{equation}
	\frac{\partial B\left(p\|q\right)}{\partial q_{j}}=q_{j}^{\lambda-2}\left(q_{j}-p_{j}\right)
\label{eq.gradBC}
\end{equation}
This gradient will be zero if $p_{i}=q_{i}\ \forall i$.\\
From the relationship (\ref{eq.fcv1bis}), the dual divergence is written:
\begin{equation}
	B\left(q\|p\right)=\frac{1}{\lambda\left(\lambda-1\right)}\sum_{i}\left[q_{i}^{\lambda}-\lambda q_{i}p_{i}^{\lambda-1}-\left(1-\lambda\right)p_{i}^{\lambda}\right]
	\label{eq.BCduale}
\end{equation}
Its gradient relative to ``$q$" will be:
\begin{equation}
		\frac{\partial B\left(q\|p\right)}{\partial q_{j}}=\frac{1}{\lambda-1}\left(q_{j}^{\lambda-1}-p_{j}^{\lambda-1}\right)
\end{equation}
Obviously, it will be zero if $p_{i}=q_{i}\ \forall i$.

\subsection{Some classical particular cases.}
If $\lambda\rightarrow 1$ we have $B\left(p\|q\right)=KL\left(p\|q\right)$.\\
This is the common point with "alpha divergences".\\
If $\lambda\rightarrow 0$ we have $B\left(p\|q\right)=IS\left(p\|q\right)$, it is the Itakura-Saito divergence \cite{itakura1968}.\\
If $\lambda=2$ we have $B\left(p\|q\right)=\sum_{i}\left(p_{i}-q_{i}\right)^{2}$, it is the Mean square deviation (Euclidean distance).\\

\section{Jensen's divergences.}
Since Jensen's divergences are measures of convexity, it is obvious that the functions
$f_{1}\left(x\right)$, $f_{2}\left(x\right)$ or $f_{c}\left(x\right)$ which differ from each other only by a linear function (i.e. which have the same second derivative) lead to the same divergence that is written:
\begin{align}
J\left(p\|q\right)=&\frac{1}{\lambda\left(\lambda-1\right)}\left\{\sum_{i}\left[\alpha p^{\lambda}_{i}+\left(1-\alpha\right)q^{\lambda}_{i}\right]\right. \nonumber \\  & \left. -\sum_{i}\left[\alpha p_{i}+\left(1-\alpha\right)q_{i}\right]^{\lambda}\right\}
\label{eq.Jensen}	
\end{align}
The corresponding gradient with respect to ``$q$" is written:
\begin{equation}
	\frac{\partial J\left(p\|q\right)}{\partial q_{j}}=\frac{1-\alpha}{\lambda-1}\left\{q^{\lambda-1}_{j}-\left[\alpha p_{j}+\left(1-\alpha\right)q_{j}\right]^{\lambda-1}\right\}
\end{equation}
This gradient will be zero if $p_{i}=q_{i}\ \forall i$.

\section{Invariance by change of scale.}
We are applying the method defined in \cite{eguchi2010}, in order to make invariant by change of scale on ``$q$" the divergences of type ``$A$", ``$B$" and ``$J$" of the preceding sections.\\
\subsection{``\textit{A}" divergences.}
The starting divergence is given by (\ref{eq.AC}):
\begin{equation}
	A\left(p\|q\right)=\frac{1}{\lambda\left(\lambda-1\right)}\sum_{i}\left(p_{i}^{\lambda}q_{i}^{1-\lambda}-\lambda p_{i}-\left(1-\lambda\right)q_{i}\right)
\end{equation}
The expression of the invariance factor $K$ is obtained explicitly, so it is the nominal value that is written:
\begin{equation}
	K_{0}=\left(\frac{\sum_{i}p_{i}^{\lambda}q_{i}^{1-\lambda}}{\sum_{i}q_{i}}\right)^{\frac{1}{\lambda}}
	\label{eq.K0AI}
\end{equation}
Then the invariant divergence by change of scale on ``$q$" is written:
\begin{equation}	AI\left(p\|q\right)=\frac{1}{\lambda-1}\left[\left(\frac{\sum_{i}p_{i}^{\lambda}q_{i}^{1-\lambda}}{\sum_{i}q_{i}}\right)^{\frac{1}{\lambda}}\sum_{i}q_{i}-\sum_{i}p_{i}\right]
\label{eq.ACIpq}	
\end{equation}
Which can also be written:
\begin{equation}	AI\left(p\|q\right)=\frac{1}{\lambda-1}\left[\sum_{i}K_{0}q_{i}-\sum_{i}p_{i}\right]
\label{eq.ACIpq2}	
\end{equation}
\subsubsection{Example}
If we consider the case $\lambda\rightarrow 1$, that is $KL\left(p\|q\right)$, we obtain:
\begin{equation}
	K_{0}=\frac{\sum_{i}p_{i}}{\sum_{i}q_{i}}
\end{equation}
\begin{equation}	KLI\left(p\|q\right)=\sum_{j}p_{j}\left[\sum_{i}\bar{p}_{i}\log\frac{\bar{p}_{i}}{\bar{q}_{i}}+\bar{q}_{i}-\bar{p}_{i}\right]=\sum_{j}p_{j}\left[\sum_{i}\bar{p}_{i}\log\frac{\bar{p}_{i}}{\bar{q}_{i}}\right]	
\end{equation}
With: $\bar{p}_{i}=\frac{{p}_{i}}{\sum_{j}p_{j}}$ and $\bar{q}_{i}=\frac{{q}_{i}}{\sum_{j}q_{j}}$.\\
That's what we would get if we do directly $\lambda\rightarrow 1$ in the expression (\ref{eq.ACIpq}).
Note that the passage to the limit $\lambda\rightarrow 1$ implies a specific calculation but effectively leads to the  Kullback-Leibler divergence invariant with respect to ``$q$".\\

The gradient of divergence (\ref{eq.ACIpq}) with respect to ``$q$" is written:
\begin{equation}
	\frac{\partial AI\left(p\|q\right)}{\partial q_{j}}=\frac{1}{\lambda}\underbrace{\left(\sum_{i}p_{i}^{\lambda}q_{i}^{1-\lambda}\right)^{\frac{1}{\lambda}}\left(\sum_{i}q_{i}\right)^{1-\frac{1}{\lambda}}}_{\textbf{T}}\,\, \left[\frac{1}{\sum_{i}q_{i}}-\frac{p_{j}^{\lambda}q_{j}^{-\lambda}}{\sum_{i}p_{i}^{\lambda}q_{i}^{1-\lambda}}\right]
	\label{eq.GACIpq}
\end{equation}
It is also possible to write this expression:
\begin{equation}
		\frac{\partial AI\left(p\|q\right)}{\partial q_{j}}=\frac{1}{\lambda}\left[K_{0}-K^{1-\lambda}_{0}p^{\lambda}_{j}q^{-\lambda}_{j}\right]
	\label{eq.GACIpqs}
\end{equation}

\textbf{Remark:} it can be observed that
\begin{equation}
	\sum_{j}q_{j}\frac{\partial AI\left(p\|q\right)}{\partial q_{j}}=0
\end{equation}
This relationship will be verified for all invariant divergences presented in this book.\\

By disregarding the constant multiplicative factor $\frac{1}{\lambda-1}$, the expression (\ref{eq.ACIpq}) is the difference of 2 terms; one can thus apply to each term of this difference the same increasing function, the ``Generalized Logarithm" for example (see appendix 1), then the extreme form, the ``Logarithm" which will lead to:\\
\begin{equation}	LAI\left(p\|q\right)=\frac{1}{\lambda}\log\sum_{i}q_{i}-\frac{1}{\lambda-1}\log\sum_{i}p_{i}+\frac{1}{\lambda\left(\lambda-1\right)}\log\sum_{i}p_{i}^{\lambda}q_{i}^{1-\lambda}
\label{eq.logalfadiv}
\end{equation}
It can be seen, which is true in all cases, that the latter divergence is not only invariant by changing scale on ``$q$" , but also by changing scale on ``$p$". This is obviously related to the use of the Logarithm, (the Generalized Logarithm is not enough).\\
For Type B divergences, the analog of (\ref{eq.logalfadiv}) will be the ``Gamma divergence".\\
The gradient with respect to ``$q$" is written:
\begin{equation}
\frac{\partial LAI\left(p\|q\right)}{\partial q_{j}}=\frac{1}{\lambda}\left[\frac{1}{\sum_{i}q_{i}}-\frac{p_{j}^{\lambda}q_{j}^{-\lambda}}{\sum_{i}p_{i}^{\lambda}q_{i}^{1-\lambda}}\right]	
\end{equation}
It is the expression(\ref{eq.GACIpq}) without the multiplicative factor ``$\textbf{T}$", and we have as previously stated:
\begin{equation}
\sum_{j}q_{j}\frac{\partial LAI\left(p\|q\right)}{\partial q_{j}}=0	
\end{equation}
We can make the same work on the dual divergence (\ref{eq.ACduale}):
\begin{equation}
	A\left(q\|p\right)=\frac{1}{\lambda\left(\lambda-1\right)}\sum_{i}\left(q_{i}^{\lambda}p_{i}^{1-\lambda}-\lambda q_{i}-\left(1-\lambda\right)p_{i}\right)
\end{equation}
By using the classical procedure for calculating the nominal invariance factor, the following results are obtained:
\begin{equation}
	K_{0}=\left(\frac{\sum_{i}p_{i}^{1-\lambda}q_{i}^{\lambda}}{\sum_{i}q_{i}}\right)^{\frac{1}{1-\lambda}}
\end{equation}
The divergence made invariant will be written:
\begin{equation}	AI\left(q\|p\right)=\frac{1}{\lambda}\left[\sum_{i}p_{i}-\left(\frac{\sum_{i}p_{i}^{1-\lambda}q_{i}^{\lambda}}{\sum_{i}q_{i}}\right)^{\frac{1}{1-\lambda}}\sum_{i}q_{i}\right]
\label{eq.ACIqp}
\end{equation}
Or also, more simply:
\begin{equation}	AI\left(q\|p\right)=\frac{1}{\lambda}\left[\sum_{i}p_{i}-\sum_{i}K_{0}q_{i}\right]
\label{eq.ACIqp2}
\end{equation}
One can then apply an increasing function on each of the terms the divergence (\ref{eq.ACIqp}), for example the Generalized Logarithm function, and go further to the Logarithm to obtain:
\begin{align}		LAI\left(q\|p\right)=&\frac{1}{\lambda}\left[\log\sum_{i}p_{i}-\frac{1}{1-\lambda}\log\sum_{i}p_{i}^{1-\lambda}q_{i}^{\lambda}\right. \nonumber \\  & \left. +\frac{\lambda}{1-\lambda}\log\sum_{i}q_{i}\right]
\end{align}
Note that this divergence is, of course, invariant with respect to ``$q$", but it is also invariant with respect to ``$p$".\\
A symmetrical divergence invariant with respect to ``$p$" and ``$q$" can be obtained by summing the (half) logarithmic invariant divergences $LAI\left(p\|q\right)$ and $LAI\left(q\|p\right)$.\\
Such a divergence is written:
\begin{equation}
LAI\left(p,q\right)=\frac{1}{\lambda(1-\lambda)}\left[\log\frac{\sum_{i}p_{i}\sum_{i}q_{i}}{\sum_{i}p^{\lambda}_{i}q^{1-\lambda}_{i}\sum_{i}q^{\lambda}_{i}p^{1-\lambda}_{i}}\right]
\label{eq.logalfadivsym}	
\end{equation}
For type B divergences, the analogue of (\ref{eq.logalfadivsym}) will be the ``symmetrical Gamma divergence ".

\subsection{``\textit{B}" divergences.}
\subsubsection {``Beta'' divergence.}
The initial divergence was (\ref{eq.BC}):
\begin{equation}
	B\left(p\|q\right)=\frac{1}{\lambda\left(\lambda-1\right)}\sum_{i}\left[p_{i}^{\lambda}-\lambda p_{i}q_{i}^{\lambda-1}-\left(1-\lambda\right)q_{i}^{\lambda}\right]
\end{equation}
The expression of the nominal invariance factor $K_{0}$ is calculated explicitly; it is written as:
\begin{equation}
	K_{0}=\frac{\sum_{i}p_{i}q_{i}^{\lambda-1}}{\sum_{i}q_{i}^{\lambda}}
	\label{eq.K0BI}
\end{equation}
This leads to the divergence:
\begin{equation}	BI\left(p\|q\right)=\frac{1}{\lambda\left(\lambda-1\right)}\left[\sum_{i}p_{i}^{\lambda}-\left(\frac{\sum_{i}p_{i}q_{i}^{\lambda-1}}{\sum_{i}q_{i}^{\lambda}}\right)^{\lambda}\sum_{i}q_{i}^{\lambda}\right]
\label{eq.BCIpq}
\end{equation}
Or also:
\begin{equation}	BI\left(p\|q\right)=\frac{1}{\lambda\left(\lambda-1\right)}\left[\sum_{i}p_{i}^{\lambda}-\sum_{i}K^{\lambda}_{0}q_{i}^{\lambda}\right]
\label{eq.BCIpq2}
\end{equation}
To recover the divergence of Basu et al. \cite{basu1998} it is then necessary, in (\ref{eq.BCIpq}), to remove the multiplicative factor $1/\lambda$ and make the parameter change $\lambda=1+\beta$.\\
The gradient of divergence (\ref{eq.BCIpq}) with respect to ``$q$" is written:
\begin{equation}
	\frac{\partial BI\left(p\|q\right)}{\partial q_{j}}=\underbrace{\left(\sum_{i}p_{i}q_{i}^{\lambda-1}\right)^{\lambda}\left(\sum_{i}q_{i}^{\lambda}\right)^{1-\lambda}}_{\textbf{S}}\ \left[\frac{q_{j}^{\lambda-1}}{\sum_{i}q_{i}^{\lambda}}-\frac{p_{j}q_{j}^{\lambda-2}}{\sum_{i}p_{i}q_{i}^{\lambda-1}}\right]
	\label{eq.GBCIpq}
\end{equation}
Or, in a simplified form:
\begin{equation}
	\frac{\partial BI\left(p\|q\right)}{\partial q_{j}}=K^{\lambda}_{0}\left[q^{\lambda-1}_{j}-\left(K_{0}\right)^{-1}p_{j}q^{\lambda-2}_{j}\right]
	\label{eq.GBCIpqs}
\end{equation}
It can be noted that:
\begin{equation}
\sum_{j}q_{j}\frac{\partial BI\left(p\|q\right)}{\partial q_{j}}=0	
\end{equation}

\subsubsection {``Gamma'' divergence.}
At this point, following the usual procedure, we can apply the ``Generalized Logarithm" on each of the terms in the brackets in (\ref{eq.BCIpq}), and then eventually we can move to the Logarithm which leads to the ``Gamma divergence" of Fujisawa and Eguchi \cite{fujisawa2008}. ($\lambda \leftrightarrow \gamma$).:\\
\begin{align}
LBI\left(p\|q\right)=&\frac{1}{\lambda\left(\lambda-1\right)}\left[\log\sum_{i}p_{i}^{\lambda}-\lambda \log\sum_{i}p_{i}q_{i}^{\lambda-1}\right. \nonumber \\  & \left. +\left(\lambda-1\right)\log\sum_{i}q_{i}^{\lambda}\right]
\label{eq.logbetadiv}	
\end{align}
This divergence is invariant by scale change not only on ``$q$", but also on ``$p$".
The gradient with respect to ``$q$" is written:
\begin{equation}
	\frac{\partial LBI\left(p\|q\right)}{\partial q_{j}}=\frac{q_{j}^{\lambda-1}}{\sum_{i}q_{i}^{\lambda}}-\frac{p_{j}q_{j}^{\lambda-2}}{\sum_{i}p_{i}q_{i}^{\lambda-1}}
\end{equation}
This expression is analogous to (\ref{eq.GBCIpq}) except for the multiplicative factor ``$\textbf{S}$" in (\ref{eq.GBCIpq}).\\
We can see that as previously stated:
\begin{equation}
\sum_{j}q_{j}\frac{\partial LBI\left(p\|q\right)}{\partial q_{j}}=0	
\end{equation}

\subsubsection{Dual ``Beta'' divergence.}
We can then come back to the dual divergence (\ref{eq.BCduale}):
\begin{equation}
	B\left(q\|p\right)=\frac{1}{\lambda\left(\lambda-1\right)}\sum_{i}\left[q_{i}^{\lambda}-\lambda q_{i}p_{i}^{\lambda-1}-\left(1-\lambda\right)p_{i}^{\lambda}\right]
\end{equation}
The expression of $K_{0}$ allowing to get the invariance is computed explicitly; we have:
\begin{equation}
	K_{0}=\left(\frac{\sum_{i}p_{i}^{\lambda-1}q_{i}}{\sum_{i}q_{i}^{\lambda}}\right)^{\frac{1}{\lambda-1}}	
\end{equation}
The corresponding invariant divergence will be written:
\begin{equation}
BI\left(q\|p\right)=\frac{1}{\lambda}\left[\sum_{i}p_{i}^{\lambda}-\left(\frac{\sum_{i}q_{i}p_{i}^{\lambda-1}}{\sum_{i}q_{i}^{\lambda}}\right)^{\frac{\lambda}{\lambda-1}}\sum_{i}q_{i}^{\lambda}\right]
\end{equation}
Or, more simply:
\begin{equation}	
BI\left(q\|p\right)=\frac{1}{\lambda}\left[\sum_{i}p_{i}^{\lambda}-\sum_{i}K^{\lambda}_{0}q_{i}^{\lambda}\right]
\end{equation}

\subsubsection {Dual ``Gamma'' divergence.}
From the previous expression, we obtain the Logarithmic form, i.e. the dual Gamma divergence:
\begin{equation}	LBI\left(q\|p\right)=\frac{1}{\lambda}\left[\log\sum_{i}p_{i}^{\lambda}-\frac{1}{1-\lambda}\log\sum_{i}q_{i}^{\lambda}+\frac{\lambda}{1-\lambda}\log\sum_{i}q_{i}p_{i}^{\lambda-1}\right]
\end{equation}
Its gradient with respect to ``$q$" is written:
\begin{equation}
	\frac{\partial LBI\left(q\|p\right)}{\partial q_{j}}=\frac{1}{\lambda-1}\left[\frac{q^{\lambda-1}_{j}}{\sum_{i}q^{\lambda}_{i}}-\frac{p^{\lambda-1}_{j}}{\sum_{i}q_{i}p^{\lambda-1}_{i}}\right]
\end{equation}

\subsubsection {Symmetrical ``Gamma'' divergence.}
It is immediately obtained as the (half) sum of the Gamma divergence and the dual Gamma divergence; it is written as follows:
\begin{equation}
LBI\left(p,q\right)=\frac{1}{\lambda-1}	\log\frac{\sum_{i}p^{\lambda}_{i}\sum_{i}q^{\lambda}_{i}}{\sum_{i}p_{i}q^{\lambda-1}_{i}{\sum_{i}q_{i}p^{\lambda-1}_{i}}}
\end{equation}
Its gradient with respect to ``$q$" is written:
\begin{equation}
	\frac{\partial LBI\left(p,q\right)}{\partial q_{j}}=\frac{\lambda}{\lambda-1}\frac{q^{\lambda-1}_{j}}{\sum_{i}q^{\lambda}_{i}}-\frac{p_{j}q^{\lambda-2}_{j}}{\sum_{i}p_{i}q^{\lambda-1}_{i}}-\frac{1}{\lambda-1}\frac{p^{\lambda-1}_{j}}{\sum_{i}q_{i}p^{\lambda-1}_{i}}
\end{equation}

\subsection{``J" divergences.}
We write:
\begin{align}
J\left(p\|Kq\right)=&\frac{1}{\lambda\left(\lambda-1\right)}\left\{\sum_{i}\left[\alpha p^{\lambda}_{i}+\left(1-\alpha\right)K^{\lambda}q^{\lambda}_{i}\right]\right. \nonumber \\  & \left. -\sum_{i}\left[\alpha p_{i}+\left(1-\alpha\right)Kq_{i}\right]^{\lambda}\right\}	
\end{align}
We deduce that:
\begin{equation}
	\frac{\partial J\left(p\|Kq\right)}{\partial K}=\frac{1-\alpha}{\lambda-1}\sum_{i}q_{i}\left\{K^{\lambda-1}q_{i}^{\lambda-1}-\left[\alpha p_{i}+\left(1-\alpha\right)K q_{i}\right]^{\lambda-1}\right\}
\end{equation}
The method proposed in \textbf{Chapter 3}, which consists in obtaining the invariance factor $K$ by solving the equation:
\begin{equation}
\frac{\partial J\left(p\|Kq\right)}{\partial K}=0	
\end{equation}
as proposed in (\ref{eq.resolK}), does not lead to an explicit solution.\\
On the other hand, if we choose to use the invariance factor:
\begin{equation}
 K=\frac{\sum_{i}p_{i}}{\sum_{i}q_{i}}	
\end{equation}
We obtain an invariant divergence which is written as follows:
\begin{align}
JI\left(p\|q\right)=&\frac{\left(\sum_{j}p_{j}\right)^{\lambda}}{\lambda\left(\lambda-1\right)}\left\{\sum_{i}\left[\alpha \overline{p}^{\lambda}_{i}+\left(1-\alpha\right)\overline{q}^{\lambda}_{i}\right]\right. \nonumber \\  & \left. -\sum_{i}\left[\alpha \overline{p}_{i}+\left(1-\alpha\right)\overline{q}_{i}\right]^{\lambda}\right\}	
\end{align}
with: $\bar{p}_{i}=\frac{{p}_{i}}{\sum_{j}p_{j}}$ and $\bar{q}_{i}=\frac{{q}_{i}}{\sum_{j}q_{j}}$.\\

Its gradient with respect to ``$q$" is written as follows:
\begin{align}
	\frac{\partial JI\left(p\|q\right)}{\partial q_{j}}&=\frac{\left(1-\alpha\right)}{\left(\lambda-1\right)}\frac{\left(\sum_{l}p_{l}^{\lambda}\right)}{\sum_{l}q_{l}}\left\{\left[\sum_{i}\bar{q}_{i}\left[\alpha \bar{p}_{i}+\left(1-\alpha\right)\bar{q}_{i}\right]^{\lambda-1}-\bar{q}^{\lambda}_{i}\right] \right. \nonumber \\ &\left.	-\left[\left[\alpha\bar{p}_{j}+\left(1-\alpha\right)\bar{q}_{j}\right]^{\lambda-1}-\bar{q}^{\lambda-1}_{j}\right]\right\}	
\end{align}
Of course, we still have the fundamental relationship:
\begin{equation}
	\sum_{j}q_{j}\frac{\partial JI\left(p\|q\right)}{\partial q_{j}}=0	
\end{equation}
\subsection{Important note on invariance.}
In the $AI\left(p\|q\right)$, $BI\left(p\|q\right)$ and $JI\left(p\|q\right)$ divergences, the invariance was established with respect to the variable ``$q$"; moreover, it was noted that these divergences are in the form of a difference of 2 positive terms.\\
 Based on this observation, an increasing function can be applied to each of the terms of the difference without changing the sign of the divergence (positive in this case); moreover, since both terms are positive, the use of the Logarithm function is possible, the use of the Generalized Logarithm being only an intermediate step.\\
We can now make an additional remark if we are considering the invariance with respect to ``$p$": we can observe that if ``$p$" is multiplied by $K$ in $AI\left(p\|q\right)$, the 2 terms of the divergence are multiplied by $K$, and only the application of a logarithmic function on each of the terms of the difference will allow the factor $K$ to vanish.\\
 Similarly, if in $BI\left(p\|q\right)$, we multiply ``$p$" by $K$, each of the terms of the difference will be multiplied by $K^{\lambda}$ and, again, only the application of the logarithmic function on each of the 2 terms of the difference will result in the invariance with respect to ``$p$", removing the $K$ factor.
 
\subsection{A particular case of the invariance factor.}
As reported in \textbf{Chapter 3}, for certain divergences, the method for calculating the invariance factor $K\left(p,q\right)$ by resolution of (\ref{eq.DiffK}) does not lead to an explicit solution, however, an expression of the form:
\begin{equation}
	K\left(p,q\right)=\frac{\sum_{j}p_{j}}{\sum_{j}q_{j}}
	\label{eq.Kbis}
\end{equation}
allows to obtain a scale change invariant form of any divergence.\\
As previously proposed for Jensen's divergences, we will introduce this factor into the ``Alpha" and ``Beta" divergences (\ref{eq.AC}) (\ref{eq.BC}) and compare the resulting expressions with the forms proposed in (\ref{eq.ACIpq}) and (\ref{eq.BCIpq}).\\
If we put (\ref{eq.Kbis}) in (\ref{eq.AC}), with the notations $\bar{p}_{i}=\frac{p_{i}}{\sum_{j}p_{j}}$ and $\bar{q}_{i}=\frac{q_{i}}{\sum_{j}q_{j}}$,, we obtain:

\begin{equation}	AI_{bis}\left(p\|q\right)=\frac{\sum_{j}p_{j}}{\lambda\left(\lambda-1\right)}\sum_{i}\left(\bar{p}_{i}^{\lambda}\bar{q}_{i}^{1-\lambda}-\lambda \bar{p}_{i}-\left(1-\lambda\right)\bar{q}_{i}\right)
\end{equation}
Which can be simplified in the following way:
\begin{equation}	AI_{bis}\left(p\|q\right)=\frac{\sum_{j}p_{j}}{\lambda\left(\lambda-1\right)}\left[\left(\sum_{i}\bar{p}_{i}^{\lambda}\bar{q}_{i}^{1-\lambda}\right)-1\right]
\end{equation}
Similarly, if we put (\ref{eq.Kbis}) in (\ref{eq.BC}), we obtain:
\begin{equation}	BI_{bis}\left(p\|q\right)=\frac{\left(\sum_{j}p_{j}\right)^{\lambda}}{\lambda\left(\lambda-1\right)}\sum_{i}\left[\bar{p}_{i}^{\lambda}-\lambda \bar{p}_{i}\bar{q}_{i}^{\lambda-1}-\left(1-\lambda\right)\bar{q}_{i}^{\lambda}\right]
\end{equation}
These 2 divergences are invariant by change of scale on ``$q$", moreover, if we disregard the multiplicative factor which depends only on ``$\sum_{j}p_{j}$", they are also invariant with respect to ``$p$".\\
Their respective gradients with respect to ``$q$" are written as:
\begin{equation}
	\frac{\partial AI_{bis}\left(p\|q\right)}{\partial q_{j}}=\frac{\sum_{l}p_{l}}{\lambda \sum_{l}q_{l}}\left[\left(\sum_{i}\bar{p}^{\lambda}_{i}\bar{q}^{1-\lambda}_{i}\right)-\bar{p}^{\lambda}_{j}\bar{q}^{\lambda}_{j}\right]
\end{equation}
and:
\begin{equation}
	\frac{\partial BI_{bis}\left(p\|q\right)}{\partial q_{j}}=\frac{\left(\sum_{l}p_{l}\right)^{\lambda}}{\sum_{l}q_{l}}\left[\left(\sum_{i}\bar{p}_{i}\bar{q}^{\lambda-1}_{i}-\bar{q}^{\lambda}_{i}\right)-\left(\bar{p}_{j}\bar{q}^{\lambda-2}_{j}-\bar{q}^{\lambda-1}_{j}\right)\right]	
\end{equation}
It can be observed that we have always the fundamental relationship typical of scale invariant divergences:
\begin{equation}
\sum_{j}q_{j}\frac{\partial AI_{bis}\left(p\|q\right)}{\partial q_{j}}=0	
\end{equation}
and
\begin{equation}
\sum_{j}q_{j}\frac{\partial BI_{bis}\left(p\|q\right)}{\partial q_{j}}=0		
\end{equation}
 
\section{Generalization in the sense of Cichocki and Ghosh.}
\subsection{Generalized divergence.}
In \cite{cichocki2011}, Cichocki proposes a general writing which allows to have a unique expression for the ``Alpha divergences " and for the ``Beta divergences "; of course, even if the writing is unique, the two remain distinct as we will see.\\
We first introduce the two parameters $a\equiv\lambda_{\alpha}$ and $b\equiv\lambda_{\beta}$ specific to the two types of divergence and the expression of the generalized divergence will be:

\begin{align}	AB\left(p\|q\right)=&-\frac{1}{a\left(b-1\right)}\left[\sum_{i}p_{i}^{a}q_{i}^{b-1}-\frac{a}{a+b-1}\sum_{i}p_{i}^{a+b-1} \right. \nonumber \\ &\left.-\frac{b-1}{a+b-1}\sum_{i}q_{i}^{a+b-1}\right]
\label{eq.AB}
\end{align}
So the Alpha and Beta divergences can be recovered in the following way:\\
- If $a+b-1=1$, that is $b-1=1-a$, we have the ``Alpha divergence".\\
- If $a=1$, that only leaves $b$, we get the ``Beta divergence".\\
To recover Cichocki's writing exactly, we have to do $a\equiv \alpha$ and $b\equiv \beta+1$.\\
Ghosh \cite{ghosh2013} suggests something similar. To retrieve our result, we just have to introduce in Ghosh's text, a multiplicative factor $\frac{1}{A+B}$, then to make the replacement: $A\equiv a$ and $B\equiv b-1$ (see appendix 3).

\subsection{Invariance by change of scale on ``$q$''.}
To make this generalized divergence invariant by change of scale, we apply on the expression (\ref{eq.AB}) the method allowing to determine the factor $K_{0}$ and we obtain:
\begin{equation}
	K_{0}=\left[\frac{\sum_{i}p_{i}^{a}q_{i}^{b-1}}{\sum_{i}q_{i}^{a+b-1}}\right]^{\frac{1}{a}}
\end{equation}
The particular cases corresponding to ``Alpha" and ``Beta" divergences are immediately revealed:\\
If $a+b-1=1$, that is if $b-1=1-a$, we have:
\begin{equation}
	K_{0,a}=\left[\frac{\sum_{i}p_{i}^{a}q_{i}^{1-a}}{\sum_{i}q_{i}}\right]^{\frac{1}{a}}
\end{equation}
If $a=1$ we obtain:
\begin{equation}
	K_{0,b}=\frac{\sum_{i}p_{i}q_{i}^{b-1}}{\sum_{i}q_{i}^{b}}
\end{equation}
This being so, the generalized divergence made invariant to the change of scale is then written down, all calculations done:
\begin{align}
	ABI\left(p\|q\right)=&\frac{1}{\left(b-1\right)\left(a+b-1\right)}\left\{\sum_{i}p_{i}^{a+b-1} \right. \nonumber \\ &\left.-\left(\sum_{i}p_{i}^{a}q_{i}^{b-1}\right)^{\frac{a+b-1}
	{a}}\left(\sum_{i}q_{i}^{a+b-1}\right)^{\frac{1-b}{a}}\right\}
	\label{eq.ABI}
\end{align}
The calculation of the gradient with respect to ``$q$" leads to:
\begin{equation}
-\frac{\partial ABI\left(p\|q\right)}{\partial q_{j}}=\frac{1}{a}\frac{\left(\sum_{i}p_{i}^{a}q_{i}^{b-1}\right)^{\frac{a+b-1}{a}}}{\left(\sum_{i}q_{i}^{a+b-1}\right)^{\frac{b-1}{a}}}\left[\frac{p_{j}^{a}q_{j}^{b-2}}{\sum_{i}p_{i}^{a}q_{i}^{b-1}}-\frac{q_{j}^{a+b-2}}{\sum_{i}q_{i}^{a+b-1}}\right]
\end{equation}
We still have the same relationship:
\begin{equation}
	\sum_{j}q_{j}\frac{\partial ABI\left(p\|q\right)}{\partial q_{j}}=0
\end{equation}
Once again, the two special cases are recovered.

\subsection{Logarithmic form.}
Disregarding the multiplicative factor, we have always a difference of two positive terms, therefore the applying the generalized logarithm or the logarithm on each of these terms leads to:
\begin{align}
		LABI\left(p\|q\right)=
		&\frac{1}{\left(b-1\right)\left(a+b-1\right)}\left\{\log\sum_{i}p_{i}^{a+b-1}-\right.\nonumber \\
		&\frac{a+b-1}{a}\log \sum_{i}p_{i}^{a}q_{i}^{b-1}-\nonumber \\
		&\left.\frac{1-b}{a}\log \sum_{i}q_{i}^{a+b-1}\right\}
		\label{eq.LABI}
\end{align}
If we compute the gradient with respect to ``$q$", we obtain:
\begin{equation}
-\frac{\partial LABI\left(p\|q\right)}{\partial q_{j}}=\frac{1}{a}\left[\frac{p_{j}^{a}q_{j}^{b-2}}{\sum_{i}p_{i}^{a}q_{i}^{b-1}}-\frac{q_{j}^{a+b-2}}{\sum_{i}q_{i}^{a+b-1}}\right]
\end{equation}
As it is always the case with invariant divergences, we have:
\begin{equation}
	\sum_{j}q_{j}\frac{\partial LABI\left(p\|q\right)}{\partial q_{j}}=0
\end{equation}
The two special cases are retrieved:\\
- If $a+b-1=1$, we have:
\begin{equation}
-\frac{\partial LABI\left(p\|q\right)}{\partial q_{j}}=\frac{1}{a}\left[\frac{p_{j}^{a}q_{j}^{-a}}{\sum_{i}p_{i}^{a}q_{i}^{1-a}}-\frac{1}{\sum_{i}q_{i}}\right]\equiv -\frac{\partial LAI\left(p\|q\right)}{\partial q_{j}}	
\end{equation}
- If $a=1$, we have:
\begin{equation}
-\frac{\partial LABI\left(p\|q\right)}{\partial q_{j}}=\frac{p_{j}q_{j}^{b-2}}{\sum_{i}p_{i}q_{i}^{b-1}}-\frac{q_{j}^{b-1}}{\sum_{i}q_{i}^{b}}\equiv -\frac{\partial LBI\left(p\|q\right)}{\partial q_{j}}	
\end{equation}

\subsubsection{Another way to implement invariance.}
From the generalized divergence (\ref{eq.AB}), we introduce the invariance factor $K\left(p,q\right)=\frac{\sum_{j}p_{j}}{\sum_{j}q_{j}}$ and we reach another form of the generalized divergence invariant by a scale change which is written:

\begin{align}	ABI_{bis}\left(p\|q\right)=&-\frac{\left(\sum_{l}p_{l}\right)^{a+b-1}}{a\left(b-1\right)}\left[\sum_{i}\bar{p}_{i}^{\left(a\right)}\bar{q}_{i}^{\left(b-1\right)}-\frac{a}{a+b-1}\sum_{i}\bar{p}_{i}^{\left(a+b-1\right)} \right. \nonumber \\ &\left.-\frac{b-1}{a+b-1}\sum_{i}\bar{q}_{i}^{\left(a+b-1\right)}\right]
\label{eq.ABIbis}
\end{align}
Its gradient with respect to ``$q$" is written:
\begin{align}
	\frac{\partial ABI_{bis}\left(p\|q\right)}{\partial q_{j}}=&\frac{(\sum_{l}p_{l})^{a+b-1}}{a \sum_{l}q_{l}}\left(\sum_{i}\left(\bar{p_{i}}\right)^{a}\left(\bar{q_{i}}\right)^{b-1}-\sum_{i}\left(\bar{q_{i}}\right)^{a+b-1} \right. \nonumber \\ &\left.+\left(\bar{q_{j}}\right)^{a+b-2}-\left(\bar{p_{j}}\right)^{a}\left(\bar{q_{j}}\right)^{b-2}\right)
\end{align}
We can verify that:
\begin{equation}
	\sum_{j}q_{j}\frac{\partial ABI_{bis}\left(p\|q\right)}{\partial q_{j}}=0
\end{equation}	

\subsection{Transition to Renyi's divergence.}
The divergence (\ref{eq.LABI}) can be written in an alternative way:
\begin{align}
		LABI\left(p\|q\right)=
		&\frac{1}{a+b-1}\left\{\frac{1}{b-1}\log\sum_{i}p_{i}^{a+b-1}-\right.\nonumber \\
		&\frac{a+b-1}{a\left(b-1\right)}\log \sum_{i}p_{i}^{a}q_{i}^{b-1}+\nonumber \\
		&\left.\frac{1}{a}\log \sum_{i}q_{i}^{a+b-1}\right\}
\end{align}
To recover Renyi's divergence, one must first make $a+b-1=1$, i.e. $b-1=1-a$, with $a\equiv\alpha$ to recover Renyi's usual notation, and one obtains something which is the invariant version of the ``$\alpha$ divergence" in logarithmic form:
\begin{align}
		LAI_{\alpha}\left(p\|q\right)=&\frac{1}{1-\alpha}\log\sum_{i}p_{i} \nonumber \\ &-\frac{1}{\alpha\left(1-\alpha\right)}\log \sum_{i}p_{i}^{\alpha}q_{i}^{1-\alpha}+\frac{1}{\alpha}\log \sum_{i}q_{i}
		\label{eq.LAIa}
\end{align}
This expression is of course also valid when the proportionality factor $K$ is equal to 1; we can remark that the relation (\ref{eq.LAIa}) can also be written:
\begin{equation}
		LABI\left(p\|q\right)=\frac{1}{\alpha\left(\alpha-1\right)}\log \sum_{i}\left(\bar{p}_{i}\right)^{\alpha}\left(\bar{q}_{i}\right)^{1-\alpha}
\end{equation}
With $\bar{p}_{i}=\frac{p_{i}}{\sum_{j}p_{j}}$ and $\bar{q}_{i}=\frac{q_{i}}{\sum_{j}q_{j}}$.\\
This is a Renyi divergence with normalized variables.\\
Now, if we assume we're dealing with probability densities, $\sum_{i}p_{i}=\sum_{i}q_{i}=1$ we get exactly the Renyi form:
\begin{equation}
		R\left(p\|q\right)=LABIS\left(p\|q\right)=\frac{1}{\alpha\left(\alpha-1\right)}\log \sum_{i}p_{i}^{\alpha}q_{i}^{1-\alpha}
\end{equation}
Once we have considered that $\sum_{i}p_{i}=\sum_{i}q_{i}$ it is of course senseless to consider a possibility of invariance, even if these sums are different from $1$.

\subsection{Dual Generalized Divergence.}
From the expression (\ref{eq.AB}), the generalized dual divergence is written:
\begin{align}
	AB\left(q\|p\right)=&-\frac{1}{a\left(b-1\right)}\left[\sum_{i}q_{i}^{a}p_{i}^{b-1} \right. \nonumber \\ &\left.-\frac{a}{a+b-1}\sum_{i}q_{i}^{a+b-1}-\frac{b-1}{a+b-1}\sum_{i}p_{i}^{a+b-1}\right]
\end{align}
From there, and according to the usual method, one can derive the invariance factor $K_{0}\left(p,q\right)$ to be inserted in $AB\left(K_{0}q\|p\right)$ to obtain an invariant form of this divergence; one thus obtains:
\begin{equation}
	K_{0}\left(p,q\right)=\left[\frac{\sum_{i}q_{i}^{a}p_{i}^{b-1}}{\sum_{i}q_{i}^{a+b-1}}\right]^{\frac{1}{b-1}}
\end{equation}
After inserting this expression in $AB\left(K_{0}q\|p\right)$, one obtains the invariant dual generalized divergence that one can write down all simplifications made:

\begin{equation}	ABI\left(q\|p\right)=\frac{1}{a\left(a+b-1\right)}\left\{\sum_{i}p_{i}^{a+b-1}-\left[\frac{\sum_{i}q_{i}^{a}p_{i}^{b-1}}{\sum_{i}q_{i}^{a+b-1}}\right]^{\frac{a}{b-1}}\sum_{i}q_{i}^{a}p_{i}^{b-1}\right\}
\end{equation}
Which can also be written in synthetic form:
\begin{equation}	ABI\left(q\|p\right)=\frac{1}{a\left(a+b-1\right)}\left\{\sum_{i}p_{i}^{a+b-1}-K^{a}_{0}\sum_{i}q_{i}^{a}p_{i}^{b-1}\right\}
\end{equation}

We can verify that this divergence is invariant with respect to ``$q$".\\
The calculation of the gradient leads to:
\begin{align}
	\frac{\partial ABI\left(q\|p\right)}{\partial q_{j}}=-\frac{1}{b-1}&\left\{\left[\frac{\sum_{i}q_{i}^{a}p_{i}^{b-1}}{\sum_{i}q_{i}^{a+b-1}}\right]^{\frac{a}{b-1}}q_{j}^{a-1}p_{j}^{b-1}\right. \nonumber \\  & \left. -\left[\frac{\sum_{i}q_{i}^{a}p_{i}^{b-1}}{\sum_{i}q_{i}^{a+b-1}}\right]^{\frac{a+b-1}{b-1}}q_{j}^{a+b-2}\right\}
	\label{eq.GABId}
\end{align}
This expression can also be written:
\begin{align}
	\frac{\partial ABI\left(q\|p\right)}{\partial q_{j}}=-\frac{1}{b-1}\frac{\left(\sum_{i}q_{i}^{a}p_{i}^{b-1}\right)^{\frac{a+b-1}{b-1}}}{\left(\sum_{i}q_{i}^{a+b-1}\right)^{\frac{a}{b-1}}}&\left[\frac{1}{\sum_{i}q_{i}^{a}p_{i}^{b-1}}q_{j}^{a-1}p_{j}^{b-1}\right. \nonumber \\  & \left. -\frac{1}{\sum_{i}q_{i}^{a+b-1}}q_{j}^{a+b-2}\right]
		\label{eq.GABId2}
\end{align}
We can note that we still have:
\begin{equation}
	\sum_{j}q_{j}\frac{\partial ABI\left(q\|p\right)}{\partial q_{j}}=0
\end{equation}

\subsubsection{Logarithmic form.}
By applying the logarithm on the two terms of the difference that appears in (\ref{eq.GABId}), we obtain the corresponding logarithmic form which is written for example:
\begin{align}	LABI\left(q\|p\right)=\frac{1}{a\left(a+b-1\right)}&\left\{\log\sum_{i}p_{i}^{a+b-1}+\frac{a}{b-1}\log\sum_{i}q_{i}^{a+b-1} \right. \nonumber \\ &\left.-\frac{a+b-1}{b-1}\log \sum_{i}q_{i}^{a}p_{i}^{b-1}\right\}
\end{align}
This divergence is invariant not only with respect to ``$q$" but also with respect to ``$p$."\\
The calculation of the gradient leads to:
\begin{equation}
	\frac{\partial LABI\left(q\|p\right)}{\partial q_{j}}=-\frac{1}{b-1}\left\{\frac{1}{\sum_{i}q_{i}^{a}p_{i}^{b-1}}q_{j}^{a-1}p_{j}^{b-1}-\frac{1}{\sum_{i}q_{i}^{a+b-1}}q_{j}^{a+b-2}\right\}
\end{equation}
It can be seen that this expression is immediately deduced from (\ref{eq.GABId2}) by deleting in the latter, the constant multiplicative factor:
\begin{equation}
	\frac{\left(\sum_{i}q_{i}^{a}p_{i}^{b-1}\right)^{\frac{a+b-1}{b-1}}}{\left(\sum_{i}q_{i}^{a+b-1}\right)^{\frac{a}{b-1}}}
\end{equation}
That observation is always true.\\
It is worth noting that we still have:
\begin{equation}
	\sum_{j}q_{j}\frac{\partial LABI\left(q\|p\right)}{\partial q_{j}}=0
\end{equation}

\subsection{Remark.}

The results of the previous section can be found from the work of Ghosh et al.\cite{ghosh2013}.\\
Indeed, we can consider that the divergence analyzed in \cite{ghosh2013} is a generalization of that of BHHJ \cite{basu1998}; it is therefore in fact the Beta divergence (type B divergence)(\ref{eq.BC}).\\
The details for making these connections are given in  \textbf{Appendix 3}.

\setcounter{table}{0}  \setcounter{equation}{0}  \setcounter{figure}{0} \setcounter{chapter}{6} \setcounter{section}{0} 
\chapter{chapter 6 -\\Other classical divergences.}  \label{chptr::chapitre6}
\section{CHI2 ($\chi^{2}$) divergences.}
\subsection{Neyman's CHI2 divergence.}
\begin{figure}[h!]
\centering
\includegraphics[width=0.7\linewidth]{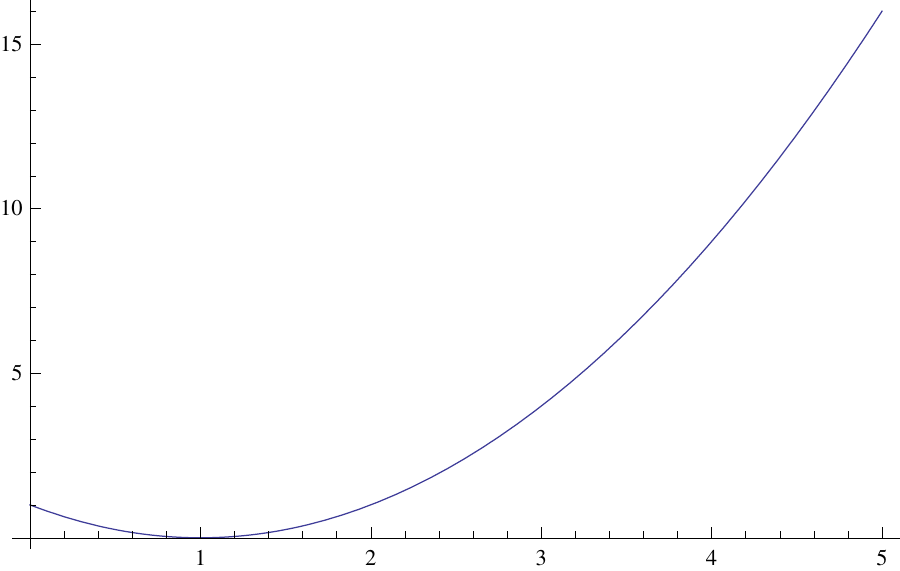}
\caption{Function $f_{c}\left(x\right)=\left(x-1\right)^{2}$}
\label{fig:CHI2Nc}
\end{figure}

\begin{figure}[h!]
\centering
\includegraphics[width=0.7\linewidth]{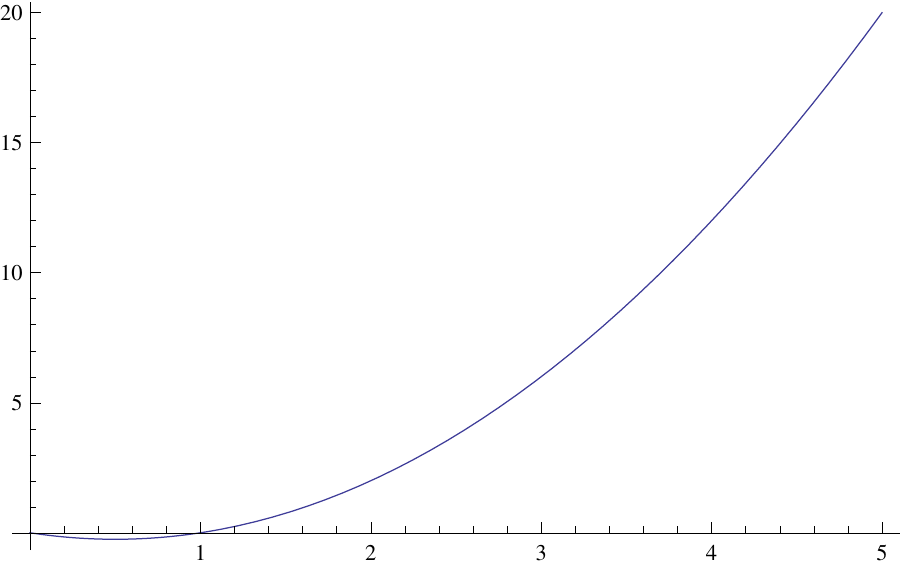}
\caption{Function $f_{1}\left(x\right)=x^{2}-x$}
\label{fig:CHI2N}
\end{figure}

This divergence is also called W Kagan's divergence \cite{basseville1989}, \cite{basseville1996}, \cite{basseville2013}; it is a Csiszär divergence based on the standard convex function v.fig.(\ref{fig:CHI2Nc}):
\begin{equation}
	f_{c}\left(x\right)=\left(x-1\right)^{2}
	\label{eq.CHI2Nc}
\end{equation}
One obtains immediately the ``$\chi_{N}^{2}$" divergence of Neymann \cite{neyman1949}:
\begin{equation}
	\chi_{N}^{2}\left(p\|q\right)=\sum_{i}\frac{\left(p_{i}-q_{i}\right)^{2}}{q_{i}}
	\label{eq.chi2Nc}
\end{equation}
The gradient with respect to ``$q$" is written:
\begin{equation}
\frac{\partial \chi_{N}^{2}\left(p\|q\right)}{\partial q_{j}}=1-\frac{p_{j}^{2}}{q_{j}^{2}}
\end{equation}
If now we introduce in (\ref{eq.chi2Nc}) a first simplification: $\sum_{i}p_{i}=\sum_{i}q_{i}$, we can obtain the divergence:
\begin{equation}
	\chi_{N}^{2}\left(p\|q\right)=\sum_{i}\frac{p_{i}^{2}}{q_{i}}-p_{i}
\end{equation}
It is derived from the simple convex function see fig.(\ref{fig:CHI2N}):
\begin{equation}
	f_{1}\left(x\right)=x^{2}-x
\end{equation}
But we can also obtain the divergence:
\begin{equation}
	\chi_{N}^{2}\left(p\|q\right)=\sum_{i}\frac{p_{i}^{2}}{q_{i}}-q_{i}
\end{equation}
which is based on the simple convex function:
\begin{equation}
	f_{2}\left(x\right)=x^{2}-1
\end{equation}
In both cases, an additional simplification $\sum_{i}p_{i}=\sum_{i}q_{i}=1$ allows to obtain:
\begin{equation}
	\chi_{N}^{2}\left(p\|q\right)=\left(\sum_{i}\frac{p_{i}^{2}}{q_{i}}\right)-1
\end{equation}
An invariant form of the divergence (\ref{eq.chi2Nc}) can be obtained; indeed, the derivation of the nominal invariance factor leads to the explicit solution:
\begin{equation}
	K_{0}\left(p,q\right)=\sqrt{\frac{\sum_{i}p^{2}_{i}q^{-1}_{i}}{\sum_{i}q_{i}}}
\end{equation}
By introducing this expression into (\ref{eq.chi2Nc}), we obtain the invariant divergence with respect to ``$q$":
\begin{equation}	\chi_{N}^{2}I\left(p\|q\right)=2\left[\left(\sum_{i}\frac{p^{2}_{i}}{q_{i}}\sum_{i}q_{i}\right)^{\frac{1}{2}}-\sum_{i}p_{i}\right]=2\left[K_{0}\sum_{i}q_{i}-\sum_{i}p_{i}\right]
\label{eq.Chi2NI}
\end{equation}
The gradient with respect to ``$q$" is written as:
\begin{equation}
	\frac{\partial \chi_{N}^{2}I\left(p\|q\right)}{\partial q_{j}}=\left(\frac{\sum_{i}\frac{p^{2}_{i}}{q_{i}}}{\sum_{i}q_{i}}\right)^{\frac{1}{2}}\left[1-\frac{\sum_{i}q_{i}}{\sum_{i}\frac{p^{2}_{i}}{q_{i}}}\frac{p^{2}_{j}}{q^{2}_{j}}\right]
\end{equation}
Starting from the expression (\ref{eq.Chi2NI}), we obtain the logarithmic form:
\begin{equation}
L\chi_{N}^{2}I\left(p\|q\right)=\log\sum_{i}\frac{p^{2}_{i}}{q_{i}}-\log\frac{\left(\sum_{i}p_{i}\right)^{2}}{\left(\sum_{i}q_{i}\right)}	
\end{equation}
Note that this divergence is not only invariant with respect to ``$q$", but also with respect to ``$p$".
The corresponding gradient with respect to ``$q$" is written:
\begin{equation}
	\frac{\partial L\chi_{N}^{2}I\left(p\|q\right)}{\partial q_{j}}=\left(\frac{1}{\sum_{i}q_{i}}\right)\left[1-\frac{\sum_{i}q_{i}}{\sum_{i}\frac{p^{2}_{i}}{q_{i}}}\frac{p^{2}_{j}}{q^{2}_{j}}\right]
\end{equation}
\subsection{Pearson's Chi2 divergence.}
\begin{figure}[h!]
\centering
\includegraphics[width=0.7\linewidth]{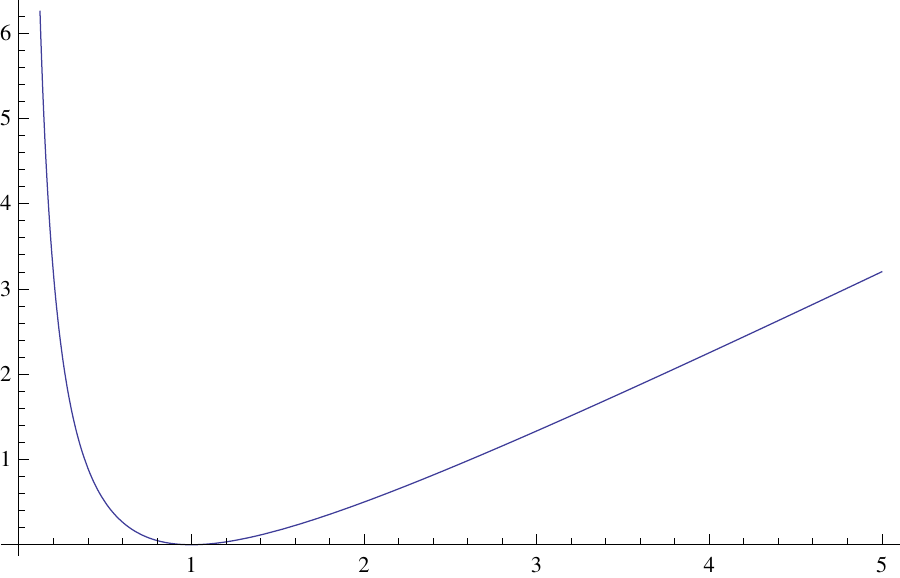}
\caption{Function $\breve{f}_{c}\left(x\right)=\frac{\left(x-1\right)^{2}}{x}$}
\label{fig:CHI2Pc}
\end{figure}

The dual divergence $\chi_{N}^{2}\left(q\|p\right)$ is based on the standard convex function ``$\breve{f}_{c}$'' see fig.(\ref{fig:CHI2Pc}), the mirror function of ``$f_{c}$'' (\ref{eq.CHI2Nc}):
\begin{equation}
	\breve{f}_{c}\left(x\right)=\frac{\left(x-1\right)^{2}}{x}
\end{equation}
This function diverges at its origin.\\
The divergence we obtained is the Pearson's $\chi_{P}^{2}$ written:
\begin{equation}
	\chi_{P}^{2}\left(p\|q\right)=\chi_{N}^{2}\left(q\|p\right)=\sum_{i}\frac{\left(q_{i}-p_{i}\right)^{2}}{p_{i}}
	\label{eq.Pearson}
\end{equation}
The gradient with respect to ``$q$" is written as:
\begin{equation}
\frac{\partial \chi_{P}^{2}\left(p\|q\right)}{\partial q_{j}}=2\left(\frac{q_{j}}{p_{j}}-1\right)
\end{equation}
If now we introduce a first simplification: $\sum_{i}p_{i}=\sum_{i}q_{i}$, we can obtain the divergence:
\begin{equation}
	\chi_{P}^{2}\left(p\|q\right)=\sum_{i}\frac{q_{i}^{2}}{p_{i}}-q_{i}
\end{equation}
This divergence is built in the sense of Csiszär on the simple convex function:
\begin{equation}
	f\left(x\right)=\frac{1}{x}-1
\end{equation}
The extreme simplification $\sum_{i}p_{i}=\sum_{i}q_{i}=$1, yields:
\begin{equation}
	\chi_{P}^{2}\left(p\|q\right)=\left(\sum_{i}\frac{q_{i}^{2}}{p_{i}}\right)-1
\end{equation}
An invariant form of the divergence (\ref{eq.Pearson}) can be obtained; indeed, the calculation of the nominal invariance factor leads to the explicit solution:
\begin{equation}
	K_{0}\left(p,q\right)=\frac{\sum_{i}q_{i}}{\sum_{i}q^{2}_{i}p^{-1}_{i}}
\end{equation}
By introducing this expression into (\ref{eq.Pearson}), we obtain the invariant divergence with respect to ``$q$":
\begin{equation}
\chi_{P}^{2}I\left(p\|q\right)=\sum_{i}p_{i}-\frac{\left(\sum_{i}q_{i}\right)^{2}}{\sum_{i}\frac{q^{2}_{i}}{p_{i}}}=\sum_{i}p_{i}-K_{0}\sum_{i}q_{i}
\label{eq.chi2PI}	
\end{equation}
The gradient with respect to ``$q$" is written:
\begin{equation}
	\frac{\partial \chi_{P}^{2}I\left(p\|q\right)}{\partial q_{j}}=\frac{\sum_{i}q_{i}}{\sum_{i}\frac{q^{2}_{i}}{p_{i}}}\left(\frac{\sum_{i}q_{i}}{\sum_{i}\frac{q^{2}_{i}}{p_{i}}}\frac{q_{j}}{p_{j}}-1\right)
\end{equation}
We derive the logarithmic form of the divergence (\ref{eq.chi2PI}):
\begin{equation}
L\chi_{P}^{2}I\left(p\|q\right)=\log\sum_{i}\frac{q^{2}_{i}}{p_{i}}-\log\frac{\left(\sum_{i}q_{i}\right)^{2}}{\sum_{i}p_{i}}	
\end{equation}
Note that this divergence is not only invariant with respect to ``$q$", but also with respect to ``$p$".\\
Its gradient with respect to ``$q$" is written:
\begin{equation}
	\frac{\partial L\chi_{P}^{2}I\left(p\|q\right)}{\partial q_{j}}=\frac{1}{\sum_{i}q_{i}}\left(\frac{\sum_{i}q_{i}}{\sum_{i}\frac{q^{2}_{i}}{p_{i}}}\frac{q_{j}}{p_{j}}-1\right)	
\end{equation}

\subsection{Extensions of Chi2 divergences.}
An extension of the $\chi^{2}$ divergences is mentioned in \cite{cichocki2009} and \cite{pardo1999}. It makes it possible, by action on the single parameter ``$\alpha$" to move gradually from Neyman's $\chi^{2}_{N}$ divergence ($\alpha=1$) to Pearson's $\chi^{2}_{P}$ divergence ($\alpha=0$).\\
This divergence that the authors attribute to Rukhin \cite{rukhin1994} is written: 
\begin{equation}
	RU_{\alpha}\left(p\|q\right)=\sum_{i}\frac{\left(p_{i}-q_{i}\right)^{2}}{\alpha q_{i}+\left(1-\alpha\right)p_{i}}
	\label{eq.rukhin}
\end{equation}
It is built in the Csiszär sense on the standard convex function:
\begin{equation}
	f_{c}\left(x\right)=\frac{\left(x-1\right)^{2}}{\alpha+\left(1-\alpha\right)x}
\end{equation}
The gradient of $RU_{\alpha}\left(p\|q\right)$, (\ref{eq.rukhin}), with respect to ``$q$" is written as follows:
\begin{equation}
	\frac{\partial RU_{\alpha}\left(p\|q\right)}{\partial q_{j}}=\frac{\left(p_{j}-q_{j}\right)\left(\alpha p_{j}-\alpha q_{j}-2p_{j}\right)}{\left[\alpha q_{j}+\left(1-\alpha\right)p_{j}\right]^{2}}
\end{equation}
The particular cases corresponding to the $\chi^{2}$ divergences of Neyman and Pearson can be found immediately.\\

\subsection{Chi2 relative divergences.}
Some ambiguity exists in the literature as to the denomination of the $\chi^{2}$ relative divergences; in any case, it is a matter of replacing ``$p$" or ``$q$" by the weighted sum $\alpha p+\left(1-\alpha\right)q\ ,\ 0\leq\alpha\leq 1$ in the $\chi^{2}$ divergences of the preceding sections.\\
For example, if in the Neyman $\chi^{2}_{N}$ Divergence (\ref{eq.chi2Nc}), we replace ``$q$" by $\alpha p+\left(1-\alpha\right)q$, we obtain an expression that is referred to in the literature as the relative Pearson Divergence; this is only a detail.

\section{Hellinger's divergence.}
It is a divergence built in the sense of Csiszär on the standard convex function:
\begin{equation}
	f_{c}\left(x\right)=\left(\sqrt{x}-1\right)^{2}
	\label{eq.fcHell}
\end{equation}
The divergence obtained is symmetrical; it is written as follows:
\begin{equation}
H\left(p\|q\right)=\sum_{i}\left(\sqrt{p_{i}}-\sqrt{q_{i}}\right)^{2}
\label{eq.DHell}	
\end{equation}
Its gradient with respect to ``$q$" is given by:
\begin{equation}
\frac{\partial H\left(p\|q\right)}{\partial q_{j}}=1-\sqrt{\frac{p_{j}}{q_{j}}}
\label{eq.gradDHell}
\end{equation}
An invariant form (with respect to ``$q$") of this divergence can be obtained by using the nominal invariance factor which is expressed as follows:
\begin{equation}
	K_{0}\left(p,q\right)=\left[\frac{\sum_{i}\sqrt{p_{i}q_{i}}}{\sum_{i}q_{i}}\right]^{2}
\end{equation}
The invariant divergence that results is written as follows:
\begin{equation}
HI\left(p\|q\right)=\sum_{i}p_{i}-\frac{\left(\sum_{i}\sqrt{p_{i}q_{i}}\right)^{2}}{\sum_{i}q_{i}}
\label{eq.DIHell}	
\end{equation}
Its gradient with respect to ``$q$" is given by:
\begin{equation}
\frac{\partial HI\left(p\|q\right)}{\partial q_{j}}=\left[\frac{\sum_{i}\sqrt{p_{i}q_{i}}}{\sum_{i}q_{i}}\right]^{2}-\frac{\sum_{i}\sqrt{p_{i}q_{i}}}{\sum_{i}q_{i}}\frac{p_{j}}{\sqrt{p_{j}q_{j}}}
\label{eq.gradDIHell}	
\end{equation}

\section{Triangular discrimination.}
This divergence mentioned by Tanéja \cite{taneja2001} is built in the sense of Csiszär on the standard convex function:
\begin{equation}
	f_{c}\left(x\right)=\frac{\left(x-1\right)^{2}}{x+1}
\end{equation}
This gives the symmetrical divergence:
\begin{equation}
	DT\left(p\|q\right)=\sum_{i}\frac{\left(p_{i}-q_{i}\right)^{2}}{p_{i}+q_{i}}
\end{equation}
Its gradient with respect to ``$q$" is given by:
\begin{equation}
\frac{\partial DT\left(p\|q\right)}{\partial q_{j}}=1-\frac{4p^{2}_{j}}{\left(p_{j}+q_{j}\right)^{2}}
\end{equation}
It is equal to zéro if $p_{j}=q_{j}\ \forall j$.\\
This divergence is to be linked to the one discussed in the following section.

\section{Harmonic mean divergence of Toussaint.}
\begin{figure}[h!]
\centering
\includegraphics[width=0.7\linewidth]{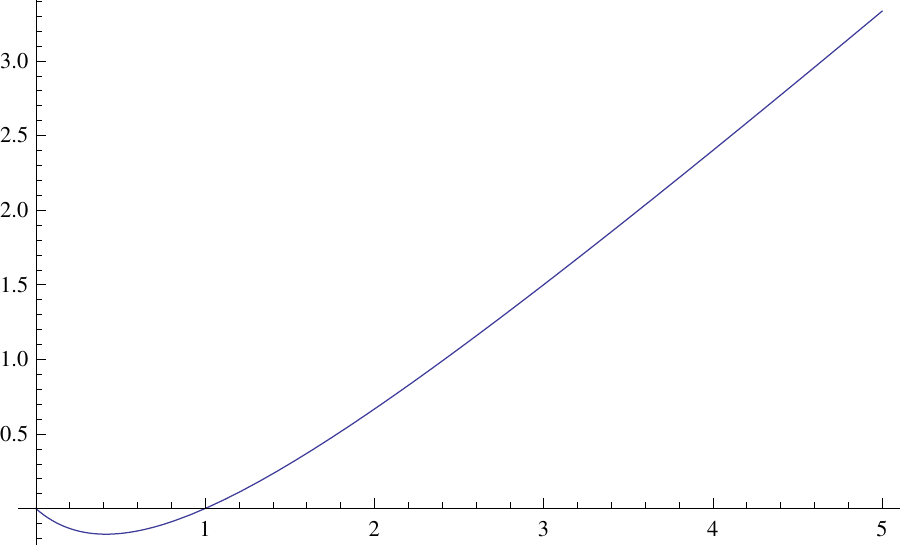}
\caption{Function $f\left(x\right)=x-\frac{2x}{1+x}$}
\label{fig:TOU}
\end{figure}

\begin{figure}[h!]
\centering
\includegraphics[width=0.7\linewidth]{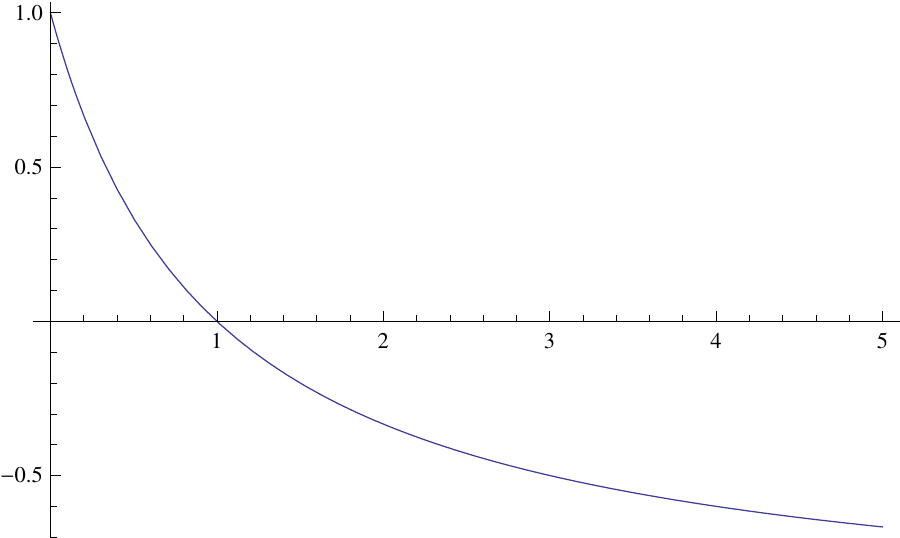}
\caption{Function $	\breve{f}\left(x\right)=1-\frac{2x}{1+x}$}
\label{fig:TOUt}
\end{figure}

She's quoted in Basseville, \cite{basseville1989},\cite{basseville1996},\cite{basseville2013}; it is close to what we would call "$M_{AH}$" \cite{taneja2001} in Chapter 7.\\
It's a Csiszär divergence based on the simple convex function see.fig.(\ref{fig:TOU}):
\begin{equation}
	f\left(x\right)=x-\frac{2x}{1+x}
\end{equation}
We obtain the divergence:
\begin{equation}
	T\left(p\|q\right)=\sum_{i}\left(p_{i}-\frac{2 p_{i}q_{i}}{p_{i}+q_{i}}\right)
\end{equation}
Its gradient with respect to ``$q$" is given by:
\begin{equation}
	\frac{\partial T\left(p\|q\right)}{\partial q_{j}}=-2\left(\frac{p_{j}}{p_{j}+q_{j}}\right)^{2}
\end{equation}
This gradient is never equal to zero.\\
The dual divergence is written as follows:
\begin{equation}
	T\left(q\|p\right)=\sum_{i}\left(q_{i}-\frac{2 p_{i}q_{i}}{p_{i}+q_{i}}\right)
\end{equation}
It is based on the mirror function of the previous one, i.e. on the simple convex function see.fig.(\ref{fig:TOUt}):
\begin{equation}
	\breve{f}\left(x\right)=1-\frac{2x}{1+x}
\end{equation}
Its gradient with respect to ``$q$" is given by:
\begin{equation}
\frac{\partial T\left(q\|p\right)}{\partial q_{j}}=1-2\left(\frac{p_{j}}{p_{j}+q_{j}}\right)^{2}	
\end{equation}
He is not zero for $p_{j}=q_{j}\ \forall j$, hence the problem previously mentioned.\\

\begin{figure}[h!]
\centering
\includegraphics[width=0.7\linewidth]{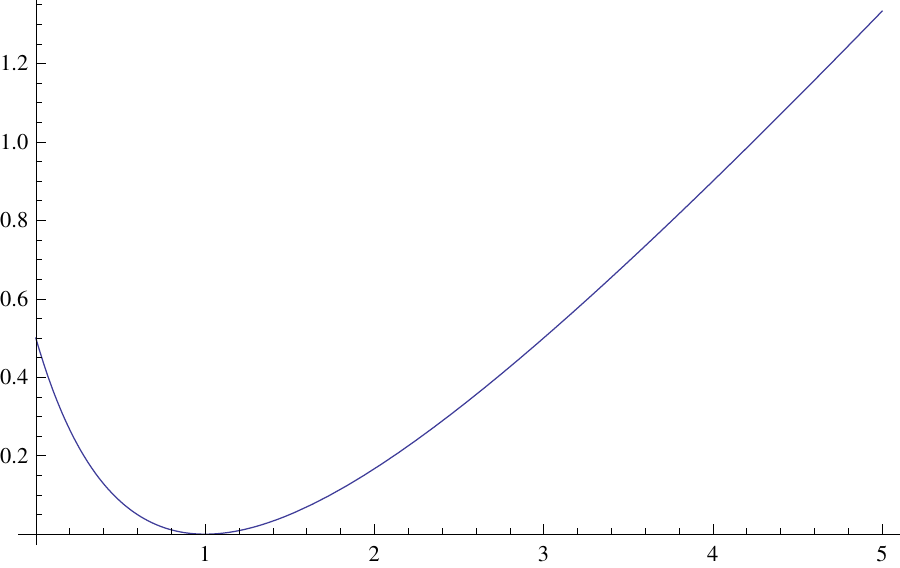}
\caption{Function $	f_{c}\left(x\right)=\frac{1}{2}\frac{\left(x-1\right)^{2}}{x+1}$}
\label{fig:TOUc}
\end{figure}

If we build the standard convex function from the function $f\left(x\right)$ or the mirror function $\breve{f}\left(x\right)$, we obtain one and the same function see.fig.(\ref{fig:TOUc}):
\begin{equation}
	f_{c}\left(x\right)=\frac{1}{2}\frac{\left(x-1\right)^{2}}{x+1}
\end{equation}
This result is indicative of the fact that the Csiszär divergence built on this function will be symmetrical.\\
It is up to a multiplicative factor, the function that allows us to construct the Triangular Discrimination in the previous section and that will enable us to build the arithmetic-harmonic mean $M_{AH}$ divergence. \cite{taneja2001}.

\section{Henze-Penrose divergence.}
\begin{figure}[h!]
\centering
\includegraphics[width=0.7\linewidth]{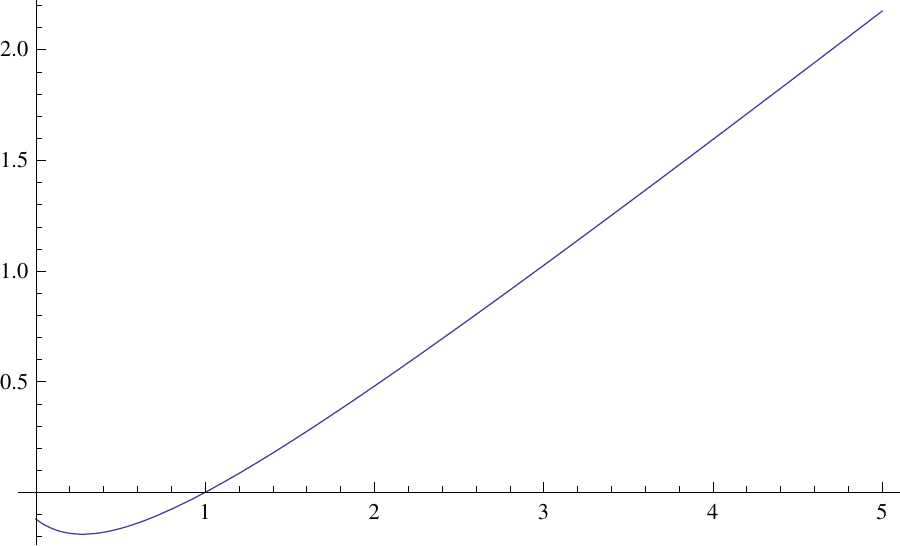}
\caption{Function $	f_{1}\left(x\right)=\frac{\beta^{2}x^{2}+\left(1-\beta\right)^{2}}{\beta x+\left(1-\beta\right)}-\left[\beta^{2}+\left(1-\beta\right)^{2}\right]$, $\beta=0.5$}
\label{fig:HP1}
\end{figure}

\begin{figure}[h!]
\centering
\includegraphics[width=0.7\linewidth]{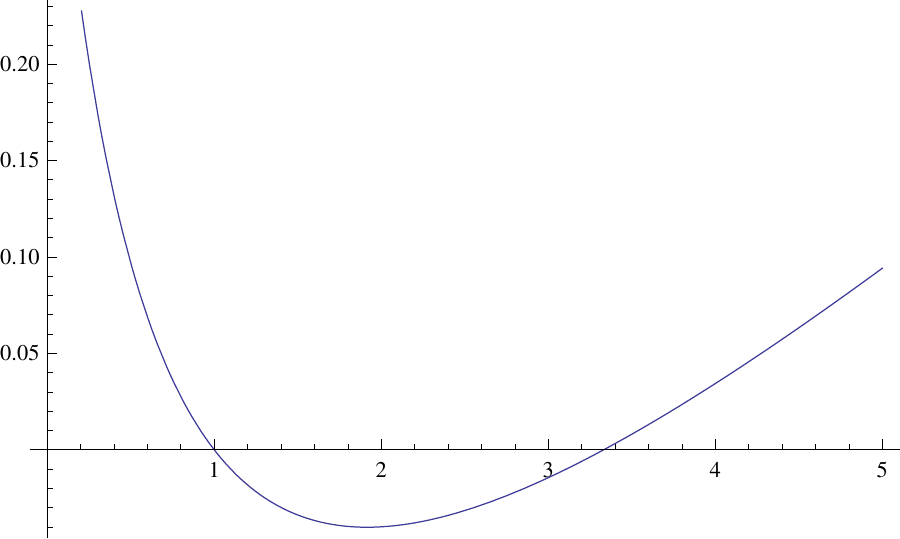}
\caption{Function $	f_{2}\left(x\right)=\frac{\beta^{2}x^{2}+\left(1-\beta\right)^{2}}{\beta x+\left(1-\beta\right)}-\left[\beta^{2}+\left(1-\beta\right)^{2}\right]x$, $\beta=0.5$}
\label{fig:HP2}
\end{figure}

This divergence proposed in \cite{henze1999} is used in \cite{neemuchwala2007} with $\beta<$1, in the form:
\begin{equation}
	HP\left(p\|q\right)=\sum_{i}\frac{\beta^{2}p_{i}^{2}+\left(1-\beta\right)^{2}q_{i}^{2}}{\beta p_{i}+\left(1-\beta\right) q_{i}}
\end{equation}
This form is oversimplified because it is not zero if $p_{i}=q_{i}\ \forall i$, even if one uses variables summed to 1.\\
The first modification consists in making it correct from this point of view, that is to say equal to zero when $p_{i}=q_{i}\;\forall i$; we obtain 2 possible forms:
\begin{equation}
		HP1\left(p\|q\right)=\sum_{i}\frac{\beta^{2}p_{i}^{2}+\left(1-\beta\right)^{2}q_{i}^{2}}{\beta p_{i}+\left(1-\beta\right) q_{i}}-\sum_{i}\left(\beta^{2}+\left(1-\beta\right)^{2}\right)q_{i}
\end{equation}
and
\begin{equation}
		HP2\left(p\|q\right)=\sum_{i}\frac{\beta^{2}p_{i}^{2}+\left(1-\beta\right)^{2}q_{i}^{2}}{\beta p_{i}+\left(1-\beta\right) q_{i}}-\sum_{i}\left(\beta^{2}+\left(1-\beta\right)^{2}\right)p_{i}
\end{equation}
The first expression is a Csiszär divergence built on the simple convex function:
\begin{equation}
	f_{1}\left(x\right)=\frac{\beta^{2}x^{2}+\left(1-\beta\right)^{2}}{\beta x+\left(1-\beta\right)}-\left[\beta^{2}+\left(1-\beta\right)^{2}\right]
\end{equation}
illustrated in the figure (\ref{fig:HP1}), for $\beta=0.5$.
\\
The second is a Csiszär divergence built on the simple convex function:
\begin{equation}
	f_{2}\left(x\right)=\frac{\beta^{2}x^{2}+\left(1-\beta\right)^{2}}{\beta x+\left(1-\beta\right)}-\left[\beta^{2}+\left(1-\beta\right)^{2}\right]x
\end{equation}
shown in the figure (\ref{fig:HP2}), for $\beta=0.5$.

\begin{figure}[h!]
\centering
\includegraphics[width=0.7\linewidth]{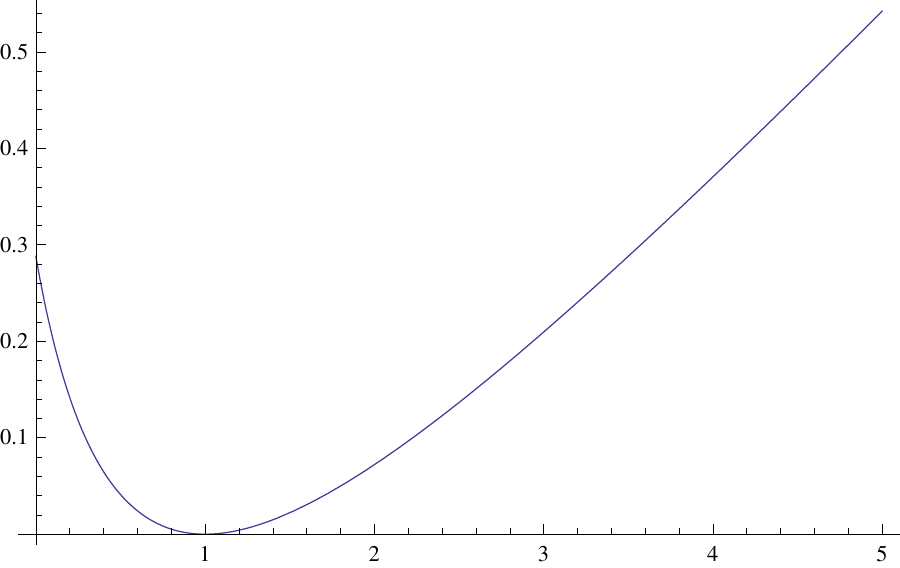}
\caption{Function $f_{c}\left(x\right)=\frac{\beta^{2}x^{2}+\left(1-\beta\right)^{2}}{\beta x+\left(1-\beta\right)}-\left[\beta^{2}+\left(1-\beta\right)^{2}\right]+\left(x-1\right)\left[2\beta^{3}-4\beta^{2}+\beta\right]$, $\beta=0.5$}
\label{fig:HPc}
\end{figure}

If we construct the standard convex function from $f_{1}$ or $f_{2}$, we obtain the same function which is written:
\begin{align}
	f_{c}\left(x\right)=&\frac{\beta^{2}x^{2}+\left(1-\beta\right)^{2}}{\beta x+\left(1-\beta\right)}-\left[\beta^{2}+\left(1-\beta\right)^{2}\right] \nonumber \\ &+\left(x-1\right)\left[2\beta^{3}-4\beta^{2}+\beta\right]
	\label{fcHP}
\end{align}
illustrated in the figure (\ref{fig:HPc}), for $\beta=0.5$.
 
By using the notation:
\begin{align}
	A=\beta^{2}+\left(1-\beta\right)^{2};\ 
	B=2\beta^{3}-4\beta^{2}+\beta
\end{align}
the Csiszär divergence built on the standard convex function (\ref{fcHP}) is written:
\begin{equation}
		HPC\left(p\|q\right)=\sum_{i}\frac{\beta^{2}p_{i}^{2}+\left(1-\beta\right)^{2}q_{i}^{2}}{\beta p_{i}+\left(1-\beta\right) q_{i}}+\sum_{i}B p_{i}-\sum_{i}\left(A+B\right) q_{i}
\end{equation}
Its gradient with respect to ``$q$" is written as:
\begin{align}
		\frac{\partial HPC\left(p\|q\right)}{\partial q_{j}}=&\left(1-\beta\right)\left\{1-\frac{2\beta^{2}p_{j}^{2}}{\left[\beta p_{j}+\left(1-\beta\right)q_{j}\right]^{2}}\right\} \nonumber \\ &-\left(2\beta^{3}-2\beta^{2}-\beta+1\right)
\end{align}
Which can also be written more simply:
\begin{equation}
		\frac{\partial HPC\left(p\|q\right)}{\partial q_{j}}=2\beta^{2}\left(1-\beta\right)\left\{1-\frac{p_{j}^{2}}{\left[\beta p_{j}+\left(1-\beta\right)q_{j}\right]^{2}}\right\}
\end{equation}
The dual divergence, built on the function $\breve{f}_{c}$ mirroring $f_{c}$, is written:
\begin{equation}
		HPC\left(q\|p\right)=\sum_{i}\frac{\beta^{2}q_{i}^{2}+\left(1-\beta\right)^{2}p_{i}^{2}}{\beta q_{i}+\left(1-\beta\right) p_{i}}+\sum_{i}B q_{i}-\sum_{i}\left(A+B\right) p_{i}
\end{equation}
Its gradient with respect to ``$q$" is written as:
\begin{equation}
\frac{\partial HPC\left(q\|p\right)}{\partial q_{j}}=2\beta\left(1-\beta\right)^{2}\left\{1-\frac{p_{j}^{2}}{\left[\beta q_{j}+\left(1-\beta\right)p_{j}\right]^{2}}\right\}	
\end{equation}
These gradients are of course zero for $p_{j}=q_{j}\ \forall j$.

\section{Polya information divergence.}
Here we consider the ``Polya information divergence" \cite{hardy1952} quoted in \cite{grendar2010}.
\subsection{Some general observations.}
\begin{figure}[h!]
\centering
\includegraphics[width=0.7\linewidth]{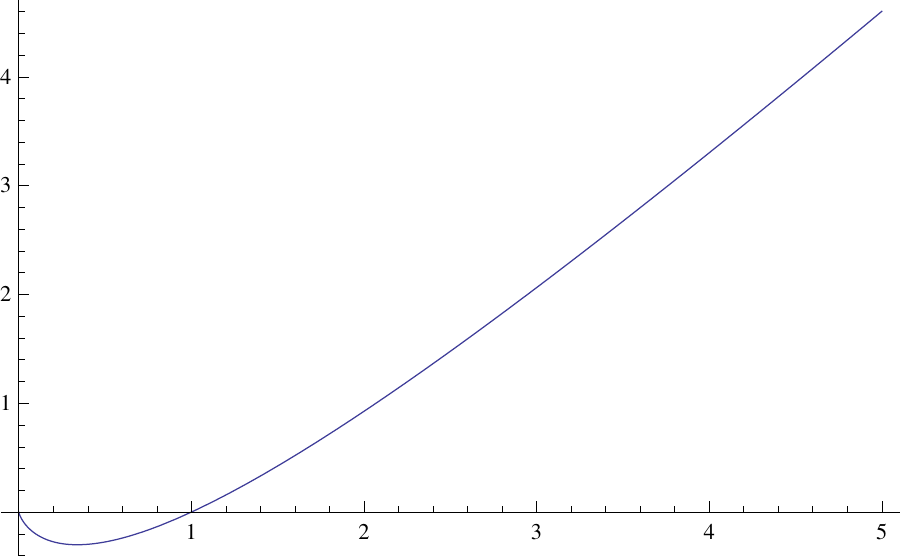}
\caption{Function $f_{1}\left(x\right)=x\log\frac{x}{1+\gamma x}+\frac{1}{\gamma}\log \frac{1}{1+\gamma x}+\frac{\left(1+\gamma\right)x}{\gamma}\log \left(1+\gamma\right)$, $\gamma=0.8$}
\label{fig:POL1}
\end{figure}

On the basis of the simplified Kullback divergence or Kullback Information expressed as: 
\begin{equation}
	IKL(a\|b)=\sum_{i}a_{i}\log\left(a_{i}/b_{i}\right)
\end{equation}
With $\gamma\geq 0$, Grendar and Niven \cite{grendar2010} are proposing the divergence:
\begin{equation}
	I_{\gamma}(p\|q)=IKL(p\|q+\gamma p)+\frac{1}{\gamma}IKL(q\|q+\gamma p)+\frac{1+\gamma}{\gamma}\log(1+\gamma)
\end{equation}
which leads to:
\begin{equation}
I_{\gamma}(p\|q)=\sum_{i}p_{i}\log\frac{p_{i}}{q_{i}+\gamma p_{i}}+\frac{1}{\gamma}\sum_{i}q_{i}\log\frac{q_{i}}{q_{i}+\gamma p_{i}}+\frac{1+\gamma}{\gamma}\log(1+\gamma)	
\end{equation}
This divergence can be used essentially with data of sum equal to 1, i.e. typically probability densities. Indeed, if in the preceding expression we do $p_{i}=q_{i}\ \forall i$, we will obtain zero only if $\sum_{i}p_{i}=\sum_{i}q_{i}=1$.\\
This means that the basic convex function that allows us to construct this divergence in the Csiszär sense can have 2 forms:
\begin{equation}
f_{1}\left(x\right)=x\log\frac{x}{1+\gamma x}+\frac{1}{\gamma}\log \frac{1}{1+\gamma x}+\left\{\frac{1+\gamma}{\gamma}\log \left(1+\gamma\right)\right\}x	
\end{equation}
shown in the figure (\ref{fig:POL1}), for $\gamma=0.8$,\\

\begin{figure}[h!]
\centering
\includegraphics[width=0.7\linewidth]{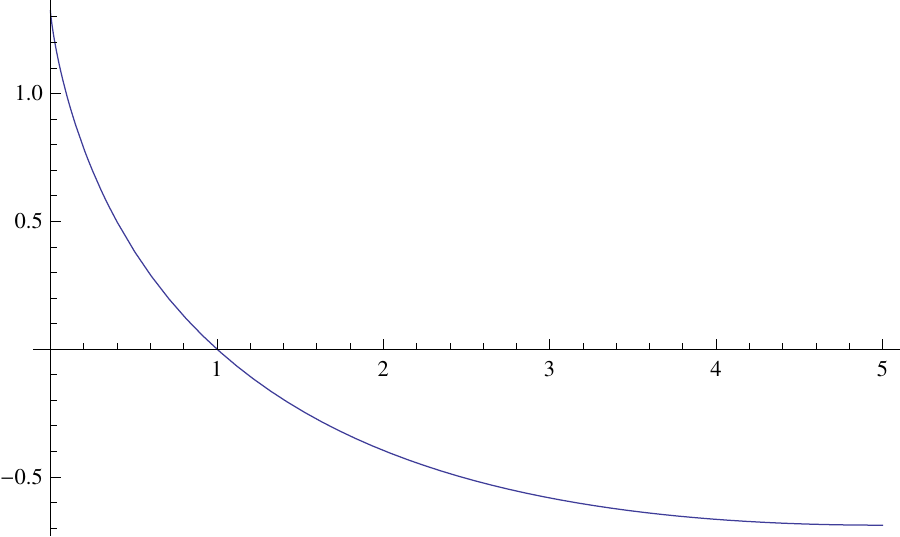}
\caption{Function $f_{2}\left(x\right)=x\log\frac{x}{1+\gamma x}+\frac{1}{\gamma}\log \frac{1}{1+\gamma x}+\frac{1+\gamma}{\gamma}\log \left(1+\gamma\right)$, $\gamma=0.8$}
\label{fig:POL2}
\end{figure}

\begin{figure}[h!]
\centering
\includegraphics[width=0.7\linewidth]{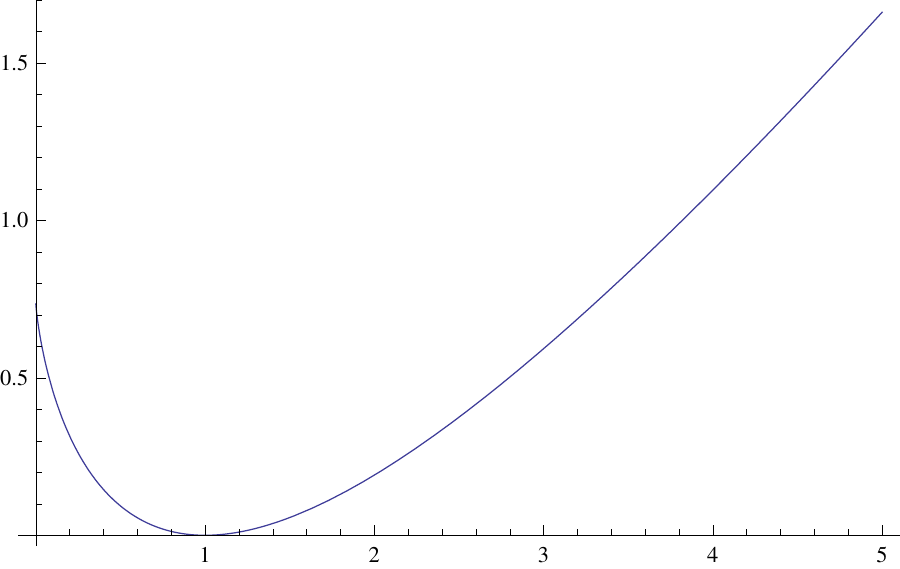}
\caption{Function $f_{c}\left(x\right)=x\log\frac{x}{1+\gamma x}+\frac{1}{\gamma}\log \frac{1}{1+\gamma x}+\frac{1+\gamma x}{\gamma}\log \left(1+\gamma\right)$, $\gamma=0.8$}
\label{fig:POLc}
\end{figure}

 or:
\begin{equation}
f_{2}\left(x\right)=x\log\frac{x}{1+\gamma x}+\frac{1}{\gamma}\log \frac{1}{1+\gamma x}+\frac{1+\gamma}{\gamma}\log \left(1+\gamma\right)	
\end{equation}
shown in the figure (\ref{fig:POL2}), for $\gamma=0.8$,\\
This being stated, these 2 functions are not standard forms because:
\begin{equation}
f'_{1}\left(x\right)=\log\frac{x}{1+\gamma x}+\frac{1+\gamma}{\gamma}\log\left(1+\gamma\right)\;\Rightarrow\;f'_{1}\left(1\right)=\frac{1}{\gamma}\log\left(1+\gamma\right)		
\end{equation}
\begin{equation}
	f'_{2}\left(x\right)=\log\frac{x}{1+\gamma x}\;\Rightarrow\;f'_{2}\left(1\right)=\log\frac{1}{1+\gamma}
\end{equation}
If we construct the corresponding standard convex function, we obtain from $f_{1}\left(x\right)$ or from $f_{2}\left(x\right)$, the same expression:
\begin{equation}
f_{c}\left(x\right)=x\log\frac{x}{1+\gamma x}+\frac{1}{\gamma}\log \frac{1}{1+\gamma x}+\frac{1+\gamma x}{\gamma}\log \left(1+\gamma\right)
\label{eq.fcpolya}	
\end{equation}
illustrated in the figure (\ref{fig:POLc}), for $\gamma=0.8$\\
This function has all of the desired properties: $f_{c}\left(1\right)=0$, $f'_{c}\left(1\right)=0$, which makes it possible to construct a Csiszär divergence that is usable in all generality:
\begin{align}
	IC_{\gamma}\left(p\|q\right)=&\sum_{i}p_{i}\log\frac{p_{i}}{q_{i}+\gamma p_{i}}+\frac{1}{\gamma}\sum_{i}q_{i}\log\frac{q_{i}}{q_{i}+\gamma p_{i}}  \nonumber \\ &+\frac{\log(1+\gamma)}{\gamma}\sum_{i}q_{i}+\gamma p_{i}
	\label{eq.divGN}
\end{align}
The gradient with respect to ``$q$" is expressed as:
\begin{equation}
	\frac{\partial IC_{\gamma}\left(p\|q\right)}{\partial q_{j}}=\frac{1}{\gamma} \log\frac{q_{j}+\gamma q_{j}}{q_{j}+\gamma p_{j}}
\end{equation}
If we are now interested in the dual divergence, it is built, in the sense of Csiszär, on the mirror function of (\ref{eq.fcpolya}), 
\begin{equation}	\breve{f_{c}}\left(x\right)=\log\frac{1}{x+\gamma}+\frac{x}{\gamma}\log\frac{x}{x+\gamma}+\frac{x+\gamma}{\gamma}\log\left(1+\gamma\right)
\end{equation}
shown in the figure (\ref{fig:POLct}), for $\gamma=0.8$:
\begin{figure}[h!]
\centering
\includegraphics[width=0.7\linewidth]{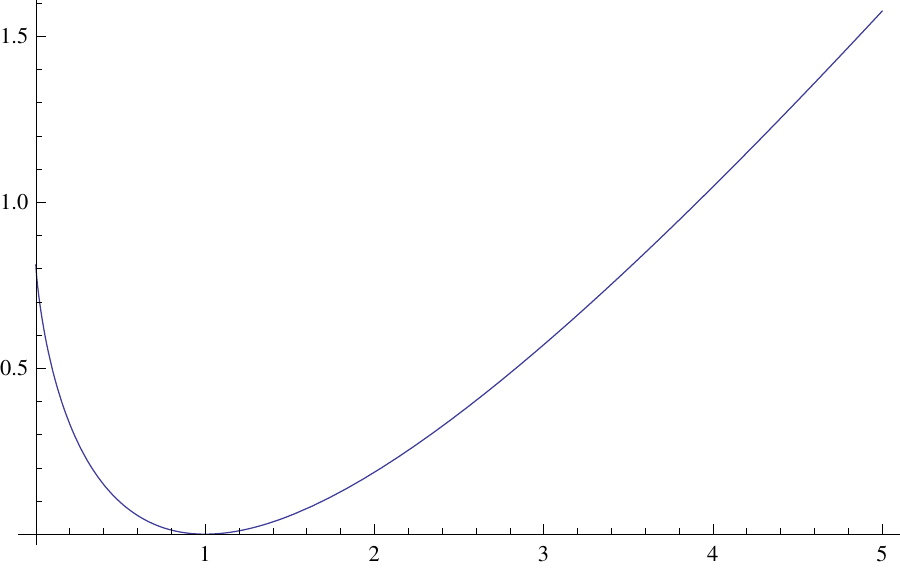}
\caption{Function $\breve{f_{c}}\left(x\right)=\log\frac{1}{x+\gamma}+\frac{x}{\gamma}\log\frac{x}{x+\gamma}+\frac{x+\gamma}{\gamma}\log\left(1+\gamma\right)$, $\gamma=0.8$}
\label{fig:POLct}
\end{figure}

This divergence is written:
\begin{align}
IC_{\gamma}\left(q\|p\right)=&\sum_{i}q_{i}\log\frac{q_{i}}{p_{i}+\gamma q_{i}}+\frac{1}{\gamma}\sum_{i}p_{i}\log\frac{p_{i}}{p_{i}+\gamma q_{i}}\nonumber \\ &+\frac{\log(1+\gamma)}{\gamma}\sum_{i}p_{i}+\gamma q_{i}	
\end{align}
The gradient with respect to ``$q$" is expressed as:
\begin{equation}
	\frac{\partial IC_{\gamma}\left(q\|p\right)}{\partial q_{j}}=\log\frac{q_{j}+\gamma q_{j}}{p_{j}+\gamma q_{j}}
\end{equation}
The symmetrical divergence can be obtained in the Jeffreys sense by making the (half) sum of the 2 preceding conjugated divergences; the divergence obtained is of course always a Csiszär divergence built on a convex function which is the (half) sum of the conjugated functions $f_{c}\left(x\right)$ and $\breve{f_{c}}\left(x\right)$.

\subsection{Invariance by change of scale.}
We can try to construct a scale-invariant version of the divergence (\ref{eq.divGN}) using the standard process.\\
First, we write:
\begin{align}
	IC_{\gamma}\left(p\|Kq\right)=&\sum_{i}p_{i}\log\frac{p_{i}}{Kq_{i}+\gamma p_{i}}+\frac{1}{\gamma}\sum_{i}Kq_{i}\log\frac{Kq_{i}}{Kq_{i}+\gamma p_{i}}\nonumber \\ &+\frac{\log(1+\gamma)}{\gamma}\sum_{i}Kq_{i}+\gamma p_{i}
	\label{eq.divGNK}
\end{align}
We then calculate the derivative of this expression with respect to $K$, and in order to determine the expression of $K$, we write that this derivative is zero.\\
Which leads us to have to solve:
\begin{equation}
\sum_{i}q_{i}\log\left(1+\frac{\gamma p_{i}}{Kq_{i}}\right)=\sum_{i}q_{i}\log\left(1+\gamma\right)	
\end{equation}
that is:
\begin{equation}
\sum_{i}q_{i}\log\frac{\left(1+\frac{\gamma p_{i}}{Kq_{i}}\right)}{\left(1+\gamma\right)}=0	
\end{equation}
There is no explicit $K$ solution. Strictly speaking, this scalar product will be zero if the 2 vectors are orthogonal, or if one of the vectors is zero; however, to illustrate the remark mentioned in Ch.3, section 3.5, we can use as invariance factor:
\begin{equation}
K^{*}=\frac{\sum_{i}p_{i}}{\sum_{i}q_{i}}
\end{equation}
By inserting this expression of $K^{*}$ in (\ref{eq.divGNK}), it comes, with all simplifications made:
\begin{align}	ICINV_{\gamma}\left(p\|q\right)=&\sum_{j}p_{j}\left\{\sum_{i}\overline{p}_{i}\log\frac{\overline{p}_{i}}{\overline{q}_{i}+\gamma \overline{p}_{i}}+\frac{1}{\gamma}\sum_{i}\overline{q}_{i}\log\frac{\overline{q}_{i}}{\overline{q}_{i}+\gamma \overline{p}_{i}}\right. \nonumber \\ &\left.+\frac{\log(1+\gamma)}{\gamma}\sum_{i}\overline{q}_{i}+\gamma \overline{p}_{i}\right\}
	\label{eq.ICINV}
\end{align}
So, if we calculate the gradient of (\ref{eq.ICINV}) with respect to ``$q$", we obtain all the calculations done:
\begin{equation}
	\frac{\partial ICINV_{\gamma}\left(p\|q\right)}{\partial q_{j}}=\frac{\sum_{l}p_{l}}{\gamma\sum_{l}q_{l}}\left[\log\frac{\overline{q}_{j}}{\overline{q}_{j}+\gamma \overline{p}_{j}}-\sum_{i}\overline{q}_{i}\log\frac{\overline{q}_{i}}{\overline{q}_{i}+\gamma \overline{p}_{i}}\right]
\end{equation}
And we can observe that, as with all the scale invariant divergences , we have:
\begin{equation}
	\sum_{j}q_{j}\frac{\partial ICINV_{\gamma}\left(p\|q\right)}{\partial q_{j}}=0
\end{equation}

\subsection{Remark.}
We can exhibit a Bregman divergence based on the function $f_{c}\left(x\right)$,\;which is written:
\begin{equation}
		BIC_{\gamma}\left(p\|q\right)=\sum_{i}p_{i}\log\frac{p_{i}}{q_{i}}-\frac{1}{\gamma}\left(1+\gamma p_{i}\right)\log\frac{1+\gamma p_{i}}{1+\gamma q_{i}}
		\label{eq.divBIC}
\end{equation}
But, of course, there is no evidence of the convexity of this divergence.

\section{Bose-Einstein and Fermi-Dirac divergences.}
\subsection{Bose-Einstein divergence.}
\begin{figure}[h!]
\centering
\includegraphics[width=0.7\linewidth]{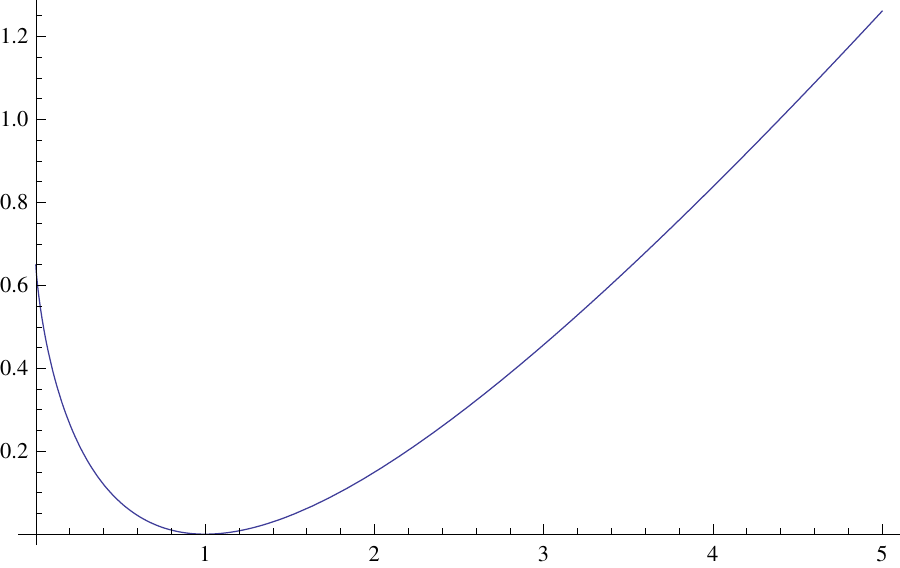}
\caption{Function $f_{c,\alpha}\left(x\right)=x \log x+\left(\alpha+x\right)\log\frac{\alpha+1}{\alpha+x}$, $\alpha=0.8$}
\label{fig:BE1c}
\end{figure}

Such divergence mentioned in Basseville \cite{basseville1989, basseville1996, basseville2013}, and in the references cited \cite{knockaert1993}, is a Csiszär divergence based on the standard convex function shown in the figure (\ref{fig:BE1c}), for $\alpha=0.8$:
\begin{equation}
	f_{c,\alpha}\left(x\right)=x \log x+\left(\alpha+x\right)\log\frac{\alpha+1}{\alpha+x}\ \ \ \ \ \alpha>0
	\label{eq.fcBE}
\end{equation}
It is written:
\begin{equation}
	BE1_{\alpha}\left(p\|q\right)=\sum_{i}p_{i}\log\frac{p_{i}}{q_{i}}+\left(\alpha q_{i}+p_{i}\right)\log\frac{\left(\alpha q_{i}+q_{i}\right)}{\left(\alpha q_{i}+p_{i}\right)}
	\label{eq.BE1}
\end{equation}
However, another divergence appears under the same name in \cite{furuichi2012}; it is written:
\begin{equation}	BE2\left(p\|q\right)=\sum_{i}p_{i}\log\frac{p_{i}}{q_{i}}+\left(1+p_{i}\right)\log\frac{\left(1+q_{i}\right)}{\left(1+p_{i}\right)}	
\end{equation}
It is not a Csiszar divergence.\\
It is deduced from a work by Furuichi \cite{furuichi2012} quoted in referring to Kapur \cite{kapur1967}.\\
It is, in fact, a Bregman divergence built on the function (\ref{eq.fcBE}) which is written:
\begin{equation}	BE2_{\alpha}\left(p\|q\right)=\sum_{i}p_{i}\log\frac{p_{i}}{q_{i}}+\left(\alpha+p_{i}\right)\log\frac{\left(\alpha+q_{i}\right)}{\left(\alpha+p_{i}\right)}
\label{eq.divBE2}	
\end{equation}
Note that except for a change in parameter, this divergence is similar to (\ref{eq.divBIC}).\\
If, in the expression (\ref{eq.divBE2}) we introduce the simplification $\alpha=1$, we recover the form proposed in \cite{furuichi2012} \cite{kapur1967} \cite{kapur1983}.\\ 
From these expressions, we can introduce the Generalized Logarithm as suggested in \cite{furuichi2012} by referring to Tsallis' entropy \cite{tsallis1988}.\\
The introduction of Tsallis entropy consists in replacing the Logarithm function by the Generalized Logarithm function denoted ``$Log_{d}$'' (see Appendix 1), which yields:
\begin{equation}
	BE1_{\alpha,d}\left(p\|q\right)=\sum_{i}p_{i}\log_{d}\frac{p_{i}}{q_{i}}+\left(\alpha q_{i}+p_{i}\right)\log_{d}\frac{\left(\alpha q_{i}+q_{i}\right)}{\left(\alpha q_{i}+p_{i}\right)}
\end{equation}
That is:
\begin{equation}
	BE1_{\alpha,d}\left(p\|q\right)=\sum_{i}p_{i}\frac{\left(\frac{p_{i}}{q_{i}}\right)^{1-d}-1}{1-d}+\left(\alpha q_{i}+p_{i}\right)\frac{\left(\frac{\alpha q_{i}+q_{i}}{\alpha q_{i}+p_{i}}\right)^{1-d}-1}{1-d}	
\end{equation}
The calculation of the gradient with respect to ``$q$" leads to:
\begin{align}
	\frac{\partial BE1_{\alpha,d}\left(p\|q\right)}{\partial q_{j}}=&-\left(\frac{p_{j}}{q_{j}}\right)^{2-d}+\frac{\alpha d}{1-d}\left(\frac{\alpha q_{j}+q_{j}}{\alpha q_{j}+p_{j}}\right)^{1-d} \nonumber \\ &+\left(\alpha+1\right)\left(\frac{\alpha q_{j}+q_{j}}{\alpha q_{j}+p_{j}}\right)^{-d}-\frac{\alpha}{1-d}
\end{align}
The limit $d\rightarrow 1$ leads to:
\begin{equation}
	\frac{\partial BE1_{\alpha,1}\left(p\|q\right)}{\partial q_{j}}=\alpha \log\frac{\alpha q_{j}+q_{j}}{\alpha q_{j}+p_{j}}	
\end{equation}
Which is the gradient of the initial Bose-Einstein divergence $BE1_{\alpha}\left(p\|q\right)$.\\
If we now consider the form $BE2_{\alpha}\left(p\|q\right)$, we obtain:
\begin{equation}	BE2_{\alpha,d}\left(p\|q\right)=\sum_{i}p_{i}\log_{d}\frac{p_{i}}{q_{i}}+\left(\alpha+p_{i}\right)\log_{d}\frac{\left(\alpha+q_{i}\right)}{\left(\alpha+p_{i}\right)}
\end{equation}
That is:
\begin{equation}	BE2_{\alpha,d}\left(p\|q\right)=\sum_{i}p_{i}\frac{\left(\frac{p_{i}}{q_{i}}\right)^{1-d}-1}{1-d}+\left(\alpha+p_{i}\right)\frac{\left(\frac{\alpha+q_{i}}{\alpha+p_{i}}\right)^{1-d}-1}{1-d}	
\end{equation}
For $\alpha=1$, this divergence is referred to as the Bose-Einstein-Tsallis divergence in \cite{furuichi2012}.\\
The gradient with respect to ``$q$" is written:
\begin{equation}
\frac{\partial BE2_{\alpha,d}\left(p\|q\right)}{\partial q_{j}}=-\left(\frac{p_{j}}{q_{j}}\right)^{2-d}+\left(\frac{\alpha+p_{j}}{\alpha+q_{j}}\right)^{d}	
\end{equation}
The limit $d\rightarrow 1$ leads to:
\begin{equation}
		\frac{\partial BE2_{\alpha,1}\left(p\|q\right)}{\partial q_{j}}=-\frac{p_{j}}{q_{j}}+\frac{\alpha+p_{j}}{\alpha+q_{j}}
\end{equation}
Which is the gradient of $BE2_{\alpha}\left(p\|q\right)$.

\subsubsection{Scale change invariance.}
From (\ref{eq.BE1}) one can write:
\begin{equation}
	BE1_{\alpha}\left(p\|Kq\right)=\sum_{i}p_{i}\log\frac{p_{i}}{Kq_{i}}+\left(\alpha Kq_{i}+p_{i}\right)\log\frac{\left(\alpha Kq_{i}+Kq_{i}\right)}{\left(\alpha Kq_{i}+p_{i}\right)}
	\label{eq.BE1K}
\end{equation}
and we have:
\begin{equation}
	\frac{\partial BE1_{\alpha}\left(p\|Kq\right)}{\partial K}=\alpha\left[\sum_{i}q_{i}\log\frac{\alpha  Kq_{i}+Kq_{i}}{\alpha Kq_{i}+p_{i}}\right]
\end{equation}
To cancel this derivative, it is necessary to solve:
\begin{equation}
	\sum_{i}q_{i}\log\frac{\alpha  Kq_{i}+Kq_{i}}{\alpha Kq_{i}+p_{i}}=0
\end{equation}
The zero scalar product implies that the corresponding vectors are orthogonal or that one of the vectors is zero.
There is no explicit solution, however, we can use:
\begin{equation}
K^{*}=\frac{\sum_{i}p_{i}}{\sum_{i}q_{i}}
\end{equation}
If we introduce this expression of $K^{*}$ in the relation (\ref{eq.BE1K}), we obtain:
\begin{equation}
BE1_{\alpha}I\left(p\|q\right)=\sum_{i}\bar{p}_{i}\log\frac{\bar{p}_{i}}{\bar{q}_{i}}+\left(\alpha \bar{q}_{i}+\bar{p}_{i}\right)\log\frac{\left(\alpha+1\right)\bar{q}_{i}}{\left(\alpha \bar{q}_{i}+\bar{p}_{i}\right)}
\label{eq.BE1I}
\end{equation}
The gradient with respect to ``$q$" expresses as:
\begin{equation}
\frac{\partial BE1_{\alpha}I\left(p\|q\right)}{\partial q_{l}}=\frac{\alpha}{\sum_{j}q_{j}}\left[\log\frac{\left(\alpha+1\right)\bar{q}_{l}}{\left(\alpha \bar{q}_{l}+\bar{p}_{l}\right)}-\sum_{i}\bar{q}_{i}\log\frac{\left(\alpha+1\right)\bar{q}_{i}}{\left(\alpha \bar{q}_{i}+\bar{p}_{i}\right)}\right]
\end{equation}

\subsection{Fermi-Dirac divergence.}
\begin{figure}[h!]
\centering
\includegraphics[width=0.7\linewidth]{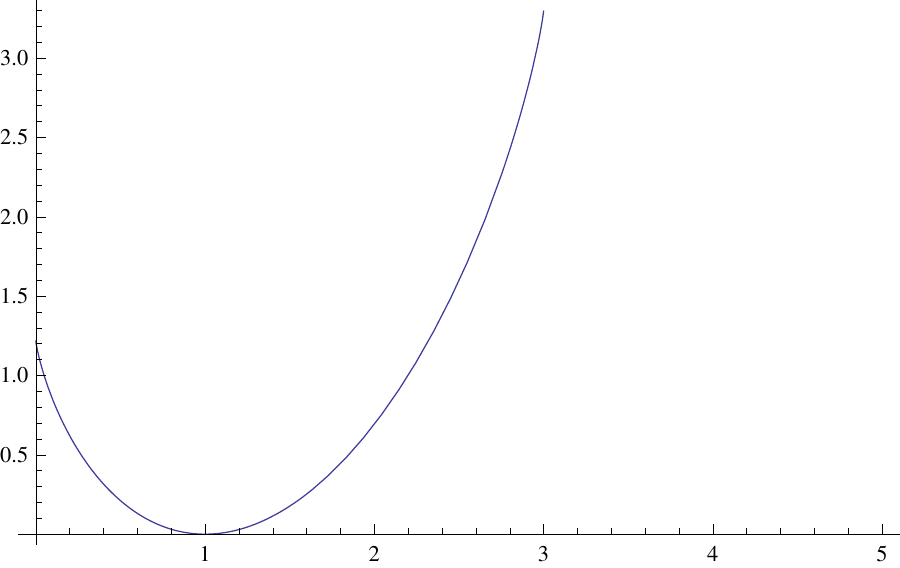}
\caption{Function $f_{c,\beta}\left(x\right)=x \log x+\left(\beta-x\right)\log\frac{\beta-x}{\beta-1}$, $\beta=3$}
\label{fig:FD1c}
\end{figure}

This divergence reported in Basseville \cite{basseville1989, basseville1996, basseville2013}, and in the mentioned references is a Csiszär divergence based on the standard convex function shown in the figure (\ref{fig:FD1c}), for $\beta=3$:
\begin{equation}
	f_{c,\beta}\left(x\right)=x \log x+\left(\beta-x\right)\log\frac{\beta-x}{\beta-1}\ \ \ \ \ \beta>1
	\label{eq.fcFD}
\end{equation}
It is written as:
\begin{equation}
	FD1_{\beta}\left(p\|q\right)=\sum_{i}p_{i}\log\frac{p_{i}}{q_{i}}+\sum_{i}\left(\beta q_{i}-p_{i}\right)\log\frac{\beta q_{i}-p_{i}}{\beta q_{i}-q_{i}}
	\label{eq.FD1}
\end{equation}
However, another divergence appears under the same name in \cite{furuichi2012} \cite{kapur1983}; it is written:
\begin{equation}	FD2\left(p\|q\right)=\sum_{i}p_{i}\log\frac{p_{i}}{q_{i}}+\left(1-p_{i}\right)\log\frac{\left(1-p_{i}\right)}{\left(1-q_{i}\right)}	
\end{equation}
This is not a Csiszâr divergence, it is in fact a Bregman divergence built on the simple convex function:
\begin{equation}
	f_{\beta}\left(x\right)=x \log x+\left(\beta-x\right)\log \left(\beta-x\right);\ \ \ x<\beta
\end{equation}
This divergence is written:
\begin{equation}
FD2_{\beta}\left(p\|q\right)=\sum_{i}p_{i}\log\frac{p_{i}}{q_{i}}+\sum_{i}\left(\beta -p_{i}\right)\log\frac{\beta -p_{i}}{\beta-q_{i}}	
\end{equation}
If in this expression, one makes $\beta=1$, one finds the expression proposed in \cite{furuichi2012}; we will note that the fact of taking $\beta=1$ is equivalent to considering that we are dealing with probability densities, i.e. $p_{i}<1\ \forall i, q_{i}<1\ \forall i$.\\
From these expressions, we can introduce the Generalized logarithm as suggested in \cite{furuichi2012} as inspired by the entropy of Tsallis \cite{tsallis1988}.\\
Using $FD1_{\beta}\left(p\|q\right)$, we obtain:
\begin{equation}
	FD1_{\beta,d}\left(p\|q\right)=\sum_{i}p_{i}\log_{d}\frac{p_{i}}{q_{i}}+\sum_{i}\left(\beta q_{i}-p_{i}\right)\log_{d}\frac{\beta q_{i}-p_{i}}{\beta q_{i}-q_{i}}
\end{equation}
That is:
\begin{equation}
	FD1_{\beta,d}\left(p\|q\right)=\sum_{i}p_{i}\frac{\left(\frac{p_{i}}{q_{i}}\right)^{1-d}-1}{1-d}+\left(\beta q_{i}-p_{i}\right)\frac{\left(\frac{\beta q_{i}-p_{i}}{\beta q_{i}-q_{i}}\right)^{1-d}-1}{1-d}	
\end{equation}
The computation of the gradient with respect to ``$q$" leads to:
\begin{align}
	\frac{\partial FD1_{\beta,d}\left(p\|q\right)}{\partial q_{j}}=&-\frac{p_{j}}{q_{j}}\left[\left(\frac{p_{j}}{q_{j}}\right)^{1-d}-\left(\frac{\beta q_{j}-p_{j}}{\beta q_{j}-q_{j}}\right)^{1-d}\right] \nonumber \\ &+\beta \frac{\left(\frac{\beta q_{j}-p_{j}}{\beta q_{j}-q_{j}}\right)^{1-d}-1}{1-d}
\end{align}
The last term of the second member is the expression of $\beta \log_{d}\left(\frac{\beta q_{i}-p_{i}}{\beta q_{i}-q_{i}}\right)$, so that when $d\rightarrow1$, we obtain:
\begin{equation}
	\frac{\partial FD1_{\beta,1}\left(p\|q\right)}{\partial q_{j}}=	\beta \log\frac{\beta q_{j}-p_{j}}{\beta q_{j}-q_{j}}
\end{equation}
Which is equal to $\frac{\partial FD1_{\beta}\left(p\|q\right)}{\partial q_{j}}$.\\
If now, we consider the divergence $FD2_{\beta}\left(p\|q\right)$ with $\beta=1$ and using the generalized logarithm, we have:
\begin{equation}	FD2_{1,d}\left(p\|q\right)=\sum_{i}p_{i}\log_{d}\frac{p_{i}}{q_{i}}+\left(1-p_{i}\right)\log_{d}\frac{\left(1-p_{i}\right)}{\left(1-q_{i}\right)}	
\end{equation}
That is:
\begin{equation}	FD2_{1,d}\left(p\|q\right)=\sum_{i}p_{i}\frac{\left(\frac{p_{i}}{q_{i}}\right)^{1-d}-1}{1-d}+\left(1-p_{i}\right)\frac{\left(\frac{1-p_{i}}{1-q_{i}}\right)^{1-d}-1}{1-d}	
\end{equation}
This divergence is referred to as the Fermi-Dirac-Tsallis divergence in \cite{furuichi2012}.\\
The calculation of the gradient with respect to ``$q$" leads to:
\begin{equation}
	\frac{\partial FD2_{1,d}\left(p\|q\right)}{\partial q_{j}}=-\left(\frac{p_{j}}{q_{j}}\right)^{2-d}+\left(\frac{1-p_{j}}{1-q_{j}}\right)^{2-d}
\end{equation}
Obviously, if $d\rightarrow 1$, we obtain:
\begin{equation}
	\frac{\partial FD2_{1,1}\left(p\|q\right)}{\partial q_{j}}=-\left(\frac{p_{j}}{q_{j}}\right)+\left(\frac{1-p_{j}}{1-q_{j}}\right)
\end{equation}
Which is equal to $\frac{\partial FD2_{\beta}\left(p\|q\right)}{\partial q_{j}}$, with $\beta=1$.

\subsubsection{Invariance by change of scale.}
From (\ref{eq.FD1}), we have:
\begin{equation}
	FD1_{\beta}\left(p\|Kq\right)=\sum_{i}p_{i}\log\frac{p_{i}}{Kq_{i}}+\sum_{i}\left(\beta Kq_{i}-p_{i}\right)\log\frac{\beta Kq_{i}-p_{i}}{\beta Kq_{i}-Kq_{i}}
	\label{eq.FD1K}
\end{equation}
Consequently:
\begin{equation}
	\frac{\partial FD1_{\beta}\left(p\|Kq\right)}{\partial K}=\beta\sum_{i}q_{i}\log\frac{\beta Kq_{i}-p_{i}}{\beta Kq_{i}-Kq_{i}}
\end{equation}
To cancel this derivative, there is no explicit solution, an alternative is to use $K^{*}$, which leads to the simplified invariant divergence:
\begin{equation}
FD1I_{\beta}\left(p\|q\right)=\sum_{i}\bar{p}_{i}\log\frac{\bar{p}_{i}}{\bar{q}_{i}}+\sum_{i}\left(\beta \bar{q}_{i}-\bar{p}_{i}\right)\log\frac{\beta \bar{q}_{i}-\bar{p}_{i}}{\left(\beta-1\right) \bar{q}_{i}}
\label{eq.FD1I}
\end{equation}
The gradient with respect to ``$q$" is expressed as:
\begin{align}
\frac{\partial FD1_{\beta}I\left(p\|q\right)}{\partial q_{l}}=&\frac{1}{\sum_{j}q_{j}}\left[2\frac{\bar{p}_{l}}{\bar{q}_{l}}-2\right. \nonumber \\  & \left.+\beta	\log\frac{\beta \bar{q}_{l}-\bar{p}_{l}}{\left(\beta-1\right) \bar{q}_{l}}-\beta\sum_{i}\bar{q}_{i}\log\frac{\beta \bar{q}_{i}-\bar{p}_{i}}{\left(\beta-1\right) \bar{q}_{i}}\right]
\label{eq.gradFD1I}
\end{align}

\section{Divergence of P.K.Bathia and S.Singh.}
\begin{figure}[h!]
\centering
\includegraphics[width=0.7\linewidth]{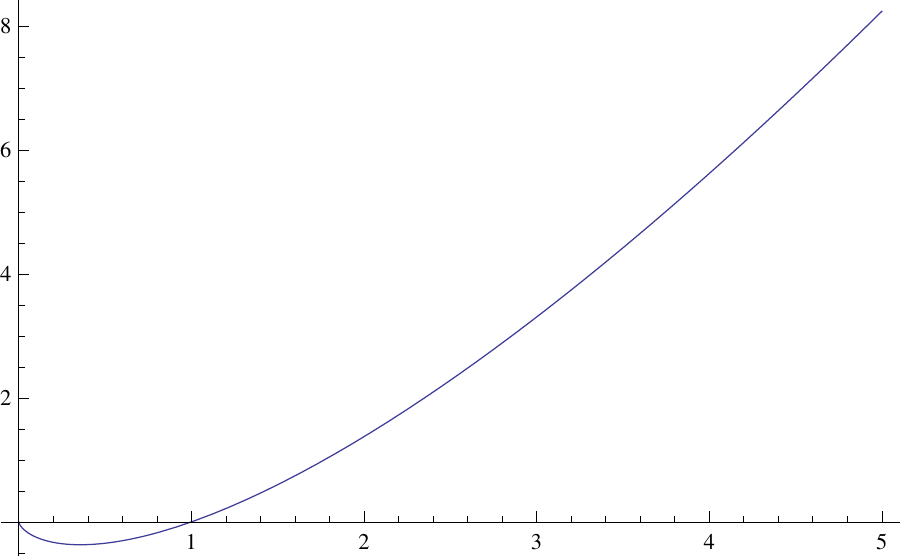}
\caption{Function $f\left(x\right)=\frac{x \sinh\left(\alpha \log x\right)}{\sinh \alpha}$, $\alpha=0.5$}
\label{fig:BS}
\end{figure}

This divergence proposed in \cite{bhatia2013} is basically a Csiszär divergence, but in this work, there are some mistakes: the basic convex function, represented on the figure (\ref{fig:BS}), for $\alpha=0.5$, is not a standard function, indeed, the authors write it:
\begin{equation}
	f\left(x\right)=\frac{x \sinh\left(\alpha \log x\right)}{\sinh \alpha}
\end{equation}

\begin{figure}[h!]
\centering
\includegraphics[width=0.7\linewidth]{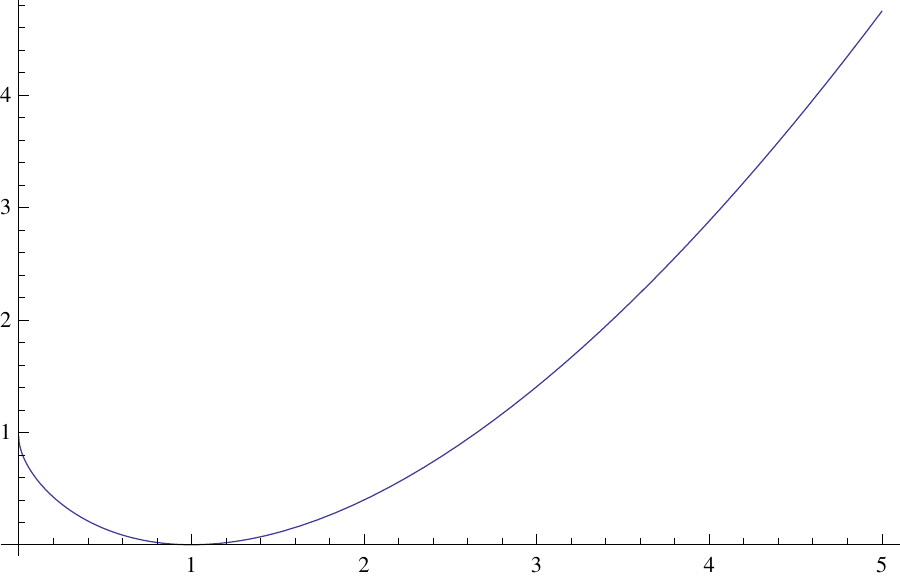}
\caption{Function $f_{c}\left(x\right)=\frac{x \sinh\left(\alpha \log x\right)-\alpha x+\alpha}{\sinh \alpha}$, $\alpha=0.5$}
\label{fig:BSc}
\end{figure}

We can note that the multiplicative factor $\frac{1}{\sinh\alpha}$ is not really helpful.\\
For this function, we do have convexity and $f\left(1\right)=0$, but we don't have $f'\left(1\right)=0$, contrary to what is stated in \cite{bhatia2013}.\\
If we still use this function to construct a Csiszär divergence, we obtain:
\begin{equation}
	SB_{\alpha}\left(p\|q\right)=\sum_{i}p_{i}\frac{\sinh\left(\alpha\log \frac{p_{i}}{q_{i}}\right)}{\sinh\alpha}
\end{equation}
The problem with this expression is classical: the gradient with respect to ``$q$" does not cancel for $p_{i}=q_{i}\ \forall i$.\\
To eliminate this kind of problem, it is necessary to introduce the corresponding standard convex function, represented in the figure (\ref{fig:BSc}), for $\alpha=0.5$, which is written as follows:
\begin{equation}
	f_{c}\left(x\right)=\frac{x \sinh\left(\alpha \log x\right)-\alpha x+\alpha}{\sinh \alpha}
	\label{eq.fcBS}
\end{equation}

We are thus led to the divergence:
\begin{equation}
	SB_{c,\alpha}\left(p\|q\right)=\frac{1}{\sinh\alpha}\left[\sum_{i}p_{i}\sinh\left(\alpha\log \frac{p_{i}}{q_{i}}\right)+\alpha\left(q_{i}-p_{i}\right]\right]
	\label{eq.SBc}
\end{equation}
Whose gradient with respect to ``$q$" is:
\begin{equation}
	\frac{\partial SB_{c,\alpha}\left(p\|q\right)}{\partial q_{j}}\approx-\frac{p_{j}}{q_{j}}\cosh\left(\alpha\log\frac{p_{j}}{q_{j}}\right)+1
\end{equation}
The dual divergence built on the mirror function of (\ref{eq.fcBS}) will be written as follows:
\begin{equation}
	SB_{c,\alpha}\left(q\|p\right)=\frac{1}{\sinh\alpha}\left[\sum_{i}q_{i}\sinh\left(\alpha\log \frac{q_{i}}{p_{i}}\right)+\alpha\left(p_{i}-q_{i}\right]\right]
\end{equation}
and its gradient with respect to ``$q$" is:
\begin{equation}
		\frac{\partial SB_{c,\alpha}\left(q\|p\right)}{\partial q_{j}}=\frac{1}{\sinh\alpha}\left[\sinh\left(\alpha\log \frac{q_{j}}{p_{j}}\right)+\alpha\cosh\left(\alpha\log \frac{q_{j}}{p_{j}}\right)-\alpha\right]
\end{equation}

\subsection{Invariance by scale change on ``$q$".}
For this divergence given by the expression (\ref{eq.SBc}), the nominal invariance factor cannot be calculated, so we use the expression $K^{*}$, which leads to the simplified invariant divergence:
\begin{equation}
	SB_{c,\alpha}\left(p\|q\right)I=\frac{1}{\sinh\alpha}\left[\sum_{i}\bar{p}_{i}\sinh\left(\alpha\log \frac{\bar{p}_{i}}{\bar{q}_{i}}\right)\right]
	\label{eq.SBI}
\end{equation}
Its gradient with respect to ``$q$" is:
\begin{equation}
	\frac{\partial SB_{c,\alpha}\left(p\|q\right)I}{\partial q_{l}}=\frac{\alpha\coth\alpha}{\sum_{j}q_{j}}\left[\sum_{i}\bar{p}_{i}\log\frac{\bar{p}_{i}}{\bar{q}_{i}}
	-\frac{\bar{p}_{l}}{\bar{q}_{l}}\log\frac{\bar{p}_{l}}{\bar{q}_{l}}\right]
	\label{eq.gradSBI}
\end{equation}

\setcounter{table}{0}  \setcounter{equation}{0}  \setcounter{figure}{0} \setcounter{chapter}{7} \setcounter{section}{0} 
\chapter{chapter 7 -\\Divergences between means \label{chptr::chapitre7}}  

The major part of this chapter can be found in the book by I.J.Taneja (online) \cite{taneja2001} and in the article by Ben-Tal et al. \cite{ben1989}.
The construction mode of these divergences is twofold: either one relies on the well-known inequalities between the classical means, or one constructs them in the sense of Csiszär by relying on the properties of the convex functions.\\
We will make this double aspect more explicit.\\
In addition, in this chapter, we will present the logarithmic forms of these divergences.\\
The invariant forms of these divergences will be developed by introducing the invariance factor $K^{*}(p,q)=\frac{\sum_{j}p_{j}}{\sum_{j}q_{j}}$ which is not necessarily the nominal invariance factor.\\
 The particular case of $M_{AG}$ for which the nominal invariance factor is explicitly calculated has been presented in the chapter dealing with "Alpha Divergences".\\
Some basic recalls on means are summarized in Appendix 2.

\section{Square root (Quadratic) - Arithmetic mean divergence.}
This divergence is founded on the inequality: $M_{S}\geq M_{A}$.

\subsection{Unweighted version.}
The most classical form of this divergence is constructed in the sense of Csiszär on the basis of the standard convex function:
\begin{equation}
	f_{c}\left(x\right)=\sqrt{\frac{x^{2}+1}{2}}-\frac{x+1}{2}
\end{equation}
We immediately obtain:
\begin{equation}	M_{SA}\left(p\|q\right)=\sum_{i}\sqrt{\frac{p_{i}^{2}+q_{i}^{2}}{2}}-\sum_{i}\frac{p_{i}+q_{i}}{2}
\end{equation}
This divergence is symmetrical.\\
Of course, if one is in a situation where one knows that $\sum_{i}p_{i}=\sum_{i}q_{i}$, (but not necessarily $=1$), this divergence becomes simpler, it becomes:
\begin{equation}
	M1_{SA}\left(p\|q\right)=\sum_{i}\sqrt{\frac{p_{i}^{2}+q_{i}^{2}}{2}}-\sum_{i}q_{i}
\end{equation}
That is to say, it is constructed in the sense of Csiszär on the simple convex function:
\begin{equation}
	f_{1}\left(x\right)=\sqrt{\frac{x^{2}+1}{2}}-1
\end{equation}
But it can also become:
\begin{equation}
	M2_{SA}\left(p\|q\right)=\sum_{i}\sqrt{\frac{p_{i}^{2}+q_{i}^{2}}{2}}-\sum_{i}p_{i}
\end{equation}
It is then built on the simple convex function:
\begin{equation}
	f_{2}\left(x\right)=\sqrt{\frac{x^{2}+1}{2}}-x
\end{equation}
We can exhibit a family of simple convex functions $f\left(x\right)$ that will give simplified divergences. Knowing $f_{c}\left(x\right)$, these functions are deduced from each other by modifying the value of $f'\left(1\right)$ in the expression:
\begin{equation}
	f\left(x\right)=f_{c}\left(x\right)+\left(x-1\right)f'\left(1\right)
\end{equation}
If we now consider the situation where $\sum_{i}p_{i}=\sum_{i}q_{i}=1$ we obtain in all cases an oversimplified form which of course, we are not always allowed to use; this form is written as follows:
\begin{equation}
		MS_{SA}\left(p\|q\right)=\sum_{i}\sqrt{\frac{p_{i}^{2}+q_{i}^{2}}{2}}-1
\end{equation}
It should be observed that this divergence cannot be constructed directly as a Csiszär divergence because the simplification $\sum_{i}p_{i}=\sum_{i}q_{i}=1$ cannot be taken into account until the divergence is exhibited.\\
Finally, it is worth noting that:
\begin{equation}
	f_{c}\left(x\right)=\frac{f1\left(x\right)+f2\left(x\right)}{2}
\end{equation}
This explains the symmetry of the divergence $M_{SA}$.

\subsection{Weighted version.}
\begin{figure}[h!]
\centering
\includegraphics[width=0.7\linewidth]{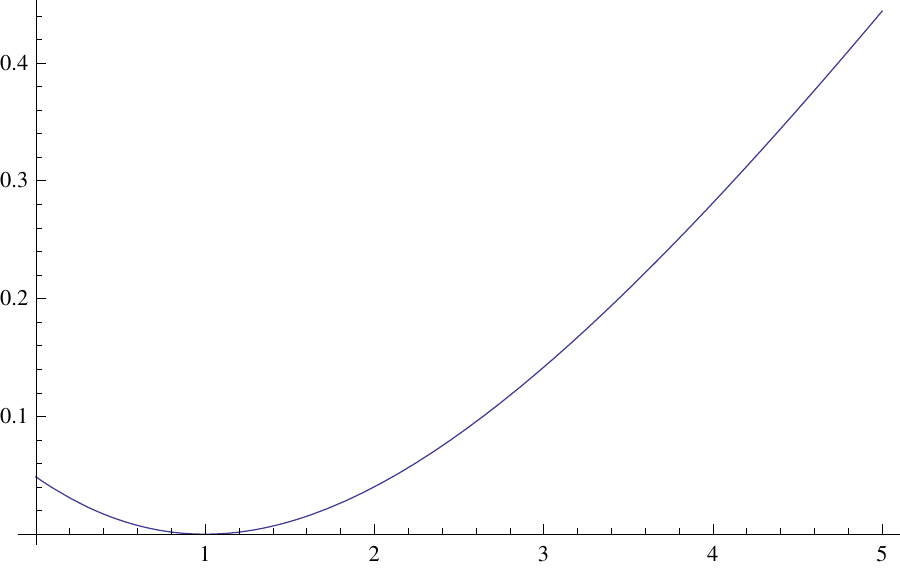}
\caption{Function $f_{c\alpha}\left(x\right)=\sqrt{\alpha x^{2}+1-\alpha}-\left(\alpha x+1-\alpha\right)$; $\alpha=0.5$}
\label{fig:MSAPc}
\end{figure}

The standard convex function allowing to construct this divergence is shown in the figure (\ref{fig:MSAPc}) for $\alpha=0.5$.
\\
It is written with $0<\alpha<1$:
\begin{equation}
	f_{c\alpha}\left(x\right)=\sqrt{\alpha x^{2}+1-\alpha}-\left(\alpha x+1-\alpha\right)
\end{equation}
This results in:
\begin{equation}
	M_{SA\alpha}\left(p\|q\right)=\sum_{i}\sqrt{\alpha p_{i}^{2}+\left(1-\alpha\right)q_{i}^{2}}-\sum_{i}\left[\alpha p_{i}+\left(1-\alpha\right)q_{i}\right]
	\label{eq.MSAa}
\end{equation}
This divergence is not symmetrical because of the weighting (unless $\alpha=1/2$).\\
Its gradient with respect to ``$q$" is written:
\begin{equation}
	\frac{\partial 	M_{SA\alpha}}{\partial q_{j}}=\left(1-\alpha\right)\left[\frac{q_{j}}{\sqrt{\alpha p_{j}^{2}+\left(1-\alpha\right)q_{j}^{2}}}-1\right]
\end{equation}

\subsubsection{Scale invariant version with respect to ``$q$".}
For this divergence, the nominal invariance factor is not obtained explicitly, so we use as the invariance factor the expression $K^{*}(p,q)$; thus we obtain the expression for the corresponding simplified invariant divergence, which is written as:
\begin{equation}
	M_{SA\alpha}I\left(p\|q\right)=\sum_{j}p_{j}\left\{\sum_{i}\sqrt{\alpha \bar{p}_{i}^{2}+\left(1-\alpha\right)\bar{q}_{i}^{2}}-\sum_{i}\left[\alpha \bar{p}_{i}+\left(1-\alpha\right)\bar{q}_{i}\right]\right\}
	\label{eq.MSAIa}
\end{equation}
After simplification and deleting the multiplying factor, we will have:
\begin{equation}
	M_{SA\alpha}I\left(p\|q\right)=\sum_{i}\sqrt{\alpha \bar{p}_{i}^{2}+\left(1-\alpha\right)\bar{q}_{i}^{2}}-1
	\label{eq.MSAIas}
\end{equation}
The gradient with respect to ``$q$" is given by:
\begin{equation}
	\frac{\partial M_{SA\alpha}I\left(p\|q\right)}{\partial q_{l}}=\frac{1-\alpha}{\sum_{j}q_{j}}\left\{\left[\alpha \bar{p}_{l}^{2}+\left(1-\alpha\right)\bar{q}_{l}^{2}\right]^{-\frac{1}{2}}\bar{q}_{l}-\sum_{i}\left[\alpha \bar{p}_{i}^{2}+\left(1-\alpha\right)\bar{q}^{2}_{i}\right]^{-\frac{1}{2}}\bar{q}^{2}_{i}\right\}
	\label{eq.gradMSAIas}
\end{equation}

\subsection{Logarithmic form.}
Since the divergence (\ref{eq.MSAa}) appears as the difference of 2 positive quantities, one can apply on each term an increasing function without changing the sign of the expression, for example the function ``Generalized Logarithm" and at the limit the ``Logarithm", which allows to write:
\begin{align}
	LM_{SA\alpha}=\log\left[\sum_{i}\sqrt{\alpha p_{i}^{2}+\left(1-\alpha\right)q_{i}^{2}}\right] -\log\left[\sum_{i}\alpha p_{i}+\left(1-\alpha\right)q_{i}\right]
\end{align}
With such a divergence, nothing is clear about the convexity, so nothing is guaranteed; however, we can calculate the gradient with respect to ``$q$": 
\begin{align}
	\frac{\partial LM_{SA\alpha}}{\partial q_{j}}=&\frac{\left(1-\alpha\right)}{\sum_{i}\sqrt{\alpha p_{i}^{2}+\left(1-\alpha\right)q_{i}^{2}}}\frac{q_{j}}{\sqrt{\alpha p_{j}^{2}+\left(1-\alpha\right)q_{j}^{2}}} \nonumber \\ &-\frac{\left(1-\alpha\right)}{\sum_{i}\alpha p_{i}+\left(1-\alpha\right)q_{i}}
\end{align}
This expression is zero if $p_{i}=q_{i}\ \forall i$.

\subsubsection{Logarithmic form invariant with respect to ``$q$".}
To obtain such a Logarithmic form we apply the function ``Log" separately on each term of (\ref{eq.MSAIas}) and we obtain:
\begin{equation}
LM_{SA\alpha}I\left(p\|q\right)=\log\sum_{i}\sqrt{\alpha \bar{p}^{2}_{i}+\left(1-\alpha\right)\bar{q}^{2_{i}}}
\label{eq.LMSAIas}	
\end{equation}
which gradient with respect to ``$q$" is:
\begin{equation}
	\frac{\partial LM_{SA\alpha}I\left(p\|q\right)}{\partial q_{l}}=\frac{1}{\sum_{i}\sqrt{\alpha \bar{p}^{2}_{i}+\left(1-\alpha\right)\bar{q}^{2}_{i}}}\frac{\partial M_{SA\alpha}I\left(p\|q\right)}{\partial q_{l}}
	\label{eq.gradLMSAIas}
\end{equation}

\section{Square root (Quadratic)- Géométric mean divergence.}
This divergence is based on the inequality: $M_{S}\geq M_{G}$.
\subsection{Unweighted version.}
It is a Csiszär divergence built on the standard convex function:
\begin{equation}
 f_c(x)=\sqrt{\frac{x^{2}+1}{2}}-\sqrt{x} 
\end{equation}
This results in the divergence:
\begin{equation}
M_{SG}\left(p\|q\right)=\sum_{i}\sqrt{\frac{p_{i}^{2}+q_{i}^{2}}{2}}-\sum_{i}\sqrt{p_{i}q_{i}}
\label{eq:MSG}
\end{equation}

\subsection{Weighted version.}
\begin{figure}[h!]
\centering
\includegraphics[width=0.7\linewidth]{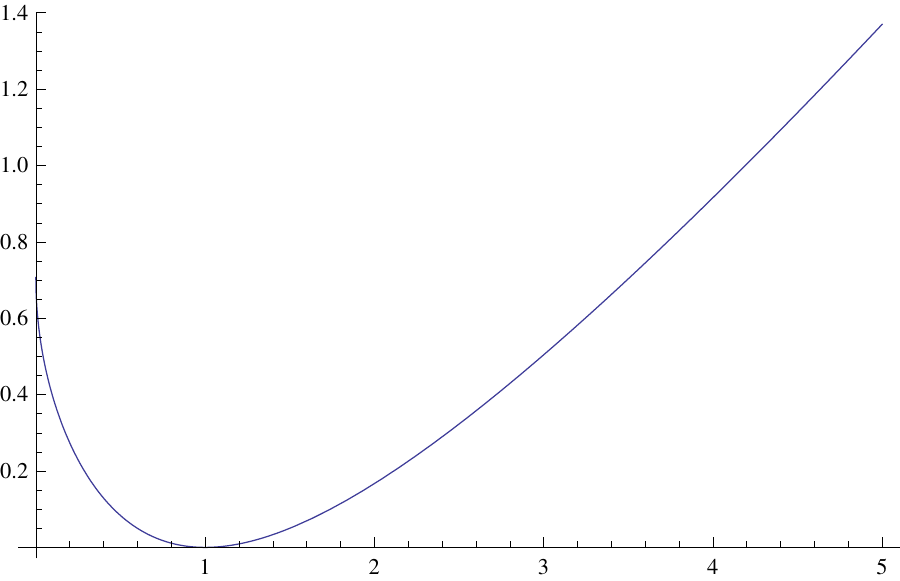}
\caption{Function $f_{c}\left(x\right)=\sqrt{\alpha x^{2}+\left(1-\alpha\right)}-x^{\alpha}$; $\alpha=0.5$}
\label{fig:MSGPc}
\end{figure}

With $0\leq \alpha \leq 1$, this divergence is based on the standard convex function: 
\begin{equation}
	f_{c}\left(x\right)=\sqrt{\alpha x^{2}+\left(1-\alpha\right)}-x^{\alpha}
\end{equation}
shown in the figure (\ref{fig:MSGPc}) for $\alpha=0.5$.
\\
The corresponding divergence will be written:
\begin{equation}
	M_{SG\alpha}\left(p\|q\right)=\sum_{i}\sqrt{\alpha p_{i}^{2}+\left(1-\alpha\right)q_{i}^{2}}-\sum_{i}p_{i}^{\alpha}q_{i}^{1-\alpha}
	\label{eq.MSGA}
\end{equation}
we deduce the gradient with respect to ``$q$":
\begin{equation}
	\frac{\partial 	M_{SG\alpha}}{\partial q_{j}}=\left(1-\alpha\right)\left[\frac{q_{j}}{\sqrt{\alpha p_{j}^{2}+\left(1-\alpha\right) q_{j}^{2}}}-p_{j}^{\alpha}q_{j}^{-\alpha}\right]
\end{equation}

\subsubsection{Scale invariant form with respect to ``$q$".}
For this divergence, the nominal invariance factor cannot be obtained explicitly, so we use the expression $K^{*}(p,q)$ as the invariance factor; thus, after simplification, we obtain the expression of the corresponding invariant divergence:
\begin{equation}
	M_{SG\alpha}I\left(p\|q\right)=\sum_{i}\sqrt{\alpha \bar{p}^{2}_{i}+\left(1-\alpha\right)\bar{q}^{2}_{i}}-\sum_{i}\bar{p}^{\alpha}_{i}\bar{q}^{1-\alpha}_{i}
	\label{eq.MSGIa}
\end{equation}
which has a gradient with respect to ``$q$" given by:
\begin{align}
\frac{\partial M_{SG\alpha}I\left(p\|q\right)}{\partial  q_{l}}=\frac{1-\alpha}{\sum_{j}q_{j}}&\left\{\frac{\bar{q}_{l}}{\sqrt{\alpha \bar{p}^{2}_{l}+\left(1-\alpha\right)\bar{q}^{2}_{l}}}-\sum_{i}\frac{\bar{q}^{2}_{i}}{\sqrt{\alpha \bar{p}^{2}_{i}+\left(1-\alpha\right)\bar{q}^{2}_{i}}}\right. \nonumber \\ & \left.-\bar{p}^{\alpha}_{l}\bar{q}^{-\alpha}_{l}+\sum_{i}\bar{p}^{\alpha}_{i}\bar{q}^{1-\alpha}_{i}\right\}	
\label{eq.gradMSGIa}
\end{align}

\subsection{Logarithmic form.}
After applying an increasing function (Log), on each of the terms of the difference appearing in (\ref{eq.MSGA}), the corresponding divergence will be written:
\begin{equation}
	LM_{SG\alpha}\left(p\|q\right)=\log\sum_{i}\sqrt{\alpha p_{i}^{2}+\left(1-\alpha\right)q_{i}^{2}}-\log\sum_{i}p_{i}^{\alpha}q_{i}^{1-\alpha}
	\label{eq.LMSGA}
\end{equation}
The gradient with respect to ``$q$" is given by:
\begin{align}
	\frac{\partial LM_{SG\alpha}}{\partial q_{j}}=\left(1-\alpha\right)&\left[\frac{1}{\sum_{i}\sqrt{\alpha p_{i}^{2}+\left(1-\alpha\right)q_{i}^{2}}}\frac{q_{j}}{\sqrt{\alpha p_{j}^{2}+\left(1-\alpha\right) q_{j}^{2}}}\right. \nonumber \\ & \left.-\frac{p_{j}^{\alpha}q_{j}^{-\alpha}}{\sum_{i}p_{i}^{\alpha}q_{i}^{1-\alpha}}\right]
	\label{eq.gradLMSGA}	
\end{align}
This expression is zero if $p_{i}=q_{i}\ \forall i$.

\subsubsection{Logarithmic form - Invariance with respect to ``$q$".}
The logarithmic version invariant with respect to "$q$" is obtained by applying the function ``Log" separately to each term of (\ref{eq.MSGIa}) and one obtains:
\begin{equation}
	LM_{SG\alpha}I\left(p\|q\right)=\log\sum_{i}\sqrt{\alpha \bar{p}^{2}_{i}+\left(1-\alpha\right)\bar{q}^{2}_{i}}-\log\sum_{i}\bar{p}^{\alpha}_{i}\bar{q}^{1-\alpha}_{i}
	\label{eq.LMSGIa}
\end{equation}
whose gradient is written as:
\begin{align}
\frac{\partial LM_{SG\alpha}I\left(p\|q\right)}{\partial  q_{l}}=&\frac{1-\alpha}{\sum_{j}q_{j}}\nonumber\\
\Bigg\{&\frac{1}{\sum_{i}\sqrt{\alpha \bar{p}^{2}_{i}+\left(1-\alpha\right)\bar{q}^{2}_{i}}}\left[\frac{\bar{q}_{l}}{\sqrt{\alpha \bar{p}^{2}_{l}+\left(1-\alpha\right)\bar{q}^{2}_{l}}}-\sum_{i}\frac{\bar{q}^{2}_{i}}{\sqrt{\alpha \bar{p}^{2}_{i}+\left(1-\alpha\right)\bar{q}^{2}_{i}}}\right]\nonumber\\&
+1-\frac{\bar{p}^{\alpha}_{l}\bar{q}^{-\alpha}_{l}}{\sum_{i}\bar{p}^{\alpha}_{i}\bar{q}^{1-\alpha}_{i}}\Bigg\}
\label{eq.gradLMSGIa}	
\end{align}

\section{Square root (Quadratic) - Harmonic mean divergence.}
This divergence is based on the inequality: $M_{S}\geq M_{H}$.
\subsection{Unweighted version.}
It is a Csiszär divergence based on the standard convex function:
\begin{equation}
	f_{c}\left(x\right)=\sqrt{\frac{x^{2}+1}{2}}-\frac{2x}{1+x}
\end{equation}
Which leads to the divergence:
\begin{equation}
	M_{SH}\left(p\|q\right)=\sum_{i}\sqrt{\frac{p_{i}^{2}+q_{i}^{2}}{2}}-\sum_{i}\frac{2p_{i}q_{i}}{p_{i}+q_{i}}
\end{equation}

\subsection{Weighted version.}
\begin{figure}[h!]
\centering
\includegraphics[width=0.7\linewidth]{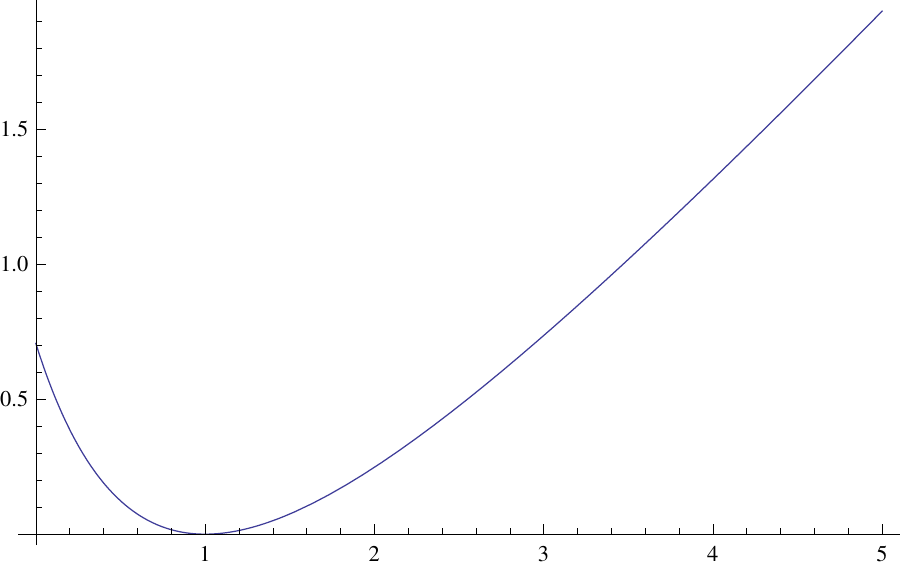}
\caption{Function $f_{c \alpha}\left(x\right)=\sqrt{\alpha x^{2}+\left(1-\alpha\right)}-\frac{x}{\left(1-\alpha\right)x+\alpha}$; $\alpha=0.5$}
\label{fig:MSHPc}
\end{figure}

With $0\leq \alpha \leq 1$, it is built in the Csiszär sense using the standard convex function :
\begin{equation}
	f_{c \alpha}\left(x\right)=\sqrt{\alpha x^{2}+\left(1-\alpha\right)}-\frac{x}{\left(1-\alpha\right)x+\alpha}
\end{equation}\\
It is shown in the figure (\ref{fig:MSHPc}) for $\alpha=0.5$.
\\
This results in the divergence:
\begin{equation}
	M_{SH \alpha}\left(p\|q\right)=\sum_{i}\sqrt{\alpha p_{i}^{2}+\left(1-\alpha\right)q_{i}^{2}}-\sum_{i}\frac{p_{i}q_{i}}{\left(1-\alpha\right)p_{i}+\alpha q_{i}}
	\label{eq.MSHA}
\end{equation}
The gradient with respect to ``$q$" is given by:
\begin{equation}
	\frac{\partial 	M_{SH\alpha}}{\partial q_{j}}=\left(1-\alpha\right)\left\{\frac{q_{j}}{\sqrt{\alpha p_{j}^{2}+\left(1-\alpha\right)q_{j}^{2}}}-\frac{p_{j}^{2}}{\left[\left(1-\alpha\right)p_{j}+\alpha q_{j}\right]^{2}}\right\}
\end{equation}
This expression is zero if $p_{i}=q_{i}\ \forall i$.

\subsubsection{Scale invariant version with respect to ``$q$".}
For this divergence, the nominal invariance factor is not obtained explicitly, so we use as the invariance factor the expression $K^{*}(p,q)$; thus, after simplification, we obtain the form of the invariant corresponding divergence:
\begin{equation}
	M_{SH\alpha}I\left(p\|q\right)=\sum_{i}\sqrt{\alpha \bar{p}^{2}_{i}+(1-\alpha)\bar{q}^{2}_{i}}-\sum_{i}\frac{\bar{p}_{i}\bar{q}_{i}}{\left(1-\alpha\right) \bar{p}_{i}+\alpha\bar{q}_{i}}
	\label{eq.MSHIa}
\end{equation}
The gradient with respect to ``$q$" is given by:

\begin{align}
\frac{\partial 	M_{SH\alpha}I\left(p\|q\right)}{\partial q_{l}}=&\frac{\left(1-\alpha\right)}{\sum_{j}q_{j}}\nonumber\\
\Bigg\{&\frac{\bar{q}_{l}}{\sqrt{\alpha \bar{p}^{2}_{l}+(1-\alpha)\bar{q}^{2}_{l}}}-\sum_{i}\frac{\bar{q}^{2}_{i}}{\sqrt{\alpha \bar{p}^{2}_{i}+(1-\alpha)\bar{q}^{2}_{i}}}\nonumber\\&
-\frac{\bar{p}^{2}_{l}}{\left[\left(1-\alpha\right) \bar{p}_{l}+\alpha\bar{q}_{l}\right]^{2}}+\sum_{i}\frac{\bar{p}^{2}_{i}\bar{q}_{i}}{\left[\left(1-\alpha\right) \bar{p}_{i}+\alpha\bar{q}_{i}\right]^{2}}\Bigg\}	
\end{align}

\subsection{Logarithmic version.}
We can write directly from (\ref{eq.MSHA}):
\begin{equation}
	LM_{SH}\left(p\|q\right)=\log\sum_{i}\sqrt{\alpha p_{i}^{2}+\left(1-\alpha\right)q_{i}^{2}}-\log\sum_{i}\frac{p_{i}q_{i}}{\left(1-\alpha\right)p_{i}+\alpha q_{i}}
\end{equation}
Its gradient with respect to ``$q$" is written as follows:
\begin{align}
\frac{\partial 	LM_{SH\alpha}}{\partial q_{j}}=\left(1-\alpha\right)&\left\{\frac{1}{\sum_{i}\sqrt{\alpha p_{i}^{2}+\left(1-\alpha\right)q_{i}^{2}}}\frac{q_{j}}{\sqrt{\alpha p_{j}^{2}+\left(1-\alpha\right)q_{j}^{2}}}\right. \nonumber \\ & \left.-\frac{1}{\sum_{i}\frac{p_{i}q_{i}}{\left(1-\alpha\right)p_{i}+\alpha q_{i}}}\frac{p_{j}^{2}}{\left[\left(1-\alpha\right)p_{j}+\alpha q_{j}\right]^{2}}\right\}	
\end{align}
This expression is zero if $p_{i}=q_{i}\ \forall i$.

\subsubsection{Logarithmic form - Invariance with respect to ``$q$".}
The logarithmic invariant version is obtained by applying the function ``Log" separately to each term of (\ref{eq.MSHIa}) and one obtains:
\begin{equation}
	LM_{SH\alpha}I\left(p\|q\right)=\log\sum_{i}\sqrt{\alpha \bar{p}^{2}_{i}+(1-\alpha)\bar{q}^{2}_{i}}-\log\sum_{i}\frac{\bar{p}_{i}\bar{q}_{i}}{\left(1-\alpha\right) \bar{p}_{i}+\alpha\bar{q}_{i}}	
\end{equation}
Its gradient with respect to ``$q$" is given by:
\begin{align}
\frac{\partial 	LM_{SH\alpha}I\left(p\|q\right)}{\partial q_{l}}=&\frac{\left(1-\alpha\right)}{\sum_{j}q_{j}}
\left\{\frac{1}{\sum_{i}\sqrt{\alpha \bar{p}^{2}_{i}+(1-\alpha)\bar{q}^{2}_{i}}}\right. \nonumber \\  & \left.\left\{\frac{\bar{q}_{l}}{\sqrt{\alpha \bar{p}^{2}_{l}+(1-\alpha)\bar{q}^{2}_{l}}}-\sum_{i}\frac{\bar{q}^{2}_{i}}{\sqrt{\alpha \bar{p}^{2}_{i}+(1-\alpha)\bar{q}^{2}_{i}}}\right\}\right. \nonumber \\  & \left.
-\frac{1}{\sum_{i}\frac{\bar{p}_{i}\bar{q}_{i}}{\left(1-\alpha\right) \bar{p}_{i}+\alpha\bar{q}_{i}}}\left\{\frac{\bar{p}^{2}_{l}}{\left[\left(1-\alpha\right) \bar{p}_{l}+\alpha\bar{q}_{l}\right]^{2}}+\sum_{i}\frac{\bar{p}^{2}_{i}\bar{q}_{i}}{\left[\left(1-\alpha\right) \bar{p}_{i}+\alpha\bar{q}_{i}\right]^{2}}\right\}\right\}	
\end{align}

\section{Arithmetic - Geometric mean divergence}.
This divergence is at the basis of the ``Alpha divergences" of Amari \cite{amari2009}; here we develop it from an elementary version. An extended version has been considered in the chapter dealing specifically with ``Alpha divergences".\\
This divergence is based on the inequality: $M_{A}\geq M_{G}$.
\subsection{Basis unweighted version.}
In the simplest version, it is a Csiszär divergence built on the standard convex function:
\begin{equation}
	f_{c}\left(x\right)=\frac{1}{2}\left(\sqrt{x}-1\right)^{2}
\end{equation}
We easily obtain:
\begin{equation}
	M_{AG}\left(p\|q\right)=\frac{1}{2}\sum_{i}\left(\sqrt{p_{i}}-\sqrt{q_{i}}\right)^{2}
\end{equation}
This symmetrical divergence is also known as Hellinger's divergence \cite{beran1977}, but it can also be written:
\begin{equation}
	M_{AG}\left(p\|q\right)=\sum_{i}\frac{p_{i}+q_{i}}{2}-\sum_{i}\sqrt{p_{i}q_{i}}
\end{equation}
So we better understand the name Arithmetic-Geometric mean divergence.\\
The simplification $\sum_{i}p_{i}=\sum_{i}q_{i}$ allows us to obtain 2 simplified divergences:
\begin{equation}
	MS1_{AG}\left(p\|q\right)=\sum_{i}q_{i}-\sum_{i}\sqrt{p_{i}q_{i}}
\end{equation}
Which is deduced in the sense of Csiszär from the simple convex function: $f_{1}\left(x\right)=1-\sqrt{x}$, and:
\begin{equation}
	MS2_{AG}\left(p\|q\right)=\sum_{i}p_{i}-\sum_{i}\sqrt{p_{i}q_{i}}
\end{equation}
deduced in the sense of Csiszär from the simple convex function: $f_{2}\left(x\right)=x-\sqrt{x}$.\\
The supplementary simplification$\sum_{i}p_{i}=\sum_{i}q_{i}=1$, leads to:
\begin{equation}
	MS_{AG}\left(p\|q\right)=1-\sum_{i}\sqrt{p_{i}q_{i}}
\end{equation}
This expression is not a divergence of Csiszär; to obtain it, one must first construct a Csiszär divergence by using one of the simple convex functions $f_{1}\left(x\right)=1-\sqrt{x}$ or $f_{2}\left(x\right)=x-\sqrt{x}$ (functions mirroring each other) deduced from $f_{c}\left(x\right)$, then introduce ``a posteriori" the simplification $\sum_{i}p_{i}=1$ or $\sum_{i}q_{i}=1$.

\subsection{Weighted version.}
\begin{figure}[h!]
\centering
\includegraphics[width=0.7\linewidth]{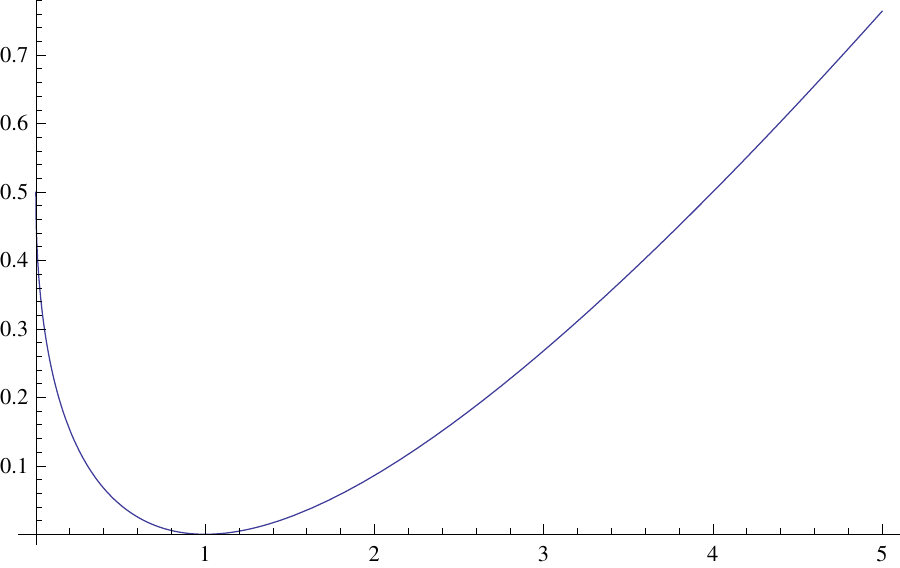}
\caption{Function $f_{c\alpha}\left(x\right)=\alpha x+\left(1-\alpha\right)-x^{\alpha}$, $\alpha=0.5$}
\label{fig:MAGPc}
\end{figure}

If one does not introduce a multiplicative factor whose sign depends on the weighting factor $0\leq\alpha\leq1$, the basic standard convex basic function will be written: 
\begin{equation}
	f_{c\alpha}\left(x\right)=\alpha x+\left(1-\alpha\right)-x^{\alpha}
\end{equation}\\
It is represented in the figure (\ref{fig:MAGPc}) for $\alpha=0.5$.
\\
The corresponding divergence is written:
\begin{equation}
	M_{AG\alpha}\left(p\|q\right)=\sum_{i}\alpha p_{i}+\left(1-\alpha\right)q_{i}-\sum_{i}p_{i}^{\alpha}q_{i}^{1-\alpha}
\end{equation}
To broaden the range of validity of this divergence to all positive values of ``$\alpha$", we will use:
\begin{equation}
	M_{AG\alpha}\left(p\|q\right)=\frac{1}{1-\alpha}\left[\sum_{i}\alpha p_{i}+\left(1-\alpha\right)q_{i}-\sum_{i}p_{i}^{\alpha}q_{i}^{1-\alpha}\right]
	\label{eq.MAGA}
\end{equation}
This is the Havrda-Charvat divergence \cite{havrda1967}, and this is exactly what is written for the ``$\alpha$" divergence.\\
There's no point in commenting on these things that are developed elsewhere.\\
The gradient with respect to ``$q$'' is written:
\begin{equation}
	\frac{\partial 	M_{AG\alpha}}{\partial q_{j}}=1-\left(\frac{p_{j}}{q_{j}}\right)^{\alpha}
\end{equation}

\subsubsection{Invariant form with respect to ``$q$".}
For such a divergence, the nominal invariance factor is obtained explicitly, the corresponding invariant form has been developed in the chapter dealing with ``Alpha Divergence".\\
 For the sake of homogeneity, we present here the case of an invariance factor of the form $K^{*}(p,q)$; we thus obtain the following expression for the corresponding simplified invariant divergence:
\begin{equation}
	M_{AG\alpha}I\left(p\|q\right)=\frac{1}{1-\alpha}\left[\sum_{i}\left[\alpha \bar{p}_{i}+\left(1-\alpha\right)\bar{q}_{i}\right]-\sum_{i}\bar{p}^{\alpha}_{i}\bar{q}^{1-\alpha}_{i}\right]
	\label{eq.MAGIa}
\end{equation}
Which simplifies into:
\begin{equation}	M_{AG\alpha}I\left(p\|q\right)=\frac{1}{1-\alpha}\left[1-\sum_{i}\bar{p}^{\alpha}_{i}\bar{q}^{1-\alpha}_{i}\right]
	\label{eq.MAGIas}
\end{equation}
The gradient with respect to ``$q$'' is written:
\begin{equation}
	\frac{\partial M_{AG\alpha}I\left(p\|q\right)}{\partial q_{l}}=\frac{1}{\sum_{j}q_{j}}\left(\bar{p}^{\alpha}_{l}\bar{q}^{\;-\alpha}_{l}-\sum_{i}\bar{p}^{\alpha}_{i}\bar{q}^{1-\alpha}_{i}\right)
	\label{eq.gradMAGIas}
\end{equation}

\subsection{Logarithmic form.}
The logarithmic form of the divergence (\ref{eq.MAGA}) can be written as follows:
\begin{equation}
	LM_{AG\alpha}\left(p\|q\right)=\frac{1}{\alpha-1}\left[\log\sum_{i}p_{i}^{\alpha}q_{i}^{1-\alpha}-\log\sum_{i}\alpha p_{i}+\left(1-\alpha\right)q_{i}\right]
	\label{eq.LMag}
\end{equation}
The gradient with respect to ``$q$'' is written:
\begin{equation}
	\frac{\partial 	LM_{AG\alpha}}{\partial q_{j}}=\frac{1}{\sum_{i}\alpha p_{i}+\left(1-\alpha\right)q_{i}}-\frac{ p_{j}^{\alpha} q_{j}^{-\alpha}}{\sum_{i} p_{i}^{\alpha} q_{i}^{1-\alpha}}
\end{equation}
If we consider that the data fields involved are probability densities, that is to say, if $\sum_{i}p_{i}=\sum_{i}q_{i}=1$, we obtain from (\ref{eq.LMag}) a simplified expression which is the classical Renyi divergence \cite{renyi1961}:
\begin{equation}
	R_{\alpha}\left(p\|q\right)=\frac{1}{\alpha-1}\left[\log\sum_{i}p_{i}^{\alpha}q_{i}^{1-\alpha}\right]
\end{equation}

\subsubsection{Logarithmic form - Invariance with respect to ``$q$".}
The invariant logarithmic version is obtained by applying the function ``Log" separately to each term of (\ref{eq.MAGIas}) and we have:
\begin{equation}
LM_{AG\alpha}I\left(p\|q\right)=\frac{1}{\alpha-1}\log\sum_{i}\bar{p}^{\alpha}_{i}\bar{q}^{1-\alpha}_{i}
\label{eq.LMAGIas}	
\end{equation}
His gradient with respect to $q$ is expressed as:
\begin{equation}
	\frac{\partial LM_{AG\alpha}I\left(p\|q\right)}{\partial q_{l}}=\frac{1}{\sum_{j}q_{j}}\left(1-\frac{\bar{p}^{\alpha}_{l}\bar{q}^{\;-\alpha}_{l}}{\sum_{i} \bar{p}^{\alpha}_{i}\bar{q}^{1-\alpha}_{i}}\right)
	\label{eq.gradLMAGIas}
\end{equation}

\section{Arithmetic - Harmonic mean divergence}
This divergence is founded on the inequality: $M_{A}\geq M_{H}$.
\subsection{Unweighted version.}
This divergence is constructed in the sense of Csiszär on the standard convex function:
\begin{equation}
	f_{c}\left(x\right)=\frac{1}{2}\frac{\left(x-1\right)^{2}}{x+1}
\end{equation}
The resulting divergence is written:
\begin{equation}
	M_{AH}\left(p\|q\right)=\frac{1}{2}\sum_{i}\frac{\left(p_{i}-q_{i}\right)^{2}}{p_{i}+q_{i}}
\end{equation}
or, more clearly:
\begin{equation}
	M_{AH}\left(p\|q\right)=\left[\sum_{i}\frac{p_{i}+q_{i}}{2}-\sum_{i}\frac{2 p_{i}q_{i}}{p_{i}+q_{i}}\right]
\end{equation}
If we consider data fields with a sum explicitly equal to $1$, we have the simplified form:
\begin{equation}
	MS_{AH}\left(p\|q\right)=\left[1-\sum_{i}\frac{2 p_{i}q_{i}}{p_{i}+q_{i}}\right]
\end{equation}

\subsection{Weighted version.}
\begin{figure}[h!]
\centering
\includegraphics[width=0.7\linewidth]{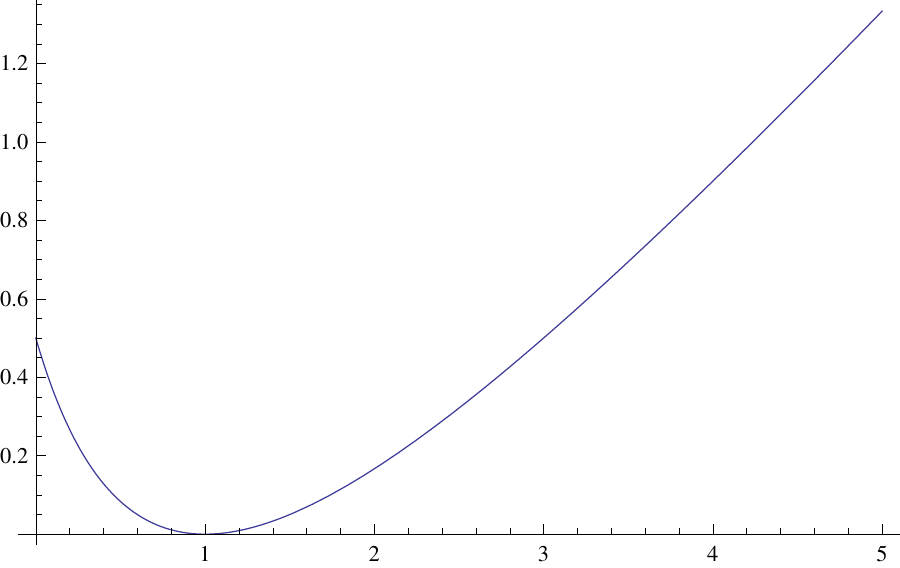}
\caption{Function $f_{c \alpha}\left(x\right)=\alpha x+\left(1-\alpha\right)-\frac{x}{\left(1-\alpha\right)x+\alpha}$, $\alpha=0.5$}
\label{fig:MAHPc}
\end{figure}

With $0\leq\alpha\leq1$, the standard convex function shown in the figure (\ref{fig:MAHPc}) for $\alpha=0.5$, is written as follows:
\begin{equation}
	f_{c \alpha}\left(x\right)=\alpha x+\left(1-\alpha\right)-\frac{x}{\left(1-\alpha\right)x+\alpha}
\end{equation}

It leads to the divergence:
\begin{equation}
	M_{AH\alpha}\left(p\|q\right)=\sum_{i}\alpha p_{i}+\left(1-\alpha\right)q_{i}-\sum_{i}\frac{p_{i}q_{i}}{\alpha q_{i}+\left(1-\alpha\right)p_{i}}
	\label{eq.MAHa}
\end{equation}
The gradient with respect to ``$q$" is expressed as:
\begin{equation}
	\frac{\partial 	M_{AH\alpha}}{\partial q_{j}}=\left(1-\alpha\right)\left\{1-\frac{p_{j}^{2}}{\left[\alpha q_{j}+\left(1-\alpha\right)p_{j}\right]^{2}}\right\}
	\label{eq.gradMAHa}
\end{equation}

\subsubsection{Scale invariant version with respect to ``$q$".}
For this divergence, the nominal invariance factor is not obtained explicitly, so as invariance factor we use the expression $K^{*}(p,q)=\frac{\sum_{j}p_{j}}{\sum_{j}q_{j}}$; we then obtain the corresponding simplified invariant divergence:
\begin{equation}
	M_{AH\alpha}I\left(p\|q\right)=\sum_{i}\left[\alpha \bar{p}_{i}+\left(1-\alpha\right)\bar{q}_{i}\right]-\sum_{i}\frac{\bar{p}_{i}\bar{q}_{i}}{\left(1-\alpha\right) \bar{p}_{i}+\alpha\bar{q}_{i}}
	\label{eq.MAHIa}
\end{equation}
Or in a simplified form:
\begin{equation}
	M_{AH\alpha}I\left(p\|q\right)=1-\sum_{i}\frac{\bar{p}_{i}\bar{q}_{i}}{\left(1-\alpha\right) \bar{p}_{i}+\alpha\bar{q}_{i}}
	\label{eq.MAHIas}
\end{equation}
The gradient with respect to ``$q$" is given by:
\begin{equation}
	\frac{\partial 	M_{AH\alpha}I\left(p\|q\right)}{\partial q_{l}}=\frac{\left(\alpha-1\right)}{\sum_{j}q_{j}}\left\{\frac{\bar{p}^{2}_{l}}{\left[\left(1-\alpha\right)\bar{p}_{l}+\alpha\bar{q}_{l}\right]^{2}}-\sum_{i}\frac{\bar{p}^{2}_{i}\bar{q}_{i}}{\left[\left(1-\alpha\right)\bar{p}_{i}+\alpha\bar{q}_{i}\right]^{2}}\right\}
	\label{eq.gradMAHIas}
\end{equation}

\subsection{Logarithmic form.}
It is written immediately from the weighted version (\ref{eq.MAHa}):
\begin{equation}
	LM_{AH\alpha}\left(p\|q\right)=\log\sum_{i}\alpha p_{i}+\left(1-\alpha\right)q_{i}-\log\sum_{i}\frac{p_{i}q_{i}}{\alpha q_{i}+\left(1-\alpha\right)p_{i}}
		\label{eq.LMAHas}
\end{equation}
The gradient with respect to ``$q$" is written as:
\begin{align}
	\frac{\partial 	LM_{AH\alpha}}{\partial q_{j}}=\left(1-\alpha\right)&\left\{\frac{1}{\sum_{i}\alpha p_{i}+\left(1-\alpha\right)q_{i}}\right. \nonumber \\ & \left.-\frac{1}{\sum_{i}\frac{p_{i}q_{i}}{\alpha q_{i}+\left(1-\alpha\right)p_{i}}}\frac{p_{j}^{2}}{\left[\alpha q_{j}+\left(1-\alpha\right)p_{j}\right]^{2}}\right\}
	\label{eq.gradLMAHas}
\end{align}

\subsubsection{Logarithmc form - Invariance with respect to ``$q$".}
The logarithmic invariant version is obtained by applying the ``Log" function separately to each term of (\ref{eq.MAHIas}) and we obtain:
\begin{equation}
	LM_{AH\alpha}I\left(p\|q\right)=-\log\sum_{i}\frac{\bar{p}_{i}\bar{q}_{i}}{\left(1-\alpha\right) \bar{p}_{i}+\alpha\bar{q}_{i}}
	\label{eq.LMAHIas}
\end{equation}
Whose gradient with respect to ``$q$" is given by:
\begin{equation}
\frac{\partial 	M_{AH\alpha}I\left(p\|q\right)}{\partial q_{l}}=\frac{\left(\alpha-1\right)}{\sum_{i}\frac{\bar{p}_{i}\bar{q}_{i}}{\left(1-\alpha\right) \bar{p}_{i}+\alpha\bar{q}_{i}}\sum_{j}q_{j}}\left[\frac{\bar{p}^{2}_{l}}{\left[\left(1-\alpha\right) \bar{p}_{l}+\alpha\bar{q}_{l}\right]^{2}}-\sum_{i}\frac{\bar{p}^{2}_{i}\bar{q}_{i}}{\left[\left(1-\alpha\right) \bar{p}_{i}+\alpha\bar{q}_{i}\right]^{2}}\right]
\label{eq.gradLMAHIas}	
\end{equation}
\setcounter{table}{0}  \setcounter{equation}{0}  \setcounter{figure}{0} \setcounter{chapter}{8} \setcounter{section}{0} 
\chapter{chapter 8 -\\Divergences between mixed fields.}  \label{chptr::chapitre8}
In this chapter, the main idea is to analyze the difference between 2 data fields ``$p$" and ``$q$" by expressing the difference in a broad sense between ``$p$" or ``$q$" and the weighted sum of ``$p$" and ``$q$"; $0\leq\alpha\leq1$ will be the corresponding weighting factor. Thus, whatever the basic divergence, one can go from a zero divergence to a divergence between the elementary fields by varying the parameter ``$\alpha$".\\
The gap between ``$p$" or ``$q$" and their weighted sum is expressed by a Kullback-Leibler divergence \cite{kullback1951} or, to generalize, by a Havrda-Charvat divergence \cite{havrda1967}, but nothing prohibits the use of other basic divergences.\\
Moreover, to broaden the picture, we will envisage taking into account in the expression of the divergences, the order of the arguments, which is understandable inasmuch as the divergences used are not generally symmetrical.\\
This argumentation describes a particular constructive process for these divergences, however, in each case examined, we will always come back to Csiszär's constructive process and the associated convex function.\\
Of course, in the literature, these studies have been proposed mainly to quantify the gap between two probability densities; here we try to leave this restrictive context and extend the applications to more general data fields.

\section{ F divergence (Jensen-Shannon relative divergence).}
\subsection{Direct form.}
\begin{figure}[h!]
\centering
\includegraphics[width=0.7\linewidth]{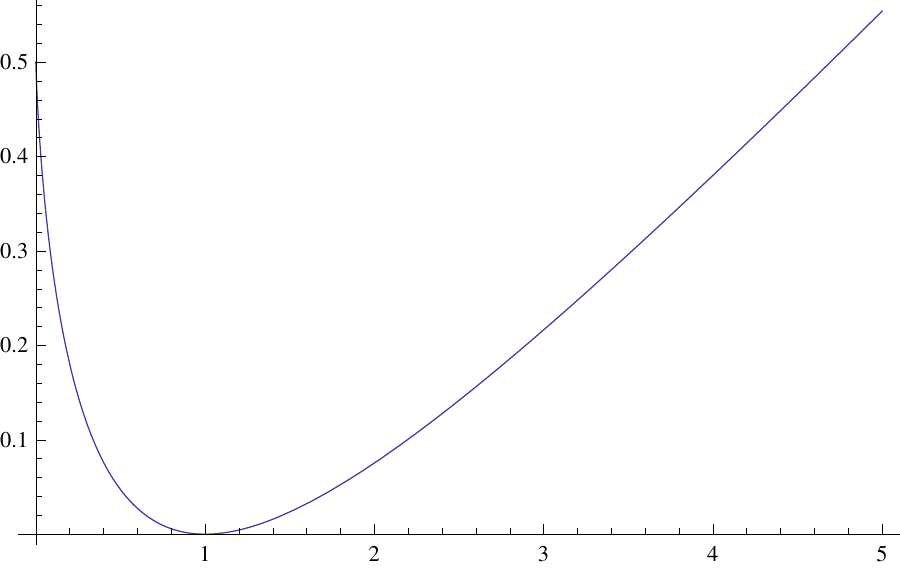}
\caption{Function $f_{c}\left(x\right)=x\log\frac{x}{\alpha x+\left(1-\alpha\right)}+\left(1-\alpha\right)\left(1-x\right)$, $\alpha=0.5$}
\label{fig:DFc}
\end{figure}

This divergence is reported by Taneja \cite{taneja2001}.\\
The constructive process can be described by:
\begin{equation}
	F\left(p\|q\right)=KL\left(p\|\alpha p+\left(1-\alpha\right)q\right)
\end{equation}
If we use as expression of K.L.'s divergence, the general form usable for any (non-normalized) variables, we obtain:
\begin{equation}
	F\left(p\|q\right)=\sum_{i}\left[p_{i}\log\frac{p_{i}}{\alpha p_{i}+\left(1-\alpha\right)q_{i}}+\left(1-\alpha\right)\left(q_{i}-p_{i}\right)\right]
\label{eq:F}	
\end{equation}
With $0\leq\alpha\leq1$, this divergence is constructed in the sense of Csiszär on the standard convex function shown in the figure (\ref{fig:DFc}) for $\alpha=0.5$, it is written as:
\begin{equation}
	f_{c}\left(x\right)=x\log\frac{x}{\alpha x+\left(1-\alpha\right)}+\left(1-\alpha\right)\left(1-x\right)
\end{equation}

\begin{figure}[h!]
\centering
\includegraphics[width=0.7\linewidth]{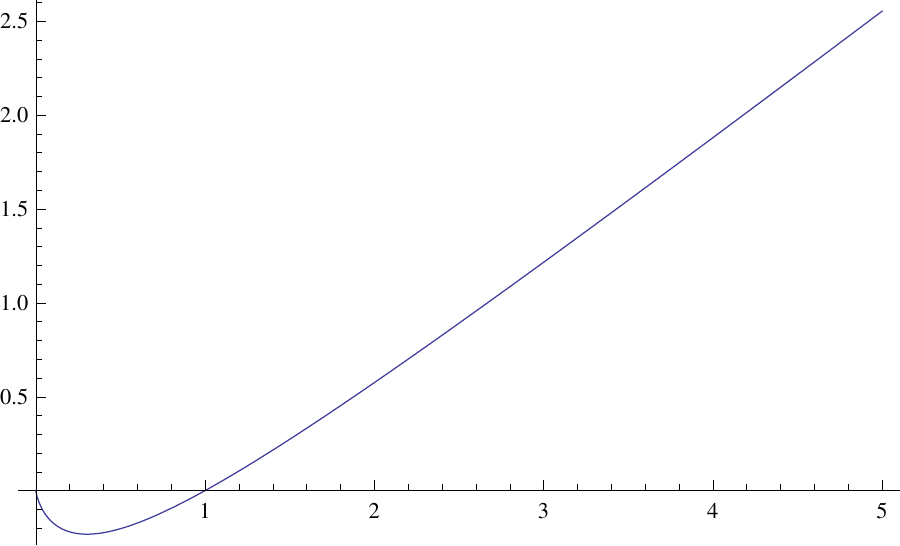}
\caption{Function $f\left(x\right)=x\log\frac{x}{\alpha x+\left(1-\alpha\right)}$, $\alpha=0.5$}
\label{fig:DF}
\end{figure}

The gradient with respect to ``$q$" is expressed as:
\begin{equation}
	\frac{\partial F\left(p\|q\right)}{\partial q_{j}}=\left(1-\alpha\right)\left(1-\frac{p_{j}}{\alpha p_{j}+\left(1-\alpha\right)q_{j}}\right)
\end{equation}
In the case where we can consider that $\sum_{i}p_{i}=\sum_{i}q_{i}$, not necessarily equal to 1, we have a simplified form, which is written as follows:
\begin{equation}
		FS\left(p\|q\right)=\sum_{i}p_{i}\log\frac{p_{i}}{\alpha p_{i}+\left(1-\alpha\right)q_{i}}
\end{equation}
With $\alpha=1/2$, this is the K divergence of Lin \cite{lin1991}.\\
With $0\leq\alpha\leq1$, it is a Csiszär divergence built on the simple convex function, shown in the figure (\ref{fig:DF}) for $\alpha=0.5$, which is expressed as follows:
\begin{equation}
	f\left(x\right)=x\log\frac{x}{\alpha x+\left(1-\alpha\right)}
\end{equation}

\subsection{Dual form.}
\begin{figure}[h!]
\centering
\includegraphics[width=0.7\linewidth]{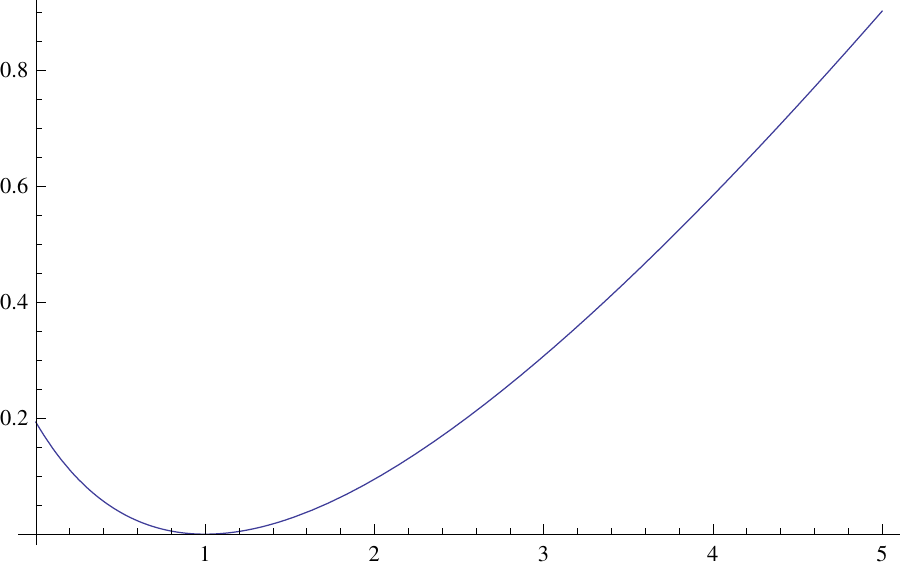}
\caption{Function $\breve{f}_{c}\left(x\right)=\log\frac{1}{\alpha+\left(1-\alpha\right)x}+\left(1-\alpha\right)\left(x-1\right)$, $\alpha=0.5$}
\label{fig:DFtc}
\end{figure}

The constructive process can be described by:
\begin{equation}
	F\left(q\|p\right)=KL\left(q\|\alpha q+\left(1-\alpha\right)p\right)
\end{equation}
If we use as expression of the K.L. divergence, the general form usable for non-normalized variables, we obtain:
\begin{equation}
	F\left(q\|p\right)=\sum_{i}\left[q_{i}\log\frac{q_{i}}{\alpha q_{i}+\left(1-\alpha\right)p_{i}}+\left(1-\alpha\right)\left(p_{i}-q_{i}\right)\right]
\end{equation}
With $0\leq\alpha\leq1$, it is a Csiszär divergence built on the standard convex function, represented in the figure (\ref{fig:DFtc}) for $\alpha=0.5$, which is expressed as follows:
\begin{equation}
	\breve{f}_{c}\left(x\right)=\log\frac{1}{\alpha+\left(1-\alpha\right)x}+\left(1-\alpha\right)\left(x-1\right)
\end{equation}

\begin{figure}[h!]
\centering
\includegraphics[width=0.7\linewidth]{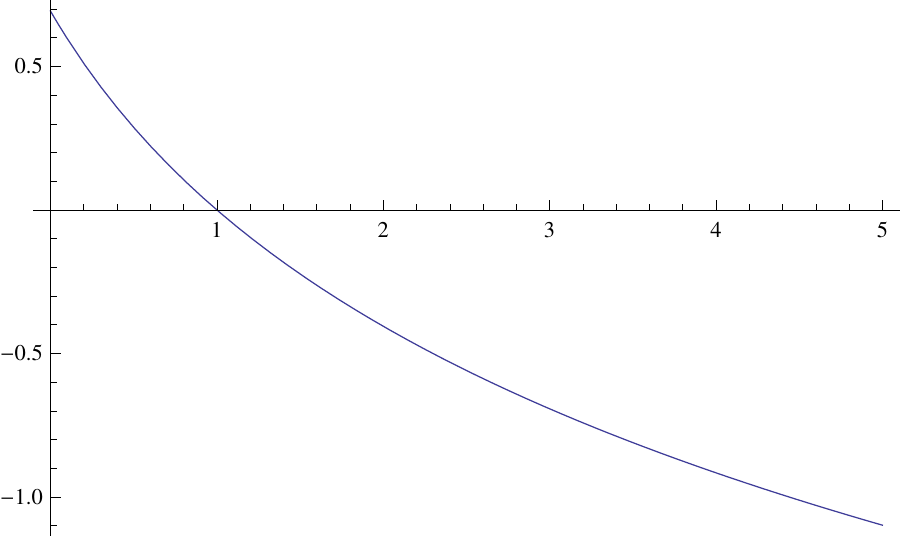}
\caption{Function $\breve{f}\left(x\right)=\log\frac{1}{\alpha+\left(1-\alpha\right)x}$, $\alpha=0.5$}
\label{fig:DFt}
\end{figure}

Its gradient with respect to ``$q$" is written as follows:
\begin{equation}
	\frac{\partial F\left(q\|p\right)}{\partial q_{j}}=\log\frac{q_{j}}{\alpha q_{j}+\left(1-\alpha\right) p_{j}}-\frac{\alpha q_{j}}{\alpha q_{j}+\left(1-\alpha\right) p_{j}}+\alpha
\end{equation}
In the case where we can consider that $\sum_{i}p_{i}=\sum_{i}q_{i}$, not necessarily equal to 1, we have a simplified form which is written as:
\begin{equation}
		FS\left(q\|p\right)=\sum_{i}q_{i}\log\frac{q_{i}}{\alpha q_{i}+\left(1-\alpha\right)p_{i}}
\end{equation}
It is a Csiszär divergence built on the simple convex function represented on the figure (\ref{fig:DFt}) for $\alpha=0.5$, written as follows:
\begin{equation}
	\breve{f}\left(x\right)=\log\frac{1}{\alpha+\left(1-\alpha\right)x}
\end{equation}

\subsection{Symmetrical form.}
It is written in the sense of Jeffreys \cite{jeffreys1946} as the (half) sum of the two previous ones:
\begin{equation}
	F\left(p,q\right)=F\left(p\|q\right)+F\left(q\|p\right)
\end{equation}
That is:
\begin{equation}
	F\left(p,q\right)=\sum_{i}p_{i}\log\frac{p_{i}}{\alpha p_{i}+\left(1-\alpha\right)q_{i}}+\sum_{i}q_{i}\log\frac{q_{i}}{\alpha q_{i}+\left(1-\alpha\right)p_{i}}	
\end{equation}
Whose gradient with respect to ``$q$" is expressed as:
\begin{align}
	\frac{\partial F\left(q,p\right)}{\partial q_{j}}=&1-\left(1-\alpha\right)\frac{p_{j}}{\alpha p_{j}+\left(1-\alpha\right)q_{j}}+\log\frac{q_{j}}{\alpha q_{j}+\left(1-\alpha\right) p_{j}}\nonumber \\ &-\frac{\alpha q_{j}}{\alpha q_{j}+\left(1-\alpha\right) p_{j}}
\end{align}

\subsection{Invariance by change of scale on ``$q$''.}
The classical method of computing the $K\left(p,q\right)$ invariance factor corresponding to the divergence $F\left(p\|q\right)$ (equation (\ref{eq:F})) involves solving in $K$ the equation:
\begin{equation}
	\frac{\partial F\left(p\|Kq\right)}{\partial K}=0
\end{equation}
That is, to solve with respect to $K$, the equation:
\begin{equation}
	\sum_{i}q_{i}\left[1-\frac{p_{i}}{\alpha p_{i}+\left(1-\alpha\right)Kq_{i}}\right]=0
\end{equation}
There's no explicit solution to this problem.\\
So we use as an invariance factor the expression $K^{*}(p,q)=\frac{\sum_{j}p_{j}}{\sum_{j}q_{j}}$.\\
In these conditions, the invariant divergence is written after simplifications with the normalized variables $\bar{p}$ and $\bar{q}$:
\begin{equation}
	FI\left(p\|q\right)=\sum_{i}\bar{p}_{i}\log\frac{\bar{p}_{i}}{\alpha \bar{p}_{i}+\left(1-\alpha\right)\bar{q}_{i}}
\label{eq:FIs}	
\end{equation}
Hence the expression of the gradient:
\begin{equation}
	\frac{\partial FI\left(p\|q\right)}{\partial q_{l}}=\frac{1-\alpha}{\sum_{j}q_{j}}\left[\sum_{i}\frac{\bar{p}_{i}\bar{q}_{i}}{\alpha \bar{p}_{i}+\left(1-\alpha\right)\bar{q}_{i}}-\frac{\bar{p}_{l}}{\alpha \bar{p}_{l}+\left(1-\alpha\right)\bar{q}_{l}}\right]
	\label{eq:gradFIs}
\end{equation}

\section{G divergence (Arithmetic-Geometric relative divergence).}
These divergences are discussed by Taneja \cite{taneja2001}, they are deduced from the previous ones by inverting the order of the compared fields; only the weighted forms are presented here.
\subsection{Direct form.}
\begin{figure}[h!]
\centering
\includegraphics[width=0.7\linewidth]{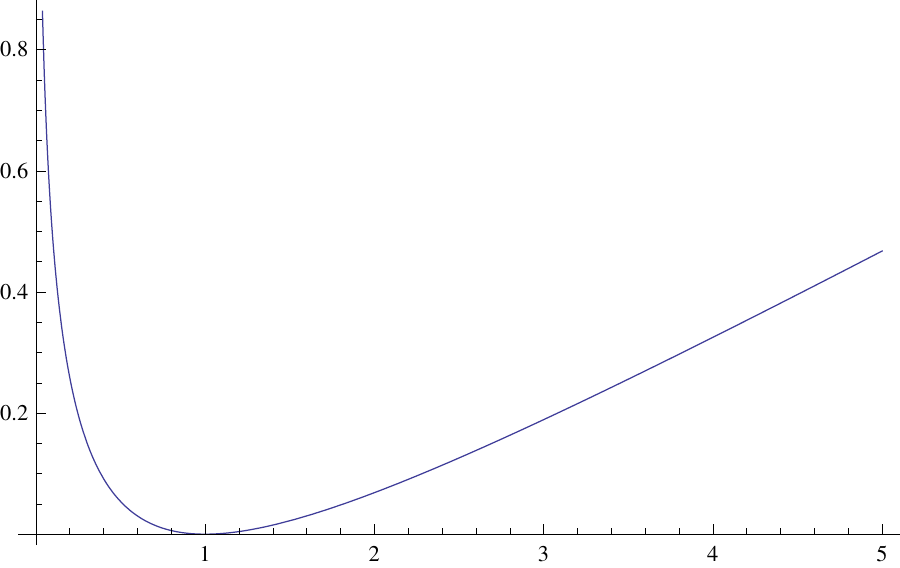}
\caption{Function $f_{c}\left(x\right)=\left[\alpha x+\left(1-\alpha\right)\right]\log\frac{\alpha x+\left(1-\alpha\right)}{x}+\left(1-\alpha\right)\left(x-1\right)$, $\alpha=0.5$}
\label{fig:DGc}
\end{figure}

The constructive process is of the following form:
\begin{equation}
	G\left(p\|q\right)=KL\left(\alpha p+\left(1-\alpha\right)q\|p\right)
\end{equation}
If we use the general form of the Kullback-Leibler divergence that can be used for non-normalized fields, we obtain:
\begin{align}
	G\left(p\|q\right)=&\sum_{i}\left[\alpha p_{i}+\left(1-\alpha\right)q_{i}\right]\log\frac{\alpha p_{i}+\left(1-\alpha\right)q_{i}}{p_{i}}\nonumber \\ &+\left(1-\alpha\right)\left(p_{i}-q_{i}\right)
	\label{eq:G}
\end{align}
It is a Csiszär's divergence built on the standard convex function, represented in the figure (\ref{fig:DGc}) for $\alpha=0.5$, which is written as follows:
\begin{equation}
	f_{c}\left(x\right)=\left[\alpha x+\left(1-\alpha\right)\right]\log\frac{\alpha x+\left(1-\alpha\right)}{x}+\left(1-\alpha\right)\left(x-1\right)
\end{equation}

\begin{figure}[h!]
\centering
\includegraphics[width=0.7\linewidth]{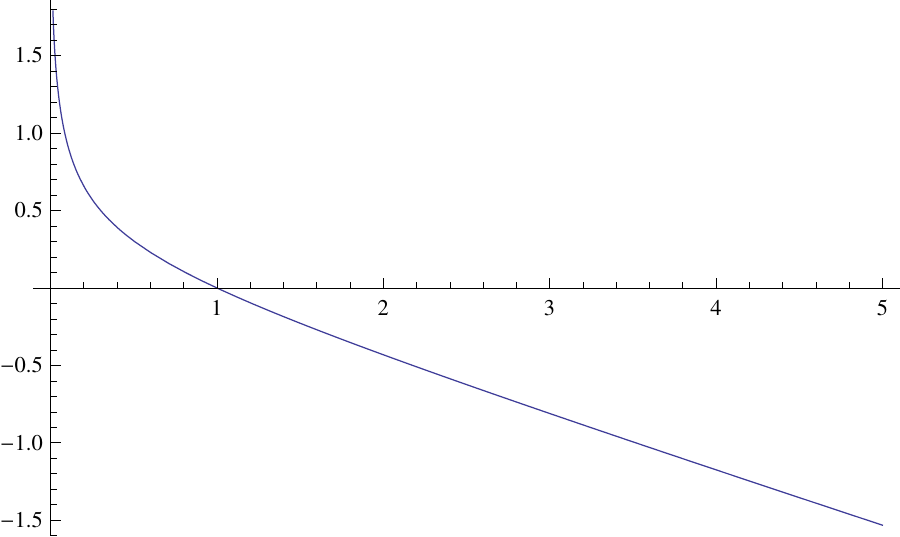}
\caption{Function $f\left(x\right)=\left[\alpha x+\left(1-\alpha\right)\right]\log\frac{\alpha x+\left(1-\alpha\right)}{x}$, $\alpha=0.5$}
\label{fig:DG}
\end{figure}

Its gradient with respect to ``$q$" is expressed as:
\begin{equation}
		\frac{\partial G\left(p\|q\right)}{\partial q_{j}}=\left(1-\alpha\right)\log\frac{\alpha p_{j}+\left(1-\alpha\right)q_{j}}{p_{j}}
\end{equation}
When we can consider that $\sum_{i}p_{i}=\sum_{i}q_{i}$, not necessarily equal to 1, a simplified form exists:
\begin{equation}
	GS\left(p\|q\right)=\sum_{i}\left[\alpha p_{i}+\left(1-\alpha\right)q_{i}\right]\log\frac{\alpha p_{i}+\left(1-\alpha\right)q_{i}}{p_{i}}
\end{equation}
It is a Csiszär divergence built on the simple convex function shown in the figure (\ref{fig:DG}) for $\alpha=0.5$, which is written:
\begin{equation}
	f\left(x\right)=\left[\alpha x+\left(1-\alpha\right)\right]\log\frac{\alpha x+\left(1-\alpha\right)}{x}
\end{equation}

\subsection{Dual form.}
\begin{figure}[h!]
\centering
\includegraphics[width=0.7\linewidth]{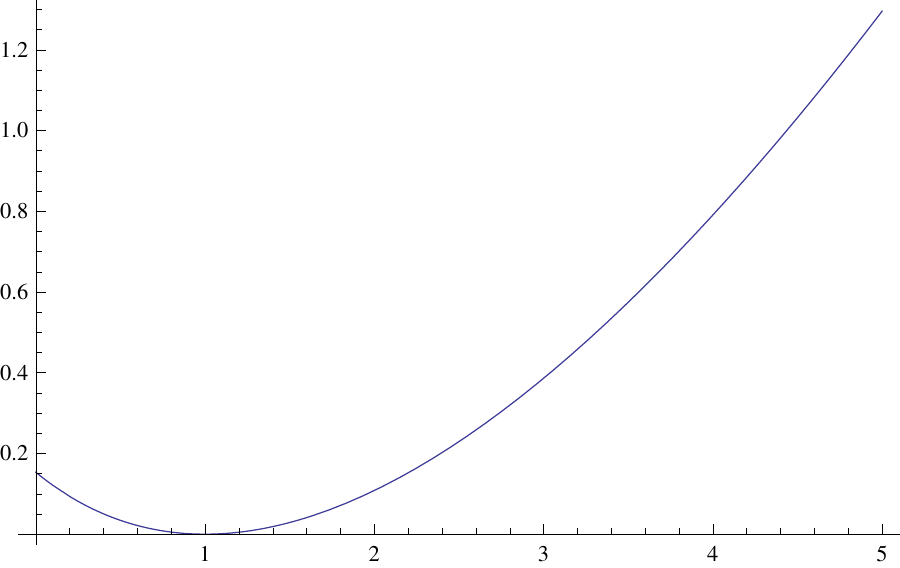}
\caption{Function $\breve{f}_{c}\left(x\right)=\left[\alpha+\left(1-\alpha\right)x\right]\log\left[\alpha+\left(1-\alpha\right)x\right]+\left(1-\alpha\right)\left(1-x\right)$, $\alpha=0.5$}
\label{fig:DGtc}
\end{figure}

The constructive process can be described by the relation:
\begin{equation}
	G\left(q\|p\right)=KL\left(\alpha q+\left(1-\alpha\right)p\|q\right)
\end{equation}
If we take the general form of the Kullback-Leibler divergence that can be used for non-normalized fields, we obtain:
\begin{align}
	G\left(q\|p\right)=&\sum_{i}\left[\alpha q_{i}+\left(1-\alpha\right)p_{i}\right]\log\frac{\alpha q_{i}+\left(1-\alpha\right)p_{i}}{q_{i}}\nonumber \\ &+\left(1-\alpha\right)\left(q_{i}-p_{i}\right)
\end{align}
It is a Csiszär divergence built on the standard convex function shown in the figure (\ref{fig:DGtc}) for $\alpha=0.5$, which is written:
\begin{equation}	\breve{f}_{c}\left(x\right)=\left[\alpha+\left(1-\alpha\right)x\right]\log\left[\alpha+\left(1-\alpha\right)x\right]+\left(1-\alpha\right)\left(1-x\right)
\end{equation}

\begin{figure}[h!]
\centering
\includegraphics[width=0.7\linewidth]{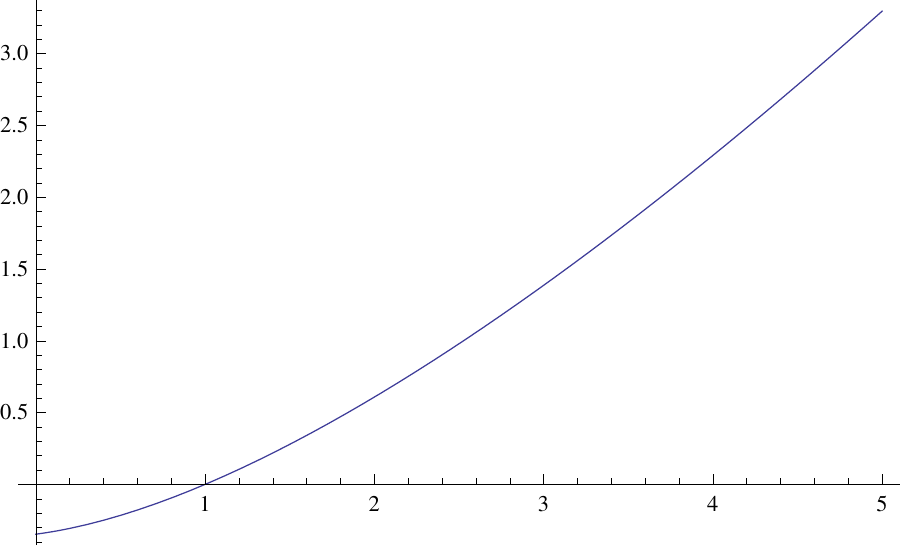}
\caption{Function $\breve{f}\left(x\right)=\left[\alpha+\left(1-\alpha\right)x\right]\log\left[\alpha+\left(1-\alpha\right)x\right]$, $\alpha=0.5$}
\label{fig:DGt}
\end{figure}

The gradient with respect to ``$q$" is given by:
\begin{equation}
	\frac{\partial G\left(q\|p\right)}{\partial q_{j}}=\alpha \log\frac{\alpha q_{j}+\left(1-\alpha\right)p_{j}}{q_{j}}+\left(1-\alpha\right)\left(1-\frac{p_{j}}{q_{j}}\right)
\end{equation}
When $\sum_{i}p_{i}=\sum_{i}q_{i}$, not necessarily equal to 1, a simplified form exists, it is written as follows:
\begin{equation}
		GS\left(q\|p\right)=\sum_{i}\left[\alpha q_{i}+\left(1-\alpha\right)p_{i}\right]\log\frac{\alpha q_{i}+\left(1-\alpha\right)p_{i}}{q_{i}}
\end{equation}
It is a Csiszär divergence built on the simple convex function represented on the figure (\ref{fig:DGt}) for $\alpha=0.5$, which is written:
\begin{equation}
	\breve{f}\left(x\right)=\left[\alpha+\left(1-\alpha\right)x\right]\log\left[\alpha+\left(1-\alpha\right)x\right]
\end{equation}

\subsection{T symmetrical divergence.}
It is classically constructed in the form of:
\begin{equation}
	T\left(p,q\right)=G\left(p\|q\right)+G\left(q\|p\right)
\end{equation}
That is:
\begin{eqnarray}
		T\left(p,q\right) &=& \sum_{i}\left[\alpha p_{i}+\left(1-\alpha\right)q_{i}\right]\log\frac{\alpha p_{i}+\left(1-\alpha\right)q_{i}}{p_{i}} \nonumber
		\\ &+& 		\sum_{i}\left[\alpha q_{i}+\left(1-\alpha\right)p_{i}\right]\log\frac{\alpha q_{i}+\left(1-\alpha\right)p_{i}}{q_{i}}
\end{eqnarray}
The corresponding gradient with respect to ``$q$" is immediately deduced from the gradients of the functions $G\left(p\|q\right)$ and $G\left(q\|p\right)$.

\subsection{Invariance by change of scale on ``$q$".}
Regarding the $G\left(p\|q\right)$ divergence given by the equation (\ref{eq:G}), the classical method does not explicitly allow to calculate a value for the invariance factor $K$.\\
By introducing the non-nominal invariance factor $K^{*}$, in the expression (\ref{eq:G}), we obtain the corresponding invariant divergence which is written after simplifications:
\begin{equation}
	GI\left(p\|q\right)=\sum_{i}\left[\alpha \bar{p}_{i}+\left(1-\alpha\right)\bar{q}_{i}\right]\log\frac{\alpha \bar{p}_{i}+\left(1-\alpha\right)\bar{q}_{i}}{\bar{p}_{i}}
	\label{eq:GI}
\end{equation}
Its gradient with respect to ``$q$'' is given by:
\begin{equation}
	\frac{\partial GI\left(p\|q\right)}{\partial q_{l}}=\frac{1-\alpha}{\sum_{j}q_{j}}\left[\log\frac{\alpha \bar{p}_{l}+\left(1-\alpha\right)\bar{q}_{l}}{\bar{p}_{l}}-\sum_{i}\bar{q}_{i}\log\frac{\alpha \bar{p}_{i}+\left(1-\alpha\right)\bar{q}_{i}}{\bar{p}_{i}}\right]
\label{eq:gradGI}
\end{equation}

\section {Generalization of the F and G divergences.}
To proceed with this generalization, the same work is done as in the previous sections, but replacing the Kullback-Leibler divergence used systematically in the constructive process, by the Havrda-Charvat divergence, i.e. something that is an "alpha (s) divergence", and relying on the fact that the K .L. divergence can be seen as a particular case of the H.C. divergence by variation of the parameter ``s" ($0\leq s\leq 1$) (we no longer use ``$\alpha$" which is already used to weight the two fields ``$p$'' and ``$q$'').\\
Indeed, the H.C. divergence is written as follows:

\begin{equation}
	HC_{s}\left(a\|b\right)=\frac{1}{s\left(s-1\right)}\sum_{i}a_{i}^{s}b_{i}^{1-s}-s a_{i}-\left(1-s\right) b_{i}
\end{equation}
This expression tends towards $KL\left(a\|b\right)$ when $s\rightarrow 1$ and towards $KL\left(b\|a\right)$ when $s\rightarrow 0$.\\
As in the previous sections, only the weighted forms will be developed; in this context the data fields that will be used are ``$p$", ``$q$", and the weighted combination ``$\alpha p+\left(1-\alpha\right)q$".

\subsection{Direct form.}
We rely on the following constructive approach:
\begin{equation}
	FG_{s,\alpha}\left(p\|q\right)=HC\left(\alpha p+\left(1-\alpha\right) q\|p\right)
\end{equation}
This results in:
\begin{align}
	FG_{s,\alpha}\left(p\|q\right)=&\frac{1}{s\left(s-1\right)}\sum_{i}\left[\alpha p_{i}+\left(1-\alpha\right) q_{i}\right]^{s}p_{i}^{1-s}-\nonumber \\ &s\left[\alpha p_{i}+\left(1-\alpha\right) q_{i}\right]-\left(1-s\right)p_{i}
	\label{eq.FGSA}	
\end{align}
It can also be said that this divergence is built, in the sense of Csiszär, on the standard convex function:
\begin{equation}
	f_{c}\left(x\right)=\frac{1}{s\left(s-1\right)}\left\{\left[\alpha x+\left(1-\alpha\right)\right]^{s}x^{1-s}-s\left[\alpha x+\left(1-\alpha\right)\right]-\left(1-s\right)x\right\}
\end{equation}
Its gradient with respect to ``$q$" is given by:
\begin{equation}
	\frac{\partial FG_{s,\alpha}\left(p\|q\right)}{\partial q_{j}}=\frac{1-\alpha}{s-1}\left\{\left[\frac{\alpha p_{j}+\left(1-\alpha\right)q_{j}}{p_{j}}\right]^{s-1}-1\right\}
\end{equation}
If we are in the case where $\sum_{i}p_{i}=\sum_{i}q_{i}$, we can use a simplified form of this divergence which is written as follows:
\begin{equation}
	FGS_{s,\alpha}\left(p\|q\right)=\frac{1}{s\left(s-1\right)}\sum_{i}\left[\alpha p_{i}+\left(1-\alpha\right) q_{i}\right]^{s}p_{i}^{1-s}-p_{i}	
\end{equation}
This is what would have been obtained if, in the mode of construction, a simplified H.C. divergence had been used, but it can also be said that this divergence is constructed in the sense of Csiszär on the simple convex function:
\begin{equation}
	f\left(x\right)=\frac{1}{s\left(s-1\right)}\left\{\left[\alpha x+\left(1-\alpha\right)\right]^{s}x^{1-s}-x\right\}
\end{equation}

\subsection{Invariant form by change of scale on ``$q$".}
Following the usual procedure, we have to solve in $K$, the equation:
\begin{equation}
	\frac{1-\alpha}{s-1}\sum_{i}q_{i}\left\{\left[\frac{\alpha p_{i}+\left(1-\alpha\right)K q_{i}}{p_{i}}\right]^{s-1}-1\right\}=0
\end{equation}
There is no explicit solution; however, it is possible to use:
\begin{equation}
	K^{*}=\frac{\sum_{i}p_{i}}{\sum_{i}q_{i}}
\end{equation}
This expression of $K^{*}$ inserted in the initial divergence (\ref{eq.FGSA}) makes it invariant by scale change and we obtain the simplified invariant divergence:
\begin{equation}
	FGI\left(p\|q\right)=\frac{1}{s\left(s-1\right)}\sum_{i}\bar{p}^{1-s}_{i}\left[\alpha \bar{p}_{i}+\left(1-\alpha\right)\bar{q_{i}}\right]^{s}
	\label{eq.FGIs}
\end{equation}
In this expression, the normalized variables $\bar{p}$ and $\bar{q}$ are used.\\
The gradient with respect to ``$q$" is written as :
\begin{equation}
	\frac{\partial FGI\left(p\|q\right)}{\partial q_{l}}=\frac{1-\alpha}{\left(s-1\right)\sum_{j}q_{j}}\left\{\left[\frac{\bar{p}_{l}}{\alpha \bar{p_{l}}+\left(1-\alpha\right)\bar{q_{l}}}\right]^{s-1}-\sum_{i}\bar{q}_{i}\left[\frac{\bar{p}_{i}}{\alpha \bar{p_{i}}+\left(1-\alpha\right)\bar{q_{i}}}\right]^{s-1}\right\}
	\label{eq.gradFGIs}	
\end{equation}

\subsection{Dual form.}
It is founded on the following constructive mode:
\begin{equation}
	FG_{s,\alpha}=HC\left(\alpha q+\left(1-\alpha\right) p\|q\right)
\end{equation}
This leads to:
\begin{align}
	FG_{s,\alpha}\left(q\|p\right)=&\frac{1}{s\left(s-1\right)}\sum_{i}\left[\alpha q_{i}+\left(1-\alpha\right) p_{i}\right]^{s}q_{i}^{1-s}\nonumber \\ &-s\left[\alpha q_{i}+\left(1-\alpha\right) p_{i}\right]-\left(1-s\right)q_{i}	
\end{align}
It can also be said that it is a Csiszär divergence built on the standard convex function:
\begin{align}	\breve{f}_{cs}\left(x\right)=\frac{1}{s\left(s-1\right)}&\left\{\left[\alpha+\left(1-\alpha\right)x\right]^{s}-\left(1-s\right)\right. \nonumber \\ & \left.-s\left[\alpha+\left(1-\alpha\right)x\right]\right\}
\end{align}
Its gradient with respect to ``$q$'' can be written as follows:
\begin{align}
	\frac{\partial FG_{s,\alpha}\left(q\|p\right)}{\partial q_{j}}=&\frac{\alpha}{s-1}\left\{\left[\frac{\alpha q_{j}+\left(1-\alpha\right) p_{j}}{q_{j}}\right]^{s-1}-1\right\}\nonumber \\ &-\frac{1}{s}\left\{\left[\frac{\alpha q_{j}+\left(1-\alpha\right) p_{j}}{q_{j}}\right]^{s}-1\right\}
\end{align}
In the case where we have $\sum_{i}p_{i}=\sum_{i}q_{i}$, we obtain the following simplified form:
\begin{equation}
	FGS_{s,\alpha}\left(q\|p\right)=\frac{1}{s\left(s-1\right)}\sum_{i}\left[\alpha q_{i}+\left(1-\alpha\right) p_{i}\right]^{s}q_{i}^{1-s}- q_{i}	
\end{equation}
It is constructed in the sense of Csiszär, on the simple convex function:
\begin{equation}
\breve{f}_{s}\left(x\right)=\frac{1}{s\left(s-1\right)}\left[\left\{\alpha+\left(1-\alpha\right)x\right\}^{s}-1\right]
\end{equation}
A further simplification can be introduced if one is dealing with probability densities, then, $\sum_{i}p_{i}=\sum_{i}q_{i}=1$.

\subsection{Symmetrical form.}
The standard convex function for constructing the symmetrical divergence in the Csiszär sense is written (with a factor of 1/2):
\begin{align}
	\widehat{f}_{cs}\left(x\right)=f_{cs}\left(x\right)+\breve{f}_{cs}\left(x\right)=\frac{1}{s\left(s-1\right)}&\left\{x\left[\frac{\alpha x+\left(1-\alpha\right)}{x}\right]^{s}\right. \nonumber \\ & \left.+\left[\alpha+\left(1-\alpha\right)x\right]^{s}-\left(1+x\right)\right\}
\end{align}
The corresponding divergence is written:
\begin{align}
	FG_{s,\alpha}\left(p,q\right)=&\frac{1}{s\left(s-1\right)}\sum_{i}\left[\alpha p_{i}+\left(1-\alpha\right)q_{i}\right]^{s}p_{i}^{1-s}\nonumber \\ &+\left[\alpha q_{i}+\left(1-\alpha\right)p_{i}\right]^{s}q_{i}^{1-s}-\left(p_{i}+q_{i}\right)
	\label{eq.FGsa}
\end{align}
The gradient with respect to ``$q$" can be written as follows:
\begin{eqnarray}
	\frac{\partial FG_{s,\alpha}\left(p,q\right)}{\partial q_{j}}&=&\frac{1-\alpha}{s-1}\left[\frac{\alpha p_{j}+\left(1-\alpha\right) q_{j}}{p_{j}}\right]^{s-1}+\frac{\alpha}{s-1}\left[\frac{\alpha q_{j}+\left(1-\alpha\right) p_{j}}{q_{j}}\right]^{s-1}
\nonumber	
	\\ &-&   \frac{1}{s}\left[\frac{\alpha q_{j}+\left(1-\alpha\right) p_{j}}{q_{j}}\right]^{s}-\frac{1}{s\left(s-1\right)}
\end{eqnarray}

\section{J  relative divergence of Dragomir.}
To broaden what has been proposed by Dragomir et al. \cite{dragomir2001}, we will develop the general weighted case.
\subsection{Direct form.}
\begin{figure}[h!]
\centering
\includegraphics[width=0.7\linewidth]{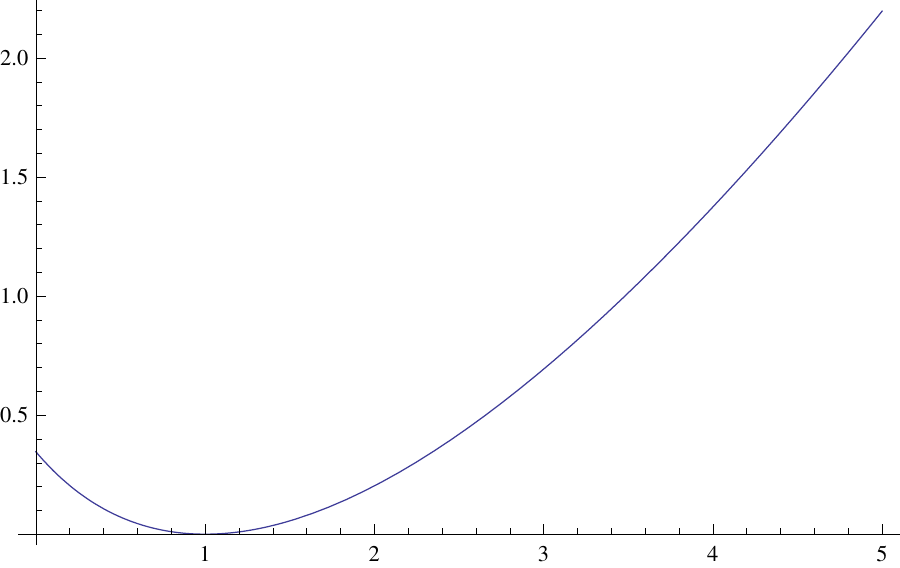}
\caption{Function $f_{c}\left(x\right)=\left(1-\alpha\right)\left(x-1\right)\log\left[\alpha+\left(1-\alpha\right)x\right]$, $\alpha=0.5$}
\label{fig:JDRc}
\end{figure}

Founding on the use of K.L.'s divergences, one can write:
\begin{equation}
	JD\left(p\|q\right)=KL\left[q\|\alpha p+\left(1-\alpha\right)q\right]+KL\left[\alpha p+\left(1-\alpha\right)q\right\|q]
\end{equation}
Using the notations in the previous sections, one can also write:
\begin{equation}
		JD\left(p\|q\right)=F\left(q\|p\right)+G\left(q\|p\right)
\end{equation}
But it can also be said that it is a Csiszär divergence built on the standard convex function represented in the figure (\ref{fig:JDRc}) for $\alpha=0.5$, written as:
\begin{equation}
	f_{c}\left(x\right)=\left(1-\alpha\right)\left(x-1\right)\log\left[\alpha+\left(1-\alpha\right)x\right]
	\label{eq.fcD}
\end{equation}

Then, we obtain the divergence:
\begin{equation}
		JD\left(p\|q\right)=\left(1-\alpha\right)\sum_{i}\left(p_{i}-q_{i}\right)\log\frac{\alpha q_{i}+\left(1-\alpha\right)p_{i}}{q_{i}}
		\label{eq.JD}
\end{equation}
The gradient with respect to ``$q$" is expressed as:
\begin{equation}
\frac{\partial JD\left(p,q\right)}{\partial q_{j}}\propto\frac{\alpha\left(p_{j}-q_{j}\right)}{\alpha q_{j}+\left(1-\alpha\right)p_{j}}-\frac{p_{j}}{q_{j}}
+\log\frac{q_{j}}{\alpha q_{j}+\left(1-\alpha\right) p_{j}}+1
\end{equation}

\subsection{Dual form.}
\begin{figure}[h!]
\centering
\includegraphics[width=0.7\linewidth]{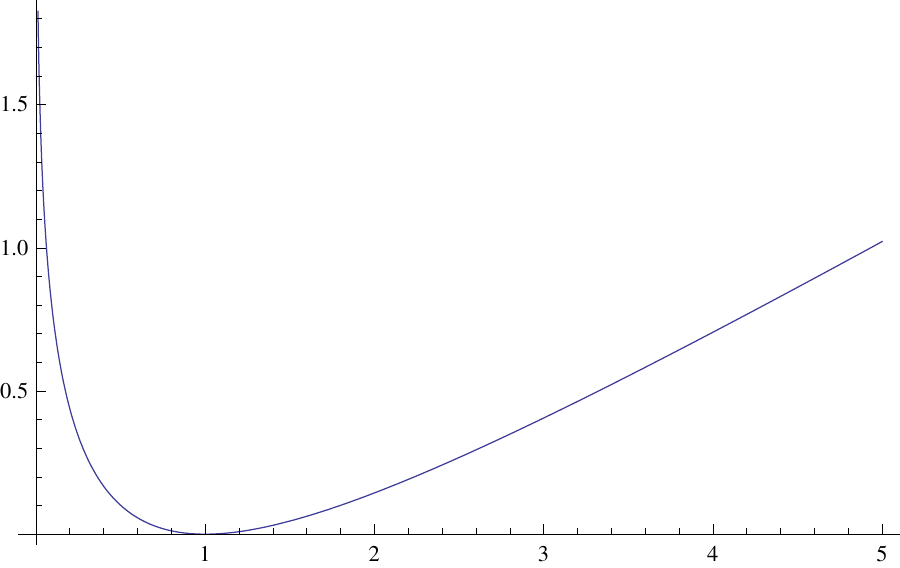}
\caption{Function $	\breve{f}_{c}\left(x\right)=\left(1-\alpha\right)\left(1-x\right)\log\left[\frac{\alpha x+\left(1-\alpha \right)}{x}\right]$, $\alpha=0.5$}
\label{fig:JDRtc}
\end{figure}

It is constructed in the sense of Csiszär on the standard convex function, mirror of the function (\ref{eq.fcD}), represented on the figure (\ref{fig:JDRtc}) for $\alpha=0.5$, that is written:
\begin{equation}
	\breve{f}_{c}\left(x\right)=\left(1-\alpha\right)\left(1-x\right)\log\left[\frac{\alpha x+\left(1-\alpha \right)}{x}\right]
\end{equation}

This results in:
\begin{equation}
	JD\left(q\|p\right)=\left(1-\alpha\right)\sum_{i}\left(q_{i}-p_{i}\right)\log\left[\frac{\alpha p_{i}+\left(1-\alpha\right)q_{i}}{p_{i}}\right]
\end{equation}
We can also say that:
\begin{equation}
	JD\left(q\|p\right)=F\left(p\|q\right)+G\left(p\|q\right)
\end{equation}
That is:
\begin{equation}
	JD\left(q\|p\right)=KL\left[p\|\alpha p+\left(1-\alpha q\right)\right]+KL\left[\alpha p+\left(1-\alpha\right)q\|p\right]
\end{equation}
The derivation of the gradient with respect to ``$q$" leads to:
\begin{equation}
	\frac{\partial JD\left(q,p\right)}{\partial q_{j}}\propto\log\left[\frac{\alpha p_{j}+\left(1-\alpha\right)q_{j}}{p_{j}}\right]+\left(1-\alpha\right)\frac{q_{j}-p_{j}}{\alpha p_{j}+\left(1-\alpha\right)q_{j}}
\end{equation}

\subsection{Invariance by change of scale on ``$q$".}
For the divergence given by the expression (\ref{eq.JD}), the nominal invariance factor cannot be computed explicitly, so we use the expression $K^{*}$, which leads to the simplified invariant divergence:
\begin{equation}
JDI\left(p\|q\right)=\left(1-\alpha\right)\sum_{i}\left(\bar{p}_{i}-\bar{q}_{i}\right)\log\frac{\alpha\bar{q}_{i}+\left(1-\alpha\right)\bar{p}_{i}}{\bar{q}_{i}}
\label{eq.JDI}
\end{equation}
The gradient with respect to ``$q$" is given by:
\begin{align}
	\frac{\partial JDI\left(p\|q\right)}{\partial q_{l}}=&\frac{1-\alpha}{\sum_{j}q_{j}}
\left[\sum_{i}\bar{q}_{i}\log\frac{\alpha\bar{q}_{i}+\left(1-\alpha\right)\bar{p}_{i}}{\bar{q}_{i}}-\log\frac{\alpha\bar{q}_{l}+\left(1-\alpha\right)\bar{p}_{l}}{\bar{q}_{l}}\right. \nonumber \\  & \left.+ \left(1-\alpha\right)\frac{\bar{p}_{l}\left(\bar{q}_{l}-\bar{p}_{l}\right)}{\bar{q}_{l}\left[\alpha\bar{q}_{l}+\left(1-\alpha\right)\bar{p}_{l}\right]}-\left(1-\alpha\right)\sum_{i}\frac{\bar{p}_{i}\left(\bar{q}_{i}-\bar{p}_{i}\right)}{\alpha\bar{q}_{i}+\left(1-\alpha\right)\bar{p}_{i}}\right]
\label{eq.gradJDI}
\end{align}

\subsection{Generalization of Dragomir's J relative divergence.}
Starting from the standard convex function (\ref{eq.fcD}), the proposed generalization consists in replacing the logarithmic function by the ``Generalized logarithm" function with the exponent ($1-d$).\\
We thus obtain the standard convex function:
\begin{equation}
	f_{c,d}\left(x\right)=\left(1-\alpha\right)\left(x-1\right)\left\{\frac{\left[\alpha+\left(1-\alpha\right)x\right]^{1-d}-1}{1-d}\right\}
\end{equation}
If we use this function to construct a Csiszär divergence, we obtain:
\begin{equation}
	JD_{d}\left(p\|q\right)=\frac{1-\alpha}{1-d}\sum_{i}\left(p_{i}-q_{i}\right)\left\{\left[\frac{\alpha q_{i}+\left(1-\alpha\right)p_{i}}{q_{i}}\right]^{1-d}-1\right\}
\end{equation}
We can see that when $d\rightarrow $2 this divergence becomes equal to:
\begin{equation}
	\left(1-\alpha\right)^{2}\sum_{i}\frac{\left(p_{i}-q_{i}\right)^{2}}{\alpha q_{i}+\left(1-\alpha\right)p_{i}}\propto M_{AH}
\end{equation}
The dual form is obtained from the mirror function of the previous one, i.e.:
\begin{equation}
	\breve{f}_{c,d}\left(x\right)=\frac{1-\alpha}{1-d}\left(1-x\right)\left\{\left[\frac{\alpha x+\left(1-\alpha\right)}{x}\right]^{1-d}-1\right\}
\end{equation}
The corresponding divergence will be written:
\begin{equation}
	JD_{d}\left(q\|p\right)=\frac{1-\alpha}{1-d}\sum_{i}\left(q_{i}-p_{i}\right)\left\{\left[\frac{\alpha p_{i}+\left(1-\alpha\right)q_{i}}{p_{i}}\right]^{1-d}-1\right\}
\end{equation}

\section {Generalizations proposed by Taneja.}
The generalization technique used by Taneja \cite{taneja1989}, consists in relying on the divergence based on the entropy of Sharma-Mittal \cite{sharma1975} in the simplified form corresponding to normalized variables (\ref{eq.SMSa}).\\
It is not much more complicated to use the general form corresponding to any variables (\ref{eq.SMas}) that is written:
\begin{equation}	SM_{\alpha,s}\left(a\|b\right)=\frac{1}{\alpha\left(s-1\right)}\left\{\left[\sum_{i}a_{i}^{\alpha}b_{i}^{1-\alpha}\right]^{\frac{s-1}{\alpha-1}}-\left[\sum_{i}\alpha a_{i}+\left(1-\alpha\right)b_{i}\right]^{\frac{s-1}{\alpha-1}}\right\}
\label{eq.SMasgen}
\end{equation}
We can remark that it is in fact the application of the Generalized Logarithm function, on the two terms of $M_{AG}\left(a\|b\right)$, with the exponent $\frac{s-1}{\alpha-1}$.\\	
We then define the 1st generalization noted by Taneja $^{1}T^{s}_{\alpha}$:
\begin{equation}
^{1}T^{s}_{\alpha}\left(p\|q\right)=\frac{1}{2}\left[SM_{\alpha,s}\left(\frac{p+q}{2}\|p\right)+SM_{\alpha,s}\left(\frac{p+q}{2}\|q\right)\right]
\label{eq.UTsa}	
\end{equation}
We can observe that for the pleasure of complicating things further, we can imagine replacing the half-sum by a weighted sum, but is it really necessary to do so ?\\
That being said, the explicit form of $^{1}T^{s}_{\alpha}$ is written:
\begin{align}
^{1}T^{s}_{\alpha}\left(p\|q\right)=\frac{1}{2\alpha\left(s-1\right)}&
\left\{\left[\sum_{i}\left(\frac{p_{i}+q_{i}}{2}\right)^{\alpha} p_{i}^{1-\alpha}
\right]^{\frac{s-1}{\alpha-1}}\right. \nonumber \\ & \left.+\left[\sum_{i}\left(\frac{p_{i}+q_{i}}{2}\right)^{\alpha}q_{i}^{1-\alpha}\right]^{\frac{s-1}{\alpha-1}} \right. \nonumber \\ & \left. -\left[\sum_{i}\alpha \left(\frac{p_{i}+q_{i}}{2}\right)+\left(1-\alpha\right)p_{i}\right]^{\frac{s-1}{\alpha-1}}\right. \nonumber \\ & \left.-\left[\sum_{i}\alpha \left(\frac{p_{i}+q_{i}}{2}\right)+\left(1-\alpha\right)q_{i}\right]^{\frac{s-1}{\alpha-1}}\right\}
\label{eq.UTsaD}
\end{align}
If we consider the simplified configuration corresponding to $\sum_{i}p_{i}=\sum_{i}q_{i}=1$, we have:
\begin{align}
^{1}TS^{s}_{\alpha}\left(p\|q\right)=\frac{1}{2\alpha\left(s-1\right)}&
\left\{\left[\sum_{i}\left(\frac{p_{i}+q_{i}}{2}\right)^{\alpha} p_{i}^{1-\alpha}
\right]^{\frac{s-1}{\alpha-1}}\right. \nonumber \\ & \left.+\left[\sum_{i}\left(\frac{p_{i}+q_{i}}{2}\right)^{\alpha}q_{i}^{1-\alpha}\right]^{\frac{s-1}{\alpha-1}} -2 \right\}
\end{align}
Now let's look at the limit expressions of $^{1}T^{s}_{\alpha}\left(p\|q\right)$ (\ref{eq.UTsaD}).\\
 If $s\rightarrow\alpha$.\\
We obtain:
\begin{align}
^{1}T^{\alpha}_{\alpha}\left(p\|q\right)=\frac{1}{\alpha\left(\alpha-1\right)}&
\left\{\sum_{i}\left(\frac{p_{i}+q_{i}}{2}\right)^{\alpha} \left(\frac{p_{i}^{1-\alpha}+q_{i}^{1-\alpha}}{2}\right)
\right. \nonumber \\ & \left.-\sum_{i}\left(\frac{p_{i}+q_{i}}{2}\right)\right\}
\label{eq.UTaa}
\end{align}
It is a Csizär divergence based on the standard convex function:
\begin{equation}	f_{c}\left(x\right)=\frac{1}{\alpha\left(\alpha-1\right)}\left[\left(\frac{x^{1-\alpha}+1}{2}\right)\left(\frac{x+1}{2}\right)^{\alpha}-\frac{x+1}{2}\right]
\end{equation}
It is similar to (\ref{eq.FGsa}).\\
Note that this expression could have been obtained by using the Havrda-Charvat divergence directly in the general definition (\ref{eq.UTsa}) instead of the divergence $SM_{\alpha,s}$.\\
Furthermore, the special case corresponding to $\sum_{i}p_{i}=\sum_{i}q_{i}=1$ is immediately written as follows:
\begin{equation}
^{1}TS^{\alpha}_{\alpha}\left(p\|q\right)=\frac{1}{\alpha\left(\alpha-1\right)}
\left\{\sum_{i}\left(\frac{p_{i}+q_{i}}{2}\right)^{\alpha} \left(\frac{p_{i}^{1-\alpha}+q_{i}^{1-\alpha}}{2}\right)
-1\right\}
\end{equation}
Of course, it cannot be constructed as a divergence of Csiszär.\\

The second limit case is to make $s\rightarrow 1$ in $^{1}T^{s}_{\alpha}\left(p\|q\right)$ (\ref{eq.UTsaD}) (we're trying to obtain something that is inspired by Renyi's divergence).\\
 We obtain:
\begin{align}
^{1}T^{1}_{\alpha}\left(p\|q\right)=\frac{1}{2\alpha\left(\alpha-1\right)}&
\left\{\log\sum_{i}\left(\frac{p_{i}+q_{i}}{2}\right)^{\alpha} p_{i}^{1-\alpha}
\right. \nonumber \\ & \left.+\log\sum_{i}\left(\frac{p_{i}+q_{i}}{2}\right)^{\alpha}q_{i}^{1-\alpha}\right. \nonumber \\ & \left. -\log\sum_{i}\alpha \left(\frac{p_{i}+q_{i}}{2}\right)+\left(1-\alpha\right)p_{i}\right. \nonumber \\ & \left.-\log\sum_{i}\alpha \left(\frac{p_{i}+q_{i}}{2}\right)+\left(1-\alpha\right)q_{i}\right\}
\end{align}
The divergence corresponding to $\sum_{i}p_{i}=\sum_{i}q_{i}=1$ will be written:
\begin{align}
^{1}TS^{1}_{\alpha}\left(p\|q\right)=\frac{1}{2\alpha\left(\alpha-1\right)}&
\left\{\log\sum_{i}\left(\frac{p_{i}+q_{i}}{2}\right)^{\alpha} p_{i}^{1-\alpha}
\right. \nonumber \\ & \left.+\log\sum_{i}\left(\frac{p_{i}+q_{i}}{2}\right)^{\alpha}q_{i}^{1-\alpha}\right\}
\end{align}

The 2nd generalization proposed by Taneja consists in starting from the expression of $^{1}T^{\alpha}_{\alpha}\left(p\|q\right)$ (\ref{eq.UTaa}) which is a difference of two positive terms, and apply to each of them the Generalized Logarithm with the power $\left(s-1/\alpha-1\right)$ which is an increasing function; thus we obtain:
\begin{align}
^{2}T^{s}_{\alpha}\left(p\|q\right)=\frac{1}{\alpha\left(s-1\right)}&
\left\{\left[\sum_{i}\left(\frac{p_{i}+q_{i}}{2}\right)^{\alpha} \left(\frac{p_{i}^{1-\alpha}+q_{i}^{1-\alpha}}{2}\right)\right]^{\frac{s-1}{\alpha-1}}
\right. \nonumber \\ & \left.-\left[\sum_{i}\left(\frac{p_{i}+q_{i}}{2}\right)\right]^{\frac{s-1}{\alpha-1}}\right\}
\label{eq.DTaa}
\end{align}
It is easy to verify that $^{2}T^{\alpha}_{\alpha}\left(p\|q\right)=^{1}T^{\alpha}_{\alpha}\left(p\|q\right)$.\\
Similarly, by making $s\rightarrow1$ in (\ref{eq.DTaa}), we obtain:
\begin{align}
^{2}T^{1}_{\alpha}\left(p\|q\right)=\frac{1}{\alpha\left(\alpha-1\right)}&
\left\{\log\sum_{i}\left(\frac{p_{i}+q_{i}}{2}\right)^{\alpha} \left(\frac{p_{i}^{1-\alpha}+q_{i}^{1-\alpha}}{2}\right)
\right. \nonumber \\ & \left.-\log\sum_{i}\left(\frac{p_{i}+q_{i}}{2}\right)\right\}
\end{align}
In any case, if we consider the simplified version $\sum_{i}p_{i}=\sum_{i}q_{i}=1$, we find the forms given by Taneja.\\

\textbf{Remark:} To fully understand the difference between these 2 generalizations, we can notice that they correspond to the following operations:\\
we consider the Havrda-Charvat divergences with a parameter ``s", noted  $HC_{s}\left(a\|b\right)$ with $a=\frac{p+q}{2}$ et $b=p$ on one hand, and $HC_{s}\left(a\|b\right)$ with $a=\frac{p+q}{2}$ et $b=q$ on the other hand.\\
* - In the case of the 1st generalization, we apply the Generalized Logarithm on the 2 terms of each of these divergences, then we make the half-sum of the divergences obtained, which gives (\ref{eq.UTsaD}) and its variants.\\
* - In the 2nd generalization, we first sum the 2 Havrda-Charvat divergences mentioned above, we obtain a difference of 2 positive terms, then, we apply the Generalized Logarithm on each of these 2 terms.\\
The results obtained by these two operating modes are different because the two operations carried out successively are not commutative.\\

To continue this small exercise, we take up an equivalent work already proposed by Taneja \cite{taneja2001} and Arndt \cite{arndt2001}; it consists in using the relation (\ref{eq.UTsa}) and inverting the order of the arguments in the divergences $SM_{\alpha,s}$, that is:
\begin{equation}
^{1}I^{s}_{\alpha}\left(p\|q\right)=\frac{1}{2}\left[SM_{\alpha,s}\left(p\|\frac{p+q}{2}\right)+SM_{\alpha,s}\left(q\|\frac{p+q}{2}\right)\right]	
\end{equation}
This can be considered as the 3rd generalization which is similar to the 1st generalization, except for the order of the arguments.\\
This being said, the explicit form of $^{1}I^{s}_{\alpha}$ is written:
\begin{align}
^{1}I^{s}_{\alpha}\left(p\|q\right)=\frac{1}{2\alpha\left(s-1\right)}&
\left\{\left[\sum_{i}\left(\frac{p_{i}+q_{i}}{2}\right)^{1-\alpha} p_{i}^{\alpha}
\right]^{\frac{s-1}{\alpha-1}}\right. \nonumber \\ & \left.+\left[\sum_{i}\left(\frac{p_{i}+q_{i}}{2}\right)^{1-\alpha}q_{i}^{\alpha}\right]^{\frac{s-1}{\alpha-1}} \right. \nonumber \\ & \left. -\left[\sum_{i}\left(1-\alpha\right) \left(\frac{p_{i}+q_{i}}{2}\right)+\alpha p_{i}\right]^{\frac{s-1}{\alpha-1}}\right. \nonumber \\ & \left.-\left[\sum_{i}\left(1-\alpha\right) \left(\frac{p_{i}+q_{i}}{2}\right)+\alpha q_{i}\right]^{\frac{s-1}{\alpha-1}}\right\}
\label{eq.UIsa}
\end{align}
In the simplified case corresponding to $\sum_{i}p_{i}=\sum_{i}q_{i}=1$, we will have the simplified forms exhibited in Arndt \cite{arndt2001}:
\begin{align}
^{1}I^{s}_{\alpha}\left(p\|q\right)=\frac{1}{2\alpha\left(s-1\right)}&
\left\{\left[\sum_{i}\left(\frac{p_{i}+q_{i}}{2}\right)^{1-\alpha} p_{i}^{\alpha}
\right]^{\frac{s-1}{\alpha-1}}\right. \nonumber \\ & \left.+\left[\sum_{i}\left(\frac{p_{i}+q_{i}}{2}\right)^{1-\alpha}q_{i}^{\alpha}\right]^{\frac{s-1}{\alpha-1}}-2\right\}
\end{align}
Similarly, if in (\ref{eq.UIsa}) we use $s=\alpha$, we will have:
\begin{align}
^{1}I^{\alpha}_{\alpha}\left(p\|q\right)=\frac{1}{\alpha\left(\alpha-1\right)}&
\left\{\sum_{i}\left(\frac{p_{i}+q_{i}}{2}\right)^{1-\alpha} \left(\frac{p_{i}^\alpha+q_{i}^\alpha}{2}\right)
\right. \nonumber \\ & \left.-\sum_{i}\left(\frac{p_{i}+q_{i}}{2}\right)\right\}
\label{eq.UIaa}
\end{align}
It is a Csizär divergence based on the standard convex function:
\begin{equation}	f_{c}\left(x\right)=\frac{1}{\alpha\left(\alpha-1\right)}\left[\left(\frac{x^\alpha+1}{2}\right)\left(\frac{x+1}{2}\right)^{1-\alpha}-\frac{x+1}{2}\right]
\end{equation}
If now, from (\ref{eq.UIsa}), we change to $s\rightarrow 1$ , we get:
\begin{align}
^{1}I^{1}_{\alpha}\left(p\|q\right)=\frac{1}{2\alpha\left(\alpha-1\right)}&
\left\{\log\sum_{i}\left(\frac{p_{i}+q_{i}}{2}\right)^{1-\alpha} p_{i}^{\alpha}
\right. \nonumber \\ & \left.+\log\sum_{i}\left(\frac{p_{i}+q_{i}}{2}\right)^{1-\alpha}q_{i}^{\alpha}\right. \nonumber \\ & \left. -\log\sum_{i}\left(1-\alpha\right) \left(\frac{p_{i}+q_{i}}{2}\right)+\alpha p_{i}\right. \nonumber \\ & \left.-\log\sum_{i}\left(1-\alpha\right) \left(\frac{p_{i}+q_{i}}{2}\right)+\alpha q_{i}\right\}
\end{align}
If we take $\sum_{i}p_{i}=\sum_{i}q_{i}=1$, we recover a simplified result from Arndt \cite{arndt2001}:
\begin{align}
^{1}IS^{1}_{\alpha}\left(p\|q\right)=\frac{1}{2\alpha\left(\alpha-1\right)}&
\left\{\log\sum_{i}\left(\frac{p_{i}+q_{i}}{2}\right)^{1-\alpha} p_{i}^{\alpha}
\right. \nonumber \\ & \left.+\log\sum_{i}\left(\frac{p_{i}+q_{i}}{2}\right)^{1-\alpha}q_{i}^{\alpha}\right\}
\end{align}
We can then define the 4th generalization which is the analogous of the 2nd one (except for the order of the arguments); to do so, we start from the relation giving $^{1}I^{\alpha}_{\alpha}\left(p\|q\right)$ (\ref{eq.UIaa}) and an extension is carried out using the Generalized Logarithm, to obtain:
\begin{align}
^{2}I^{s}_{\alpha}\left(p\|q\right)=\frac{1}{\alpha\left(s-1\right)}&
\left\{\left[\sum_{i}\left(\frac{p_{i}+q_{i}}{2}\right)^{1-\alpha} \left(\frac{p_{i}^\alpha+q_{i}^\alpha}{2}\right)\right]^{\frac{s-1}{\alpha-1}}
\right. \nonumber \\ & \left.-\left[\sum_{i}\left(\frac{p_{i}+q_{i}}{2}\right)\right]^{\frac{s-1}{\alpha-1}}\right\}
\end{align}
This divergence was exhibited by Taneja in \cite{taneja1989}.\\
If we perform the simplification $\sum_{i}p_{i}=\sum_{i}q_{i}=1$, we retrieve an expression given by Arndt \cite{arndt2001}:\\
On this relationship, moving to the $s\rightarrow1$ limit gives us something close to a Renyi divergence:
\begin{align}
^{2}I^{1}_{\alpha}\left(p\|q\right)=\frac{1}{\alpha\left(\alpha-1\right)}&
\left\{\log\sum_{i}\left(\frac{p_{i}+q_{i}}{2}\right)^{1-\alpha} \left(\frac{p_{i}^\alpha+q_{i}^\alpha}{2}\right)
\right. \nonumber \\ & \left.-\log\sum_{i}\left(\frac{p_{i}+q_{i}}{2}\right)\right\}
\end{align}
Finally, making the simplification $\sum_{i}p_{i}=\sum_{i}q_{i}=$1, we retrieve Arndt's result \cite{arndt2001}:
\begin{equation}
^{2}IS^{1}_{\alpha}\left(p\|q\right)=\frac{1}{\alpha\left(\alpha-1\right)}
\left\{\log\sum_{i}\left(\frac{p_{i}+q_{i}}{2}\right)^{1-\alpha} \left(\frac{p_{i}^\alpha+q_{i}^\alpha}{2}\right)
\right\}
\end{equation}
On generalizations 3 and 4, we can make the same remarks as on generalizations 1 and 2 as to the order of the operations which made it possible to obtain them; the non-commutativity of these operations leads to \textit{``Generalization 3 $\neq$ Generalization 4"}, as we had obtained \textit{``Generalization 1 $\neq$ Generalization 2"}.\\
Another generalization (say the 5th generalization) is proposed by Taneja in \cite{taneja1989} and may be illustrative of the fact that:\\

 ``\textbf{Every inequality can give rise to a divergence}".\\
 
 Indeed, let us consider the inequality between weighted arithmetic mean $\left(MA\right)_{i}=\beta p_{i}+\left(1-\beta\right) q_{i}$ and the weighted geometric mean $\left(MG\right)_{i}=p^{\beta}_{i}q^{1-\beta}_{i}$; 
\begin{equation}
\left(MG\right)_{i}-\left(MA\right)_{i}\leq0\ \ \ \Rightarrow\ \ \left(MG\right)_{i}^{1-r}-\left(MA\right)_{i}^{1-r}\leq0	\ \ ; \ \ r\leq1
\end{equation}
And then:
\begin{equation}
\left(MG\right)_{i}^{1-r}\left(MA\right)_{i}^{r}-\left(MA\right)_{i}\leq0	
\end{equation}
The extension to the entire field leads to:
\begin{equation}
\sum_{i}\left(MG\right)_{i}^{1-r}\left(MA\right)_{i}^{r}-\sum_{i}\left(MA\right)_{i}\leq0	
\end{equation}
And finally:
\begin{equation}
\frac{1}{r-1}\left[\sum_{i}\left(MG\right)_{i}^{1-r}\left(MA\right)_{i}^{r}-\sum_{i}\left(MA\right)_{i}\right]\geq0	
\end{equation}
From there, we can apply on each of the terms of the difference, an increasing function without changing the sign of the inequality, why not a generalized logarithmic function with the power $\left(s-1/r-1\right)$, to obtain the divergence proposed by Taneja:
\begin{equation}
^{3}T^{s}_{r}\left(p\|q\right)=\frac{1}{s-1}
\left\{\left[\sum_{i}\left(MG\right)_{i}^{1-r}\left(MA\right)_{i}^{r}\right]^{\frac{s-1}{r-1}}-\left[\sum_{i}\left(MA\right)_{i}\right]^{\frac{s-1}{r-1}}\right\}	
\end{equation}
In fact, Taneja \cite{taneja1989} proposes something less general than that, because he limits it to the case of variables having a sum equal to 1, in which case the 2nd term of the difference will be equal to 1.\
The gradient with respect to ``$q_{j}$" is written as follows:
\begin{align}
	\frac{\partial \left[^{3}T^{s}_{r}\left(p\|q\right)\right]}{\partial q_{j}}=\frac{1-\beta}{r-1}&
\left\{\left[\sum_{i}\left(MG\right)_{i}^{1-r}\left(MA\right)_{i}^{r}\right]^{\frac{s-r}{r-1}}\left[\left(1-r\right)\frac{\left(MA\right)^{r}_{i}}{\left(MG\right)^{r}_{i}}\frac{p^{\beta}_{j}}{q^{\beta}_{j}}+r\frac{\left(MA\right)^{r-1}_{i}}{\left(MG\right)^{r-1}_{i}}\right]\right. \nonumber \\ & \left.
	-\left[\sum_{i}\left(MA\right)_{i}\right]^{\frac{s-r}{r-1}}\right\}
	\label{eq.^{3}T^{s}_{r}}
\end{align}
The 2 classical limit cases can be considered:\\

* $s=r$ leads to:
\begin{equation}
^{3}T^{r}_{r}\left(p\|q\right)=\frac{1}{r-1}\left\{\sum_{i}\left(MG\right)_{i}^{1-r}\left(MA\right)_{i}^{r}-\sum_{i}\left(MA\right)_{i}\right\}	
\end{equation}
* $s\rightarrow 1$ leads to the Logarithmic form:
\begin{equation}
^{3}T^{1}_{r}\left(p\|q\right)=\frac{1}{r-1}\left\{\log\left[\sum_{i}\left(MG\right)_{i}^{1-r}\left(MA\right)_{i}^{r}\right]-\log\left[\sum_{i}\left(MA\right)_{i}\right]\right\}	
\end{equation}
The corresponding gradients with respect to ``$q$" are obtained immediately as a limit case of (\ref{eq.^{3}T^{s}_{r}}).\\

This last divergence is referred to in Taneja \cite{taneja1989} and Neemuchwala \cite{neemuchwala2007} with the name ``$\alpha GA$" which can be noted here as ``$rGA$", and for which they only consider the case of standardized variables.\\
Of course, this fifth generalization is in relation with the divergences based on inequalities between the generalized means developed in \textbf{chapter 7}.\\

If we want to extend this procedure to other means, we rely on the inequalities between weighted and unweighted means  $MH\leq MG\leq MA\leq MQ$.\\
 With:
\begin{equation}
	\left(MH\right)_{i}=\frac{p_{i}q_{i}}{\beta p_{i}-\left(1-\beta\right)q_{i}}\ \ \ \ ; \ \ \ \ 	\left(MQ\right)_{i}=\sqrt{\beta p^{2}_{i}+\left(1-\beta\right)q^{2}_{i}}
\end{equation}
We can write:
\begin{equation}
^{3}THG^{s}_{r}\left(p\|q\right)=\frac{1}{s-1}\left\{\left[\sum_{i}\left(MH\right)_{i}^{1-r}\left(MG\right)_{i}^{r}\right]^{\frac{s-1}{r-1}}-\left[\sum_{i}\left(MG\right)_{i}\right]^{\frac{s-1}{r-1}}\right\}	
\end{equation}\\
\begin{equation}
^{3}THA^{s}_{r}\left(p\|q\right)=\frac{1}{s-1}\left\{\left[\sum_{i}\left(MH\right)_{i}^{1-r}\left(MA\right)_{i}^{r}\right]^{\frac{s-1}{r-1}}-\left[\sum_{i}\left(MA\right)_{i}\right]^{\frac{s-1}{r-1}}\right\}	
\end{equation}\\
\begin{equation}
^{3}THQ^{s}_{r}\left(p\|q\right)=\frac{1}{s-1}\left\{\left[\sum_{i}\left(MH\right)_{i}^{1-r}\left(MQ\right)_{i}^{r}\right]^{\frac{s-1}{r-1}}-\left[\sum_{i}\left(MQ\right)_{i}\right]^{\frac{s-1}{r-1}}\right\}	
\end{equation}\\
The divergence:
\begin{equation}
^{3}TGA^{s}_{r}\left(p\|q\right)=\:^{3}T^{s}_{r}\left(p\|q\right)
\end{equation}
has been previously expressed;
\begin{equation}
^{3}TGQ^{s}_{r}\left(p\|q\right)=\frac{1}{s-1}\left\{\left[\sum_{i}\left(MG\right)_{i}^{1-r}\left(MQ\right)_{i}^{r}\right]^{\frac{s-1}{r-1}}-\left[\sum_{i}\left(MQ\right)_{i}\right]^{\frac{s-1}{r-1}}\right\}	
\end{equation}\\
\begin{equation}
^{3}TAQ^{s}_{r}\left(p\|q\right)=\frac{1}{s-1}\left\{\left[\sum_{i}\left(MA\right)_{i}^{1-r}\left(MQ\right)_{i}^{r}\right]^{\frac{s-1}{r-1}}-\left[\sum_{i}\left(MQ\right)_{i}\right]^{\frac{s-1}{r-1}}\right\}	
\end{equation}\\
Logarithmic forms can be inferred immediately.\\
Simplified forms of these divergences are developed in \textbf{chapter 7} dealing with divergences between means.
\setcounter{table}{0}  \setcounter{equation}{0}  \setcounter{figure}{0} \setcounter{chapter}{9} \setcounter{section}{0} 
\chapter{chapter 9 -\\Smoothness constraint regularization.}  \label{chptr::chapitre9}
We consider here, the application of scale invariant divergences to Tikhonov-type smoothness constraint regularization \cite{tikhonov1974methods}.\\
This analysis has an obvious interest when one must introduce this type of regularization in a problem for which there is a constraint on the sum of the unknown parameters and, in particular, when it is associated with a non-negativity constraint.
Indeed, when regularizing an inverse problem by a smoothness constraint, we are led to minimize in ``$x$", an objective function of the form:
\begin{equation}
J\left(x\right)=D_{1}\left(y\|Hx\right)+\gamma D_{2}\left(x\|x_{d}\right)	
\end{equation}
The first term $D_{1}$ is a divergence which represents the attachment to the data (data consistency), i.e. a discrepancy between the measurements and the model (here linear), whereas the second term $D_{2}$ is a divergence which represents a discrepancy between the solution ``$x$" and a default solution ``$x_{d}$"; this default solution must obviously satisfy the constraints imposed on the solution, for example, $\left[x_{d}\right]_{i}\geq 0\;\forall i$, and if we have a sum constraint, $\sum_{i}\left[x_{d}\right]_{i}=\sum_{i}x_{i}$.\\
The coefficient ``$\gamma$" is the regularization factor that adjusts the relative importance of the 2 terms of $J\left(x\right)$.\\
We consider here, the case of divergences which are invariant by change of scale on ``$q$" constructed on the basis of ``$\alpha$" and ``$\beta$" divergences, and we consider more specifically the logarithmic forms which are not only invariant with respect to ``$q$", but also with respect to ``$p$".\\
On the other hand, it was emphasized that when a divergence is made invariant with respect to ``$q$" using the invariance factor $K^{*}=\frac{\sum_{j}p_{j}}{\sum_{j}q_{j}}$, the resulting divergence  is also invariant with respect to ``$p$"; this type of invariant divergence will also be considered.\\
The invariant divergences with respect to both arguments have the advantage of applying in all cases, whether the default solution is a constant, or a linear function of the solution itself (Laplacian case for example).\\
In the minimization problem, we will use the expression of the gradient:
\begin{equation}
\frac{\partial J\left(x\right)}{\partial x_{l}}=\frac{\partial D_{1}\left(y\|Hx\right)}{\partial x_{l}}+\gamma \frac{\partial D_{2}\left(x\|x_{d}\right)}{\partial x_{l}}	
\end{equation}
If the divergences $D_{1}$ and $D_{2}$ are invariant divergences, the characteristic property of such divergences will result in:
\begin{equation}
	\sum_{l}x_{l}\frac{\partial J\left(x\right)}{\partial x_{l}}=0
\end{equation}

\section{Invariant divergences derived from the ``$\alpha$" and ``$\beta$" divergences.}
\subsection{Constant default solution = C.}
In this case, all invariant divergences can be used, whether logarithmic ($LAI, LBI$) or not ($AI, BI$); we will just develop the case for logarithmic forms.\\
In the following discussions, for the two-dimensional case, the tables under consideration are written in lexicographical order.

\subsubsection{* With ``$LAI$".}
We have:
\begin{equation}	LAI\left(c\|x\right)=\frac{1}{a}\log\sum_{i}x_{i}-\frac{1}{a-1}\log\sum_{i}c_{i}+\frac{1}{a\left(a-1\right)}\log\sum_{i}c_{i}^{a}x_{i}^{1-a}
\end{equation}
Taking into account the fact that all $c_{i}$ are equals whatever ``$i$'', we have the gradient:
\begin{equation}
	\frac{\partial 	LAI\left(c\|x\right)}{\partial x_{j}}=\frac{1}{a}\left[\frac{1}{\sum_{i}x_{i}}-\frac{x_{j}^{-a}}{\sum_{i}x_{i}^{1-a}}\right]
\end{equation}

\subsubsection{* With ``$LBI$".}
We then have:
\begin{equation}	LBI\left(c\|x\right)=\frac{1}{b\left(b-1\right)}\log\sum_{i}c_{i}^{b}+\frac{1}{b}\log\sum_{i}x_{i}^{b}-\frac{1}{b-1}\log\sum_{i}c_{i}x_{i}^{b-1}
\end{equation}
Which gradient is:
\begin{equation}
	\frac{\partial 	LBI\left(c\|x\right)}{\partial x_{j}}=\frac{x_{j}^{b-1}}{\sum_{i}x_{i}^{b}}-\frac{x_{j}^{b-2}}{\sum_{i}x_{i}^{b-1}}
\end{equation}

\subsubsection{* With the general form ``$LABI$''.}
The divergence considered is written as follows:
\begin{align}
		LABI\left(c\|x\right)=
		&\frac{1}{a+b-1}\left\{\frac{1}{b-1}\log\sum_{i}c_{i}^{a+b-1}-\right.\nonumber \\
		&\frac{a+b-1}{a\left(b-1\right)}\log \sum_{i}c_{i}^{a}x_{i}^{b-1}+ \nonumber \\
&\left.\frac{1}{a}\log \sum_{i}x_{i}^{a+b-1}\right\}
\end{align}
After some simple calculations, the gradient is written as follows:
\begin{equation}
	\frac{\partial LABI\left(c\|x\right)}{\partial x_{j}}=\frac{1}{a}\left[\frac{x_{j}^{a+b-2}}{\sum_{i}x_{i}^{a+b-1}}-\frac{x_{j}^{b-2}}{\sum_{i}x_{i}^{b-1}}\right]
\end{equation}
The two previous particular cases can be found without difficulty by respectively making $a+b-1=$1, i.e. $b-1=1-a$ (LAI) and $a=1$ (LBI).\\

\subsection{Use of the Laplacian operator.}
In this case, the default solution ``$x_{d}$" is a smoothed version of the solution ``$x$", i.e. in a trivial way, a solution which does not contain high frequencies or anyway they are reduced; in the classical sense of Tikhonov \cite{tikhonov1974methods}, the regularization term is written as follows:
\begin{equation}
	D_{2}\left(x\|x_{d}\right)=\left\|Lx\right\|^{2}
\end{equation}
Here we must distinguish between 2 situations:\\

\textbf{* 1 - One-dimensional case.}\\
In this case, the regularization term ``$Lx$" corresponds to the convolution of ``$x$" by the mask $\left[-\frac{1}{2}\  ;\  1 \ ;\  -\frac{1}{2}\right]$.\\
We can therefore write the result in the form ``$x-Tx$" where the term ``$Tx$" written in the sense of the matrix vector product, represents the convolution of ``$x$" by the mask $\left[\frac{1}{2}\ ;\  0 \ ;\  \frac{1}{2}\right]$.\\
 
\textbf{* 2 - Two-dimensional case.}\\
In this case, the regularization term ``$Lx$" corresponds to the convolution of the table ``$x$" by the mask: 
$
\begin{pmatrix}0&-1/4&0\\ 
-1/4&1&-1/4\\ 
0&-1/4&0\\ 
\end{pmatrix}.  
$
We can therefore write the result in the form ``$x-Tx$" where the term ``$Tx$" written in the sense of the matrix vector product, represents the convolution of ``$x$" by the mask:
$
\begin{pmatrix}0&1/4&0\\ 
1/4&0&1/4\\ 
0&1/4&0\\ 
\end{pmatrix}.
$\\

At the end of this operation, the table ``$Tx$" is ordered lexicographically, as well as the table ``$x$".\\

Note that in these conditions, the matrix $T$ is symmetrical, its columns are normed to $1$, that is $\sum_{i}t_{ij}=1\ \forall j$, then $\sum_{i}\left(Tx\right)_{i}=\sum_{i}x_{i}$ .\\
Consequently, as a regularization term one must write a divergence invariant by scale change between $p\equiv x$ and $q\equiv Tx$ for example; since ``$p$" and ``$q$" both depend on the ``$x$" solution, it is understandable that we must use the Logarithmic form which is invariant with respect to both arguments (we can verify that the non Logarithmic forms do not allow us to develop algorithms that spontaneously maintain the flow, i.e. which do not spontaneously respect the sum constraint).

\subsubsection{* With ``$LAI$''.}
Using the notations $p\equiv x$ and $q\equiv Tx$, the divergence is expressed as:
\begin{equation}	
LAI\left(x\|Tx\right)=\frac{1}{a}\log\sum_{i}\left(Tx\right)_{i}-\frac{1}{a-1}\log\sum_{i}x_{i}+\frac{1}{a\left(a-1\right)}\log\sum_{i}x_{i}^{a}\left(Tx\right)_{i}^{1-a}
\label{eq.LAIlap}
\end{equation}
The associated gradient is written:
\begin{align}
\frac{\partial LAI\left(x\|Tx\right)}{\partial x_{j}}=
&\frac{1}{a}\left[\frac{\sum_{i}t_{i,j}}{\sum_{i}\left(Tx\right)_{i}}-\frac{\sum_{i}t_{i,j}x_{i}^{a}\left(Tx\right)_{i}^{-a}}{\sum_{i}x_{i}^{a}\left(Tx\right)_{i}^{1-a}}\right]+\nonumber \\
&\frac{1}{a-1}\left[\frac{x_{j}^{a-1}\left(Tx\right)_{j}^{1-a}}{\sum_{i}x_{i}^{a}\left(Tx\right)_{i}^{1-a}}-\frac{1}{\sum_{i}x_{i}}\right]
\label{eq.gradLAIlap}
\end{align}
We may observe that, as always with this type of divergence:\\
$\sum_{j}x_{j}\frac{\partial }{\partial x_{j}}=0$.

\subsubsection{* With ``$LBI$".}
With the notations in the previous paragraph, the divergence considered is written as follows:
\begin{equation}	
LBI\left(x\|Tx\right)=\frac{1}{b\left(b-1\right)}\log\sum_{i}x_{i}^{b}+\frac{1}{b}\log\sum_{i}\left(Tx\right)_{i}^{b}-\frac{1}{b-1}\log\sum_{i}x_{i}\left(Tx\right)_{i}^{b-1}
\label{eq.LBIlap}
\end{equation}
The gradient will be:
\begin{align}
\frac{\partial LBI\left(x\|Tx\right)}{\partial x_{j}}=
&\frac{1}{b-1}\left[\frac{x_{j}^{b-1}}{\sum_{i}x_{i}^{b}}-\frac{\left(Tx\right)_{j}^{b-1}}{\sum_{i}x_{i}\left(Tx\right)_{i}^{b-1}}\right]+\nonumber\\
&\left[\frac{\sum_{i}t_{i,j}\left(Tx\right)_{i}^{b-1}}{\sum_{i}\left(Tx\right)_{i}^{b}}-\frac{\sum_{i}t_{i,j}x_{i}\left(Tx\right)_{i}^{b-2}}{\sum_{i}x_{i}\left(Tx\right)_{i}^{b-1}}\right]
\label{eq.gradLBIlap}
\end{align}
Again, as always with this type of divergence we have:\\
$\sum_{j}x_{j}\frac{\partial }{\partial x_{j}}=0$.

\subsubsection{* With the general form ``$LABI$".}
If we use the general form of the divergence, the two previous cases can be recovered without any problem, indeed, the divergence is written as follows:
\begin{align}
		LABI\left(x\|Tx\right)=
		&\frac{1}{a+b-1}\left\{\frac{1}{b-1}\log\sum_{i}x_{i}^{a+b-1}- \right. \nonumber \\
		&\frac{a+b-1}{a\left(b-1\right)}\log \sum_{i}x_{i}^{a}\left(Tx\right)_{i}^{b-1}+ \nonumber \\
		&\left. \frac{1}{a}\log \sum_{i}\left(Tx\right)_{i}^{a+b-1}\right\}
\end{align}
We can calculate the gradient, which is expressed as follows:
\begin{align}
\frac{\partial LABI\left(x\|Tx\right)}{\partial x_{j}}= &\frac{1}{a}\left[\frac{\sum_{i}t_{i,j}\left(Tx\right)_{i}^{a+b-2}}{\sum_{i}\left(Tx\right)_{i}^{a+b-1}}-\frac{\sum_{i}t_{i,j}x_{i}^{a}\left(Tx\right)_{i}^{b-2}}{\sum_{i}x_{i}^{a}\left(Tx\right)_{i}^{b-1}}\right]+ \nonumber \\ 
&\frac{1}{b-1}\left[\frac{x_{j}^{a+b-2}}{\sum_{i}x_{i}^{a+b-1}}-\frac{x_{j}^{a-1}\left(Tx\right)_{j}^{b-1}}{\sum_{i}x_{i}^{a}\left(Tx\right)_{i}^{b-1}}\right]
\end{align}

The expressions of the two preceding paragraphs are found immediately by using respectively: $a+b-1=1$, i.e. $b-1=1-a$ (which gives $LAI$) and $a=1$ (which gives $LBI$).

\subsection{Important Note.}
If, in the Laplacian regularization, we try to use the non logarithmic forms of the invariant divergences, that is ``$AI$'' or ``$BI$'', we loose the property of spontaneous maintaining of the flow, in fact these divergences are invariant only with respect to the second argument; to show this, it is enough to evaluate the gradient of $AI\left(x\|Tx\right)$ for example, and to note that the relation which allows the maintaining of the flow is not fulfilled.\\
In fact, if we use the expression of ``$AI$'', where $p\equiv x$ and $q\equiv Tx$, it comes:

\begin{equation}	AI\left(x\|Tx\right)=\frac{1}{a-1}\left[\left(\sum_{i}\left(Tx\right)_{i}\right)^{\frac{a-1}{a}}\left(\sum_{i}x_{i}^{a}\left(Tx\right)_{i}^{1-a}\right)^{\frac{1}{a}}-\sum_{i}x_{i}\right]
\end{equation}
And according to the expression of the gradient we have:
\begin{equation}
\sum_{j}x_{j}\frac{\partial AI\left(x\|Tx\right)}{\partial x_{j}}=	\frac{1}{1-a}\left\{\sum_{j}x_{j}-\left[\sum_{j}x_{j}^{a}\left(Tx\right)_{j}^{1-a}\right]^{\frac{1}{a}}\left[\sum_{j}\left(Tx\right)_{j}\right]^{\frac{a-1}{a}}\right\}
\end{equation}
This expression could be equal to zero at convergence $(x\rightarrow Tx)$, but certainly not at each iteration.

\section{Effect of the invariance factor $K^{*}=\frac{\sum_{j}p_{j}}{\sum_{j}q_{j}}$.}
In the study of the invariance factors K, the particular role of $K^{*}=\frac{\sum_{j}p_{j}}{\sum_{j}q_{j}}$ has been noted; indeed, the introduction of this factor in any divergence $D\left(p\|q\right)$ , in order to make it invariant with respect to ``$q$", led to a divergence $DI\left(p\|q\right)$ which is invariant with respect to this argument, obtained easily from the divergence $D\left(p\|q\right)$ by replacing ``$p$" and ``$q$" by normalized variables $\bar{p_{i}}=\frac{p_{i}}{\sum_{j}p_{j}}$ and $\bar{q_{i}}=\frac{q_{i}}{\sum_{j}q_{j}}$ respectively, with a multiplying factor depending only on $\sum_{j}p_{j}$ that can be omitted.\\
Disregarding this multiplicative factor, we obtain a divergence that is not only invariant with respect to ``$q$", but also with respect to ``$p$".\\
This observation allows us to consider the use of such divergences in the case of smoothness constraint regularization using the Laplacian of the solution, as developed in the previous sections; it is a possible alternative to the use of logarithmic forms ``$LAI$" et``$LBI$".\\
Some examples to illustrate this point are given in Appendix 6.\\

\subsection{Applications of some divergences to regularization.}

We analyze in this paragraph the regularization in the sense of Tikhonov \cite{tikhonov1974methods} with the Laplacian, as previously considered for the divergences $LAI\left(p\|q\right)$ and $LBI\left(p\|q\right)$.\\
We will focus on the cases of a few classical divergences that are rendered invariant by the use of the invariance factor $K^{*}=\frac{\sum_{j}p_{j}}{\sum_{j}q_{j}}$; note that for the K.L. divergence, this factor is the nominal invariance factor $K_{0}$.\\
For the application under consideration here, we have:
\begin{equation}
\bar{p}_{i}=\frac{x_{i}}{\sum_{j}x_{j}},\ \ \ \ \ \bar{q}_{i}=\frac{\left(Tx\right)_{i}}{\sum_{j}\left(Tx\right)_{j}}
\label{eq.varred}                               
\end{equation}
With $\sum_{j}t_{jl}=1\ \forall l$ (the columns of $T$ are normalized to $1$), it comes:
\begin{equation}
\frac{\partial \bar{p}_{i}}{\partial x_{l}}=\frac{1}{\sum_{j}x_{j}}\left(\delta_{i,l}-\bar{p}_{i}\right),\ \ \ \ \frac{\partial \bar{q}_{i}}{\partial x_{l}}=\frac{1}{\sum_{j}\left(Tx\right)_{j}}\left(t_{il}-\bar{q}_{i}\right)
\label{eq.deviv part}	
\end{equation}

\subsubsection{1 - Mean square deviation.}
The invariant divergence we obtain is written as follows (see (\ref{eq.EQMInorm}) in Appendix 6):
\begin{equation}
	EQMI\left(p\|q\right)=\sum_{i}\left(\bar{p}_{i}-\bar{q}_{i}\right)^{2}
\end{equation}
Then:
\begin{equation}
	\frac{\partial EQMI\left(p\|q\right)}{\partial x_{l}}=\sum_{i}\left(\bar{p}_{i}-\bar{q}_{i}\right)\left(\frac{\partial \bar{p}_{i}}{\partial x_{l}}-\frac{\partial \bar{q}_{i}}{\partial x_{l}}\right)
\end{equation}
Taking into account the fact that $\sum_{j}x_{j}=\sum_{j}\left(Tx\right)_{j}$, we obtain after a few simple calculations:
\begin{equation}
\frac{\partial EQMI\left(p\|q\right)}{\partial x_{l}}=\frac{1}{\sum_{j}x_{j}}\left[\left(\bar{p}_{l}-\bar{q}_{l}\right)-\sum_{i}\left(\bar{p}_{i}-\bar{q}_{i}\right)^{2}-\sum_{i}t_{il}\left(\bar{p}_{i}-\bar{q}_{i}\right)\right]	
\end{equation}
or also:
\begin{equation}
\frac{\partial EQMI\left(p\|q\right)}{\partial x_{l}}=\frac{1}{\sum_{j}x_{j}}\left\{\left[\left(I-T^{T}\right)\left(\bar{p}-\bar{q}\right)\right]_{l}-\sum_{i}\left(\bar{p}_{i}-\bar{q}_{i}\right)^{2}\right\}	
\end{equation}
Coming back to the true variables (\ref{eq.varred}), it comes:                               
\begin{equation}
\frac{\partial EQMI\left(p\|q\right)}{\partial x_{l}}=\frac{1}{\left(\sum_{j}x_{j}\right)^{2}}\left\{\left[\left(I-T^{T}\right)\left(x-Tx\right)\right]_{l}-\frac{1}{\sum_{j}x_{j}}\sum_{i}\left[x_{i}-\left(Tx\right)_{i}\right]^{2}\right\}	
\end{equation}
We can verify that, as always:
\begin{equation}
	\sum_{l}x_{l}\frac{\partial EQMI\left(p\|q\right)}{\partial x_{l}}=0
\end{equation}

\subsubsection{2 - Kullback-Leibler divergence.}

Taking into account the invariance factor that is used, and the possible simplifications, the following is obtained (see (\ref{eq.KLInorm}) in Appendix 6):
\begin{equation}
	KLI\left(p\|q\right)=\sum_{i}\bar{p}_{i}\log\frac{\bar{p}_{i}}{\bar{q}_{i}}
\end{equation}
Then:
\begin{equation}
	\frac{\partial KLI\left(p\|q\right)}{\partial x_{l}}=\sum_{i}\left(\log \frac{\bar{p}_{i}}{\bar{q}_{i}}+1\right)\frac{\partial \bar{p}_{i}}{\partial x_{l}}-\sum_{i}\frac{\bar{p}_{i}}{\bar{q}_{i}}\frac{\partial \bar{q}_{i}}{\partial x_{l}}
\end{equation}
With the partial derivative expressions previously mentioned (\ref{eq.deviv part}), and taking into account that $\sum_{j}x_{j}=\sum_{j}\left(Tx\right)_{j}$, it comes all calculations made:
\begin{equation}
	\frac{\partial KLI\left(p\|q\right)}{\partial x_{l}}=\frac{1}{\sum_{j}x_{j}}\left[\log\frac{\bar{p}_{l}}{\bar{q}_{l}}+1-\sum_{i}t_{il}\frac{\bar{p}_{i}}{\bar{q}_{i}}-\sum_{i}\bar{p}_{i}\log\frac{\bar{p}_{i}}{\bar{q}_{i}}\right]	
\end{equation}
With $\sum_{i}t_{il}=1$, we can write:
\begin{equation}
	\frac{\partial KLI\left(p\|q\right)}{\partial x_{l}}=\frac{1}{\sum_{j}x_{j}}\left\{\log\frac{\bar{p}_{l}}{\bar{q}_{l}}-\sum_{i}\bar{p}_{i}\log\frac{\bar{p}_{i}}{\bar{q}_{i}}+\left[T^{T}\left(1-\frac{\bar{p}}{\bar{q}}\right)\right]_{l}\right\}	
\end{equation}
Coming back to the true variables (\ref{eq.varred}), it comes:
\begin{equation}
\frac{\partial KLI\left(p\|q\right)}{\partial x_{l}}=\frac{1}{\sum_{j}x_{j}}\left\{\log\frac{x_{l}}{\left(Tx\right)_{l}}-\sum_{i}\frac{x_{i}}{\sum_{j}x_{j}}\log\frac{x_{i}}{\left(Tx\right)_{i}}+\left[T^{T}\left(1-\frac{x}{\left(Tx\right)}\right)\right]_{l}\right\}	
\end{equation}
We can verify that, as always:
\begin{equation}
	\sum_{l}x_{l}\frac{\partial KLI\left(p\|q\right)}{\partial x_{l}}=0
\end{equation}

\subsubsection{3 - Neyman's Chi2 divergence.}

The invariant divergence is written (see (\ref{eq.chi2Nnorm}) in Appendix 6):
\begin{equation}
	\chi^{2}_{N}I\left(p\|q\right)=\sum_{i}\frac{\left(\bar{p}_{i}-\bar{q}_{i}\right)^{2}}{\bar{q}_{i}}
\end{equation}
Then:
\begin{equation}
\frac{\partial \chi^{2}_{N}I\left(p\|q\right)}{\partial x_{l}}=2\sum_{i}\left(\frac{\bar{p_{i}}}{\bar{q_{i}}}-1\right)\frac{\partial \bar{p_{i}}}{\partial x_{l}}+\sum_{i}\left[1-\left(\frac{\bar{{p}_{i}}}{\bar{q}_{i}}\right)^{2}\right]\frac{\partial \bar{q_{i}}}{\partial x_{l}}
\end{equation}
Considering the expressions of the partial derivatives of $\bar{p_{i}}$ and $\bar{q_{i}}$ with respect to ``$x_{l}$" (\ref{eq.deviv part}), it comes all calculations made, and with the simplifications already indicated:
\begin{equation}
\frac{\partial \chi^{2}_{N}I\left(p\|q\right)}{\partial x_{l}}=\frac{1}{\sum_{j}x_{j}}\left[2\frac{\bar{p_{l}}}{\bar{q_{l}}}-\sum_{i}t_{il}\left(\frac{\bar{p_{i}}}{\bar{q_{i}}}\right)^{2}-\sum_{i}\frac{\bar{p}^{2}_{i}}{\bar{q_{i}}}\right]	
\end{equation}
Turning back to the true variables, we have:
\begin{equation}
\frac{\partial \chi^{2}_{N}I\left(p\|q\right)}{\partial x_{l}}=\frac{1}{\sum_{j}x_{j}}\left\{2\frac{x_{l}}{\left(Tx\right)_{l}}-\left[T^{T}\frac{x^{2}}{\left(Tx\right)^{2}}\right]_{l}-\frac{1}{\sum_{j}x_{j}}\sum_{i}\frac{x^{2}_{i}}{\left(Tx\right)}_{i}\right\}	
\end{equation}
And, then:
\begin{equation}
	\sum_{l}x_{l}\frac{\partial \chi^{2}_{N}I\left(p\|q\right)}{\partial x_{l}}=0
\end{equation}

\subsubsection{4 - Pearson's Chi2 divergence.}
The invariant divergence is written (see (\ref{chi2Pnorm}) in Appendix 6):
\begin{equation}
	\chi^{2}_{P}I\left(p\|q\right)=\sum_{i}\frac{\left(\bar{p}_{i}-\bar{q}_{i}\right)^{2}}{\bar{p}_{i}}
\end{equation}
Then:
\begin{equation}
\frac{\partial \chi^{2}_{P}I\left(p\|q\right)}{\partial x_{l}}=\sum_{i}\left[1-\left(\frac{\bar{q_{i}}}{\bar{p_{i}}}\right)^{2}\right]\frac{\partial \bar{p_{i}}}{\partial x_{l}}+2\sum_{i}\left[\left(\frac{\bar{{q}_{i}}}{\bar{p}_{i}}\right)-1\right]\frac{\partial \bar{q_{i}}}{\partial x_{l}}
\end{equation}
Considering the expressions of the partial derivatives of $\bar{p_{i}}$ and $\bar{q_{i}}$ with respect to ``$x_{l}$" (\ref{eq.deviv part}), it comes with all the calculations made, with the simplifications already indicated:
\begin{equation}
\frac{\partial \chi^{2}_{P}I\left(p\|q\right)}{\partial x_{l}}=\frac{1}{\sum_{j}x_{j}}\left[-\left(\frac{\bar{q_{l}}}{\bar{p_{l}}}\right)^{2}+2\sum_{i}t_{il}\frac{\bar{{q}_{i}}}{\bar{p}_{i}}-\sum_{i}\frac{\bar{q}^{2}_{i}}{\bar{p}_{i}}\right]
\end{equation}
And, turning back to the true variables:
\begin{equation}
\frac{\partial \chi^{2}_{P}I\left(p\|q\right)}{\partial x_{l}}=\frac{1}{\sum_{j}x_{j}}\left\{-\left(\frac{\left(Tx\right)_{l}}{x_{l}}\right)^{2}+2\left[T^{T}\frac{\left(Tx\right)}{x}
\right]_{l}-\frac{1}{\sum_{j}\left(Tx\right)_{j}}\sum_{i}\frac{\left(Tx\right)^{2}_{i}}{x_{i}}\right\}
\end{equation}
Then, we have:
\begin{equation}
	\sum_{l}x_{l}\frac{\partial \chi^{2}_{P}I\left(p\|q\right)}{\partial x_{l}}=0
\end{equation}
\setcounter{table}{0}  \setcounter{equation}{0}  \setcounter{figure}{0} \setcounter{chapter}{10} \setcounter{section}{0} 
\chapter{chapitre 10 -\\Algorithmic point of view.}  \label{chptr::chapitre10}
\section{Some preliminary remarks.}
In this chapter, different algorithmic methods for minimizing divergences, under non-negativity constraint and under sum constraint, are described. The divergences considered are assumed to be convex and differentiable.\\
We will first consider the minimization problem under the non-negativity constraint.\\
We will distinguish three cases corresponding to specific problems:\\
* Simplified divergences that can be used for fields of explicitly equal sums (possibly equal to 1).\\
* General divergences which contain no simplification related to the sum of the data fields.\\
* Divergences invariant by change of scale.\\
In a second step, one will introduce in addition, the sum constraint on the unknown parameters.\\

We're essentially interested in the separable divergences that will be noted:
\begin{equation}
	D\left(p\|q\right)=\sum_{i}D\left(p_{i}\|q_{i}\right)=\sum_{i}d_{i}
\end{equation}

The non-separable case is easily deduced from this.\\

Moreover, we will only consider the case of a linear model:
\begin{equation}
	q=Hx\ \ \ \Rightarrow\ \ \ q_{i}=\sum_{j}h_{ij}x_{j}\ \ \ \Rightarrow\ \ \ \frac{\partial q_{i}}{\partial x_{l}}=h_{il}
\end{equation}
And we will have:
\begin{equation}
	\frac{\partial D\left(p\|q\right)}{\partial x_{l}}=\sum_{i}\frac{\partial d_{i}}{\partial q_{i}}\frac{\partial q_{i}}{\partial x_{l}}
	\label{eq.dDsdq}
\end{equation}
Consequently, with the notation:
\begin{equation}
	\frac{\partial d_{i}}{\partial q_{i}}=\left[\frac{\partial D}{\partial q}\right]_{i}
\end{equation}
We will have:
\begin{equation}
		\frac{\partial D\left(p\|q\right)}{\partial x_{l}}=\sum_{i}h_{il}\frac{\partial d_{i}}{\partial q_{i}}=\left[H^{T}\frac{\partial D}{\partial q}\right]_{l}
\end{equation}
This justifies ``\textit{a posteriori}'', the fact that for the divergences considered in this work, the expression of the gradient with respect to the variable ``$q$" has been given almost systematically.
\section{Split Gradient Method (S.G.M.).}
In this section, we present the method of minimization under non-negativity constraint, widely developed elsewhere \cite{lanteri2001}, \cite{lanteri2002}, which allows us to understand the origin of multiplicative algorithms and we can easily see why these algorithms can suffer from convergence problems.\\
This method applies to divergences $D\left(p\|q\right)$ having a minimum for $p_{i}=q_{i}\ \;\forall i$; it is quite obvious that this excludes simplified divergences which do not have this property. We will always consider that the ``$p_{i}$" are the ``$y_{i}$" measures and that the model ``$q_{i}$"  is linearly related to the true unknowns ``$x_{i}$" by the relation $q_{i}=\left(Hx\right)_{i}$.\\
In a later section dedicated to the introduction of the sum constraint, it will be shown that, subject to a change of variables, the S.G.M. method can be applied on simplified divergences which have the property of being zero for $p_{i}=q_{i}\ \;\forall i$.\\
In the following sections we will note:
\begin{equation}
	D\left(p\|q\right)=D\left(y\|Hx\right)\equiv D\left(x\right)
\end{equation}

\subsection{Non-negativity constraint.}
For a strictly convex divergence $D\left(x\right)$, we consider the problem:\\
Minimize with respect to ``$x$'': $D\left(x\right)$\\
Subject to the constraint $x\geq0\ \  \equiv\ \ x_{i}\geq0\ \ \forall i$.\\
Let $L\left(x,\lambda\right)$ the Lagrangian associated with this problem, which is expressed in the form of:
\begin{equation}
L\left(x,\lambda\right)=D\left(x\right)-\left\langle \lambda,g\left(x\right)\right\rangle	
\end{equation}
``$\lambda$'' is the vector of Lagrange multipliers whose components are positive or zero ($\lambda_{i}\geq0\ \ \forall i$); $\left\langle \lambda,g\left(x\right)\right\rangle$ is for the scalar product; $g\left(x\right)\equiv g\left(x_{i}\right)\equiv g_{i}\left(x\right)\ \ \forall i$ is a function for expressing constraints; this function must be positive when the constraints are inactive ($x_{i}>$0) and zero when the constraints are active ($x_{i}=$0), moreover, if ``$x^{*}$'' denotes the optimum, the zeros of $g_{i}\left(x^{*}\right)/\left[\nabla g\left(x^{*}\right)\right]_{i}$ must be the same as those of $g_{i}\left(x^{*}\right)$.\\

We propose to build a constrained minimization algorithm based on Karush, Kuhn, Tucker (KKT) conditions \cite{boyd2004}, \cite{bertsekas1995},\cite{culioli2012introduction} which must be satisfied at the optimum of the problem, i.e. at the solution $(x^{*},\lambda^{*})$.\\
These conditions are written:
\begin{equation}
	\nabla_{x}L\left(x^{*},\lambda^{*}\right)=0\ \Rightarrow\ \lambda_{i}^{*}\left[\nabla g\left(x^{*}\right)\right]_{i}=\left[\nabla D\left(x^{*}\right)\right]_{i}\ \Rightarrow\ \lambda_{i}^{*}=\frac{\left[\nabla D\left(x^{*}\right)\right]_{i}}{\left[\nabla g\left(x^{*}\right)\right]_{i}}
\end{equation}
\begin{equation}
	g\left(x^{*}\right)\geq0\ \ \Rightarrow\ \ g_{i}\left(x^{*}\right)\geq0\ \ \forall i
	\label {eq.cg}
\end{equation}
\begin{equation}
	\lambda^{*}\geq 0\ \ \Rightarrow\ \ \lambda_{i}^{*}\geq 0\ \ \forall i
\label {eq.cla}	
\end{equation}
\begin{equation}
	\left\langle \lambda^{*} g\left(x^{*}\right)\right\rangle = 0
\end{equation}
Note that this last condition, which is a scalar product, results in fact, taking into account (\ref {eq.cg}) and (\ref {eq.cla}), in a set of complementary conditions:
\begin{equation}
	\lambda_{i}^{*} g_{i}\left(x^{*}\right)=0\ \ \Rightarrow\ \ \frac{\left[\nabla D\left(x^{*}\right)\right]_{i}}{\left[\nabla g\left(x^{*}\right)\right]_{i}}g_{i}\left(x^{*}\right)=0 \ \ \forall i
\end{equation}
Taking into account the properties of $g\left(x\right)$, this last condition can be written as follows:
\begin{equation}
	\left[\nabla D\left(x^{*}\right)\right]_{i}g_{i}\left(x^{*}\right)=0 \ \ \forall i
	\label{eq.base} 
\end{equation}
To express it in more explicit way, we will have:
\begin{enumerate}
	\item Either:
\begin{equation}
g_{i}\left(x^{*}\right)>0\ \ \Rightarrow \ \ The\ constraint\ is\ inactive
\end{equation}
then (\ref{eq.base}) is fulfilled if:
\begin{equation}
\lambda_{i}^{*}=0\ \ \Leftrightarrow\ \ \left[\nabla D\left(x^{*}\right)\right]_{i}=0
\label{eq.kktlz}
\end{equation}
then ``$x^{*}$'' is the unconstrained solution of the problem.\\
\item Or:
\begin{equation}
\lambda_{i}^{*}>0\ \ \Leftrightarrow\ \ \left[\nabla D\left(x^{*}\right)\right]_{i}>0
\label{eq.kktlp}
\end{equation}
then (\ref{eq.base}) is fulfilled if:
\begin{equation}
g_{i}\left(x^{*}\right)=0\ \ \Rightarrow \ \ The\ constraint\ is\ active
\end{equation}

\item or, in an extreme case, we have simultaneously:
\begin{equation}
\lambda_{i}^{*}=0\ \ \Rightarrow\ \ \left[\nabla D\left(x^{*}\right)\right]_{i}=0\ \ et\ \ g_{i}\left(x^{*}\right)=0	
\end{equation}

\end{enumerate}
Then, (\ref{eq.base}) is ``a fortiori" fulfilled, this is the case where the minimum corresponding to $\nabla D\left(x\right)=0$ is located on the constraint.\\

\begin{figure}[ht!]
\centering
\includegraphics[height=\linewidth, angle=-90]{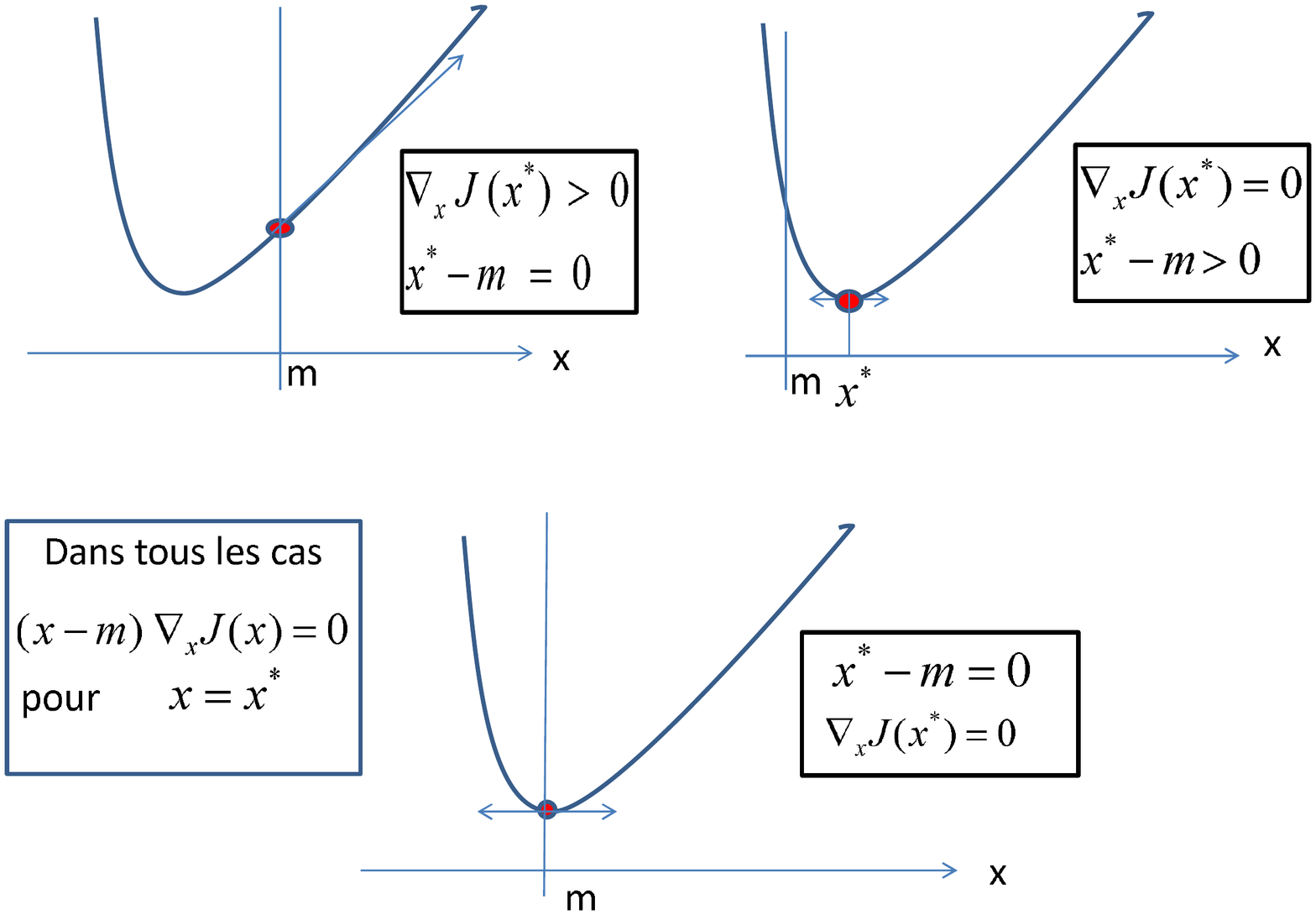}
\caption{Karush, Kuhn, Tucker conditions for a constraint of the form $x\geq m$}
\label{fig:KKT}
\end{figure}

These situations are illustrated (in the one-dimensional case) in the figure (\ref{fig:KKT}) for a constraint $x_{i}\geq m\ \ \forall i$, by expressing the constraint using the function $g\left(x\right)=x-m$.

The form of the function $g\left(x\right)$ used to express the constraints is of fundamental importance since it directly influences the final form and the behavior of the algorithm.\\
For a non-negativity constraint, in the simplest case, one chooses $g\left(x\right)=x$ that is $g_{i}\left(x\right)=x_{i}\ \ \forall i$.

\subsubsection{Algorithm.}
To solve (\ref{eq.base}), we rely on the fact that $-\nabla D\left(x\right)$ is a descent direction for the unconstrained problem, and that if $M$ is a positive definite matrix, then any vector of the form $M\left[-\nabla D\left(x\right)\right]$ is also a descent direction for this problem; this is especially true if $M$ is a diagonal matrix with positive terms.\\
We can also similarly rely on the successive substitution method \cite{hildebrand1987} to write an iterative algorithm in the form:
\begin{equation}
x_{i}^{k+1}=x_{i}^{k}+\alpha_{i}^{k}x_{i}^{k}f_{i}\left(x^{k}\right)\left[-\nabla D\left(x^{k}\right)\right]_{i}
\label{eq.algofond}	
\end{equation}
In this expression $\alpha_{i}^{k}>0$ is the descent step size which will ensure both the fulfillment of the constraint and the convergence of the algorithm (in the sense of contracting applications, we can say that the proper choice of the step $\alpha_{i}^{k}$ ensures that the algorithm (\ref{eq.algofond}) of the type $x^{k+1}=F\left(x^{k}\right)$ is a contraction).\\
The function $f_{i}\left(x\right)$ is a function having positive values when we are in the constraint domain; it depends on the forms of $D\left(x\right)$ and $g\left(x\right)$ as well as on the wished properties for the algorithm, in particular if we want to obtain a purely multiplicative form.\\
We now analyze the properties of the algorithm, and in particular, whether or not the solution fulfills the KKT conditions.\\
The non-negativity constraint imposes restrictions on the descent step size. Indeed, to remain within a general framework, we consider that even if $D\left(x\right)$ is not defined everywhere, its definition domain is entirely contained within the constraint domain.\\

\textbf{ Therefore, the condition $x_{i}\geq 0 \ \forall i$ must be imposed first.}\\

 Thus, we define an  interior points method (all successive estimates will fulfill the constraints).\\
So, as a first step, we need to ensure that at each iteration we have:
\begin{equation}
 x_{i}^{k}\geq 0 \ \ \Rightarrow\ \  x_{i}^{k+1}\geq 0\ \ \forall i	
\end{equation}
Which leads, from (\ref{eq.algofond}), to the condition:
\begin{equation}
	1+\alpha_{i}^{k}f_{i}\left(x^{k}\right)\left[-\nabla D\left(x^{k}\right)\right]_{i}\geq 0\ \ \forall i
	\label{eq.condpos}
\end{equation}
For $\left[-\nabla D\left(x^{k}\right)\right]_{i}> 0$, the condition (\ref{eq.condpos}) is always satisfied, in this case, for the relevant ``$i$" components, there are no restrictions on the descent stepsize as far as the fulfillment of the constraints is concerned.\\
On the other hand, for $\left[-\nabla D\left(x^{k}\right)\right]_{i}< 0$ the condition (\ref{eq.condpos}) leads for the corresponding ``$i$" components to:
\begin{equation}
	\alpha_{i}^{k}\leq\frac{1}{f_{i}\left(x^{k}\right)\left[\nabla D\left(x^{k}\right)\right]_{i}}
	\label{eq.pasmaxi}
\end{equation}
This results in a set of descent step size values that guarantee the non-negativity of each component separately.\\
At iteration ``$k$", the maximum step size to ensure the non-negativity of all the components of the solution will therefore be :
\begin{equation}
\alpha_{M}^{k}=\min_{i}\frac{1}{f_{i}\left(x^{k}\right)\left[\nabla D\left(x^{k}\right)\right]_{i}}\ \ \forall i/\ \ \left[\nabla D\left(x^{k}\right)\right]_{i}> 0,\ \ x_{i}>0
\label{eq.pasmaxM}	
\end{equation}
This maximum step size value being obtained, the value of the descent step size independent of the components, ensuring the convergence of the algorithm, must be calculated by a one-dimensional minimisation method (possibly a simplified method of the Armijo type \cite{armijo1966}), within the interval $\left[0,\alpha_{M}^{k}\right]$ in the direction:
\begin{equation}
	d_{k}=diag\left[f_{i}\left(x^{k}\right)\right]\ diag\left[x_{i}^{k}\right]\left[-\nabla D\left(x^{k}\right)\right]
	\label{eq.dirde}
\end{equation}
This direction, which is no longer that of the opposite of the gradient, remains nevertheless a descent direction for $D\left(x\right)$.\\
The general algorithm is therefore written:
\begin{equation}
	x^{k+1}=	x^{k}+\alpha^{k}\;d^{k}
	\label{eq.algoref}
\end{equation}
We can now check that the conditions of KKT (\ref{eq.kktlz}) and (\ref{eq.kktlp}) are fulfilled at optimum, indeed: 
\begin{itemize}
	\item if $x_{i}^{*}>0$, according to (\ref{eq.algofond}), we have necessarily $\left[\nabla D\left(x^{*}\right)\right]_{i}=0$, that is $\lambda^{*}_{i}=0$; the optimum is the unconstrained minimum, the constraint is inactive. 
	\item if $x_{i}\approx 0$ ($x_{i}=\epsilon>0$), according to (\ref{eq.algofond}), to attain the constraint it will be necessary to have $\left[\nabla D\left(x\right)\right]_{i}>0$ that is $\lambda^{*}_{i}>0$, then, the optimum is located on the constraint which is active.
	\item of course, one can have simultaneously $x_{i}^{*}=0$ and $\left[\nabla D\left(x^{*}\right)\right]_{i}=0$; then, the optimum is located on the constraint which is active and it also corresponds to the unconstrained minimum. 
\end{itemize}

\subsection{Standard Algorithms - Non multiplicative}
The algorithms of this form that we will be led to use, strictly correspond to the previous method, in which we take:
\begin{equation}
	f_{i}\left(x\right)=1 \ \ \forall i
\end{equation}

\subsection{Multiplicative Algorithms.}
In the algorithm defined by (\ref{eq.dirde}) and (\ref{eq.algoref}), a particular choice of the function $f_{i}\left(x\right)$ will allow to obtain purely multiplicative algorithms.\\
To characterize it for this specific purpose, we observe that if at the convergence of the algorithm we have $\nabla D\left(x^{*}\right)=$0, then, necessarily the gradient can be written as the difference of 2 terms of the same sign (we can even say positive without taking anything away from the generality).\\
We will therefore write:
\begin{equation}
-\nabla D\left(x^{k}\right)=U\left(x^{k}\right)-V\left(x^{k}\right)\ \ \ ;\ \ \ U\left(x^{k}\right)>0\ \ \;\ \  V\left(x^{k}\right)>0	
\end{equation}
If we define now:
\begin{equation}
	f_{i}\left(x^{k}\right)=\frac{1}{\left[V\left(x^{k}\right)\right]_{i}}
	\label{eq.fibase}
\end{equation}
The algorithm (\ref{eq.algofond}), (\ref{eq.dirde}), (\ref{eq.algoref}) can be written:
\begin{equation}
x_{i}^{k+1}=x_{i}^{k}+\alpha_{i}^{k}x_{i}^{k}\frac{1}{\left[V\left(x^{k}\right)\right]_{i}}\left[U\left(x^{k}\right)-V\left(x^{k}\right)\right]_{i}	
\end{equation}
Or also:
\begin{equation}
x_{i}^{k+1}=x_{i}^{k}+\alpha_{i}^{k}x_{i}^{k}\left\{\frac{\left[U\left(x^{k}\right)\right]_{i}}{\left[V\left(x^{k}\right)\right]_{i}}-1\right\}
\label{eq.algomrelax}	
\end{equation}
We can say, even if it is somewhat trivial, that we have substituted the convergence condition $\left[U\left(x^{k}\right)\right]_{i}-\left[V\left(x^{k}\right)\right]_{i}=0$ by the condition $\frac{\left[U\left(x^{k}\right)\right]_{i}}{\left[V\left(x^{k}\right)\right]_{i}}=1$.\\
According to (\ref{eq.pasmaxi}), we will have:
\begin{equation}
	\alpha_{i}^{k}\leq\frac{1}{1-\frac{\left[U\left(x^{k}\right)\right]_{i}}{\left[V\left(x^{k}\right)\right]_{i}}}\ \ ; if \ \left[-\nabla D\left(x^{k}\right)\right]_{i}< 0\ \ \Leftrightarrow\ \ \left[U\left(x^{k}\right)\right]_{i}<\left[V\left(x^{k}\right)\right]_{i}
\end{equation}
The values for $\alpha_{i}^{k}$ will all necessarily be greater than 1, so does $\alpha_{M}^{k}$, and we have to remember that these are limit values ensuring that the solution is positive or zero.\\
If we choose a descent step size $\alpha^{k}=1\ \forall k$, the algorithm takes the form "Purely Multiplicative (P.M.)" which is written:
\begin{equation}
x_{i}^{k+1}=x_{i}^{k}\left\{\frac{\left[U\left(x^{k}\right)\right]_{i}}{\left[V\left(x^{k}\right)\right]_{i}}\right\}
\label{eq.algomult}	
\end{equation}
On the contrary, the algorithm (\ref{eq.algomrelax}) can be interpreted as a relaxed version of (\ref{eq.algomult}).\\ 
It is quite obvious that in (\ref{eq.algomult}) the successive estimates will be positive if the initial estimate is positive, and that moreover, if at a given iteration, a component of the solution reaches the value zero (i.e. is on the constraint), this component will remain at zero during the following iterations.\\
However, it should be noted that in general, nothing guarantees the convergence of a purely multiplicative algorithm of type (\ref{eq.algomult}), only a specific analysis of the divergence $D\left(x\right)$  considered will make it possible to ensure it.\\
Such an analysis has been performed for 2 particular divergences, the root mean square (RMS) divergence that led to the ISRA algorithm \cite{daube1986}, \cite{de1987} and the Kullback-Leibler divergence that led to the E.M. algorithms \cite{dempster1977},\cite{lange1984}  or Richardson-Lucy \cite{richardson1972}, \cite{lucy1974}.\\
Moreover, we must note that even if the algorithm (\ref{eq.algomult}) converges, its convergence speed is fixed, because there is no longer any adjustment parameter available. Indeed, for a fixed descent direction, this role is played by the descent step size.\\
The low convergence speed is one of the drawbacks (frequently mentioned in the literature), of purely multiplicative algorithms of the type (\ref{eq.algomult}).

\subsection{Effect of the function expressing the constraints.}
We will show that by modifying the function expressing the constraints $g_{i}\left(x\right)$, we can achieve a change in the speed of convergence of the algorithms exhibited in the preceding sections (as long as they converge).\\
For this aim we choose for example $g_{i}\left(x\right)=x_{i}^{1/n}$, then, the relation (\ref{eq.base}) is written:
\begin{equation}
	\left[x_{i}^{k}\right]^{1/n}\left\{\left[U\left(x^{k}\right)\right]_{i}-\left[V\left(x^{k}\right)\right]_{i}\right\}=0
	\label{eq.mult}
\end{equation}
or also:	
\begin{equation}
	x_{i}^{k}\left\{\left[U\left(x^{k}\right)\right]_{i}^{n}-\left[V\left(x^{k}\right)\right]_{i}^{n}\right\}=0
	\label{eq.multbis}
\end{equation}

\subsubsection{Multiplicative Algorithms.}
The equation (\ref{eq.multbis}) can be written as follows:
\begin{equation}	\frac{x_{i}^{k}}{\left[V\left(x^{k}\right)\right]_{i}^{n}}\left\{\left[U\left(x^{k}\right)\right]_{i}^{n}-\left[V\left(x^{k}\right)\right]_{i}^{n}\right\}=0\ \ \Rightarrow\ \ x_{i}^{k}\left\{\frac{\left[U\left(x^{k}\right)\right]_{i}^{n}}{\left[V\left(x^{k}\right)\right]_{i}^{n}}-1\right\}=0
	\label{eq.tercor}	
\end{equation}
Then, the algorithm (\ref{eq.algofond})(\ref{eq.algoref}) can be written as:
\begin{equation}	x_{i}^{k+1}=x_{i}^{k}+\alpha^{k}x_{i}^{k}\left\{\underbrace{\frac{\left[U\left(x^{k}\right)\right]_{i}^{n}}{\left[V\left(x^{k}\right)\right]_{i}^{n}}-1}_{A}\right\}
	\label{eq.algoacc}	
\end{equation}
To illustrate the change in the behavior of the algorithm (\ref{eq.algoacc}) with respect to the algorithm (\ref{eq.algomrelax}), in connection with the form of $g_{i}\left(x^{k}\right)$, we are looking at the corrective term ``\textbf{A}" which can be developed in the form:
\begin{equation}	\frac{1}{\left[V\left(x^{k}\right)\right]_{i}^{n}}\left\{\sum_{p=0}^{n-1}\left[U\left(x^{k}\right)\right]_{i}^{n-p-1}\left[V\left(x^{k}\right)\right]_{i}^{p}\right\}\left\{\left[U\left(x^{k}\right)\right]_{i}-\left[V\left(x^{k}\right)\right]_{i}\right\}=0	
\end{equation}
Which means that the $f_{i}\left(x^{k}\right)$ factor here takes the form:
\begin{equation}	f_{i}\left(x^{k}\right)=\frac{1}{\left[V\left(x^{k}\right)\right]_{i}^{n}}\left\{\sum_{p=0}^{n-1}\left[U\left(x^{k}\right)\right]_{i}^{n-p-1}\left[V\left(x^{k}\right)\right]_{i}^{p}\right\}
\end{equation}
The expression of $f_{i}\left(x^{k}\right)=\frac{1}{\left[V\left(x^{k}\right)\right]_{i}}$ (\ref{eq.fibase}) here is multiplied by:
\begin{equation}
B=\frac{1}{\left[V\left(x^{k}\right)\right]_{i}^{n-1}}\left\{\sum_{p=0}^{n-1}\left[U\left(x^{k}\right)\right]_{i}^{n-p-1}\left[V\left(x^{k}\right)\right]_{i}^{p}\right\}	
\end{equation}
Apart from the obvious change in the direction of descent, we observe that when we are close to convergence, i.e. when $\left[U\left(x^{k}\right)\right]_{i}\rightarrow \left[V\left(x^{k}\right)\right]_{i}$ this change results in $B\rightarrow n$; so the descent step size is progressively multiplied by a factor that tends towards ``$n$"; so we can expect a gain of a factor of about ``$n$" on the convergence speed.\\
In a relaxed algorithm, this effect is counterbalanced by the computation of the descent step size in the modified direction, so there's no risk of divergence for the algorithm, but then the gain in speed is partially suppressed.\\
On the other hand, if we consider the algorithm written in the form (\ref{eq.algoacc}) and we choose a descent step equal to 1, we will obtain a purely multiplicative algorithm which is written as follows:
\begin{equation}
	x_{i}^{k+1}=x_{i}^{k}\frac{\left[U\left(x^{k}\right)\right]_{i}^{n}}{\left[V\left(x^{k}\right)\right]_{i}^{n}}
	\label{eq.algoaccmult}	
\end{equation}
This is the form frequently proposed in the litterature \cite{shepp1982}, \cite{zaccheo1996}, to achieve a convergence speed gain of a ``$n$" factor, but again, as is often the case in fixed step size algorithms, nothing guarantees convergence, and, to put it roughly: if you want to go too fast, you often go to the wall....

\subsubsection{Standard Algorithm - Non multiplicative.}
In this case, the equation (\ref{eq.multbis}) can be written:
\begin{equation}
x_{i}^{k}\left\{\sum_{p=0}^{n-1}\left[U\left(x^{k}\right)\right]_{i}^{n-p-1}\left[V\left(x^{k}\right)\right]_{i}^{p}\right\}\left\{\left[U\left(x^{k}\right)\right]_{i}-\left[V\left(x^{k}\right)\right]_{i}\right\}=0		
\end{equation}
The fonction $f_{i}\left(x^{k}\right)=1$ is here replaced by the expression:
\begin{equation}
f_{i}\left(x^{k}\right)=\sum_{p=0}^{n-1}\left[U\left(x^{k}\right)\right]_{i}^{n-p-1}\left[V\left(x^{k}\right)\right]_{i}^{p}	
\end{equation}
It should be noted, however, that with such algorithms, for convergence problems, the calculation of the descent step size is indispensable, which somewhat cancels out the effect of acceleration, if any.

\subsection{Specific case of simplified divergences.}
The application of the SGM method and, in general, of the methods described in the previous section, implies that the divergences considered have a minimum (possibly zero) when $p_{i}=q_{i}\ \forall i$.\\
This being the case, in the expressions of certain divergences, simplifications of the type $\sum_{i}p_{i}=\sum_{i}q_{i}$ have been explicitly introduced, or in the case of comparisons between probability densities, $\sum_{i}p_{i}=\sum_{i}q_{i}=1$.\\
When this type of simplification has been introduced in a divergence, although the divergence considered is zero for $p_{i}=q_{i}\ \forall i$, the gradient with respect to ``$q$" may either never cancel out or cancel out for $p_{i}\neq q_{i}$, as has been shown in many examples.\\
These divergences are thus ``a priori" not usable in the context which is ours; however one can nevertheless try to exploit them in spite of the simplifications carried out, by introducing variables which are explicitly of the same sum, for example normalized variables.\\
 We will notice that this operation, whose objective is to enable the use of the simplified divergences, allows to take into account the sum constraint at the same time as the non-negativity constraint.\\
 This is fully understandable because when one tries to make a divergence invariant with respect to ``$q$" by introducing an invariance factor $K^{*}\left(p,q\right)=\frac{\sum_{j}p_{j}}{\sum_{j}q_{j}}$, one obtains a divergence which has the same expression as the initial divergence, but with normalized variables $\bar{p_{i}}=\frac{p_{i}}{\sum_{j}p_{j}}$ and $\bar{q_{i}}=\frac{q_{i}}{\sum_{j}q_{j}}$, thus, the normalized variables are used ``de facto'', so it is not surprising that the use of such variables allows to take into account the sum constraint; this is a kind of justification for the introduction of normalized variables that may seem somewhat artificial.\\
In order to explicitly account for these normalizations in the simplified $D\left(p\|q\right)$ divergences, we will write them down as $D\left(\bar{p}\|\bar{q}\right)$.\\ 
If we denote the measurements by ``$y_{i}$'' and the linear model by $\left(Hx\right)_{i}$, the variables used in the expressions of the simplified divergences will be:
\begin{equation}
 \bar{p_{i}}=\frac{p_{i}}{\sum_{j}p_{j}}C=\frac{y_{i}}{\sum_{j}y_{j}}C\ \ \ ;\ \ \ \bar{q_{i}}=\frac{q_{i}}{\sum_{j}q_{j}}C=\frac{\left(Hx\right)_{i}}{\sum_{j}\left(Hx\right)_{j}}C \ \ \ ; \ \ \  C>0
 \label{eq.varN}	
\end{equation}
Two cases occur in a classical way (not restrictive): either the data fields are typically probability densities, then $C=1$, or we deal with deconvolution problems such as those encountered in image deconvolution, then $C=\sum_{i}y_{i}$.\\
Of course, the problem of convexity of the divergence thus modified arises; indeed, taking into account the normalization introduced, if the divergence considered $D\left(\bar{p}\|\bar{q}\right)$ is convex with respect to ``$\:\bar{q}\:$'', is it still convex with respect to ``$q$", therefore with respect to ``$x$"?\\
One fact is certain: if $\sum_{j}\left(q\right)_{j}=\sum_{j}\left(Hx\right)_{j}$ is maintained constant during the minimization process, $D\left(\bar{p}\|\bar{q}\right)$ will be convex with respect to ``$x$".\\
We will use the notation:
\begin{equation}
	D\left(\bar{p}\|\bar{q}\right)=\sum_{i}d\left(\bar{p_{i}}\|\bar{q_{i}}\right)=\sum_{i}d_{i}
\end{equation}
Then, we have:
\begin{equation}
	\frac{\partial D\left(\bar{p}\|\bar{q}\right)}{\partial q_{m}}=\sum_{i}\frac{\partial d\left(\bar{p_{i}}\|\bar{q_{i}}\right)}{\partial \bar{q_{i}}}\frac{\partial \bar{q_{i}}}{\partial q_{m}}
\end{equation}
With:
\begin{equation}
\frac{\partial \bar{q_{i}}}{\partial q_{m}}=\frac{C}{\sum_{j}q_{j}}\left(\delta_{i m}-\frac{\bar{q_{i}}}{C}\right)	
\end{equation}
So it comes:
\begin{equation}
\frac{\partial D\left(\bar{p}\|\bar{q}\right)}{\partial q_{m}}=\frac{C}{\sum_{j}q_{j}}\left[\frac{\partial d\left(\bar{p_{m}}\|\bar{q_{m}}\right)}{\partial \bar{q_{m}}}-\sum_{i}\frac{\bar{q_{i}}}{C}\frac{d\left(\bar{p_{i}}\|\bar{q_{i}}\right)}{\partial \bar{q_{i}}}\right]	
\end{equation}
If we now take into account the fact that we are using a linear model, we have:
\begin{equation}
	\frac{\partial D\left(\bar{p}\|\bar{q}\right)}{\partial x_{l}}=\sum_{j}\frac{\partial D\left(\bar{p}\|\bar{q}\right)}{\partial q_{j}}\frac{\partial q_{j}}{\partial x_{l}}=\sum_{j}h_{j l}\frac{\partial D\left(\bar{p}\|\bar{q}\right)}{\partial q_{j}}
\end{equation}
so:
\begin{equation}
	-\frac{\partial D\left(\bar{p}\|\bar{q}\right)}{\partial x_{l}}=\frac{C}{\sum_{j}q_{j}}\left\{\sum_{i}h_{i,l}\left[-\frac{\partial d_{i}\left(\bar{p}\|\bar{q}\right)}{\partial \bar{q_{i}}}\right]-\left(\sum_{j}h_{j,l}\right)\sum_{i}\frac{\bar{q_{i}}}{C}\left[-\frac{\partial d_{i}\left(\bar{p}\|\bar{q}\right)}{\partial \bar{q_{i}}}\right]\right\}	
\end{equation}
Or also:
\begin{equation}
	-\frac{\partial D\left(\bar{p}\|\bar{q}\right)}{\partial x_{l}}=\frac{C}{\sum_{j}\left(Hx\right)_{j}}\left\{\left[H^{T}-\left(\sum_{j}h_{j,l}\right)\ \frac{\bar{Q}}{C}\right]\left[-\frac{\partial D\left(\bar{p}\|\bar{q}\right)}{\partial \bar{q}}\right]\right\}_{l}
\end{equation}
In this expression, $\overline{Q}$ is a matrix in which all the lines are identical and equal to $\left[\bar{{q}_{1}}\ \bar{q_{2}}\ .....\bar{q_{N}}\right]$ with $\sum_{i}\bar{q_{i}}=C$ and $\left[\frac{\partial D}{\partial \bar{q}}\right]$ is the vector of partial derivatives $\left[\frac{\partial d_{i}}{\partial \bar{q_{i}}}\right]$.\\
Another writting may be more appealing, considering Richardson's  \cite{richardson1972} and Lucy's \cite{lucy1974} deconvolution algorithms (for example):
\begin{equation}
	-\frac{\partial D\left(\bar{p}\|\bar{q}\right)}{\partial x_{l}}=C\frac{\sum_{j}h_{j,l}}{\sum_{j}q_{j}}\left\{\left[\overline{H}^{T}- \frac{\overline{Q}}{C}\right]\left[-\frac{\partial D\left(\bar{p}\|\bar{q}\right)}{\partial \bar{q}}\right]\right\}_{l}
\end{equation}
Here, $\overline{H}$ is a matrix whose columns are of sum 1.\\

\textbf{Algorithm.}\\

From the above considerations, we can write an SGM-type algorithm in the form:
\begin{equation}	
x^{k+1}_{l}=x^{k}_{l}+\delta^{k}x^{k}_{l}\left[	-\frac{\partial D\left(\bar{p}\|\bar{q}\right)}{\partial x_{l}}\right]
\end{equation}
that is:
\begin{equation}	
x^{k+1}_{l}=x^{k}_{l}+\delta^{k}x^{k}_{l}C\frac{\sum_{j}h_{j,l}}{\sum_{j}q^{k}_{j}}\left\{\left[\overline{H}^{T}- \frac{Q^{k}}{C}\right]\left[-\frac{\partial D\left(\bar{p}\|\bar{q}\right)}{\partial \bar{q}^{k}}\right]\right\}_{l}
\label{eq.suivant1}
\end{equation}
or also:
\begin{equation}	x^{k+1}_{l}=x^{k}_{l}+\delta^{k}\left\{x^{k}_{l}C\frac{\sum_{j}h_{j,l}}{\sum_{j}q^{k}_{j}}\left[\overline{H}^{T}\left[-\frac{\partial D\left(\bar{p}\|\bar{q}\right)}{\partial \bar{q}^{k}}\right]\right]_{l}- x^{k}_{l}\frac{\sum_{j}h_{j,l}}{\sum_{j}q^{k}_{j}}\left[\overline{Q}^{k}\left[-\frac{\partial D\left(\bar{p}\|\bar{q}\right)}{\partial \bar{q}^{k}}\right]\right]_{l}\right\}
\label{eq.suivant2}
\end{equation}
With this algorithm, we have:
\begin{equation}
	\sum_{l}x^{k+1}_{l}=\sum_{l}x^{k}_{l}
\end{equation}
To show this, we calculate: $\sum_{l}x^{k+1}_{l}$ and we look at what are becoming the terms between the braces of (\ref{eq.suivant2}).\\
 For the first of these terms, we have:
\begin{equation}
	\sum_{l}x^{k}_{l}C\frac{\sum_{j}h_{j,l}}{\sum_{j}q^{k}_{j}}\left[\overline{H}^{T}\left[-\frac{\partial D}{\partial \bar{q}^{k}}\right]\right]_{l}=\frac{C}{\sum_{j}q^{k}_{j}}\sum_{l}q^{k}_{l}\left[-\frac{\partial D}{\partial \bar{q}^{k}}\right]_{l}=\sum_{l}\bar{q}^{k}_{l}\left[-\frac{\partial D}{\partial \bar{q}^{k}}\right]_{l}
\end{equation}
For the second term of the brace, given the properties of the matrix $\overline{Q}$ in which all the lines are identical, the term:
\begin{equation}
\left[\overline{Q}^{k}\left[-\frac{\partial D}{\partial \bar{q}^{k}}\right]\right]_{l}=\left[\sum_{i}\bar{q}^{k}_{i}\left[-\frac{\partial D}{\partial \bar{q}^{k}}\right]_{i}\right]_{l}	
\end{equation}
is independent of ``$l$"; therefore, we can write:
\begin{equation}
\sum_{l}x^{k}_{l}\frac{\sum_{j}h_{j,l}}{\sum_{j}q^{k}_{j}}\left[\overline{Q}^{k}\left[-\frac{\partial D}{\partial \bar{q}^{k}}\right]\right]_{l}=\frac{\sum_{j}\sum_{l}h_{j,l}x^{k}_{l}}{\sum_{j}q^{k}_{j}}\sum_{i}\bar{q}^{k}_{i}\left[-\frac{\partial D}{\partial \bar{q}^{k}}\right]_{i}=\sum_{i}\bar{q}^{k}_{i}\left[-\frac{\partial D}{\partial \bar{q}^{k}}\right]_{i}	
\end{equation}
To summarize, if we are returning to the equation (\ref{eq.suivant2}), the term between braces will give in the sum, a null contribution, and we will have:
\begin{equation}
\sum_{l}x^{k+1}_{l}=\sum_{l}x^{k}_{l}
\label{eq.sommex}	
\end{equation}
In conclusion, the algorithm (\ref{eq.suivant1}),(\ref{eq.suivant2}) ensures that the sum constraint is satisfied. This is a consequence of the introduction of the variables in the form (\ref{eq.varN}).\\
However, it should be noted that the relation (\ref{eq.sommex}) is absolutely not equivalent to:
\begin{equation}
\sum_{l}\left(Hx^{k}\right)_{l}=\sum_{l}\left(Hx^{k+1}\right)_{l}
\label{eq.sommeHx}	
\end{equation}
so that the convexity of the initial divergence from the ``$x$" variable is still not assured unless the $H$ matrix is column-normalized.\\
If this is the case, the relations (\ref{eq.sommex}) and (\ref{eq.sommeHx}) become equivalent; this property will be reused in the following section.

\subsubsection{Purely multiplicative form of the algorithm.}
 Given the properties of $\overline{H}$ and of the matrix $\frac{\overline{Q}}{C}$ (for each line, the sum of the terms is equal to 1), in the relation (\ref{eq.suivant1}), we can make $-\frac{\partial D}{\partial \bar{q}}\geq0$ without changing the expression of the gradient, by adding a constant ($-\min_{i}\left[-\frac{\partial d_{i}}{\partial \bar{q}_{i}}\right]$).\\
We denote $\left[-\frac{\partial D}{\partial \bar{q}}\right]_{d}$ the shifted vector.\\
\begin{equation}
\left[-\frac{\partial D}{\partial \bar{q}}\right]_{d}=-\frac{\partial D}{\partial \bar{q}}-\min_{i}\left[-\frac{\partial d_{i}}{\partial \bar{q}_{i}}\right]+\epsilon	
\end{equation}
We can therefore envisage a decomposition of the gradient into a difference of 2 positive terms.\\ 
From this remark, in the case of the simplified divergences considered here, one can write an algorithm in the form: 
\begin{equation}
	x^{k+1}_{l}=x^{k}_{l}+\delta^{k}x^{k}_{l}\left\{\left[\overline{H}^{T}- \frac{\overline{Q}^{k}}{C}\right]\left[-\frac{\partial D}{\partial \bar{q}^{k}}\right]_{d}\right\}_{l}
\end{equation}
Or also:
\begin{equation}
	x^{k+1}_{l}=x^{k}_{l}+\delta^{k} x^{k}_{l}\left\{\left[\overline{H}^{T}\left(-\frac{\partial D}{\partial \bar{q}^{k}}\right)_{s}\right]_{l}-\sum_{i}\frac{\bar{q}^{k}_{i}}{C}\left(-\frac{\partial d_{i}}{\partial \bar{q}^{k}_{i}}\right)_{d}\right\}
\end{equation}
It has already been noted that the second term of the expression between braces is a constant independent of the component ``$l$".\\
As long as $\sum_{i}\bar{q}^{k}_{i}\left(-\frac{\partial d_{i}}{\partial \bar{q}^{k}_{i}}\right)_{d}\neq 0$, which is generally the case, and with a descent step size $\delta_{l}=1$, we can exhibit a purely multiplicative form:
\begin{equation}
	x^{k+1}_{l}=x^{k}_{l}\frac{\left[\overline{H}^{T}\left[-\frac{\partial D}{\partial \bar{q}^{k}}\right]_{d}\right]_{l}}{\sum_{i}\frac{\bar{q}^{k}_{i}}{C}\left(-\frac{\partial d_{i}}{\partial \bar{q}^{k}_{i}}\right)_{d}}=x^{k}_{l}\frac{C\left[H^{T}\left[-\frac{\partial D}{\partial \bar{q}^{k}}\right]_{d}\right]_{l}}{\sum_{j}h_{j,l}\ \ \sum_{i}\bar{q}^{k}_{i}\left(-\frac{\partial d_{i}}{\partial \bar{q}^{k}_{i}}\right)_{d}}
\end{equation}

At this point, if we calculate $\sum_{l} x^{k+1}_{l}$, we obtain nothing but a weighted sum that is difficult to interpret, indeed, $\sum_{j}h_{j,l}$ depends on the ``$l$" component, on the other hand, if we make the assumption of a matrix $H$ normalized to 1 in columns, i.e. $\sum_{j}h_{j,l}=1\;\forall l$, (which is the case in classical deconvolution problems),  we arrive at something simpler which is written as follows:
\begin{equation}
	\sum_{l}x^{k+1}_{l}=\sum_{l}\left(Hx^{k}\right)_{l}=\sum_{l}x^{k}_{l}
\end{equation}
Consequently, by introducing the variables in the form (\ref{eq.varN}), we also obtained a purely multiplicative algorithm making it possible to fulfill the sum constraint, at least in the case where the $H$ matrix is normalized in columns, i.e. in the case of the usual convolution with a positive kernel with integral equal to 1.

\section{Introduction of the sum constraint.}
In this section, we will rely on the developments of the previous section and introduce an additional constraint which is the sum constraint on the unknowns ``$x_{i}$". Several methods are possible to introduce this constraint; their use depends essentially on the properties of the divergences considered.

\subsection{Simplified divergences.}
The case of these divergences has been addressed in the previous paragraph.\\
Indeed, the introduction of an invariance factor equal to
\begin{equation}
	K^{*}\left(p,q\right)=\frac{\sum_{j}p_{j}}{\sum_{j}q_{j}}
\end{equation}
allows to show simplified divergences in which the variables are explicitly of the same sum.
Moreover, we observe that (with the exception of a multiplicative factor which depends only on $\sum_{j}p_{j}$ and which can be omitted), these divergences are not only invariant with respect to the variable ``$q$", but also with respect to the variable ``$p$".\\
 This property will make it possible to use this type of divergence for regularization problems in which the both arguments of the divergence depends on the true variable ``$x$".

\subsection{Non-simplified divergences.}
This paragraph concerns the Csiszär divergences built on standard convex functions as well as the Bregman divergences and the Jensen divergences.\\
With these divergences, we will have $p_{i}=y_{i}$ and $q_{i}=\left(Hx\right)_{i}$.\\
The algorithms for minimizing these divergences under non-negativity constraint are derived from KKT conditions and can result in multiplicative algorithms.\\
The sum constraint $\sum_{i}x_{i}=C,\ \ C>0$ can be taken into account by introducing the variable change $x_{i}=\frac{u_{i}}{\sum_{j}u_{j}}C$.\\
The problem that arises then is that of the convexity of the divergence to be minimized with respect to the new variable ``$u$".\\
 An answer to this problem can be the following: the objective function (i.e. the considered divergence) being supposed to be convex with respect to ``$q$", therefore with respect to ``$x$" since $q=Hx$, it will be convex with respect to ``$u$" if $\sum_{j}u_{j}$ is constant during the minimization process, i.e. in the course of the iterations.\\
We must therefore write an algorithm of minimization with respect to ``$u$" of the divergence $D\left(p\|q\right)=D\left(y\|Hx\right)$ with: $x_{i}=\frac{u_{i}}{\sum_{j}u_{j}}C$; this algorithm must be such that during the iterations $\sum_{j}u^{k}_{j}=Cst$.\\
Using the notation:
\begin{equation}
	D\left(p\|q\right)=\sum_{i}d\left(p_{i}\|q_{i}\right)
\end{equation}
we will have:
\begin{equation}
	\frac{\partial D}{\partial x_{m}}=\sum_{i}\frac{\partial d\left(p_{i}\|q_{i}\right)}{\partial q_{i}}\frac{\partial q_{i}}{\partial x_{m}}=\left[H^{T}\frac{\partial D\left(p\|q\right)}{\partial q}\right]_{m}
\end{equation}
Then the gradient of $D$ with respect to ``$u$" is written as follows:
\begin{equation}
\left[\frac{\partial D}{\partial u}\right]_{l}=\sum_{m}\frac{\partial D}{\partial x_{m}}\frac{\partial x_{m}}{\partial u_{l}}	
\end{equation}
With:
\begin{equation}
\frac{\partial x_{m}}{\partial u_{l}}=\frac{1}{\sum_{j}u_{j}}\left(C\;\delta_{l\:m}-x_{l}\right)	
\end{equation}
we have:
\begin{equation}
\left[\frac{\partial D}{\partial u}\right]_{l}=\frac{1}{\sum_{j}u_{j}}\left[C\;\frac{\partial D}{\partial x_{l}}-\sum_{m}x_{m}\frac{\partial D}{\partial x_{m}}\right]	
\end{equation}
and we have, immediately:
\begin{equation}
\sum _{l}u_{l}\left[-\frac{\partial D}{\partial u}\right]_{l}=0
\label{eq:pcle}	
\end{equation}
\textbf{This relationship is the key point of the affair}.\\

Indeed, the iterative minimization algorithm with respect to ``$u$", founded on the KKT conditions, ensuring non-negativity is written in the form (\ref{eq.algofond}):
\begin{equation}
u^{k+1}_{l}=u^{k}_{l}+\delta^{k}u^{k}_{l}\left[-\frac{\partial D}{\partial u^{k}}\right]_{l}
\label{eq.algoU}	
\end{equation}
Then, with (\ref{eq:pcle}) we will have:
\begin{equation}
\sum_{l}u^{k+1}_{l}=\sum_{l}u^{k}_{l}	
\end{equation}
In the course of the iterations, the trajectory of the solutions will thus be contained in the convexity domain of the objective function; one can thus return to the initial variable ``$x$" by dividing the 2 members of (\ref{eq.algoU}) by $\sum_{l}u^{k+1}_{l}=\sum_{l}u^{k}_{l}$, and we will have:
\begin{equation}
x^{k+1}_{l}=x^{k}_{l}+\delta^{k}\frac{1}{\sum_{j}u_{j}^{k}}x^{k}_{l}\left[C\;\left(-\frac{\partial D}{\partial x^{k}}\right)_{l}-\sum_{m}x^{k}_{m}\left(-\frac{\partial D}{\partial x^{k}}\right)_{m}\right]	
\end{equation}
Or also:
\begin{equation}
x^{k+1}_{l}=x^{k}_{l}+\zeta^{k}x^{k}_{l}\left[C\;\left(-\frac{\partial D}{\partial x^{k}}\right)_{l}-\sum_{m}x^{k}_{m}\left(-\frac{\partial D}{\partial x^{k}}\right)_{m}\right]
\label{eq.algocontraint}	
\end{equation}
Note that in this algorithm the second term in square brackets is a constant that is independent of the component considered; moreover, the quantity:
\begin{equation}
\zeta^{k}=\delta^{k}\frac{1}{\sum_{j}u_{j}^{k}}	
\end{equation}
represents globally a descent step size which is calculated at each iteration by a one-dimensional minimization procedure as indicated for the S.G.M. method.\\ 
Taking into account the relation $\sum_{i}x_{i}=C$, we will thus have at each iteration, independently of the descent step size:
\begin{equation}
\sum_{l}x^{k+1}_{l}=\sum_{l}x^{k}_{l}	
\end{equation}
 To fulfill the sum constraint at each iteration, it is thus only necessary to set an initial estimate ``$x^{0}$" such that $\sum_{i}x^{0}_{i}=C$.

\subsubsection{* Multiplicative form.}
From the algorithm, a multiplicative form can be derived.\\
 To do this, we shift the components of the opposite of the gradient to make them all positive: 
\begin{equation}
	\left[-\frac{\partial D}{\partial x^{k}_{j}}\right]_{d}=\left[-\frac{\partial D}{\partial x^{k}_{j}}\right]-\min_{l}\left[-\frac{\partial D}{\partial x^{k}_{l}}\right]+\epsilon
\end{equation}
The introduction of this offset doesn't change the algorithm.\\
We can then propose an algorithm that will be the basis of the purely multiplicative form, which is written:
\begin{equation}
x^{k+1}_{l}=x^{k}_{l}+\zeta^{k}x^{k}_{l}\left[\frac{C\;\left(-\frac{\partial D}{\partial x^{k}_{l}}\right)_{d}}{\sum_{m}x^{k}_{m}\left(-\frac{\partial D}{\partial x^{k}_{m}}\right)_{d}}-1\right]
\label{eq.algocontraintmult}	
\end{equation}
With this algorithm, the sum-holding property is preserved over iterations, regardless of the descent step size.\\
The purely multiplicative form will, as always, be obtained by using a descent step size $\zeta^{k}=1\ \forall k$, but then convergence is no longer guaranteed. 
\subsection{Divergences invariant by change of scale.}
For these divergences, we rely on the method developed in the previous sections to build a minimization algorithm that takes into account the non-negativity constraint.\\
As far as the summation constraint is concerned, we have two methods at our disposition:\\
* Either we normalize at each iteration which does not lead to any variation of the objective function taking into account the invariance property.\\
* Either we rely on a specific property of invariant divergences we previously established.\\

\textbf{* Property:} Considering a divergence $D\left(p\|q\right)$, which we render invariant by scale change using any invariance factor, whether it corresponds to it (``\textit{nominal}" invariance factor) or not, we obtain a divergence $DI\left(p\|q\right)$, for which we always have:
\begin{equation}
	\sum_{l}q_{l}\frac{\partial DI\left(p\|q\right)}{\partial q_{l}}=0
	\label{eq.prid}
\end{equation}
\textbf{This property is analogous to(\ref{eq:pcle}).}\\

\textbf{Example: Mean square deviation}.\\

The basic divergence is:
\begin{equation}
	MC=\sum_{i}\left(p_{i}-q_{i}\right)^{2}
\end{equation}
By introducing the invariance factor $K$, we have:
\begin{equation}
	MCI=\sum_{i}\left(p_{i}-Kq_{i}\right)^{2}
\end{equation}
The gradient with respect to ``$q_{j}$" is written, all simplifications made:
\begin{equation}
\frac{\partial MCI}{\partial q_{j}}=-2K\left(p_{j}-Kq_{j}\right)-2\frac{\partial K}{\partial q_{j}}\sum_{i}q_{i}\left(p_{i}-Kq_{i}\right)	
\label{eq.gradMCI}
\end{equation}
Then, we have:
\begin{equation}
	\sum_{j}q_{j}\frac{\partial MCI}{\partial q_{j}}=-2\left[\sum_{j}q_{j}\left(p_{j}-Kq_{j}\right)\right]\left(K+\sum_{j}q_{j}\frac{\partial K}{\partial q_{j}}\right)
\end{equation}
Then, if we use the nominal invariance factor, i.e.:
\begin{equation}
	K_{0}\left(p,q\right)=\frac{\sum_{l}p_{l}q_{l}}{\sum_{l}q^{2}_{l}}
\end{equation}
which is a solution to the differential equation:
\begin{equation}
K+\sum_{j}q_{j}\frac{\partial K}{\partial q_{j}}=0
\label{eq.diffK}	
\end{equation}
the expression of the gradient (\ref{eq.gradMCI}) becomes simpler because the second term is zero.\\
We still have the property (\ref{eq.prid}).\\
Similarly, if we use another expression of $K$, for example: 
\begin{equation}
	K^{*}\left(p,q\right)=\frac{\sum_{l}p_{l}}{\sum_{l}q_{l}}
\end{equation}
which is another solution of (\ref{eq.diffK}), we also obtain an invariant divergence, the expression of the gradient (\ref{eq.gradMCI}) does not simplify, but we still have the property (\ref{eq.prid}).\\

\subsubsection{Some comments on the use of $K^{*}\left(p,q\right)$.}
We recall that when a divergence $D(p\|q)$ is made invariant by using the invariance factor $K^{*}(p,q)$, the divergence obtained is similar to the initial divergence in which the variables ``$p$" and ``$q$" are replaced by normalized variables $\bar{p_{i}}=\frac{p_{i}}{\sum_{j}p_{j}}$ and $\bar{q_{i}}=\frac{q_{i}}{\sum_{j}q_{j}}$, possibly with a multiplicative scalar factor which depends only on the measures ``$p_{i}$" and which can be omitted; simplifications can then appear.\\
Therefore, the invariant divergence $DI(p\|q)$ obtained with this particular invariance factor is denoted as $D(\bar{p}\|\bar{q})$ with:
\begin{equation}
D(\bar{p}\|\bar{q})=\sum_{i}d(\bar{p}_{i}\|\bar{q}_{i})	
\end{equation}
 and the gradient with respect to ``$q$" is written as follows:
\begin{equation}
\frac{\partial DI(p\|q)}{\partial q_{l}}=\frac{\partial D(\bar{p}\|\bar{q})}{\partial q_{l}}	
\end{equation}
It is calculated as follows:
\begin{equation}
\frac{\partial D(\bar{p}\|\bar{q})}{\partial q_{l}}=\sum_{i}\frac{\partial d(\bar{p_{i}}\|\bar{q_{i}})}{\partial q_{l}}=\sum_{i}\frac{\partial d(\bar{p_{i}}\|\bar{q_{i}})}{\partial \bar{q}_{i}}\frac{\partial \bar{q}_{i}}{\partial q_{l}}
\end{equation}
With:
\begin{equation}
	\frac{\partial \bar{q}_{i}}{\partial q_{l}}=\frac{\delta_{il}}{\sum_{j}q_{j}}-\frac{\bar{q}_{i}}{\sum_{j}q_{j}}
\end{equation}
So:
\begin{equation}
\frac{\partial D(\bar{p}\|\bar{q})}{\partial q_{l}}=\frac{1}{\sum_{j}q_{j}}\left[\frac{\partial d(\bar{p_{l}}\|\bar{q_{l}})}{\partial \bar{q}_{l}}-	\sum_{i}\bar{q}_{i}\frac{\partial d(\bar{p_{i}}\|\bar{q_{i}})}{\partial \bar{q}_{i}}\right]
\end{equation}
We can observe that the classical relationship for invariant divergences:
\begin{equation}
	\sum_{l}q_{l}\frac{\partial DI(p\|q)}{\partial q_{l}}=\sum_{l}q_{l}\frac{\partial D(\bar{p}\|\bar{q})}{\partial q_{l}}=0
\end{equation}
is fulfilled.\\

\textbf{Algorithm:}\\
In the context of the use of a linear model $q=Hx$, the calculation of the gradient with respect to ``$x$", true unknowns of the problem, follows by the relation:
\begin{equation}
\frac{\partial DI(p\|q)}{\partial x_{l}}=\left[H^{T}\frac{\partial DI(p\|q)}{\partial q}\right]_{l}
\label{eq.graddixl}	
\end{equation}
Given this relationship, the algorithms for minimizing scale invariant divergences under non-negativity and sum-of-unknowns constraints can be written directly according to the method indicated above, in the form of:
\begin{equation}
	x^{k+1}_{l}=x^{k}_{l}+\alpha^{k}x^{k}_{l}\left[-H^{T}\frac{\partial DI(p\|q^{k})}{\partial q^{k}}\right]_{l}
	\label{eq.algo flux const}
\end{equation}
Non-negativity and convergence are ensured by a computation of the descent step size, but whatever its value, we will always have:
\begin{equation}
	\sum_{l}x^{k+1}_{l}=\sum_{l}x^{k}_{l}
\end{equation}
Indeed:
\begin{equation}
	\sum_{l}x^{k+1}_{l}=\sum_{l}x^{k}_{l}+\alpha^{k}\sum_{l}\left\{x^{k}_{l}\left[H^{T}\left(-\frac{\partial DI\left(p\|q\right)}{\partial q^{k}}\right)\right]_{l}\right\}
\end{equation}
Which can be written:
\begin{equation}
	\sum_{l}x^{k+1}_{l}=\sum_{l}x^{k}_{l}+\alpha^{k}\sum_{l}\left\{\left(Hx^{k}\right)_{l}\left(-\frac{\partial DI\left(p\|q\right)}{\partial q^{k}}\right)_{l}\right\}
\end{equation}
And then:
\begin{equation}
	\sum_{l}x^{k+1}_{l}=\sum_{l}x^{k}_{l}+\alpha^{k}\sum_{l}\left\{q^{k}_{l}\left(-\frac{\partial DI\left(p\|q\right)}{\partial q^{k}}\right)_{l}\right\}
\end{equation}
Taking into account the property (\ref{eq.prid}) we have:
\begin{equation}
	\sum_{l}x^{k+1}_{l}=\sum_{l}x^{k}_{l}
\end{equation}
It is therefore enough to fix the sum of the components of the initial estimate to maintain this sum during the iterative process.\\
 Moreover, as for all the algorithms proposed in the preceding sections, \textbf{the maximum step size $\alpha^{k}_{M}$ ensuring non-negativity is calculated first}, then the descent step size ensuring the convergence of the algorithm is calculated at each iteration by a one-dimensional search method in the interval $[0,\alpha^{k}_{M}]$.\\
On the other hand, it should be noted that a purely multiplicative algorithm deduced from the SGM method would not spontaneously possess this property of maintaining the sum.\\
However, nothing prevents us from proposing such a purely multiplicative form. Indeed, for all the divergences considered, we can observe (although there is no formal demonstration of this property) that the opposite of the gradient is always written as the difference of 2 positive terms, i.e.:
\begin{equation}
\left[-H^{T}\frac{\partial DI(p\|q^{k})}{\partial q^{k}}\right]_{l}=U^{k}_{l}-V^{k}_{l}\ \ 	;U^{k}_{l}>0\ \ ;V^{k}_{l}>0
\end{equation}
Then, following the S.G.M. method, we can propose a multiplicative algorithm obtained using a descent step size $\alpha^{k}=1 \ \ \forall k$ written as:
\begin{equation}
x^{k+1}_{l}=x^{k}_{l}\left[\frac{U^{k}}{V^{k}}\right]_{l}	
\end{equation}
Of course, there is no evidence of convergence of such algorithms for all the divergences, since the descent step size is fixed, however, the non-negativity constraint is satisfied for a positive initial estimate.\\
On the other hand, the sum constraint, which is not automatic as in algorithms of the form (\ref{eq.algo flux const}), can be ensured here by a normalization at each iteration; indeed, taking into account the invariance property of the divergences considered, this operation does not modify the value of the objective function i.e. of the divergence.

\setcounter{table}{0}  \setcounter{equation}{0}  \setcounter{figure}{0} \setcounter{chapter}{11} \setcounter{section}{0} 
\chapter[chapter 11 -\\Applications to the N.M.F.]{chapter 11 -\\Applications to the Non Negative Matrix Factorization.} \label{chptr::chapitre11}

 \section{Introduction.}
In this chapter, there is no attempt to repeat the work dealing with factorization in non-negative matrices (N.M.F.); a very extensive bibliography on this subject can be found in \cite{cichocki2009} and in \cite{berry2007algorithms}.\\
Instead, we will develop algorithms that take into account the sum constraints involved in NMF.\\
In addition, we will demonstrate the interest of using scale invariant divergences; indeed, some properties of these divergences allow to take into account easily sum constraints. These properties will be established, and their influence on the corresponding algorithms will be shown.\\
Regularized algorithms are also considered, whether the divergences are invariant or not.\\
Some problems related to purely multiplicative regularized algorithms as proposed in the literature have been reported in \textbf{Appendix 9} and it will be explained why such problems do not appear in the methods proposed here. 

 \subsection{Linear unmixing.}
 The problem is presented in the context of hyperspectral data.\\
 At the simplest level, which is the \textbf{linear unmixing problem}, we have the measurement of a spectrum (1 only). The measured data are therefore the intensity values at different wavelengths, arranged as a vector ``$y$"; this can be seen as the decomposition over wavelengths of the intensity in 1 pixel of an image. The sum of these intensities thus represents the total intensity ``$a$" in the pixel considered $\sum_{i}y_{i}=a$.\\
 On the other hand, we have a certain number of spectra of simple elements (endmembers); each of these spectra is arranged in a column (vector), the juxtaposition of these columns forming a table that can be considered as a matrix $H$.\\
 The problem is to find the weights ``$x$" (positive or zero, of course) such that the observed spectrum is described as a weighted sum of the basic spectra, that is, of endmembers.\\
 So we have to solve with respect to ``$x$" a problem that is written in a matrix form:
\begin{equation}
	y=Hx
\end{equation}
  Furthermore, we would want that these unknown coefficients sum to 1, in order to obtain percentages of each of the elementary spectra in the measured spectrum.\\
 Without further details on the elementary spectra available (for example the integrals of these spectra), the problem seems to me to be unsolvable; for a solution to be envisaged, all the spectra considered, both the measured and the reference spectra, would have to be of the same integral (i.e. the corresponding vectors would have to be of the same sum).\\
One possibility is to use reference spectra normalised to 1, and to impose as a constraint on the sum of the weighting coefficients $\sum_{i}x_{i}=\sum_{i}y_{i}$, the percentages being then easily obtained.\\
Another solution is to use reference spectra normalized to 1, and to normalize the measured spectrum to 1, then the constraint on the sum of the weighting coefficients at 1 is immediate.\\
In both cases, the common point is the normalization of the reference spectra, i.e. the sum to 1 of the columns of the $H$ matrix.\\ 
However, if we perform a simulation, that is, if we generate the measured spectrum as a weighted sum of elementary spectra with weighting coefficients of sum 1, we can hope to find a solution to the inverse problem, whatever the integrals of the reference spectra, but then we are very far from the real problem.\\
The real case that comes closest to this situation is the one where the ``endmembers" with the ``right" properties are extracted from the data prior to solving the inverse problem.\\ 
Finally, in a real case, an additional problem may arise, because if we have a very large set of reference spectra, even if we assume that all integrals problems are solved, we may think that all reference spectra will not necessarily contribute to the measured spectrum, then the corresponding weighting coefficients will be zero, so we have to consider that the solution of the problem i.e. the vector of weighting coefficients has a ``\textbf{sparse}" structure, which introduces an additional constraint.  

\subsection{Multispectral case.}
Here, the problem becomes a little more complicated.\\
 Suppose we have an image, the result of an observation; the intensity in each of the pixels of the image is decomposed according to wavelengths as in the previous section, so we will have as many measured spectra as there are pixels in the image; each of these spectra (vector) constitutes a column of a table that represents the measurements. For the moment, it is only a table (not yet a matrix).\\
Note that the sum of the components in a column represents the total intensity in the corresponding pixel and there is no reason why all these values should be equal.\\
Knowing as in the previous section a set of reference spectra (column vectors) juxtaposed to form a matrix $H$, the problem is thus to find for each measured spectrum, the vector (column) of the weighting coefficients.\\
This is just a succession of problems similar to the one described in the previous section.\\
If, now, to see the problem as a whole, the (non-negative) weights are arranged in column vectors, the juxtaposition of these columns forms a table $X$ which can be called a matrix and treated as such.\\
 The problem can now be written in matrix form:
\begin{equation}
Y=H\;X 
\end{equation}
Where $Y$ is of dimension $\left(L*C\right)$; ``$L$" is the number of wavelengths, ``$C$" is the number of pixels in the observed image,\\
$H$ is of dimension $\left(L*M\right)$; ``$M$" is the number of reference spectra,\\
$X$ is $\left(M*C\right)$.\\
Obviously, ``$L$" and ``$C$" are fixed by the experimental conditions from which the measured spectra are obtained, while ``$M$" depends on the number of reference spectra that are available.\\
If the columns of the matrix $H$ are normalized to $1$, then we have $\sum_{i}Y_{ij}=\sum_{i}X_{ij},\ \forall j$; this last inequality will constitute a constraint of the problem.\\

\subsection{Case of the NMF.}
Here, things get even more complicated, indeed, in this case, we only have the $Y$ measurements and we're looking for both the reference spectra matrix $H$ and the weighting coefficients matrix $X$.\\
The first difficulty which appears in an obvious way is the multiplicity of the solutions, indeed, the relation $Y=HX$ can be written $Y=HDD^{-1}X$, with $D$ invertible with non-negative terms, and any decomposition of the form $\widetilde{H}=HD$, $\widetilde{X}=D^{-1}X$ also provides a solution.\\
From this point of view, the introduction of constraints such as $\sum_{i}H_{ij}=1\ \forall j$ and $\sum_{i}X_{ij}=\sum_{i}Y_{ij}\ \forall j$, helps to reduce the ambiguity.\\
However, an appropriate choice of ``$M$" is generally critical, this choice is of course problem dependent; nevertheless, ``$M$" is often chosen such that $M<<\min(L,C)$, which is singularly imprecise.
In Lee and Seung's article \cite{lee2001algorithms} another rule is proposed, it is written $M<\frac{LC}{L+C}$, which corresponds to the rule \textit{``number of data $>$ number of unknowns"} which can be understood in the context of systems of linear equations.\\
Finally, it's clear that the product $HX$ is only an approximation of $Y$ to the rank ``$M$".

\section{Generalities.}
All the problems previously mentioned imply looking for the solution of a minimization problem with respect to the unknown parameters, a discrepancy between the measures $Y\equiv P$ and the linear model $HX\equiv Q$. 
By noting $D\left(P\|Q\right)$ the divergence between the tables (matrices) $P$ and $Q$:
\begin{equation}
D\left(P\|Q\right)=\sum_{ij}D\left(P_{ij}\|Q_{ij}\right)=\sum_{ij}D_{ij}	
\end{equation}
With for our applications:
\begin{equation}
	P_{ij}=Y_{ij};\ \ \ Q_{ij}=\left[HX\right]_{ij}=\sum_{l}H_{_{il}}X_{lj}
\end{equation}
The minimization techniques proposed in the literature and used in this book require the calculation of the
 gradient of this divergence with respect to the elements of the matrix $H$ as well as with respect to the elements of the matrix $X$.\\

\subsubsection{* Gradient with respect to the elements of the matrix $X$.}
By introducing the matrix $\left[A\right]$ such that:
\begin{equation}
	\left[A\right]_{ij}=\left[\frac{\partial D}{\partial Q}_{ij}\right]
	\label{eq.matA}
\end{equation}
We can write in matrix form:
\begin{equation}
\left[\frac{\partial D\left(P\|Q\right)}{\partial X}\right]_{nm}=\frac{\partial D\left(P\|Q\right)}{\partial X_{nm}}=\left[H^{T}A\right]_{nm}=\sum_{l}H^{T}_{nl}\left[\frac{\partial D}{\partial Q}_{lm}\right]
\label{eq.gradX}	
\end{equation}

\subsubsection{- Details of the calculation.}
The gradient is derived as follows:
\begin{equation}
	\frac{\partial D\left(P\|Q\right)}{\partial X_{nm}}=\sum_{ij}\frac{D\left(Y_{ij}\|\left[HX\right]_{ij}\right)}{\partial \left[HX\right]_{ij}}\frac{\partial\left[HX\right]_{ij}}{\partial X_{nm}}
\end{equation}
With:
\begin{equation}
	\frac{\partial \left[HX\right]_{ij}}{\partial X_{nm}}=\sum_{l}H_{_{il}}\frac{\partial X_{lj}}{\partial X_{nm}}=\sum_{l}H_{il}\delta_{ln}\delta_{jm}
\end{equation}
We have:
\begin{equation}
	\frac{\partial \left[HX\right]_{ij}}{\partial X_{nm}}=H_{in}\delta_{jm}	
\end{equation}
Then:
\begin{equation}
\frac{\partial D\left(P\|Q\right)}{\partial X_{nm}}=\sum_{ij}\frac{\partial D\left(Y_{ij}\|\left[HX\right]_{ij}\right)}{\partial \left[HX\right]_{ij}}H_{in}\delta_{jm}	
\end{equation}
And finally:
\begin{equation}
\frac{\partial D\left(P\|Q\right)}{\partial X_{nm}}=\sum_{i}\frac{\partial D\left(Y_{im}\|\left[HX\right]_{im}\right)}{\partial \left[HX\right]_{im}}H_{in}
\end{equation}
Or also:
\begin{equation}
\frac{\partial D\left(P\|Q\right)}{\partial X_{nm}}=\sum_{i}H^{T}_{ni}\frac{\partial D\left(Y_{im}\|\left[HX\right]_{im}\right)}{\partial \left[HX\right]_{im}}=\sum_{i}H^{T}_{ni}\left[\frac{\partial D}{\partial Q}_{im}\right]
\end{equation}

\subsubsection{* Gradient with respect to the elements of the matrix $H$.}
Using the matrix $\left[A\right]$ (\ref{eq.matA}), we will have in matrix form:
\begin{equation}
\left[\frac{\partial D\left(P\|Q\right)}{\partial H}\right]_{nm}=\frac{\partial D\left(P\|Q\right)}{\partial H_{nm}}=\left[A X^{T}\right]_{nm}=\sum_{j}\left[\frac{\partial D}{\partial Q_{nj}}\right]X^{T}_{jm}
\label{eq.gradH}	
\end{equation}

\subsubsection{- Details of the calculation.}
The calculation is as follows:
\begin{equation}
	\frac{\partial D\left(P\|Q\right)}{\partial H_{nm}}=\sum_{ij}\frac{\partial D\left(Y_{ij}\|\left[HX\right]_{ij}\right)}{\partial \left[HX\right]_{ij}}\frac{\partial\left[HX\right]_{ij}}{\partial H_{nm}}
\end{equation}
We have:
\begin{equation}
	\frac{\partial \left[HX\right]_{ij}}{\partial H_{nm}}=\sum_{l}\frac{\partial H_{il}}{\partial H_{nm}}X_{lj}=\sum_{l}\delta_{in}\delta_{lm}X_{lj}
\end{equation}
Then:
\begin{equation}
	\frac{\partial \left[HX\right]_{ij}}{\partial H_{nm}}=\delta_{in}	X_{mj}
\end{equation}
So:
\begin{equation}
\frac{\partial D\left(P\|Q\right)}{\partial H_{nm}}=\sum_{ij}\frac{\partial D\left(Y_{ij}\|\left[HX\right]_{ij}\right)}{\partial \left[HX\right]_{ij}}\delta_{in}	X_{mj}
\end{equation}
And finally:
\begin{equation}
\frac{\partial D\left(P\|Q\right)}{\partial H_{nm}}=\sum_{j}\frac{\partial D\left(Y_{nj}\|\left[HX\right]_{nj}\right)}{\partial \left[HX\right]_{nj}}X_{mj}
\end{equation}
Or also:
\begin{equation}
\frac{\partial D\left(P\|Q\right)}{\partial H_{nm}}=\sum_{j}\frac{\partial D\left(Y_{nj}\|\left[HX\right]_{nj}\right)}{\partial \left[HX\right]_{nj}}X^{T}_{jm}=\sum_{j}\left[\frac{\partial D}{\partial Q_{nj}}\right]X^{T}_{jm}
\end{equation}

\section{Algorithmic.}
\subsection{Principle of the method.}
The unknowns of the problem are the matrices $H$ and $X$ of the product $HX\equiv Q$, the data are the matrix $Y\equiv P$. We look for $H$ and $X$ by minimizing a divergence between the matrices $Y$ and $HX$ which is written as follows:
\begin{equation}
	X;\ H; = \underbrace{Arg\;min\;}_{X,H}D\left(Y\|\left[HX\right]\right)=\underbrace{Arg\;min\;}_{X,H}\sum_{j}\sum_{i}D\left(Y_{ij}\|\left[HX\right]_{ij}\right)
\end{equation}
We can notice that we so write a divergence between the columns of the two matrices, and that we then sum over all the columns.\\
The method of minimization generally proposed is an iterative method in which one operates alternately on the two unknowns which can be summarized for example according to the diagram:\\
1 - Iteration on $H$: $X^{k}\;,\;H^{k}\rightarrow X^{k}\;,\;H^{k+1}$\\
2 - Iteration on $X$: $X^{k}\;,\;H^{k+1}\rightarrow X^{k+1}\;,\;H^{k+1}$\\
It is this flowchart that will be considered in the rest of the computations, but the order of update is not mandatory.\\
In this diagram, one point remains to be clarified: is the updating done element by element, or column by column?\\
From my point of view, given the constraints imposed on the columns of $H$ and $X$, it is obvious to operate column by column.\\
Indeed, the constraints imposed on the unknowns are as follows:\\
* $H_{ij}\geq 0\ \ \forall i,j$.\\
* $X_{ij}\geq 0\ \ \forall i,j$.\\
* $\sum_{i}H_{ij}=1\ \ \forall j$.\\ 
* $\sum_{i}X_{ij}=\sum_{i}Y_{ij}\ \ \forall j$.\\ 
If the columns of the matrix $Y$ are normalized to 1, the sum constraint on $X$ becomes:\\
* $\sum_{i}X_{ij}=1\ \ \forall j$.\\ 
The main difficulty encountered in this problem is that even if the divergence considered is separately convex in $H$ and $X$, it is not jointly convex with respect to the two unknowns; consequently, there is no guarantee that the absolute minimum is reached by the proposed algorithmic method, one can only say that a local minimum of the objective function is reached, therefore a solution that is not necessarily optimal.  

\subsubsection{* An alternative approach.}
In the problems of searching for saddle points in constrained optimization based on Lagrangian methods, we are led to operate on primary variables in a minimization step and on dual variables in a maximization step.\\
Two algorithmic methods are proposed \cite{culioli2012introduction} (p.131-132, fig.3.10-3.11):\\
1 - The Arrow-Hurwitz method, in which an iteration step (minimization) is performed on the primal variables followed by an iteration step (maximization) on the dual variables.\\
The procedure described in the previous section has a clear analogy with this method, it is of course understood that in the case of NMF, we must successively minimize with respect to the two unknowns.\\
2 -  The Uzawa method in which one performs a mininization (several iterative steps) on the primal variables until a stop criterion is met, followed by a maximization (several iterative steps) on the dual variables until a stop criterion is met, and so on.\\
The Uzawa method allows us to consider a procedure for MMF that would consist of first minimizing $H$ (for example), doing several iteration steps (until a stop criterion is fulfilled), then repeating the same procedure on $X$, and so on.\\
Such a procedure has never been considered to my knowledge for NMF, however, the example shown in \cite{culioli2012introduction} seems to indicate better results with the Uzawa method in terms of convergence speed.\\
Of course, this is just an example applied to saddle point research, i.e. in a very different context, and therefore nothing is well established....\\
This alternative method is not used here; it has been used for blind deconvolution in \cite{lanteri1994blind}, \cite{lanteri1995comparison}.

\subsection{Iterations with respect to $H$.}
\textbf{Whatever the type of divergence used, whether it is scale invariant with respect to $Q$ or non-invariant, taking into account the sum constraint on the columns of $H$ implies to proceed by means of a variable change method.}\\

\textbf{- Important remark.}\\
\textbf{Multiplicative algorithms can always be obtained by introducing a shift in the terms constituting the gradient of the divergence considered; this point will be developed in the following section.}\\

\textbf{**In this section, we note $Q=H^{k}X^{k}$.}
\subsubsection{1 - Non-multiplicative algorithm.}
In order to take into account the constraint $\sum_{i}H_{ij}=1\ \forall j$, we proceed to the change of variables:
\begin{equation}
	H_{ij}=\frac{z_{ij}}{\sum_{l}z_{lj}}
	\label{eq.chgtvarH}
\end{equation}
First, we build an iterative algorithm on ``$z$", then we go back to an algorithm on $H$.\\
To do that, we have:
\begin{equation}
	\frac{\partial D}{\partial z_{nm}}=\sum_{i}\frac{\partial D}{\partial H_{im}}\frac{\partial H_{im}}{\partial z_{nm}}
\end{equation}
With:
\begin{equation}
	\frac{\partial H_{im}}{\partial z_{nm}}=\frac{1}{\sum_{l}z_{lm}}\left(\delta_{in}-H_{im}\right)
\end{equation}
We have:
\begin{equation}
	\frac{\partial D}{\partial z_{nm}}=\sum_{i}\frac{\partial D}{\partial H_{im}}\frac{1}{\sum_{l}z_{lm}}\left(\delta_{in}-H_{im}\right)
\end{equation}
And finally:
\begin{equation}
	\frac{\partial D}{\partial z_{nm}}=\frac{1}{\sum_{l}z_{lm}}\left[\frac{\partial D}{\partial H_{nm}}-\sum_{i}H_{im}\frac{\partial D}{\partial H_{im}}\right]
	\label{eq.derivDznm}
\end{equation}
One can immediately verify that we have:
\begin{equation}
	\sum_{n}z_{nm}\frac{\partial D}{\partial z_{nm}}=0
	\label{eq.propgradD}
\end{equation}
We are reasoning on the column ``$m$" and we write in a first step, an algorithm on the ``$n$" component:
\begin{equation}
	z^{k+1}_{nm}=z^{k}_{nm}+\alpha^{k}_{nm}z^{k}_{nm}\left[-\frac{\partial D}{\partial z^{k}_{nm}}\right]
\end{equation}
Taking into account the considerations on the determination of the descent step size outlined in \textbf{Chapter 10}, an algorithm that ensures the non-negativity and convergence of the algorithm for all ``$n$" components can be written in the form:
\begin{equation}
	z^{k+1}_{nm}=z^{k}_{nm}+\alpha^{k}_{m}z^{k}_{nm}\left[-\frac{\partial D}{\partial z^{k}_{nm}}\right]
\end{equation}
and finally, taking into account (\ref{eq.derivDznm}):
\begin{equation}
	z^{k+1}_{nm}=z^{k}_{nm}+\alpha^{k}_{m}\frac{z^{k}_{nm}}{\sum_{l}z^{k}_{lm}}\left[\left(-\frac{\partial D}{\partial H^{k}_{nm}}\right)-\sum_{i}H^{k}_{im}\left(-\frac{\partial D}{\partial H^{k}_{im}
	}\right)\right]
	\label{eq.algoA}
\end{equation}
From this expression, according to (\ref{eq.chgtvarH}) and (\ref{eq.propgradD}) we can easily verify that we have, $\forall m$:
\begin{equation}
	\sum_{n}z^{k+1}_{nm}=\sum_{n}z^{k}_{nm}
	\label{eq.sommefixe}
\end{equation}
With such an algorithm, we move in a solution space such that $\sum_{n}z^{k}_{nm}=Cte\ \forall k$, so if the divergence considered is convex with respect to $H$, the convexity is maintained during the change of variables.\\
Therefore, going back to the algorithm (\ref{eq.algoA}), and dividing the whole by $\sum_{n}z^{k}_{nm}$, it comes:
\begin{equation}
		H^{k+1}_{nm}=H^{k}_{nm}+\alpha^{k}_{m}\frac{H^{k}_{nm}}{\sum_{l}z^{k}_{lm}}\left[\left(-\frac{\partial D}{\partial H^{k}_{nm}}\right)-\sum_{i}H^{k}_{im}\left(-\frac{\partial D}{\partial H^{k}_{im}
	}\right)\right]
\end{equation}
or also, with $\delta^{k}_{m}=\frac{\alpha^{k}_{m}}{\sum_{l}z^{k}_{lm}}$:
\begin{equation}
		H^{k+1}_{nm}=H^{k}_{nm}+\delta^{k}_{m}H^{k}_{nm}\left[\left(-\frac{\partial D}{\partial H^{k}_{nm}}\right)-\sum_{i}H^{k}_{im}\left(-\frac{\partial D}{\partial H^{k}_{im}
	}\right)\right]
	\label{eq.algonmH}
\end{equation}
Which can be written in condensed form:
\begin{equation}
		H^{k+1}_{nm}=H^{k}_{nm}+\delta^{k}_{m}H^{k}_{nm}\left\{\left(-\frac{\partial D}{\partial H^{k}_{nm}}\right)-\left[\left(H^{k}\right)^{T}\left(-\frac{\partial D}{\partial H^{k}
	}\right)\right]_{mm}\right\}
	\label{eq.algonmHcond}
\end{equation}
Note that the second of the terms between braces is a constant for the column ``$m$" considered.\\
It can then be observed that regardless of the descent step size $\delta^{k}_{m}$, we have $\forall m$:
\begin{equation}
	\sum_{n}H^{k+1}_{nm}=\sum_{n}H^{k}_{nm}=1
	\label{eq.sommefixeH}
\end{equation}
\textbf{The algorithm (\ref{eq.algonmH}) (\ref{eq.algonmHcond}) represents the algorithm corresponding to the iterations on $H$ for the NMF.\\
It is initialized with:\\
*  $H^{0}$ such that $\sum_{n}H^{0}_{nm}=1\ \forall m$.\\
*  $X^{0}$ such that $\sum_{n}X^{0}_{nm}=\sum_{n}Y_{nm}\ \forall m$.\\
Note that in this algorithm, one can update the $H$ matrix, column by column.}\\
Using the expression of the gradient relative to $H$ (\ref{eq.gradH}), we can also write the algorithm (\ref{eq.algonmH}) in the form:
\begin{align}
		H^{k+1}_{nm}=H^{k}_{nm}+\delta^{k}_{m}H^{k}_{nm}&\left[\left(\sum_{l}\left[-\frac{\partial D}{\partial Q_{nl}}\right]X^{T}_{lm}\right)\right. \nonumber \\  & \left.-\sum_{i}H^{k}_{im}\left(\sum_{l}\left[-\frac{\partial D}{\partial Q_{il}}\right]X^{T}_{lm}\right)\right]
	\label{eq.algonmHbis}
\end{align}
Or else in condensed form:
\begin{equation}
		H^{k+1}_{nm}=H^{k}_{nm}+\delta^{k}_{m}H^{k}_{nm}\left\{\left[\left(-\frac{\partial D}{\partial Q}\right)X^{T}\right]_{nm}-\left[\left(H^{k}\right)^{T}\left(-\frac{\partial D}{\partial Q}\right)X^{T}\right]_{mm}\right\}
	\label{eq.algonmHbiscond}
\end{equation}
With, of course, in that expression, $X\equiv X^{k}$ and $Q=H^{k}X^{k}$.

\subsubsection{2 - Multiplicative algorithm.}
If we want to obtain a multiplicative algorithm while preserving the non-negativity and sum constraints, we must first introduce an shift of the components of the opposite of the gradient $\left[-\frac{\partial D}{\partial H^{k}_{ij}}\right]$ in such a way as to make them all non-negative; to do this, one writes:
\begin{equation}
	\left[-\frac{\partial D}{\partial H^{k}_{ij}}\right]_{d}=\left[-\frac{\partial D}{\partial H^{k}_{ij}}\right]-\min_{i}\left[-\frac{\partial D}{\partial H^{k}_{ij}}\right]+\epsilon
\end{equation}
Taking into account the definition of $H_{ij}$ (\ref{eq.chgtvarH}), the algorithm (\ref{eq.algonmH}) can be rewritten:
\begin{equation}
		H^{k+1}_{nm}=H^{k}_{nm}+\delta^{k}_{m}H^{k}_{nm}\left[\left(-\frac{\partial D}{\partial H^{k}_{nm}}\right)_{d}-\sum_{i}H^{k}_{im}\left(-\frac{\partial D}{\partial H^{k}_{im}
	}\right)_{d}\right]
	\label{eq.algonmHd}
\end{equation}
With this algorithm, the property (\ref{eq.sommefixeH}) is maintained and using the SGM method presented in \textbf{Chapter 10}, we can write an algorithm which is the basis of the multiplicative forms, as follows:
\begin{equation}
		H^{k+1}_{nm}=H^{k}_{nm}+\delta^{k}_{m}H^{k}_{nm}\left[\frac{\left(-\frac{\partial D}{\partial H^{k}_{nm}}\right)_{d}}{\sum_{i}H^{k}_{im}\left(-\frac{\partial D}{\partial H^{k}_{im}}\right)_{d}}-1 \right]
\end{equation}
As is always the case, a purely multiplicative algorithm can be obtained if one chooses for all iterations a descent step equal to 1.\\
The algorithm is then written:
\begin{equation}
		H^{k+1}_{nm}=H^{k}_{nm}\left[\frac{\left(-\frac{\partial D}{\partial H^{k}_{nm}}\right)_{d}}{\sum_{i}H^{k}_{im}\left(-\frac{\partial D}{\partial H^{k}_{im}}\right)_{d}}\right]
\label{eq.algomH}
\end{equation}
Using the expression of (\ref{eq.gradH}), we can also write, with $X\equiv X^{k}$ et $Q=H^{k}X^{k}$:
\begin{equation}
		H^{k+1}_{nm}=H^{k}_{nm}\left[\frac{\left(\sum_{l}\left[-\frac{\partial D}{\partial Q_{nl}}\right]X^{T}_{lm}\right)_{d}}{\sum_{i}H^{k}_{im}\left(\sum_{l}\left[-\frac{\partial D}{\partial Q_{il}}\right]X^{T}_{lm}\right)_{d}}\right]
\label{eq.algomHbis}
\end{equation}
 But of course, if the non-negativity and sum constraints are fulfilled, there is no guarantee that such an algorithm will be convergent.

\subsection{Iterations on $X$.}
For these unknowns, in order to impose the constraint $\sum_{i}X_{ij}=\sum_{i}Y_{ij}\ \forall j$ (possibly =1), two cases must be distinguished depending on whether the divergence to be minimized is scale invariant or not.\\
If a non-invariant divergence is used, a change in variables similar to that used for $H$ must be introduced.\\
If an invariant divergence is used, the implementation of the general property of these divergences deduced from (\ref{eq.Pfond}), allows to ensure that the sum constraint is satisfied, without it being necessary to introduce a change of variables; this aspect will be developed in the section devoted to the scale invariant divergences.\\

\textbf{In this section, we develop the case of non-invariant divergences; we therefore proceed by change of variables.}\\

\textbf{The use of non-invariant divergences associated with the variable change method will lead to multiplicative algorithms for $X$ in a similar way to what has been done for $H$.}\\

\textbf{We have obtained $H^{k+1}$, we know $X^{k}$, we are looking for $X^{k+1}$, so here, $Q=H^{k+1}X^{k}$.}

\subsubsection{1 - Non-multiplicative algorithm.}
Taking into account the sum constraints imposed on the columns of the matrix $X$, the change of variables corresponding to $X$ is written as follows:
\begin{equation}
	X_{ij}=\frac{t_{ij}}{\sum_{l}t_{lj}}\sum_{l}Y_{lj}\ \ \Leftrightarrow\ \ \sum_{i}X_{ij}=\sum_{l}Y_{lj}
	\label{eq.chgtvarX}
\end{equation}
For a divergence $D\left(P\|Q\right)$, the gradient with respect to ``$t$" is written as:
\begin{equation}
	\frac{\partial D\left(P\|Q\right)}{\partial t_{nm}}=\sum_{i}\frac{\partial D\left(P\|Q\right)}{\partial X_{im}}\frac{\partial X_{im}}{\partial t_{nm}}
\end{equation}
With:
\begin{equation}
\frac{\partial X_{im}}{\partial t_{nm}}=\frac{\sum_{l}Y_{lm}}{\sum_{l}t_{lm}}\left[\delta_{in}\right]-\frac{X_{im}}{\sum_{l}t_{lm}}=\frac{\sum_{l}Y_{lm}}{\sum_{l}t_{lm}}\left[\delta_{in}\right]-\frac{X_{im}}{\sum_{l}t_{lm}}	
\end{equation}
We have:
\begin{equation}
\frac{\partial D\left(P\|Q\right)}{\partial t_{nm}}=\frac{\sum_{l}Y_{lm}}{\sum_{l}t_{lm}}\frac{\partial D\left(P\|Q\right)}{\partial X_{nm}}-\frac{1}{\sum_{l}t_{lm}}\sum_{i}X_{im}\frac{\partial D\left(P\|Q\right)}{\partial X_{im}}
\end{equation}
From this we deduce the iterative algorithm corresponding to the component ``$n$" of the column ``$m$":
\begin{align}
	t^{k+1}_{nm}=t^{k}_{nm}+\alpha^{k}_{nm}t^{k}_{nm}&\left[\frac{\sum_{l}Y_{lm}}{\sum_{l}t^{k}_{lm}}\left(-\frac{\partial D\left(P\|Q\right)}{\partial X^{k}_{nm}}\right)\right. \nonumber \\  & \left.-\frac{1}{\sum_{l}t^{k}_{lm}}\sum_{i}X^{k}_{im}\left(-\frac{\partial D\left(P\|Q\right)}{\partial X^{k}_{im}}\right)\right]
\end{align}
Taking into account the considerations on the determination of the descent step size as specified in \textbf{Chapter 10}, an algorithm that will guarantee the non-negativity and convergence for all the components ``$n$" can be written in the form: 
\begin{align}
	t^{k+1}_{nm}=t^{k}_{nm}+\alpha^{k}_{m}t^{k}_{nm}&\left[\frac{\sum_{l}Y_{lm}}{\sum_{l}t^{k}_{lm}}\left(-\frac{\partial D\left(P\|Q\right)}{\partial X^{k}_{nm}}\right)\right. \nonumber \\  & \left.-\frac{1}{\sum_{l}t^{k}_{lm}}\sum_{i}X^{k}_{im}\left(-\frac{\partial D\left(P\|Q\right)}{\partial X^{k}_{im}}\right)\right]
	\label{eq.algoB}
\end{align}
We can then easily verify that we have:
\begin{equation}
	\sum_{n}t^{k+1}_{nm}=\sum_{n}t^{k}_{nm}
\end{equation}
Consequently, if the divergence considered is convex with respect to $X$, the change of variables will not modify this property.\\
By introducing the descent step size $\delta^{k}_{m}=\frac{\alpha^{k}_{m}}{\sum_{l}t^{k}_{lm}}$, we can therefore go back to the initial variables by dividing the two members of (\ref{eq.algoB}) by $\sum_{n}t^{k}_{nm}$, which leads to:
\begin{align}
	X^{k+1}_{nm}=X^{k}_{nm}+\delta^{k}_{m}X^{k}_{nm}&\left[\left(\sum_{l}Y_{lm}\right)\left(-\frac{\partial D\left(P\|Q\right)}{\partial X^{k}_{nm}}\right)\right. \nonumber \\  & \left.-\sum_{i}X^{k}_{im}\left(-\frac{\partial D\left(P\|Q\right)}{\partial X^{k}_{im}}\right)\right]
	\label{eq.algoBbis}
\end{align}
Which can also be written in condensed form:
\begin{align}
	X^{k+1}_{nm}=X^{k}_{nm}+\delta^{k}_{m}X^{k}_{nm}&\left\{\left(\sum_{l}Y_{lm}\right)\left(-\frac{\partial D\left(P\|Q\right)}{\partial X^{k}}\right)_{nm}\right. \nonumber \\  & \left.-\left[\left(X^{k}\right)^{T}\left(-\frac{\partial D\left(P\|Q\right)}{\partial X^{k}}\right)\right]_{mm}\right\}
	\label{eq.algoBbiscond}
\end{align}
Note that the second term in the braces is a constant for the whole ``$m$" column. \\
With an initialization such as $\sum_{n}X^{0}_{nm}=\sum_{n}Y_{nm}\ \forall m$, taking into account the change of variables (\ref{eq.chgtvarX}), we can check without difficulty that the sum constraint on the columns of $X$ is satisfied and that we have:
\begin{equation}
	\sum_{n}X^{k+1}_{nm}=\sum_{n}X^{k}_{nm}=\sum_{n}Y_{nm}\ \ \forall m
\end{equation}
In this algorithm, we will use the expression of the gradient with respect to $X$:
\begin{equation}
\left[\frac{\partial D\left(P\|Q\right)}{\partial X_{nm}}\right]=\sum_{l}H^{T}_{nl}\left[\frac{\partial D}{\partial Q}\right]_{lm}	
\end{equation}
With $H\equiv H^{k+1}$ and $Q=H^{k+1}X^{k}$.
	
\subsubsection{2 - Multiplicative algorithm.}
If we wish to obtain a multiplicative algorithm while preserving the non-negativity and sum constraints, we must first introduce a shift of the components of the opposite of the gradient $\left[-\frac{\partial D}{\partial X^{k}_{ij}}\right]$ in such a way as to make them all non-negative; to do this, one writes:
\begin{equation}
	\left[-\frac{\partial D}{\partial X^{k}_{ij}}\right]_{d}=\left[-\frac{\partial D}{\partial X^{k}_{ij}}\right]-\min_{i}\left[-\frac{\partial D}{\partial X^{k}_{ij}}\right]+\epsilon
\end{equation}
We can then rewrite the algorithm (\ref{eq.algoBbis}) in the form:
\begin{align}
	X^{k+1}_{nm}=X^{k}_{nm}+\delta^{k}_{m}X^{k}_{nm}&\left[\left(\sum_{l}Y_{lm}\right)\left(-\frac{\partial D\left(P\|Q\right)}{\partial X^{k}_{nm}}\right)_{d}\right. \nonumber \\  & \left.-\sum_{i}X^{k}_{im}\left(-\frac{\partial D\left(P\|Q\right)}{\partial X^{k}_{im}}\right)_{d}\right]	
\end{align}
According to the SGM method presented in \textbf{Chapter 10}, we can deduce the algorithm:
\begin{equation}
	X^{k+1}_{nm}=X^{k}_{nm}+\delta^{k}_{m}X^{k}_{nm}\left[\frac{\left(\sum_{l}Y_{lm}\right)\left(-\frac{\partial D\left(P\|Q\right)}{\partial X^{k}_{nm}}\right)_{d}}{\sum_{i}X^{k}_{im}\left(-\frac{\partial D\left(P\|Q\right)}{\partial X^{k}_{im}}\right)_{d}}-1\right]
\end{equation}
Then, using a descent step size equal to 1, we obtain the purely multiplicative algorithm:
\begin{equation}
	X^{k+1}_{nm}=X^{k}_{nm}\left[\frac{\left(\sum_{l}Y_{lm}\right)\left(-\frac{\partial D\left(P\|Q\right)}{\partial X^{k}_{nm}}\right)_{d}}{\sum_{i}X^{k}_{im}\left(-\frac{\partial D\left(P\|Q\right)}{\partial X^{k}_{im}}\right)_{d}}\right]
	\label{eq.algomultX}
\end{equation}
We can easily verify that we have: $\sum_{n}X^{k+1}_{nm}=\sum_{n}Y_{nm}$.\\
Nevertheless, for this algorithm, as for all algorithms with a fixed descent step size, there is no guarantee of convergence, even if all the constraints are fulfilled.\\

\textbf{* Important remark :} We will see in the following sections that, in the particular case of \textbf{``Alpha divergences"}, given the specific expression of the gradient, the iterative algorithm on $H$ (\ref{eq.algonmH}) and the iterative algorithm on $X$ (\ref{eq.algoBbis}) have a simplified expression that makes the shift unnecessary. Thus, a multiplicative algorithm can be written directly using the SGM method presented in \textbf{Chapter 10}.

\subsection{Non-invariant divergences - Overview of the algorithms.}
By performing alternately the iterations, first on $H$, then on $X$, and with the initialisations:\\
* $H^{0}$ such that $\sum_{n}H^{0}_{nm}=1\ \forall m$\\
* $X^{0}$ such that $\sum_{n}X^{0}_{nm}=\sum_{n}Y_{nm}\ \forall m$;\\

we first have:\\  
** $X\equiv X^{k}$, $H\equiv H^{k}$ et $Q=H^{k}X^{k}$\\

The algorithm obtained using the change of variables on $H$ is written as follows (\ref{eq.algonmH}):
\begin{equation}
		H^{k+1}_{nm}=H^{k}_{nm}+\delta^{k}_{m}H^{k}_{nm}\left[\left(-\frac{\partial D}{\partial H^{k}_{nm}}\right)-\sum_{i}H^{k}_{im}\left(-\frac{\partial D}{\partial H^{k}_{im}
	}\right)\right]
\end{equation}
or also (\ref{eq.algonmHbis}):
\begin{align}
		H^{k+1}_{nm}=H^{k}_{nm}+\delta^{k}_{m}H^{k}_{nm}&\left[\left(\sum_{l}\left[-\frac{\partial D}{\partial Q_{nl}}\right]X^{T}_{lm}\right)\right. \nonumber \\  & \left.-\sum_{i}H^{k}_{im}\left(\sum_{l}\left[-\frac{\partial D}{\partial Q_{il}}\right]X^{T}_{lm}\right)\right]
\end{align}
The algorithm can equally be obtained in a purely multiplicative form (\ref{eq.algomHbis}):
\begin{equation}
		H^{k+1}_{nm}=H^{k}_{nm}\left[\frac{\left(\sum_{l}\left[\frac{-\partial D}{\partial Q_{nl}}\right]X^{T}_{lm}\right)_{d}}{\sum_{i}H^{k}_{im}\left(\sum_{l}\left[-\frac{\partial D}{\partial Q_{il}}\right]X^{T}_{lm}\right)_{d}}\right]
\end{equation}
At this point, we have:\\

**  $H\equiv H^{k+1}$, $X\equiv X^{k}$ and $Q=H^{k+1}X^{k}$.\\

With $\delta^{k}_{m}=\frac{\alpha^{k}_{m}}{\sum_{l}t^{k}_{lm}}$, the algorithm obtained using the change of variables on $X$ (\ref{eq.algoBbis}) is written:\\	
\begin{align}
	X^{k+1}_{nm}=X^{k}_{nm}+\delta^{k}_{m}X^{k}_{nm}&\left[\left(\sum_{l}Y_{lm}\right)\left(-\frac{\partial D\left(P\|Q\right)}{\partial X^{k}_{nm}}\right)\right. \nonumber \\  & \left.-\sum_{i}X^{k}_{im}\left(-\frac{\partial D\left(P\|Q\right)}{\partial X^{k}_{im}}\right)\right]
\end{align}
In this algorithm, we will use the expression of the gradient with respect to $X$:\\
\begin{equation}
\left[\frac{\partial D\left(P\|Q\right)}{\partial X_{nm}}\right]=\sum_{l}H^{T}_{nl}\left[\frac{\partial D}{\partial Q}\right]_{lm}	
\end{equation}
The corresponding algorithm in purely multiplicative form (\ref{eq.algomultX}) is written as follows: 
\begin{equation}
	X^{k+1}_{nm}=X^{k}_{nm}\left[\frac{\left(\sum_{l}Y_{lm}\right)\left(-\frac{\partial D\left(P\|Q\right)}{\partial X^{k}_{nm}}\right)_{d}}{\sum_{i}X^{k}_{im}\left(-\frac{\partial D\left(P\|Q\right)}{\partial X^{k}_{im}}\right)_{d}}\right]
\end{equation}\\
At this point, a cycle has been run on the 2 variables; we have therefore obtained $H^{k+1}$ and $X^{k+1}$.

\section{Applications to some typical divergences.}
In this section, the algorithms corresponding to the ``Alpha" and ``Beta" divergences are discussed.

\subsection {``Alpha'' divergences.}
The basic divergence (\ref{eq.AC}) when applied to tables is written as follows:
\begin{equation}
	A\left(P\|Q\right)=\frac{1}{\lambda\left(\lambda-1\right)}\sum_{j}\sum_{i}\left[P^{\lambda}_{ij}Q^{1-\lambda}_{ij}-\lambda P_{ij}-\left(1-\lambda\right)Q_{ij}\right]
	\label{eq.ACbis}
\end{equation}
The expression of the gradient with respect to $Q$ is written as follows:
\begin{equation}
	\frac{\partial A\left(P\|Q\right)}{\partial Q_{nm}}=\frac{1}{\lambda}\left[1-\left(\frac{P_{nm}}{Q_{nm}}\right)^{\lambda}\right]
	\label{eq.gradA/Q}
\end{equation}

\textbf{We proceed alternately to iterations, first on $H$, then on $X$.}
\subsubsection {* Itérations on $H$.}

 ** First, we have: $X\equiv X^{k}$, $H\equiv H^{k}$, $P=Y$ and $Q=H^{k}X^{k}$\\

From (\ref{eq.gradA/Q}), we have:
\begin{equation}
	\frac{\partial A\left(P\|Q\right)}{\partial H_{nm}}=\sum_{j}\frac{\partial A\left(P\|Q\right)}{\partial Q_{nj}}X^{T}_{jm}=\frac{1}{\lambda}\sum_{j}\left[1-\left(\frac{P_{nj}}{Q_{nj}}\right)^{\lambda}\right]X^{T}_{jm}
	\label{eq.gradA/H}
\end{equation}
The iterative algorithm on $H$ is then written from (\ref{eq.algonmHbis}):
\begin{align}
		H^{k+1}_{nm}=H^{k}_{nm}+\frac{\delta^{k}_{m}}{\lambda}H^{k}_{nm}&\left[\left(\sum_{l}\left[\left(\frac{P_{nl}}{Q_{nl}}\right)^{\lambda}-1\right]X^{T}_{lm}\right)\right. \nonumber \\  & \left.-\sum_{i}H^{k}_{im}\left(\sum_{l}\left[\left(\frac{P_{il}}{Q_{il}}\right)^{\lambda}-1\right]X^{T}_{lm}\right)\right]
			\label{eq.AlgoA/H}
\end{align}
This algorithm can be used as it is, or a multiplicative algorithm can be obtained by performing an offset as defined in the previous section.\\
 However, in the particular case of the Alpha divergences, a simplification occurs and the algorithm (\ref{eq.AlgoA/H}) is rewritten:
\begin{align}
		H^{k+1}_{nm}=H^{k}_{nm}+\frac{\delta^{k}_{m}}{\lambda}H^{k}_{nm}&\left[\left(\sum_{l}\left(\frac{P_{nl}}{Q_{nl}}\right)^{\lambda}X^{T}_{lm}\right)\right. \nonumber \\  & \left.-\sum_{i}H^{k}_{im}\left(\sum_{l}\left(\frac{P_{il}}{Q_{il}}\right)^{\lambda}X^{T}_{lm}\right)\right]
			\label{eq.AlgoA/Hbis}
\end{align}
Consequently, the multiplicative algorithm will be obtained without the need for an offset, by taking a descent step size equal to $\frac{\delta^{k}_{m}}{\lambda}=1$, from the expression:
\begin{equation}
		H^{k+1}_{nm}=H^{k}_{nm}+\frac{\delta^{k}_{m}}{\lambda}H^{k}_{nm}\left[\frac{\sum_{l}\left(\frac{P_{nl}}{Q_{nl}}\right)^{\lambda}X^{T}_{lm}}{\sum_{i}H^{k}_{im}\left(\sum_{l}\left(\frac{P_{il}}{Q_{il}}\right)^{\lambda}X^{T}_{lm}\right)}-1\right]
			\label{eq.AlgoA/Hmult}
\end{equation}
\subsubsection {* Itérations on $X$.}

 ** First, we have: $X\equiv X^{k}$, $H\equiv H^{k+1}$, $P=Y$ and $Q=H^{k+1}X^{k}$\\

From the expression (\ref{eq.gradA/Q}), we have:
\begin{equation}
	\frac{\partial A\left(P\|Q\right)}{\partial X_{nm}}=\sum_{j}H^{T}_{nj}\frac{\partial A\left(P\|Q\right)}{\partial Q_{jm}}
\end{equation}
Or also:
\begin{equation}
	\frac{\partial A\left(P\|Q\right)}{\partial X_{nm}}=\frac{1}{\lambda}\left[\sum_{j}H^{T}_{nj}-\sum_{j}H^{T}_{nj}\left(\frac{P_{jm}}{Q_{jm}}\right)^{\lambda}\right]
	\label{eq.gradA/X}
\end{equation}
Considering the properties of $H$, it comes:
\begin{equation}
	\frac{\partial A\left(P\|Q\right)}{\partial X_{nm}}=\sum_{j}H^{T}_{nj}\frac{\partial A\left(P\|Q\right)}{\partial Q_{jm}}=\frac{1}{\lambda}\left[1-\sum_{j}H^{T}_{nj}\left(\frac{P_{jm}}{Q_{jm}}\right)^{\lambda}\right]
	\label{eq.gradA/Xbis}
\end{equation}
From (\ref{eq.algoBbis}) the iterative algorithm on $X$ is then written:
\begin{align}
	X^{k+1}_{nm}=X^{k}_{nm}+\frac{\delta^{k}_{m}}{\lambda}X^{k}_{nm}&\left\{\sum_{l}Y_{lm}\left[\sum_{j}H^{T}_{nj}\left(\frac{P_{jm}}{Q_{jm}}\right)^{\lambda}-1\right]\right. \nonumber \\  & \left.-\sum_{i}X^{k}_{im}\left[\sum_{j}H^{T}_{ij}\left(\frac{P_{jm}}{Q_{jm}}\right)^{\lambda}-1\right]\right\}
	\label{eq.AlgoA/X}
\end{align}
Which can be simplified by:
\begin{align}
	X^{k+1}_{nm}=X^{k}_{nm}+\frac{\delta^{k}_{m}}{\lambda}X^{k}_{nm}& \left\{\sum_{l}Y_{lm}\left[\sum_{j}H^{T}_{nj}\left(\frac{P_{jm}}{Q_{jm}}\right)^{\lambda}\right]\right. \nonumber \\  & \left.-\sum_{i}X^{k}_{im}\left[\sum_{j}H^{T}_{ij}\left(\frac{P_{jm}}{Q_{jm}}\right)^{\lambda}\right]\right\}
	\label{eq.AlgoA/Xbis}
\end{align}
The multiplicative algorithm can be obtained from (\ref{eq.AlgoA/X}) by introducing an offset to render the bracketed terms positive.\\
However, due to the simplified form of the Alpha divergences (\ref{eq.AlgoA/Xbis}), a multiplicative algorithm can also be obtained without the need for an offset.\\
Indeed, such an algorithm is obtained by taking a descent step $\frac{\delta^{k}_{m}}{\lambda}$ equal to $1\;\forall k$ from the expression:
\begin{equation}
	X^{k+1}_{nm}=X^{k}_{nm}+\frac{\delta^{k}_{m}}{\lambda}X^{k}_{nm}\left\{\frac{\sum_{l}Y_{lm}\left[\sum_{j}H^{T}_{nj}\left(\frac{P_{jm}}{Q_{jm}}\right)^{\lambda}\right]}{\sum_{i}X^{k}_{im}\left[\sum_{j}H^{T}_{ij}\left(\frac{P_{jm}}{Q_{jm}}\right)^{\lambda}\right]}-1\right\}
	\label{eq.AlgoA/Xmult}
\end{equation}

\subsection {``Beta'' divergences.}
The basis divergence (\ref{eq.BC}) applied to tables is written as follows:
\begin{equation}
	B\left(P\|Q\right)=\frac{1}{\lambda\left(\lambda-1\right)}\sum_{j}\sum_{i}\left[P^{\lambda}_{ij}-\lambda P_{ij}Q^{\lambda-1}_{ij}-\left(1-\lambda\right)Q^{\lambda}_{ij}\right]
	\label{eq.BCbis}
\end{equation}
The expression of the gradient with respect to $Q$ is written as follows:
\begin{equation}
	\frac{\partial B\left(P\|Q\right)}{\partial Q_{nm}}=Q^{\lambda-1}_{nm}-Q^{\lambda-2}_{nm}P_{nm}
	\label{eq.gradB/Q}
\end{equation}
\textbf{We proceed alternately to iterations, first on $H$, then on $X$.}

\subsubsection {* Itérations on $H$.}

 ** First, we have: $X\equiv X^{k}$, $H\equiv H^{k}$, $P=Y$ and $Q=H^{k}X^{k}$\\

From the expression (\ref{eq.gradB/Q}), we have:
\begin{equation}
	\frac{\partial B\left(P\|Q\right)}{\partial H_{nm}}=\sum_{j}\frac{\partial B\left(P\|Q\right)}{\partial Q_{nj}}X^{T}_{jm}=\sum_{j}\left[Q^{\lambda-1}_{nj}-Q^{\lambda-2}_{nj}P_{nj}\right]X^{T}_{jm}
	\label{eq.gradB/H}
\end{equation}
The iterative algorithm on $H$ is then written from (\ref{eq.algonmHbis}):
\begin{align}
		H^{k+1}_{nm}=H^{k}_{nm}+\delta^{k}_{m}H^{k}_{nm}& \left[\sum_{j}\left(Q^{\lambda-2}_{nj}P_{nj}-Q^{\lambda-1}_{nj}\right)X^{T}_{jm}\right. \nonumber \\  & \left.-\sum_{i}H^{k}_{im}\sum_{j}\left(Q^{\lambda-2}_{ij}P_{ij}-Q^{\lambda-1}_{ij}\right)X^{T}_{jm}\right]
			\label{eq.AlgoB/H}
\end{align}
To obtain a multiplicative algorithm that maintains the flux, there is no other solution than to shift the terms in "$\sum_{j}$" of the previous expression to render them positive.

\subsubsection {* Itérations on $X$.}
 ** First, we have: $X\equiv X^{k}$, $H\equiv H^{k+1}$, $P=Y$ and $Q=H^{k+1}X^{k}$\\

From the expression (\ref{eq.gradB/Q}), we have:
\begin{equation}
	\frac{\partial B\left(P\|Q\right)}{\partial X_{nm}}=\sum_{j}H^{T}_{nj}\frac{\partial B\left(P\|Q\right)}{\partial Q_{jm}}=\sum_{j}H^{T}_{nj}\left(Q^{\lambda-1}_{jm}-Q^{\lambda-2}_{jm}P_{jm}\right)
	\label{eq.gradB/X}
\end{equation}
From (\ref{eq.algoBbis}) the iterative algorithm on $X$ is written as follows:
\begin{align}
	X^{k+1}_{nm}= X^{k}_{nm}+\delta^{k}_{m}X^{k}_{nm}& \left[\left(\sum_{l}Y_{lm}\right)\sum_{j}H^{T}_{nj}\left(Q^{\lambda-2}_{jm}P_{jm}-Q^{\lambda-1}_{jm}\right)\right. \nonumber \\  & \left.-\sum_{i}X^{k}_{im}\sum_{j}H^{T}_{ij}\left(Q^{\lambda-2}_{jm}P_{jm}-Q^{\lambda-1}_{jm}\right)\right]
	\label{eq.AlgoB/X}
\end{align}
Considering the properties of the columns in the $X$ matrix, the flux will be maintained throughout the iterations.\\
The multiplicative algorithm can be obtained from (\ref{eq.AlgoB/X}) by introducing an offset in order to make positive the terms of the form "$\sum_{j}$" of the previous expression.\\

\textbf{Note: Unlike the case of Alpha divergences, a grouping of terms with the same sign in (\ref{eq.AlgoB/H}) and (\ref{eq.AlgoB/X}) can lead to a multiplicative algorithm; however, the algorithm thus obtained will not have the properties of preserving the sum constraints during iterations.\\}

\section{Fundamental property of scale invariant divergences.}

Two cases are considered in this paragraph:\\
In the first case, the divergence we will construct will be scale invariant when all the elements of the $Q$ matrix are multiplied by the same positive constant. The divergences obtained will have interesting properties in the case of blind deconvolution for example.\\
In the second case the scale invariance property relates to each column of the matrix $Q$, i.e. the divergence obtained will be invariant when each column of $Q$ is multiplied by a positive scalar which is specific to it. The divergences obtained will have properties that will be specifically adapted to the case of NMF.

\subsection{General case - All components of $Q$ are concerned simultaneously.}
Here, starting with a divergence $D\left(P\|Q\right)$ , we first construct a divergence $DI\left(P\|Q\right)$  that remains scale invariant when the elements of the $Q$ matrix are multiplied by the same positive constant.\\
 In this case, the nominal invariance factor, valid for the $Q$ matrix as a whole, is calculated as follows:
\begin{equation}
	K_{0}\left(P,Q\right)=Arg\;min_{K}\sum_{ij}\frac{\partial D\left(P_{ij}\|KQ_{ij}\right)}{\partial K} 
\end{equation}
\subsubsection{* Fundamental property.}
\textbf{The general property of the scale invariant divergences with respect to the $Q$ matrix as a whole is written as follows:}
\begin{equation}
\sum_{nm}Q_{nm}\frac{\partial D\left(P\|KQ\right)}{\partial Q_{_{nm}}}=\sum_{nm}Q_{nm}\frac{\partial DI\left(P\|Q\right)}{\partial Q_{_{nm}}}=0
\label{eq.propgen}	
\end{equation}
Obviously, this is strictly equivalent to the calculus developed in Chapter 3, equation (\ref{eq.Pfond}), as long as the matrices $Q$ and $\frac{\partial DI\left(P\|Q\right)}{\partial Q}$ are written lexicographically.\\

\textbf{Note that with the matrix relationship $Q\equiv HX$, the invariance on $Q$ translates into an invariance on $X$, but not on $H$.}
\subsubsection{- Detail of the calculation.}
The general property of these divergences can be established as follows.\\
We have:
\begin{equation}
	\frac{\partial D\left(P\|KQ\right)}{\partial Q_{_{nm}}}=\sum_{ij}\frac{\partial D\left(P_{ij}\|KQ_{ij}\right)}{\partial \left(KQ_{ij}\right)}\frac{\partial \left(KQ_{ij}\right)}{\partial Q_{nm}}
\end{equation}
but:
\begin{equation}
	\frac{\partial \left(KQ_{ij}\right)}{\partial Q_{nm}}=\frac{\partial \left(KQ_{ij}\right)}{\partial K}\frac{\partial K}{\partial Q_{nm}}+\frac{\partial \left(KQ_{ij}\right)}{\partial Q_{_{ij}}}\frac{\partial Q_{_{ij}}}{\partial Q_{nm}}
\end{equation}
That is:
\begin{equation}
	\frac{\partial \left(KQ_{ij}\right)}{\partial Q_{nm}}=Q_{ij}\frac{\partial K}{\partial Q_{nm}}+K\delta_{in}\delta_{jm}
\end{equation}
which results in:
\begin{equation}
\frac{\partial D\left(P\|KQ\right)}{\partial Q_{_{nm}}}=\sum_{ij}Q_{ij}\frac{\partial D\left(P_{ij}\|KQ_{ij}\right)}{\partial \left(KQ_{ij}\right)}\frac{\partial K}{\partial Q_{nm}}+K \frac{\partial D\left(P_{nm}\|KQ_{nm}\right)}{\partial \left(KQ_{_{nm}}\right)}
\end{equation}
then, with $Q_{ij}=\frac{\partial \left(KQ_{ij}\right)}{\partial K}$, we obtain:
\begin{equation}
\sum_{nm}Q_{nm}\frac{\partial D\left(P\|KQ\right)}{\partial Q_{_{nm}}}=\left[\sum_{ij}\frac{\partial D\left(P_{ij}\|KQ_{ij}\right)}{\partial K}\right]\left[K+\sum_{nm}Q_{nm}\frac{\partial K}{\partial Q_{nm}}\right]	
\end{equation}
The first term of the second member is zero by definition of the nominal invariance factor, otherwise the second term of the second member is zero if the invariance factor satisfies the differential equation, whether or not the invariance factor is the nominal invariance factor.\\
So we have: 
\begin{equation}
\sum_{nm}Q_{nm}\frac{\partial DI\left(P\|Q\right)}{\partial Q_{_{nm}}}=0	
\end{equation}
Therefore, this relation is equivalent to the relation (\ref{eq.Pfond}) which is established in the \textbf{Chapter 3} dealing with scale invariant divergences.

\subsection{The invariance is derived for each column of $Q$.}

This situation is particularly adapted to NMF; it corresponds to the case where a  sum constraint on the component occurs on each column of the $X$ matrix (in the product $\left[HX\right]\equiv Q$).\\
It is quite clear that the multiplication of the column $X_{\textbf{.}j}$ by a constant, results in the multiplication of the homologous column $\left[HX\right]_{\textbf{.}j}\equiv Q_{\textbf{.}j}$ by the same constant.\\
A divergence must therefore be constructed with a scale invariance property when multiplying the columns of $Q$ by various coefficients; consequently, an invariance factor must be calculated for each column of $Q$.\\
By denoting the column index as ``$j$" and the row index as ``$i$", the basic divergence is written as follows:
\begin{equation}
	D\left(P\|Q\right)=\sum_{j}\sum_{i}D\left(P_{ij}\|Q_{ij}\right)=\sum_{j}\sum_{i}D_{ij}
\end{equation}
The nominal invariance factor corresponding to the column ``$m$" of the matrix $Q$ is defined by:
\begin{equation}
	K_{0m}=Arg\;min_{K}\sum_{i}D\left(P_{im}\|KQ_{im}\right)
\end{equation}
When the derivation is possible, the different nominal invariance factors are therefore obtained by calculating for every column ``$m$", $K_{0m}$ such that:
\begin{equation}
	\sum_{i}\frac{\partial D\left(P_{im}\|KQ_{im}\right)}{\partial K}=0
\end{equation}
Obviously, the scalar factor $K_{0m}$ will only contain terms corresponding to the columns labeled ``$m$" in the $P$ and $Q$ matrices.\\
The resulting scale invariant divergence will be written:
\begin{equation}
	DI\left(P\|Q\right)=\sum_{j}\sum_{i}D\left(P_{ij}\|K_{0j}Q_{ij}\right)
\end{equation}

\subsubsection{* Property.}
\textbf{The general property of the scale invariant divergences, when the invariance factor is calculated separately for each column ``$m$" of $Q$ is written:}
\begin{equation}
	\sum_{n}Q_{nm}\frac{\partial DI\left(P\|Q\right)}{\partial Q_{nm}}=0
	\label{eq.propgenbis}
\end{equation}

\subsubsection{- Détail of the calculation.}
The derivation of the gradient with respect to $Q$ is as follows:
\begin{equation}
	\frac{\partial DI\left(P\|Q\right)}{\partial Q_{nm}}=\sum_{j}\sum_{i}\frac{\partial D\left(P_{ij}\|K_{0j}Q_{ij}\right)}{\partial\left(K_{0j}Q_{ij}\right)}\frac{\partial \left(K_{0j}Q_{ij}\right)}{\partial Q_{nm}}
\end{equation}
In this expression, the term $\frac{\partial \left(K_{0j}Q_{ij}\right)}{\partial Q_{nm}}$ is zero, unless $j=m$, so it leaves:
\begin{equation}
	\frac{\partial DI\left(P\|Q\right)}{\partial Q_{nm}}=\sum_{i}\frac{\partial D\left(P_{im}\|K_{0m}Q_{im}\right)}{\partial\left(K_{0m}Q_{im}\right)}\frac{\partial \left(K_{0m}Q_{im}\right)}{\partial Q_{nm}}
	\label{eq.gradnm}
\end{equation}
The second term of the second member is written:
\begin{equation}
\frac{\partial \left(K_{0m}Q_{im}\right)}{\partial Q_{nm}}=\frac{\partial \left(K_{0m}Q_{im}\right)}{\partial K_{0m}}\frac{\partial K_{0m}}{\partial Q_{nm}}+\frac{\partial \left(K_{0m}Q_{im}\right)}{\partial Q_{im}}\frac{\partial Q_{im}}{\partial Q_{nm}}
\end{equation}
That is:
\begin{equation}
\frac{\partial \left(K_{0m}Q_{im}\right)}{\partial Q_{nm}}=Q_{im}	\frac{\partial K_{0m}}{\partial Q_{nm}}+K_{0m}\delta_{in}
\end{equation}
Coming back then to the expression of the gradient (\ref{eq.gradnm}), we obtain:
\begin{equation}
	\frac{\partial DI\left(P\|Q\right)}{\partial Q_{nm}}=\sum_{i}Q_{im}\frac{\partial D\left(P_{im}\|K_{0m}Q_{im}\right)}{\partial\left(K_{0m}Q_{im}\right)}\frac{\partial K_{0m}}{\partial Q_{nm}}+K_{0m}\frac{\partial D\left(P_{nm}\|K_{0m}Q_{nm}\right)}{\partial\left(K_{0m}Q_{nm}\right)}
\end{equation}
And finally:
\begin{equation}
	\sum_{n}Q_{nm}\frac{\partial DI\left(P\|Q\right)}{\partial Q_{nm}}=\left[\sum_{n}\frac{\partial D\left(P_{nm}\|K_{0m}Q_{nm}\right)}{\partial\left(K_{0m}Q_{nm}\right)}\right]\left[K_{0m}+\sum_{n}Q_{nm}\frac{\partial K_{0m}}{\partial Q_{nm}}\right]
\end{equation}
The first term of the second member is zero if the invariance factor is the nominal invariance factor for the column under consideration (``$m$"), while the second term of the second member is zero for all invariance factors satisfying the differential equation:
\begin{equation}
K_{0m}+\sum_{n}Q_{nm}\frac{\partial K_{0m}}{\partial Q_{nm}}=0	
\end{equation}
Consequently:
\begin{equation}
	\sum_{n}Q_{nm}\frac{\partial DI\left(P\|Q\right)}{\partial Q_{nm}}=0\ \forall m
	\label{eq.propgenter}
\end{equation}
This relation is similar to the property (\ref{eq.Pfond}) applied on each column of the matrix $Q$.

\section{Iterative algorithm for invariant divergences.}
For these divergences denoted $DI$, by operating column by column, we rely on the property (\ref{eq.propgenbis}) (\ref{eq.propgenter}).

\subsection{Iterations on $H$.}

\textbf{Since the invariance on $Q$ translates into an invariance on $X$, but not on $H$, there is no particular thing to expect regarding iterations on $H$.}

Consequently, the iterative algorithm on $H$ will be deduced from the variable change method; it is given in (\ref{eq.algonmH}) (\ref{eq.algonmHbis}) and is written for an invariant divergence $DI$, with $H=H^{k},\ X=X^{k},\ Q=H^{k}X^{k}$:
\begin{align}
		H^{k+1}_{nm}=H^{k}_{nm}+\delta^{k}_{m}H^{k}_{nm}&\left[\left(\sum_{l}\left[-\frac{\partial DI}{\partial Q_{nl}}\right]X^{T}_{lm}\right)\right. \nonumber \\  & \left.-\sum_{i}H^{k}_{im}\left(\sum_{l}\left[-\frac{\partial DI}{\partial Q_{il}}\right]X^{T}_{lm}\right)\right]
	\label{eq.algonmHinv}
\end{align}
\textbf{A multiplicative algorithm can be obtained by proceeding in an similar way to what has been previously developed for non-invariant divergences.}\\

At this point, we have obtained $H^{k+1}$.
\subsection{Iterations on $X$.}
\subsubsection{1 - Non-multiplicative algorithm.}
The iterative algorithm on $X$ allowing to exploit the properties of invariant divergences, can be written, with the method developed in \textbf{Chapter 10}, as follows:
\begin{equation}
	X^{k+1}_{nm}=X^{k}_{nm}+\alpha^{k}_{m}X^{k}_{nm}\left[-\frac{\partial DI}{\partial X^{k}_{nm}}\right]
	\label{eq.algoXI}
\end{equation}
With this algorithm, the constraint $\sum_{n}X^{k}_{nm}=Cte\;\forall k$ is ensured as shown in the following calculation.\\
In the previous sections, we have shown (\ref{eq.matA}) that with:
\begin{equation}
	\left[A\right]_{ij}=\left[\frac{\partial D}{\partial Q}\right]_{ij}	
\end{equation}
One could write (\ref{eq.gradX}) in matrix form:
\begin{equation}
\left[\frac{\partial DI\left(P\|Q\right)}{\partial X}\right]_{nm}=\frac{\partial DI\left(P\|Q\right)}{\partial X_{nm}}=\left[H^{T}A\right]_{nm}=\sum_{l}H^{T}_{nl}\frac{\partial DI\left(P\|Q\right)}{\partial Q_{lm}}	
\end{equation}
Taking into account the fact that $H\equiv H^{k+1}$ and $Q=H^{k+1}X^{k}$, the algorithm (\ref{eq.algoXI}) can therefore be rewritten:
\begin{equation}
	X^{k+1}_{nm}=X^{k}_{nm}+\alpha^{k}_{m}X^{k}_{nm}\left[\sum_{i}H^{T}_{ni}\left(-\frac{\partial DI\left(P\|Q\right)}{\partial Q_{im}}\right)\right]
	\label{eq.algoXIbis}
\end{equation}
Then:
\begin{equation}
	\sum_{n}X^{k+1}_{nm}=\sum_{n}X^{k}_{nm}+\alpha^{k}_{m}\sum_{n}\left\{X^{k}_{nm}\left[\sum_{i}H^{T}_{ni}\left(-\frac{\partial DI\left(P\|Q\right)}{\partial Q_{im}}\right)\right]\right\}
\end{equation}
By reversing the order of the summations, it comes:
\begin{equation}
	\sum_{n}X^{k+1}_{nm}=\sum_{n}X^{k}_{nm}+\alpha^{k}_{m}\sum_{i}\left(-\frac{\partial DI\left(P\|Q\right)}{\partial Q_{im}}\right)\sum_{n}H^{T}_{ni}X^{k}_{nm}
\end{equation}
or also:
\begin{equation}
	\sum_{n}X^{k+1}_{nm}=\sum_{n}X^{k}_{nm}+\alpha^{k}_{m}\sum_{i}\left(HX^{k}\right)_{im}\left(-\frac{\partial DI\left(P\|Q\right)}{\partial Q_{im}}\right)
\end{equation}
So, with (\ref{eq.propgenbis}), it comes:
\begin{equation}
	\sum_{n}X^{k+1}_{nm}=\sum_{n}X^{k}_{nm}
\end{equation}
It is therefore sufficient to choose an initial estimate with the correct sum (i.e. $\sum_{n}X^{0}_{nm}=\sum_{n}Y_{nm}$), to ensure that this sum is maintained during the iterations.\\

\textbf{The algorithm (\ref{eq.algoXI}) (\ref{eq.algoXIbis}) represents the iterative algorithm on $X$ when minimizing an invariant divergence.}

\subsubsection{2 - Multiplicative algorithm.}
A multiplicative algorithm can be obtained from (\ref{eq.algoXI}) (\ref{eq.algoXIbis}) by operating as follows:\\
The SGM method is strictly applied; the opposite of the gradient is first split into a difference of 2 positive terms:
\begin{equation}
	-\frac{\partial DI}{\partial X^{k}_{nm}}=U^{k}_{nm}-V^{k}_{nm}\ ;\ \ U^{k}_{nm}>0\ ;\ \ V^{k}_{nm}>0
\end{equation}
This allows in a first step (with the usual restrictions on the descent step size), to write a multiplicative algorithm in the form:
\begin{equation}
	\widetilde{X}^{k+1}_{nm}=X^{k}_{nm}\left[\frac{U^{k}_{nm}}{V^{k}_{nm}}\right]
	\label{eq.algomultXinv}
\end{equation}
But then the property $\sum_{n}\widetilde{X}^{k+1}_{nm}=\sum_{n}X^{k}_{nm}$ is lost, unlike what we had with algorithms (\ref{eq.algoXI}) (\ref{eq.algoXIbis}).\\
However, we can recover this property by introducing an extra step of normalization by writing:
\begin{equation}
	X^{k+1}_{nm}=\frac{\widetilde{X}^{k+1}_{nm}}{\sum_{n}\widetilde{X}^{k+1}_{nm}}\sum_{n}Y_{nm}
	\label{eq.normalisation}
\end{equation}
Due to the property of invariance with respect to $X$, this normalization does not lead to any change in the invariant divergence considered.\\
The association of (\ref{eq.algomultXinv}) and (\ref{eq.normalisation}) is thus a multiplicative algorithm having the summation property on the columns of $X$. This is subject of course to the convergence of (\ref{eq.algomultXinv}) in relation to the choice of the descent step size (equal to 1).

\section{Applications to some particular invariant divergences.}
In the previous sections, it was shown that the proposed algorithms used, for a given divergence, the gradient expression:
\begin{equation}
	\frac{\partial 	D\left(P\|Q\right)}{\partial Q_{nm}}
\end{equation}
Consequently, we give for some invariant divergences, the expression of this gradient and the iterative algorithm on $X$ that we can obtain.

\subsection{Invariant ``Alpha'' divergences.}
From the basic divergence (\ref{eq.ACbis}), the invariance on each column of the table (matrix) $Q$ is obtained by deriving the nominal invariance factor (which is possible); this gives, for the column ``$j$":
\begin{equation}
	K_{0j}=\left(\frac{\sum_{i}P^{\lambda}_{ij}Q^{1-\lambda}_{ij}}{\sum_{i}Q_{ij}}\right)^{\frac{1}{\lambda}}
	\label{eq.K0Aj}
\end{equation}
The divergence we obtain, invariant with respect to the columns of $Q$, is written as follows:
\begin{equation}
	AI_{Q}\left(P\|Q\right)=\frac{1}{\lambda-1}\sum_{j}\left[\left(\frac{\sum_{l}P^{\lambda}_{lj}Q^{1-\lambda}_{lj}}{\sum_{l}Q_{lj}}\right)^{\frac{1}{\lambda}}\sum_{i}Q_{ij}-\sum_{i}P_{ij}\right]
	\label{eq.AIQ}
\end{equation}
That is:
\begin{equation}
	AI_{Q}\left(P\|Q\right)=\frac{1}{\lambda-1}\sum_{j}\left[K_{0j}\sum_{i}Q_{ij}-\sum_{i}P_{ij}\right]
	\label{eq.AIQs}
\end{equation}
The construction of the algorithms implies the calculation of the gradient with respect to the elements of the column ``$m$" of the matrix $Q$ which is written after some calculations:
\begin{align}
	\frac{\partial 	AI_{Q}\left(P\|Q\right)}{\partial Q_{nm}}=\frac{1}{\lambda}& \left[\left(\frac{\sum_{l}P^{\lambda}_{lm}Q^{1-\lambda}_{lm}}{\sum_{l}Q_{lm}}\right)^{\frac{1}{\lambda}}\right. \nonumber \\  & \left.-\left(\frac{\sum_{l}P^{\lambda}_{lm}Q^{1-\lambda}_{lm}}{\sum_{i}Q_{lm}}\right)^{\frac{1-\lambda}{\lambda}}\frac{P^{\lambda}_{nm}}{Q^{\lambda}_{nm}}\right]
	\label{eq.gradAI}
\end{align}
That is:
\begin{equation}
	\frac{\partial 	AI_{Q}\left(P\|Q\right)}{\partial Q_{nm}}=\frac{1}{\lambda}\left[K_{0m}-K_{0m}^{1-\lambda}\frac{P^{\lambda}_{nm}}{Q^{\lambda}_{nm}}\right]
	\label{eq.gradAIQ}
\end{equation}
We can verify that the fundamental relation (\ref{eq.propgenbis}) for scale invariant divergences with respect to the columns of $Q$ is satisfied and that we have:
\begin{equation}
	\sum_{n}Q_{nm}\frac{\partial 	AI_{Q}\left(P\|Q\right)}{\partial Q_{nm}}=0
\end{equation}
In the particular case ($\lambda=1$) which corresponds to the Kullback-Leibler divergence, we obtain from the expression of the gradient (\ref{eq.gradAIQ}), the gradient of the K.L. divergence , invariant with respect to the columns of $Q$; such divergence is denoted $KLI_{Q}$:
\begin{equation}
	\frac{\partial 	KLI_{Q}\left(P\|Q\right)}{\partial Q_{nm}}=\frac{\sum_{l}P_{lm}}{\sum_{l}Q_{lm}}-\frac{P_{nm}}{Q_{nm}}
	\label{eq.gradKLIQ}
\end{equation}

\textbf{* Itérations on $H$.}\\

With (\ref{eq.gradAIQ}), it comes:
\begin{equation}
	\frac{\partial AI_{Q}\left(P\|Q\right)}{\partial H_{nm}}=\sum_{l}\frac{\partial AI_{Q}\left(P\|Q\right)}{\partial Q_{nl}}X^{T}_{lm}=\frac{1}{\lambda}\sum_{l}\left[K_{0l}-K_{0l}^{1-\lambda}\frac{P^{\lambda}_{nl}}{Q^{\lambda}_{nl}}\right]X^{T}_{lm}
\end{equation}
With $X=X^{k},\ H=H^{k},\ Q=H^{k}X^{k}$ and taking into account the fact that the invariance factor $K_{0l}$ given by (\ref{eq.K0Aj}) varies with the iteration because it depends on $Q$, the iterative algorithm on $H$ deduced from (\ref{eq.algonmHinv}) is then written after simplification:
\begin{equation}
	H^{k+1}_{nm}=H^{k}_{nm}+\frac{\delta^{k}_{m}}{\lambda}H^{k}_{nm}\left[\sum_{l}K_{0l}^{1-\lambda}\left(\frac{P^{\lambda}_{nl}}{Q^{\lambda}_{nl}}-\sum_{i}H^{k}_{im}\frac{P^{\lambda}_{il}}{Q^{\lambda}_{il}}\right)X^{T}_{lm}\right]
\end{equation}
\textbf{* Itérations on $X$.}\\

At this point in the calculation, we have $H\equiv H^{k+1}$, $X\equiv X^{k}$, $Q\equiv H^{k+1}X^{k}$.\\
It has been shown that:
\begin{equation}
\frac{\partial 	AI_{Q}\left(P\|Q\right)}{\partial X_{nm}}= \sum_{j}H^{T}_{nj}\frac{\partial AI_{Q}\left(P\|Q\right)}{\partial Q_{jm}}
\end{equation}
That is, with (\ref{eq.gradAIQ}):
\begin{equation}
\frac{\partial 	AI_{Q}\left(P\|Q\right)}{\partial X_{nm}}=\frac{1}{\lambda}\sum_{j}H^{T}_{nj}\left[K_{0m}-K_{0m}^{1-\lambda}\frac{P^{\lambda}_{jm}}{Q^{\lambda}_{jm}}\right]
\label{eq.gradAIX}
\end{equation}
Considering that $\sum_{j}H^{T}_{nj}=1$, the iterative algorithm on $X$ given by (\ref{eq.algoXI}) is written:
\begin{equation}
	X^{k+1}_{nm}=X^{k}_{nm}+\frac{\alpha^{k}_{m}}{\lambda}X^{k}_{nm}\left[\left(K^{k}_{0m}\right)^{1-\lambda}\sum_{j}H^{T}_{nj}\left(\frac{P^{\lambda}_{jm}}{Q^{\lambda}_{jm}}\right)-K^{k}_{0m}\right]
\end{equation}
In this expression, we have taken into account the fact that the invariance factor depends on iteration, indeed, we recall that:
\begin{equation}
	K_{0m}=\left(\frac{\sum_{i}P^{\lambda}_{im}Q^{1-\lambda}_{im}}{\sum_{i}Q_{im}}\right)^{\frac{1}{\lambda}}
\end{equation}
with $Q\equiv H^{k+1}X^{k}$; to avoid any ambiguity, we have denoted his expression: $K^{k}_{0m}$.\\ 
Since $K^{k}_{0m}>$0 and considering the very particular expression of the gradient, this also leads to the expression which is at the basis of the multiplicative algorithm:
\begin{equation}
	X^{k+1}_{nm}=X^{k}_{nm}+\frac{K^{k}_{0m}\alpha^{k}_{m}}{\lambda}X^{k}_{nm}\left[\left(K^{k}_{0m}\right)^{-\lambda}\sum_{j}H^{T}_{nj}\left(\frac{P^{\lambda}_{jm}}{Q^{\lambda}_{jm}}\right)-1\right]
	\label{eq.AlgoAIX}
\end{equation}
By taking a descent step size $\frac{K^{k}_{0m}\alpha^{k}_{m}}{\lambda}$ equal to 1, we obtain the purely multiplicative form:
\begin{equation}
	X^{k+1}_{nm}=X^{k}_{nm}\left(K^{k}_{0m}\right)^{-\lambda}\sum_{j}H^{T}_{nj}\left(\frac{P^{\lambda}_{jm}}{Q^{\lambda}_{jm}}\right)
\end{equation}
Using the expression of $K^{k}_{0m}$, it can be shown after some simple calculations that $\sum_{n}X^{k+1}_{nm}=\sum_{n}X^{k}_{nm}$.
\subsubsection{Logarithmic form.}

As proposed for Alpha divergences, we can give a logarithmic form for the invariant divergence $AI_{Q}$, it is written as follows:
\begin{align}
	LAI_{Q}\left(P\|Q\right)=\frac{1}{\lambda-1}&\left\{\log\left[\sum_{j}\sum_{i}\left(\frac{\sum_{l}P^{\lambda}_{lj}Q^{1-\lambda}_{lj}}{\sum_{l}Q_{lj}}\right)^{\frac{1}{\lambda}}Q_{ij}\right]\right. \nonumber \\  & \left.-\log\left[\sum_{j}\sum_{i}P_{ij}\right]\right\}
	\label{eq.LAIQ}
\end{align}
It can be observed that, as with all logarithmic forms, this divergence is invariant not only with respect to the columns of $Q$, but also with respect to the columns of $P$.\\
The gradient with respect to $Q$ is written as:
\begin{align}
	\frac{\partial 	LAI_{Q}\left(P\|Q\right)}{\partial Q_{nm}}=\frac{A}{\lambda}&\left[\left(\frac{\sum_{l}P^{\lambda}_{lm}Q^{1-\lambda}_{lm}}{\sum_{l}Q_{lm}}\right)^{\frac{1}{\lambda}}\right. \nonumber \\  & \left.-\left(\frac{\sum_{l}P^{\lambda}_{lm}Q^{1-\lambda}_{lm}}{\sum_{i}Q_{lm}}\right)^{\frac{1-\lambda}{\lambda}}\frac{P^{\lambda}_{nm}}{Q^{\lambda}_{nm}}\right]
	\label{eq.gradLAI}
\end{align}
That is:
\begin{equation}
\frac{\partial 	LAI_{Q}\left(P\|Q\right)}{\partial Q_{nm}}=A\frac{\partial 	AI_{Q}\left(P\|Q\right)}{\partial Q_{nm}}	
\end{equation}
with:
\begin{equation}
	A=\frac{1}{\sum_{i}\sum_{j}\left(\frac{\sum_{l}P^{\lambda}_{lj}Q^{1-\lambda}_{lj}}{\sum_{l}Q_{lj}}\right)^{\frac{1}{\lambda}}Q_{ij}}=\frac{1}{\sum_{i}\sum_{j}K_{0j}Q_{ij}}
\end{equation}
The gradient with respect to $Q$ can then be written in a simpler way:
\begin{equation}
	\frac{\partial 	LAI_{Q}\left(P\|Q\right)}{\partial Q_{nm}}=A\frac{K_{0m}}{\lambda}\left(1-K^{-\lambda}_{0m}\frac{P^{\lambda}_{nm}}{Q^{\lambda}_{nm}}\right)
\end{equation}
We then deduce the expression of the gradient with respect to $X$:
\begin{equation}
\frac{\partial 	LAI_{Q}\left(P\|Q\right)}{\partial X_{nm}}=\sum_{j}H^{T}_{nj}\frac{\partial LAI_{Q}\left(P\|Q\right)}{\partial Q_{jm}}	
\end{equation}
That is:
\begin{equation}
\frac{\partial 	LAI_{Q}\left(P\|Q\right)}{\partial X_{nm}}=A\frac{K_{0m}}{\lambda}\left(1-K^{-\lambda}_{0m}\sum_{j}H^{T}_{nj}\frac{P^{\lambda}_{jm}}{Q^{\lambda}_{jm}}\right)	
\end{equation}
Hence the iterative algorithm on $X$:
\begin{equation}
	X^{k+1}_{nm}=X^{k}_{nm}+A\frac{K^{k}_{0m}\alpha^{k}_{m}}{\lambda}X^{k}_{nm}\left[\left(K^{k}_{0m}\right)^{-\lambda}\sum_{j}H^{T}_{nj}\left(\frac{P^{\lambda}_{jm}}{Q^{\lambda}_{jm}}\right)-1\right]
	\label{eq.AlgoLAIX}
\end{equation}
As with the $AI(P\|Q)$ divergence, we will take into account the fact that $H\equiv H^{k+1}$, $Q=H^{k+1}X^{k}$, and therefore $K_{0m}\equiv K^{k}_{0m}$ depends on the iteration.\\
Note that this algorithm is analogous to (\ref{eq.AlgoAIX}), except for the multiplicative factor ``$A$" that occurs in the correction term and can be included in the descent step size; therefore, the sum constraint (flux-holding) property highlighted in (\ref{eq.AlgoAIX}) is maintained.\\
Due to the particular expression of the Gradient, a multiplicative flux-holding algorithm is obtained as for the $AI(P\|Q)$ divergence.

\subsection{Invariant ``Beta'' divergences .}
The basis divergence (\ref{eq.BC}) applied to tables is written as follows:
\begin{equation}
	B\left(p\|q\right)=\frac{1}{\lambda\left(\lambda-1\right)}\sum_{j}\sum_{i}\left[P_{ij}^{\lambda}-\lambda P_{ij}Q_{ij}^{\lambda-1}-\left(1-\lambda\right)Q_{ij}^{\lambda}\right]
\end{equation}
The invariance on each column of the table (matrix) $Q$ is obtained by calculating the nominal invariance factor (which is possible); thus, for the column ``$j$" we obtain:
\begin{equation}
	K_{0j}=\frac{\sum_{i}P_{ij}Q^{\lambda-1}_{ij}}{\sum_{i}Q^{\lambda}_{ij}}
	\label{eq.K0Bj}
\end{equation}
The divergence generated, invariant with respect to the columns of $Q$, is written as follows:
\begin{equation}
	BI_{Q}\left(P\|Q\right)=\frac{1}{\lambda\left(\lambda-1\right)}\sum_{j}\left[\sum_{i}P_{ij}^{\lambda}-\left(\frac{\sum_{l}P_{lj}Q^{\lambda-1}_{lj}}{\sum_{l}Q^{\lambda}_{lj}}\right)^{\lambda}\sum_{i}Q_{ij}^{\lambda}\right]
	\label{eq.BIQ}
\end{equation}
Or also:
\begin{equation}
	BI_{Q}\left(P\|Q\right)=\frac{1}{\lambda\left(\lambda-1\right)}\sum_{j}\left[\sum_{i}P_{ij}^{\lambda}-K^{\lambda}_{0j}\sum_{i}Q_{ij}^{\lambda}\right]
	\label{eq.BIQs}
\end{equation}
The calculation of the gradient of this divergence with respect to $Q$ leads to:
\begin{align}
	\frac{\partial BI_{Q}\left(P\|Q\right)}{\partial Q_{nm}}=&\left[\left(\frac{\sum_{l}P_{lm}Q^{\lambda-1}_{lm}}{\sum_{l}Q^{\lambda}_{lm}}\right)^{\lambda}Q^{\lambda-1}_{nm}\right. \nonumber \\  & \left.-\left(\frac{\sum_{l}P_{lm}Q^{\lambda-1}_{lm}}{\sum_{l}Q^{\lambda}_{lm}}\right)^{\lambda-1}P_{nm}Q^{\lambda-2}_{nm}\right]
	\label{eq.gradBIQ}
\end{align}
That is to say more simply:
\begin{equation}
	\frac{\partial BI_{Q}\left(P\|Q\right)}{\partial Q_{nm}}=\left[K^{\lambda}_{0m}Q^{\lambda-1}_{nm}-K^{\lambda-1}_{0m}P_{nm}Q^{\lambda-2}_{nm}\right]
	\label{eq.gradBIQbis}
\end{equation}
One can verify that the fundamental relationship for divergences invariant with respect to the columns of $Q$ is satisfied.\\
The case ($\lambda=1$) corresponding to the Kullback-Leibler divergence is obtained without difficulty; its gradient has been given by (\ref{eq.gradKLIQ}), whereas in the particular case ($\lambda=2$) corresponding to the mean square deviation $MCI_{Q}$, the gradient obtained from (\ref{eq.gradBIQ}) is written as follows:
 \begin{equation}
	\frac{\partial MCI_{Q}\left(P\|Q\right)}{\partial Q_{nm}}=\left(\frac{\sum_{l}P_{lm}Q_{lm}}{\sum_{l}Q^{2}_{lm}}\right)^{2}Q_{nm}-\left(\frac{\sum_{l}P_{lm}Q_{lm}}{\sum_{l}Q^{2}_{lm}}\right)P_{nm}
	\label{eq.gradMCIQ}
\end{equation}\\
\textbf{* Itérations on $H$.}\\

With (\ref{eq.gradBIQbis}), we have:
\begin{equation}
	\frac{\partial BI_{Q}\left(P\|Q\right)}{\partial H_{nm}}=\sum_{l}\frac{\partial BI_{Q}\left(P\|Q\right)}{\partial Q_{nl}}X^{T}_{lm}=\sum_{l}\left[K^{\lambda}_{0l}Q^{\lambda-1}_{nl}-K^{\lambda-1}_{0l}P_{nl}Q^{\lambda-2}_{nl}\right]X^{T}_{lm}	
\end{equation}
With $X=X^{k},\ H=H^{k},\ Q=H^{k}X^{k}$ and taking into account that the invariance factor $K_{0l}$ given by (\ref{eq.K0Bj}) varies with the iteration via $Q$, the iterative algorithm on $H$ deduced from (\ref{eq.algonmHinv}) is written as follows:
\begin{align}
	H^{k+1}_{nm}=H^{k}_{nm}+\delta^{k}_{m}H^{k}_{nm}&\left[\sum_{l}K_{0l}^{\lambda-1}\left(P_{nl}Q^{\lambda-2}_{nl}-\sum_{i}H^{k}_{im}P_{il}Q^{\lambda-2}_{il}\right)X^{T}_{lm}\right. \nonumber \\  & \left.-\sum_{l}K_{0l}^{\lambda}\left(\sum_{i}H^{k}_{im}Q^{\lambda-1}_{il}-Q^{\lambda-1}_{nl}\right)X^{T}_{lm}\right]
\end{align}
\textbf{* Itérations on $X$.}\\

At this point of the calculus, we have $H\equiv H^{k+1}$, $X\equiv X^{k}$, $Q\equiv H^{k+1}X^{k}$.\\
Taking into account that:
\begin{equation}
\frac{\partial 	BI_{Q}\left(P\|Q\right)}{\partial X_{nm}}=\sum_{j}H^{T}_{nj}\frac{\partial BI_{Q}\left(P\|Q\right)}{\partial Q_{jm}}
\end{equation}
It comes from the expression (\ref{eq.gradBIQbis}):
\begin{equation}
\frac{\partial 	BI_{Q}\left(P\|Q\right)}{\partial X_{nm}}=\sum_{j}H^{T}_{nj}\left[K^{\lambda}_{0m}Q^{\lambda-1}_{jm}-K^{\lambda-1}_{0m}P_{jm}Q^{\lambda-2}_{jm}\right]
\label{eq.gradBIX}
\end{equation}
To take into account the fact that the invariance factor varies with iteration, we will denote $K^{k}_{0m}$; it is given by the expression (\ref{eq.K0Bj}), where $Q\equiv H^{k+1}X^{k}$.\\
The iterative algorithm on $X$ given by (\ref{eq.algoXI}) is then written as follows:
\begin{equation}
	X^{k+1}_{nm}=X^{k}_{nm}+\alpha^{k}_{m}\left(K^{k}_{0m}\right)^{\lambda-1}X^{k}_{nm}\left[\sum_{j}H^{T}_{nj}\left(P_{jm}Q^{\lambda-2}_{jm}-K^{k}_{0m}Q^{\lambda-1}_{jm}\right)\right]
	\label{eq.AlgoBIX}
\end{equation}
One can check that with this algorithm, the property $\sum_{n}X^{k+1}_{nm}=\sum_{n}X^{k}_{nm}$ is ensured, however, a multiplicative form deduced from this expression will not have this property any more.\\
To obtain a multiplicative form having the property of maintaining the flow, it will be necessary to proceed in 2 steps as indicated in \textbf{(11.6.2)}.

\subsubsection{Logarithmic form.}
As proposed for the Alpha divergences, we can give a logarithmic form for the $BI_{Q}$ divergence; it is written as follows:
\begin{align}
	LBI_{Q}\left(P\|Q\right)=\frac{1}{\lambda\left(\lambda-1\right)}&\left\{\log\left[\sum_{j}\sum_{i}P_{ij}^{\lambda}\right]\right. \nonumber \\  & \left.-\log\left[\sum_{j}\sum_{i}\left(\frac{\sum_{l}P_{lj}Q^{\lambda-1}_{lj}}{\sum_{l}Q^{\lambda}_{lj}}\right)^{\lambda}Q_{ij}^{\lambda}\right]\right\}
	\label{eq.LBIQ}
\end{align}
It can be observed that, as with all logarithmic forms, this divergence is invariant not only with respect to columns of $Q$, but also with respect to columns of $P$.\\
The gradient with respect to $Q$ is written as follows:
\begin{align}
	\frac{\partial LBI_{Q}\left(P\|Q\right)}{\partial Q_{nm}}=B & \left[\left(\frac{\sum_{l}P_{lm}Q^{\lambda-1}_{lm}}{\sum_{l}Q^{\lambda}_{lm}}\right)^{\lambda}Q^{\lambda-1}_{nm}\right. \nonumber \\  & \left.-\left(\frac{\sum_{l}P_{lm}Q^{\lambda-1}_{lm}}{\sum_{l}Q^{\lambda}_{lm}}\right)^{\lambda-1}P_{nm}Q^{\lambda-2}_{nm}\right]
	\label{eq.gradLBIQ}
\end{align}
That is:
\begin{equation}
\frac{\partial LBI_{Q}\left(P\|Q\right)}{\partial Q_{nm}}=B\ \frac{\partial BI_{Q}\left(P\|Q\right)}{\partial Q_{nm}}
\end{equation}
with:
\begin{equation}
	B=\frac{1}{\sum_{j}\sum_{i}Q_{ij}^{\lambda}\left(\frac{\sum_{l}P_{lj}Q^{\lambda-1}_{lj}}{\sum_{l}Q^{\lambda}_{lj}}\right)^{\lambda}}=\frac{1}{\sum_{i}\sum_{j}K^{\lambda}_{0j}Q_{ij}^{\lambda}}
\end{equation}
The gradient with respect to $Q$ is written in a simpler way:
\begin{equation}
	\frac{\partial 	LBI_{Q}\left(P\|Q\right)}{\partial Q_{nm}}=BK^{\lambda-1}_{0m}\left(K_{0m}Q^{\lambda-1}_{nm}-P_{nm}Q^{\lambda-2}_{nm}\right)
\end{equation}
We then deduce the expression of the gradient with respect to $X$:
\begin{equation}
\frac{\partial 	LBI_{Q}\left(P\|Q\right)}{\partial X_{nm}}=\sum_{j}H^{T}_{nj}\frac{\partial LBI_{Q}\left(P\|Q\right)}{\partial Q_{jm}}	
\end{equation}
That is:
\begin{equation}
\frac{\partial 	LBI_{Q}\left(P\|Q\right)}{\partial X_{nm}}=BK^{\lambda-1}_{0m}\left[\sum_{j}H^{T}_{nj}\left(K_{0m}Q^{\lambda-1}_{jm}-P_{jm}Q^{\lambda-2}_{jm}\right)\right]	
\end{equation}
Hence the iterative algorithm on $X$:
\begin{equation}
	X^{k+1}_{nm}=X^{k}_{nm}+\alpha^{k}_{m}BK^{\lambda-1}_{0m}X^{k}_{nm}\left[\sum_{j}H^{T}_{nj}\left(P_{jm}Q^{\lambda-2}_{jm}-K^{k}_{0m}Q^{\lambda-1}_{jm}\right)\right]
	\label{eq.AlgoLBIX}
\end{equation}
As with the $BI(P\|Q)$ divergence, we will take into account the fact that $H\equiv H^{k+1}$, $Q=H^{k+1}X^{k}$ and that, consequently, $K_{0m}\equiv K^{k}_{0m}$ is dependent on iteration.\\
Note that this algorithm is analogous to (\ref{eq.AlgoBIX}), except for the multiplicative factor ``$B$" which is involved in the corrective term and can be included in the descent step size.\\
We can verify that with this algorithm, the property $\sum_{n}X^{k+1}_{nm}=\sum_{n}X^{k}_{nm}$ is ensured, however, a multiplicative form deduced from this expression will no longer possess this property .\\
To obtain a multiplicative form having the property of flow holding, it will be necessary to proceed in 2 steps as indicated in \textbf{(11.6.2)}.

\section{Régularisation. Non-invariant divergences }
In the context of NMF, a smoothness constraint regularization can only be performed on the columns of the $H$ matrix, in fact, each column represents the spectrum of an elementary component.\\
With respect to the $X$ matrix, whose columns contain weighting factors, such a smoothness constraint regularization has no meaning.\\
On the other hand, a much more logical form of regularization is to introduce a certain degree of ``\textit{sparsity}" on the elements in the columns of $X$.\\
In any case, solving the regularized problem requires minimizing with respect to $H$ and $X$, a functional of the form:
\begin{equation}
	J\left(H,X\right)=D\left(Y\|HX\right)+\gamma DH\left(H\right)+\mu DX\left(X\right)
\end{equation}
The minimization is carried out under the non-negativity and sum constraints already specified.\\
The ``$\gamma$" and ``$\mu$" coefficients are the positive regularization factors.\\
Since the iterative process of minimization is performed alternately on $H$ then on $X$ ( for example), one works first concerning $H$ on a composite divergence:
\begin{equation}
DCH=D\left(Y\|HX\right)+\gamma DH\left(H\right)	
\end{equation}
then, concerning $X$, on a composite divergence: 
\begin{equation}
DCX=D\left(Y\|HX\right)+\mu DX\left(X\right)	
\end{equation} 

\subsection{Regularization of the columns of the matrix $H$.}

The regularization is performed separately on each of the columns of $H$.\\
The basic algorithm is given by (\ref{eq.algonmH}) applied to the composite divergence $DCH$; it is written:
\begin{align}
	H^{k+1}_{nm}=H^{k}_{nm}+\delta^{k}_{m}H^{k}_{nm}&\left\{-\frac{\partial \left(DCH\right)}{\partial H^{k}_{nm}}\right. \nonumber \\  & \left.+\sum_{i}H^{k}_{im}\frac{\partial \left(DCH\right)}{\partial H^{k}_{im}}\right\}
	\label{eq.algoHreg}
\end{align}
We are now considering 2 particular cases of regularization by smoothness constraint:

\subsubsection{1 - Euclidean norm of the solution.}
\begin{equation}
	DH\left(H\right)=DH\left(H\|C\right)=\frac{1}{2}\sum_{j}\sum_{i}\left(H_{ij}-C\right)^{2}
\label{eq.regticste}
\end{equation}
For columns of $H$ with ``N" components, taking into account the constraint $\sum_{i}H_{ij}=1\ \forall j$, we will have $C=\frac{1}{N}$.\\
\begin{equation}
\frac{\partial DH}{\partial H_{nm}}=\left(H_{nm}-C\right)
\label{eq.gradregticste}	
\end{equation}

\subsubsection{2 - Euclidean norm of the Laplacian of the solution.}
One can as well use a regularization term expressed by the Euclidean norm of the Laplacian, which will be expressed by:
\begin{equation}
	DH\left(H\right)=DH\left(H\|TH\right)=\frac{1}{2}\sum_{j}\sum_{i}\left[H_{ij}-\left(TH_{.j}\right)_{i}\right]^{2}
\end{equation}
Where the matrix operation $\left(TH_{.j}\right)$ corresponds to the convolution of the column ``$j$" of $H$ by the mask $\left[0.5\ ; \  0\ ; \ 0.5\right]$.\\
The regularization is performed independently for each column of $H$.\\
In this case, we have:
\begin{equation}
	\frac{\partial DH}{\partial H_{nm}}=H_{nm}-\left[\left(T+T^{T}\right)H_{.m}\right]_{n}+\left[\left(T^{T}T\right)H_{.m}\right]_{n}
\end{equation}
Of course, taking into account the symmetry of the convolution mask, $T^{T}=T$.

\subsection{Regularization of the columns of the matrix $X$.}
The regularized algorithm is developed on the basis of the algorithm (\ref{eq.algoBbis}), taking into account the fact that the divergence considered is the composite divergence $DCX$;  
\begin{align}
	X^{k+1}_{nm}=X^{k}_{nm}+\delta^{k}_{m}X^{k}_{nm}&\left[\left(\sum_{l}Y_{lm}\right)\left(-\frac{\partial DCX}{\partial X^{k}_{nm}}\right)\right. \nonumber \\  & \left.-\sum_{i}X^{k}_{im}\left(-\frac{\partial DCX}{\partial X^{k}_{im}}\right)\right]
	\label{eq.algoBregbase}
\end{align}
\subsubsection{* Expression of the penalty term.}

Here again, the regularization applies to each column of the $X$ matrix, but, contrary to the $H$ matrix, the regularization does not relate to the smoothness of the solution, but to its ``\textit{sparsity}".\\
The measure of ``\textit{sparsity}" used here was proposed by Hoyer \cite{hoyer2004}; by denoting $X_{j}$ the column ``j" of $X$, if ``$N$" is the number of elements of $X_{j}$, the sparsity factor for this column is expressed by:
\begin{equation}
	s_{j}=\frac{\sqrt{N}-\frac{\left\|X_{j}\right\|_{1}}{\left\|X_{j}\right\|_{2}}}{\sqrt{N}-1}
\end{equation}
This relationship induces a relationship between $\left\|X_{j}\right\|_{1}$ and $\left\|X_{j}\right\|_{2}$ which we will express in the form:
\begin{equation}
	\frac{\left\|X_{j}\right\|^{2}_{2}}{\left\|X_{j}\right\|^{2}_{1}}=A_{j}^{2}\ \ \ ;\ \ \ A_{j}=\frac{1}{\sqrt{N}-s_{j}\left(\sqrt{N}-1\right)}
\end{equation}
The complete ``\textit{sparsity}" corresponds to $s_{j}=1$, that is $A_{j}=1$, while the absence of ``\textit{sparsity}" corresponds to $s_{j}=0$, that is $A_{j}=\frac{1}{\sqrt{N}}$.\\
Furthermore, the constraints imposed during the minimization process for a column ``$j$" of $X$, which are the non-negativity and the fixed sum, will result in $\left\|X_{j}\right\|_{1}=\left\|Y_{j}\right\|_{1}$, at each iteration.\\
We will express the corresponding penalty term as:
\begin{equation}
	DX\left(X\right)=\sum_{j}\left[\frac{1}{2}\left\|X_{j}\right\|^{2}_{2}-\frac{A_{j}^{2}}{2}\left\|X_{j}\right\|^{2}_{1}\right]^{2}
	\label{eq.sparhoyer}
\end{equation}
The gradient with respect to $X$ will be expressed by:
\begin{equation}
	\frac{\partial DX\left(X\right)}{\partial X_{nm}}=X_{nm}\left\|X_{m}\right\|^{2}_{2}+A_{m}^{4}\left\|X_{m}\right\|^{3}_{1}-A_{m}^{2}X_{nm}\left\|X_{m}\right\|^{2}_{1}-A_{m}^{2}\left\|X_{m}\right\|_{1}\left\|X_{m}\right\|^{2}_{2}
	\label{eq.gradspar}
\end{equation}
Or also, by taking into account the sum constraint:
\begin{equation}
	\frac{\partial DX\left(X\right)}{\partial X_{nm}}=X_{nm}\left\|X_{m}\right\|^{2}_{2}+A_{m}^{4}\left\|Y_{m}\right\|^{3}_{1}-A_{m}^{2}X_{nm}\left\|Y_{m}\right\|^{2}_{1}-A_{m}^{2}\left\|Y_{m}\right\|_{1}\left\|X_{m}\right\|^{2}_{2}
	\label{eq.gradsparsimpl}
\end{equation}

\section{Some examples of regularized algorithms.}

Regarding the ``\textit{data consistency}" term, we will deal with the ``$\alpha$" divergence and the ``$\beta$" divergence.\\
The regularization term on $H$ will be the Euclidean norm of the solution.\\
The regularization on $X$ will be made in the sense of Hoyer.

\subsection{Regularized ``Alpha'' divergence.}
\subsubsection{* Régularization on $H$.}
The basis algorithm is given by (\ref{eq.algoHreg}).\\
The gradients involved in the basic algorithm are given respectively by (\ref{eq.gradA/H}) and (\ref{eq.gradregticste}).\\
Taking into account the simplifications related to the property $\sum_{i}H^{k}_{im}=1$, and with $X\equiv X^{k}$, the regularized algorithm is written according to (\ref{eq.algonmH}):
\begin{align}
		H^{k+1}_{nm}=H^{k}_{nm}&+\frac{\delta^{k}_{m}}{\lambda}H^{k}_{nm}\left[\left(\sum_{l}\left(\frac{Y_{nl}}{(H^{k}X)_{nl}}\right)^{\lambda}X^{T}_{lm}\right)-\gamma H^{k}_{nm}\right. \nonumber \\  & \left.-\sum_{i}H^{k}_{im}\left(\sum_{l}\left(\frac{Y_{il}}{(H^{k}X)_{il}}\right)^{\lambda}X^{T}_{lm}-\gamma H^{k}_{im}\right)\right]
			\label{eq.AlgoA/Hreg}
\end{align}
It can also be written in a condensed form:
\begin{align}
		H^{k+1}_{nm}=H^{k}_{nm}&+\frac{\delta^{k}_{m}}{\lambda}H^{k}_{nm}\left\{\left[\left(\frac{Y}{(H^{k}X)}\right)^{\lambda}X^{T}-\gamma H^{k}\right]_{nm}\right. \nonumber \\  & \left.-\left[\left(H^{k}\right)^{T}\left(\frac{Y}{(H^{k}X)}\right)^{\lambda}X^{T}-\gamma \left(H^{k}\right)^{T}H^{k}\right]_{mm}\right\}
			\label{eq.AlgoA/Hregcond}
\end{align}
We can easily verify that the $\lambda=1$ case corresponding to the Kullback-Leibler divergence can be found as well by using the ``$\beta$" divergences.\\

\textbf{- Multiplicative form.}\\

By denoting:
\begin{equation}
	U^{k}_{nm}=\left[\left(\frac{Y}{(H^{k}X)}\right)^{\lambda}X^{T}-\gamma H^{k}\right]_{nm}
\end{equation}
The algorithm (\ref{eq.AlgoA/Hreg}) is written:
\begin{equation}
	H^{k+1}_{nm}=H^{k}_{nm}+\frac{\delta^{k}_{m}}{\lambda}H^{k}_{nm}\left\{U^{k}_{nm}-\sum_{i}\left(H^{k}_{mi}\right)^{T}U^{k}_{im}\right\}
\end{equation}
We proceed to the shift:
\begin{equation}
\left[U^{k}_{nm}\right]_{d}=U^{k}_{nm}-\min_{n}\left(U^{k}_{nm}\right)+\epsilon>0	
\end{equation}
Taking into account the properties of $H$, the algorithm is written as follows:
\begin{equation}
	H^{k+1}_{nm}=H^{k}_{nm}+\frac{\delta^{k}_{m}}{\lambda}H^{k}_{nm}\left\{\left[U^{k}_{nm}\right]_{d}-\sum_{i}\left(H^{k}_{mi}\right)^{T}\left[U^{k}_{im}\right]_{d}\right\}
\end{equation}
This makes it possible to write an algorithm which will be at the basis of the purely multiplicative algorithm:
\begin{equation}
	H^{k+1}_{nm}=H^{k}_{nm}+\frac{\delta^{k}_{m}}{\lambda}H^{k}_{nm}\left\{\frac{\left[U^{k}_{nm}\right]_{d}}{\sum_{i}\left(H^{k}_{mi}\right)^{T}\left[U^{k}_{im}\right]_{d}}-1\right\}
	\label{eq.AlgoA/Hregcondm}
\end{equation}
The purely multiplicative algorithm will be obtained by choosing a descent step size  $\frac{\delta^{k}_{m}}{\lambda}=1$.

\subsubsection{* Régularization on $X$.}
The basic algorithm is given by (\ref{eq.algoBregbase}) and the expressions of the involved gradients are given by (\ref{eq.gradA/Xbis}) and (\ref{eq.gradsparsimpl}).\\
Considering these expressions and the constraint $\sum_{l}X^{k}_{lm}=\sum_{l}Y_{lm}$, some simplifications occur in the expression of the basic algorithm which can be rewritten in the simplified form:
\begin{equation}
	X^{k+1}_{nm}=X^{k}_{nm}+\delta^{k}_{m}X^{k}_{nm}\left[\left(\sum_{l}Y_{lm}\right)U^{k}_{nm}-\sum_{i}X^{k}_{im}U^{k}_{im}\right]
	\label{eq.AlgoA/Xregspar}
\end{equation}
With:
\begin{align}
	U^{k}_{nm}= &\frac{1}{\lambda}\sum_{j}H^{T}_{nj}\left(\frac{P_{jm}}{Q_{jm}}\right)^{\lambda}\nonumber \\  & +\mu \left(X^{k}_{nm}-A^{2}_{m}\|Y_{m}\|_{1}\right)\left(A^{2}_{m}\|Y_{m}\|^{2}_{1}-\|X_{m}\|^{2}_{2}\right)
\end{align}
In this expression we take into account the fact that:
\begin{equation}
Q=H^{k+1}X^{k}	
\end{equation}

\subsubsection{- Multiplicative form.}
To obtain such a form, one carries out, as already proposed, the shift which allows to make the expression of $U^{k}_{nm}$ positive whatever the component ``$n$" considered:
\begin{equation}
	\left[U^{k}_{nm}\right]_{d}=U^{k}_{nm}-\min_{n}\left(U^{k}_{nm}\right)+\epsilon
\end{equation}
The algorithm (\ref{eq.AlgoA/Xregspar}) will be written identically with the shifted quantities.\\
We can then propose a pseudo multiplicative algorithm which is written as follows:
\begin{equation}
	X^{k+1}_{nm}=X^{k}_{nm}+\delta^{k}_{m}X^{k}_{nm}\left[\frac{\left(\sum_{l}Y_{lm}\right)\left[U^{k}_{nm}\right]_{d}}{\sum_{i}X^{k}_{im}\left[U^{k}_{im}\right]_{d}}-1\right]
\end{equation}
This makes it possible to obtain, by taking a descent step size equal to 1, the purely multiplicative form:
\begin{equation}
	X^{k+1}_{nm}=X^{k}_{nm}\left[\frac{\left(\sum_{l}Y_{lm}\right)\left[U^{k}_{nm}\right]_{d}}{\sum_{i}X^{k}_{im}\left[U^{k}_{im}\right]_{d}}\right]
\end{equation}

\subsection{Regularized ``Beta'' divergence.}
\subsubsection{* Régularization on $H$.}

The basis algorithm is given by (\ref{eq.algoHreg}).\\
The gradients involved in this algorithm are given respectively by (\ref{eq.gradB/H}) and (\ref{eq.gradregticste}).\\
Using a procedure analogous to that of the previous paragraph, taking into account the simplifications related to the property $\sum_{i}H^{k}_{im}=1$, and with $X\equiv X^{k}$, we introduce the expression:
\begin{equation}
	U^{k}_{nm}=\sum_{j}\left[\left(H^{k}X\right)^{\lambda-2}_{nj}Y_{nj}-\left(H^{k}X\right)^{\lambda-1}_{nj}\right]X^{T}_{jm}-\gamma H^{k}_{nm}
\end{equation}
Then, the regularized algorithm is written according to (\ref{eq.algonmH}):
\begin{equation}
		H^{k+1}_{nm}=H^{k}_{nm}+\delta^{k}_{m}H^{k}_{nm} \left(U^{k}_{nm}-\sum_{i}H^{k}_{im}U^{k}_{im}\right)
			\label{eq.AlgoB/Hreg}
\end{equation}
Or also in condensed form:
\begin{equation}
		H^{k+1}_{nm}=H^{k}_{nm}+\delta^{k}_{m}H^{k}_{nm} \left\{U^{k}_{nm}-\left[\left(H^{k}\right)^{T}U^{k}\right]_{mm}\right\}
			\label{eq.AlgoB/Hregbis}
\end{equation}
The case $\lambda=1$ corresponding to the Kullback-Leibler divergence allows us to recover the result already obtained for the ``$\alpha$" divergence, while $\lambda=2$ corresponds to the mean square deviation.\\
In order to obtain a multiplicative algorithm, we perform the shift:
\begin{equation}
\left[U^{k}_{nm}\right]_{d}=U^{k}_{nm}-\min_{n}\left(U^{k}_{nm}\right)+\epsilon	
\end{equation}
The multiplicative algorithm then takes a form similar to (\ref{eq.AlgoA/Hregcondm}), subject to choosing a descent step size equal to $\delta^{k}_{m}=1\ \forall k,m$.

\subsubsection{* Régularization on $X$.} 
The basis algorithm is given by (\ref{eq.algoBregbase}) and the expressions for the involved gradients are given by (\ref{eq.gradB/X}) and (\ref{eq.gradsparsimpl}).\\
Taking into account these expressions and the constraint $\sum_{l}X^{k}_{lm}=\sum_{l}Y_{lm}$, some simplification occurs in the expression of the basis algorithm which can be rewritten in the simplified form:
\begin{equation}
	X^{k+1}_{nm}=X^{k}_{nm}+\delta^{k}_{m}X^{k}_{nm}\left[\left(\sum_{l}Y_{lm}\right)U^{k}_{nm}-\sum_{i}X^{k}_{im}U^{k}_{im}\right]
	\label{eq.AlgoB/Xregspar}
\end{equation}
With:
\begin{align}
	U^{k}_{nm}=&\sum_{j}H^{T}_{nj}\left[P_{jm}\left(Q_{jm}\right)^{\lambda-2}-\left(Q_{jm}\right)^{\lambda-1}\right]\nonumber \\  & +\mu \left(X^{k}_{nm}-A^{2}_{m}\|Y_{m}\|_{1}\right)\left(A^{2}_{m}\|Y_{m}\|^{2}_{1}-\|X_{m}\|^{2}_{2}\right)
\end{align}
In this expression we take into account the fact that:
\begin{equation}
H\equiv H^{k+1}\ \ ;\ \  Q=H^{k+1}X^{k}	
\end{equation}

\subsubsection{- Multiplicative form.}
To obtain such a form, one performs, as already proposed, the shift that allows to make the expression of $U^{k}_{nm}$ positive whatever the component ``$n$" considered:
\begin{equation}
	\left[U^{k}_{nm}\right]_{d}=U^{k}_{nm}-\min_{n}\left(U^{k}_{nm}\right)+\epsilon
\end{equation}
The algorithm (\ref{eq.AlgoA/Xregspar}) will be written identically with the shifted quantities.\\
We can then propose a pseudo multiplicative algorithm which is written as follows:
\begin{equation}
	X^{k+1}_{nm}=X^{k}_{nm}+\delta^{k}_{m}X^{k}_{nm}\left[\frac{\left(\sum_{l}Y_{lm}\right)\left[U^{k}_{nm}\right]_{d}}{\sum_{i}X^{k}_{im}\left[U^{k}_{im}\right]_{d}}-1\right]
\end{equation}
This allows us to obtain, by taking a descent step size equal to 1, the purely multiplicative form:
\begin{equation}
	X^{k+1}_{nm}=X^{k}_{nm}\left[\frac{\left(\sum_{l}Y_{lm}\right)\left[U^{k}_{nm}\right]_{d}}{\sum_{i}X^{k}_{im}\left[U^{k}_{im}\right]_{d}}\right]
\end{equation}

\section{ Regularization. Scale invariant divergences.}
In this case, the resolution of the regularized problem implies the minimization with respect to $H$ and $X$ of a functional founded on the use of invariant divergences, which is written:
\begin{equation}
	JI\left(H,X\right)=DI\left(Y\|HX\right)+\gamma DIH\left(H\right)+\mu DIX\left(X\right)
\end{equation}
The minimization is carried out under the non-negativity and sum constraints already specified.\\
The ``$\gamma$" and ``$\mu$" coefficients are the positive regularization factors.\\
The iterative process of minimization being alternately on $H$ and then on $X$ (for example).\\
As far as the $H$ is concerned, we operate on a composite divergence:
\begin{equation}
DCIH=DI\left(Y\|HX\right)+\gamma DIH\left(H\right)	
\end{equation}
\textbf{Regarding the ``\textit{data consistency}" term $DI\left(Y\|HX\right)$, since the invariance with respect to $Q=HX$ results in an invariance with respect to $X$, but not with respect to $H$, there's nothing particular to be expected from using an invariant divergence concerning the iterations on $H$.}\\

Thus, regardless of the choice of ``\textit{data consistency}" term, the minimization with respect to $H$ will necessarily involve a change in variables as already indicated; more generally, concerning $H$, one can think that it is not even necessary to introduce an invariant divergence inasmuch as the $D\left(Y\|HX\right)$ divergence reflecting the ``\textit{data consistency}" and the corresponding invariant divergence $DI\left(Y\|HX\right)$ have the same minimum.\\ 
As a consequence, a regularization on the columns of $H$ using an invariant form of the penalty term is not absolutely necessary, but, for the sake of homogeneity, it is still possible to operate by a variable change on $H$ in $DCIH$.\\ 

Concerning $X$, we operate on a composite divergence: 
\begin{equation}
DCIX=DI\left(Y\|HX\right)+\mu DIX\left(X\right)	
\end{equation}
For these divergences denoted $DI$, we rely on the properties (\ref{eq.propgen})(\ref{eq.propgenbis}), and more precisely, in our particular problem, on (\ref{eq.propgenbis})(\ref{eq.propgenter}).\\

An iterative algorithm on $X$ can be written, with the method developed in \textbf{Chapter 10}, in the form given by (\ref{eq.algoXI}):
\begin{equation}
	X^{k+1}_{nm}=X^{k}_{nm}+\alpha^{k}_{m}X^{k}_{nm}\left[-\frac{\partial DCIX}{\partial X^{k}_{nm}}\right]
	\label{eq.algoDCIX}
\end{equation}
Of course, this writing implies the use of a penalty term $DIX$ that is invariant with respect to $X$.\\
We can easily show that with this algorithm, the constraint $\sum_{n}X^{k}_{nm}=Cte\;\forall k$ is fulfilled.\\

\textbf{The algorithm (\ref{eq.algoDCIX}) is the iterative algorithm on $X$ when minimizing a regularized invariant divergence.}

 \subsubsection{* Penalty term on $X$. Invariant form.}
The regularization with respect to $X$ reflects the ``\textit{sparsity}" of the components in a column of the $X$ matrix.\\
An invariant form of the penalty term, inferred from (\ref{eq.sparhoyer}) is written:
\begin{equation}
	DIX\left(X\right)=\frac{1}{2}\sum_{j}\left[\frac{\|X_{j}\|^{2}_{2}}{\|X_{j}\|^{2}_{1}}-A^{2}\right]^{2}
\end{equation}
The gradient with respect to $X$ will be written as follows:
\begin{equation}
	\frac{\partial DIX\left(X\right)}{\partial X_{nm}}=\left(A^{2}-\frac{\|X_{m}\|^{2}_{2}}{\|X_{m}\|^{2}_{1}}\right)\left(\frac{\|X_{m}\|^{2}_{2}}{\|X_{m}\|^{2}_{1}}-\frac{X_{nm}}{\|X_{m}\|_{1}}\right)\frac{1}{\|X_{m}\|_{1}}
	\label{eq.gradDIX}
\end{equation}
We can easily verify that:
\begin{equation}
	\sum_{n}X_{nm}\frac{\partial DIX\left(X\right)}{\partial X_{nm}}=0
\end{equation}

\subsection{ Regularized invariant ``Alpha'' divergence.}
The ``\textit{data attachment}" term is the Invariant Alpha divergence  (\ref{eq.AIQ}) or (\ref{eq.AIQs}).\\
With $P\equiv Y$, $Q\equiv HX$ and $K_{0m}$ given by (\ref{eq.K0Aj}), of which we reproduce the expression;
\begin{equation}
	K_{0m}=\left(\frac{\sum_{n}P^{\lambda}_{nm}Q^{1-\lambda}_{nm}}{\sum_{n}Q_{nm}}\right)^{\frac{1}{\lambda}}
		\label{eq.K0Ajbis}
\end{equation}
The expressions of the involved gradients are given by (\ref{eq.gradAIX}) and (\ref{eq.gradDIX}); they lead to:
\begin{align}
	\frac{\partial DCIX}{\partial X_{nm}}=&\frac{1}{\lambda}\sum_{j}H^{T}_{nj}\left[K_{0m}-K_{0m}^{1-\lambda}P^{\lambda}_{jm}Q^{-\lambda}_{jm}\right]\nonumber \\  & +\mu \left(A^{2}-\frac{\|X_{m}\|^{2}_{2}}{\|X_{m}\|^{2}_{1}}\right)\left(\frac{\|X_{m}\|^{2}_{2}}{\|X_{m}\|^{2}_{1}}-\frac{X_{nm}}{\|X_{m}\|_{1}}\right)\frac{1}{\|X_{m}\|_{1}}
\end{align}
Taking into account the fact that the invariance factor depends on the iteration, the iterative algorithm on $X$ (\ref{eq.algoDCIX}), is written:
\begin{align}
	X^{k+1}_{nm}=X^{k}_{nm}&+\alpha^{k}_{m}X^{k}_{nm}\Bigg\{\frac{1}{\lambda}\sum_{j}H^{T}_{nj}\left[\left(K^{k}_{0m}\right)^{1-\lambda}P^{\lambda}_{jm}\left(Q^{k}_{jm}\right)^{-\lambda}-K^{k}_{0m}\right] \nonumber \\  & -\mu \left(A^{2}-\frac{\|X^{k}_{m}\|^{2}_{2}}{\|X^{k}_{m}\|^{2}_{1}}\right)\left(\frac{\|X^{k}_{m}\|^{2}_{2}}{\|X^{k}_{m}\|^{2}_{1}}-\frac{X^{k}_{nm}}{\|X^{k}_{m}\|_{1}}\right)\frac{1}{\|X^{k}_{m}\|_{1}}\Bigg\}
\end{align}
In this expression, the following relationships must be taken into account:\\
* $H\equiv H^{k+1}$\\
* $Q^{k}\equiv H^{k+1}X^{k}$\\
* $\|X^{k}_{m}\|_{1}=\|Y_{m}\|_{1}$ by initialization and properties of the algorithm.\\
* the invariance factor $K^{k}_{0m}$ is deduced from (\ref{eq.K0Ajbis}) by introducing the expression of $Q^{k}$ mentioned above.\\

\subsection{ Regularized invariant ``Beta'' divergence.}
The ``\textit{data attachment}" term is the  invariant Beta divergence (\ref{eq.BIQ}) or (\ref{eq.BIQs}).\\
With $P\equiv Y$, $Q\equiv HX$ and $K_{0m}$ given by (\ref{eq.K0Bj})
\begin{equation}
	K_{0m}=\frac{\sum_{n}P_{nm}Q^{\lambda-1}_{nm}}{\sum_{n}Q^{\lambda}_{nm}}
	\label{eq.K0Bm}
\end{equation}
The expressions of the involved gradients are given by (\ref{eq.gradBIX}) and (\ref{eq.gradDIX}); they lead to:\\
\begin{align}
	\frac{\partial DCIX}{\partial X_{nm}}=&\sum_{j}H^{T}_{nj}\left[K^{\lambda}_{0m}Q^{\lambda-1}_{jm}-K^{\lambda-1}_{0m}P_{jm}Q^{\lambda-2}_{jm}\right]\nonumber \\  & +\mu \left(A^{2}-\frac{\|X_{m}\|^{2}_{2}}{\|X_{m}\|^{2}_{1}}\right)\left(\frac{\|X_{m}\|^{2}_{2}}{\|X_{m}\|^{2}_{1}}-\frac{X_{nm}}{\|X_{m}\|_{1}}\right)\frac{1}{\|X_{m}\|_{1}}
\end{align}\\
Taking into account the variation of the invariance factor during the iteration, we will write:
\begin{equation}
	K^{k}_{0m}=\frac{\sum_{n}P_{nm}\left(Q^{k}_{nm}\right)^{\lambda-1}}{\sum_{n}\left(Q^{k}_{nm}\right)^{\lambda}}
\end{equation}\\
 The iterative algorithm on $X$ (\ref{eq.algoDCIX}), is then written:
\begin{align}
	X^{k+1}_{nm}=X^{k}_{nm}&+\alpha^{k}_{m}X^{k}_{nm}\Bigg\{\sum_{j}H^{T}_{nj}\left[\left(K^{k}_{0m}\right)^{\lambda-1}P_{jm}\left(Q_{jm}\right)^{{\lambda-2}}-\left(K^{k}_{0m}\right)^{\lambda}\left(Q_{jm}\right)^{\lambda-1}\right] \nonumber \\  & -\mu \left(A^{2}-\frac{\|X^{k}_{m}\|^{2}_{2}}{\|X^{k}_{m}\|^{2}_{1}}\right)\left(\frac{\|X^{k}_{m}\|^{2}_{2}}{\|X^{k}_{m}\|^{2}_{1}}-\frac{X^{k}_{nm}}{\|X^{k}_{m}\|_{1}}\right)\frac{1}{\|X^{k}_{m}\|_{1}}\Bigg\}
\end{align}\\
In this expression, the following relationships must be taken into account:\\
* $H\equiv H^{k+1}$\\
* $Q\equiv H^{k+1}X^{k}$\\
* $\|X^{k}_{m}\|_{1}=\|Y_{m}\|_{1}$ by initialization and properties of the algorithm.\\

\section{ Overview of the techniques implemented.}

\subsection{ Iterations on $H$.}
For any kind of divergence, with or without regularization:\\
a) Method of variables change.\\
b) The obtention of multiplicative algorithms requires shifting the components of the gradient.

\subsection{ Iterations on $X$.}
\subsubsection{1 - Non-invariant divergences, with or without a non-invariant regularization term.}
a) Method of variables change.\\
b) The obtention of multiplicative algorithms requires shifting the components of the gradient.

\subsubsection{2 - Invariant divergences, with or without an invariant regularization term.}
a) Use of the specific property of invariant divergences.\\
b) Obtaining multiplicative algorithms with fulfilled sum constraint implies the application of the SGM method, followed by normalization.

\setcounter{table}{0}  \setcounter{equation}{0}  \setcounter{figure}{0} \setcounter{chapter}{12} \setcounter{section}{0} 
\chapter{chapter 12 -\\Application to Blind Deconvolution.}  \label{chptr::chapitre12}

\section{Introduction.}
In this chapter, we do not attempt to revisit the well-known difficulties of signal deconvolution problems; these aspects have been widely developed in the literature dedicated to inverse problems \cite{demoment1985deconvolution} \cite{bertero1998}, and more precisely to the deconvolution of images and spectra \cite{andrews1977digital}, \cite{jansson1984deconvolution}, \cite{jansson2014deconvolution2}.\\
In a general context, the use of scale invariant divergences to obtain multiplicative algorithms that simultaneously ensure non-negativity and the sum constraint, implies, at each iteration, to carry out a normalization that will not modify the value of the invariant divergence considered.\\
However, with regard to particular the problem of blind deconvolution, and more precisely with regard to the respect of the non-negativity constraints, and especially the sum constraints, the use of scale invariant divergences and of their properties leads to particularly well adapted reconstruction algorithms, without the help of a normalization procedure.\\
On the other hand, whatever the type of divergence used, the method of  variables change makes it possible to obtain multiplicative algorithms taking into account the sum constraints.\\

\textbf{Consequently, the first part of this chapter is devoted to the use of invariant divergences and to their minimization within the framework of blind deconvolution problems, then the second part is devoted to non-invariant divergences and the associated multiplicative algorithms.\\
In both cases, we will discuss the smoothness constraint regularization of the solution.}\\

To establish a parallel with the NMF problem of the previous chapter, we will take advantage of the fact that the convolution product is commutative as opposed to the matrix product.\\
Indeed the basic relationship is given by a Fredholm integral of the first kind:
\begin{equation}
	\textbf{Y}\left(r,s\right)=\int^{+\infty}_{-\infty}\int^{+\infty}_{-\infty}\textbf{H}\left(r,s,\xi,\eta\right)\textbf{X}\left(\xi,\eta\right)d\xi d\eta
\end{equation}
In the case of convolution, the kernel $\textbf{H}$ is invariant under translation:
\begin{equation}
\textbf{H}\left(r,s,\xi,\eta\right)=\textbf{H}\left(r-\xi,s-\eta\right)	
\end{equation}
We denote:
\begin{equation}
	\textbf{H}\otimes \textbf{X}=\int^{+\infty}_{-\infty}\int^{+\infty}_{-\infty}\textbf{H}\left(r-\xi,s-\eta\right)\textbf{X}\left(\xi,\eta\right)d\xi d\eta
\end{equation}
In imaging, in the case of classical deconvolution $\textbf{Y}\left(r,s\right)$ is the measured image, usually corrupted by noise,  $\textbf{H}\left(r,s,\xi,\eta\right)$ is the supposedly known, non-negative system Point Spread Function (PSF), $\textbf{X}\left(\xi,\eta\right)$ is the unknown non-negative image we're trying to reconstruct.\\
Moreover, the  Point Spread Function of the system, i.e. the kernel of the integral equation, has generally an integral equal to 1.
\begin{equation}
	\int\int \textbf{H}\left(u,v\right)du dv=1
	\label{eq.inth}
\end{equation}
As a result of which:
\begin{equation}
	\int\int \textbf{Y}\left(r,s\right)dr ds=\int\int \textbf{X}\left(\xi,\eta\right)d\xi d\eta
\label{eq.intyx}	
\end{equation}
The ``simple" deconvolution problem is similar to the linear unmixing problem discussed in the previous chapter.\\
We can notice that in this problem, the non-negativity of the $\textbf{X}$ solution and the relation (\ref{eq.intyx}) are the constraints.\\
On the other hand, in the blind deconvolution problem, we dispose only of the measured image $\textbf{Y}$, and the unknowns are the kernel $\textbf{H}$ and the true object $\textbf{X}$.\\
The constraints are the non-negativity of $\textbf{H}$ and $\textbf{X}$ on the one hand, and the relations (\ref{eq.inth}) and (\ref{eq.intyx}) on the other hand.\\
Of course, despite the constraints, the problems associated with multiple solutions are not completely solved.\\
This problem is similar to the NMF problem in the previous chapter, with the already mentioned advantage of the commutativity of the convolution product.
In the convolution problem, the measured image, the unknown object and the impulse response are both pixelated and contained in arrays of the same dimensions.\\
The basic relationship:
\begin{equation}
 \textbf{Y}=\textbf{H}\otimes \textbf{X}\ \Leftrightarrow\ \textbf{Y}\left(r,s\right)=\int^{+\infty}_{-\infty}\int^{+\infty}_{-\infty}\textbf{H}\left(r-\xi,s-\eta\right)\textbf{X}\left(\xi,\eta\right)d\xi d\eta
\end{equation} 
is written after discretization, in the vector form:
\begin{equation}
	y=Hx
	\label{eq.relyx}
\end{equation}
In this matrix form, the measured image of size $\left(N*N\right)$ is lexicographically ordered as a vector ``$y$" of length $\left(N^{2}\right)$; the same applies to the unknown object ``$x$", then, "$H$" is a circulating Toeplitz block matrix of dimension $\left(N^{2}*N^{2}\right)$.\\
Given the commutativity of the convolution product, we also have:
\begin{equation}
\textbf{Y}=\textbf{X}\otimes \textbf{H}\Leftrightarrow\\ \textbf{Y}\left(r,s\right)=\int^{+\infty}_{-\infty}\int^{+\infty}_{-\infty}\textbf{X}\left(r-\xi,s-\eta\right)\textbf{H}\left(\xi,\eta\right)d\xi d\eta
\end{equation}
which results in:
\begin{equation}
	y=Xh
		\label{eq.relyh}
\end{equation}
Here, vectors ``$h$" and ``$y$" of length $\left(N^{2}\right)$ correspond to the lexicographical organization of the arrays $\textbf{H}$ and $\textbf{Y}$, then, $X$ is a circulating Toeplitz block matrix of size $\left(N^{2}*N^{2}\right)$.\\
By using the notations already introduced in the analysis of the divergences, we will have:
\begin{equation}
	p\equiv y\ \ ;\ \ q\equiv Hx\equiv Xh
	\label{eq.relxypq}
\end{equation}
Referring to the notations of relations (\ref{eq.relyx}) (\ref{eq.relyh}) (\ref{eq.relxypq}), the resolution of the blind deconvolution problem consists in recovering the PSF ``$h$" and the object ``$x$" by minimizing a divergence $D\left(p\|q\right)$ under the constraints:
\begin{equation}
	x_{i}\geq 0\ \forall i\;\;;h_{i}\geq0\ \forall i\;\;;\sum_{i}h_{i}=1\;\;;\sum_{i}x_{i}=\sum_{i}y_{i}
\end{equation}
\textbf{Due to the lexicographical organization of the images (tables), all considerations and operations concerning divergences between vectors are directly applicable in deconvolution or blind deconvolution problems; the image aspects (tables) will only appear in calculations using Fourier transforms.}\\

Note that during calculations, expressions of the form $\left(H^{T}x\right)$ and $\left(X^{T}h\right)$ will be used; these operations expressed in the context of continuous functions correspond to:
\begin{equation}
	H^{T}u\equiv \int\int \textbf{H}\left(-r,-s,\xi,\eta\right)\textbf{U}\left(\xi,\eta\right)d\xi d\eta
\end{equation}
with, in the case of convolution:
\begin{equation}
	 \textbf{H}\left(-r,-s,\xi,\eta\right)=\textbf{H}\left(-r-\xi,-s-\eta\right)
\end{equation}
and
\begin{equation}
	X^{T}u\equiv \int\int \textbf{X}\left(-r,-s,\xi,\eta\right)\textbf{U}\left(\xi,\eta\right)d\xi d\eta
\end{equation}
with, in the case of convolution:
\begin{equation}
	\textbf{X}\left(-r,-s,\xi,\eta\right)=\textbf{X}\left(-r-\xi,-s-\eta\right)
\end{equation}

\subsection{Some details related to the calculus.}
To perform the calculus of the convolution products, it is interesting to use the Fourier Transforms noted $\mathcal{F}$ ; indeed, the classical property of the convolution product allows to write \cite{roddier1971distributions}:
\begin{equation}
\int^{+\infty}_{-\infty}\int^{+\infty}_{-\infty}\textbf{H}\left(r-\xi,s-\eta\right)\textbf{X}\left(\xi,\eta\right)d\xi d\eta=\mathcal{F}^{-1}\left[\mathcal{F}\left(\textbf{H}\right)\mathcal{F}\left(\textbf{X}\right)\right]	
\end{equation}
and moreover, we have:
\begin{equation}
	\textbf{H}\left(-r,-s\right)=\mathcal{F}\left[\mathcal{F}\left(\textbf{H}\left(r,s\right)\right)\right]\ \ ;\ \ \textbf{X}\left(-r,-s\right)=\mathcal{F}\left[\mathcal{F}\left(\textbf{X}\left(r,s\right)\right)\right] 
\end{equation}
In doing so, we will only handle tables of the size of the object, image, and PSF.\\
As a result, in the reconstruction iterative process, the unknown images (tables) are updated as a whole.

\section{Algorithmic method.}
The method described applies to all divergences, whether invariant or not.\\
Leaving aside the regularization aspects, the problem is posed as a constrained minimization problem of a divergence between the measures $p\equiv y$ and the convolutional model $q\equiv Hx\equiv Xh\equiv \left[h\otimes x\right]$, which reflects the ``\textit{data attachment}''.\\	
$\underbrace{\min}_{x,h} DI\left(y\|\left[h\otimes x\right]\right)$\\
 s.t. $\ \ x_{i}>0\ \forall i\ ;\  \sum_{i}x_{i}=\sum_{i}y_{i}\ $ et $\ h_{i}>0\ \forall i\ ;\  \sum_{i}h_{i}=1$.\\
	
The proposed blind deconvolution method is structurally inspired by the one initially proposed in \cite{ayers1988iterative} and extended by the use of iterative algorithms in \cite{lanteri1994blind} \cite{lanteri1995comparison} for example.\\
This method, in which we operate alternately on the two unknowns is analogous to the one developed for MMF in \textbf{chapter 11}.\\

It can be synthesized, for example, according to the following scheme:\\
1 - $x^{k}\;,\;h^{k}\rightarrow x^{k}\;,\;h^{k+1}$\\
2 - $x^{k}\;,\;h^{k+1}\rightarrow x^{k+1}\;,\;h^{k+1}$\\

This is the scheme that will be considered in further calculations.\\
Note that we might as well iterate first on ``$x$" and then on ``$h$".\\
In blind deconvolution, as in simple deconvolution, the updating is done on all the components of ``$x$" simultaneously; similarly, the updating to ``$h$" is done on all the components simultaneously.\\
This last point differs from what is done in the case of the NMF where the updating $H^{k}\rightarrow H^{k+1}$ is carried out successively column by column as well as the updating $X^{k}\rightarrow X^{k+1}$ .\\
The alternative iterative method inspired by Uzawa's algorithm suggested in \textbf{Chapter 11} for NMF, can also be considered for blind deconvolution; it was used in \cite{lanteri1994blind}, \cite{lanteri1995comparison}. Such a method will not be developed here. 

\section{Application to scale invariant divergences.}

\subsection{Property of scale invariant divergences with respect to ``$q$".}

The use of invariant divergences with respect to ``$q$" is particularly interesting for the blind deconvolution problem, because their use associated with the commutativity of the convolution product allows to take into account simply the sum constraints on the two unknowns, i.e. $\sum_{i}h_{i}=1\;\;;\sum_{i}x_{i}=\sum_{i}y_{i}$.\\
Let's recall the property of invariant divergences arising here (\ref{eq.Pfond}):
\begin{equation}
	\sum_{m}q_{m}\left[\frac{\partial DI\left(p\|q\right)}{\partial q}\right]_{m}=0
\end{equation}

\subsubsection{* Initialization.}
The algorithm is initialized with $x^{0}\ /\ \sum_{i}x^{0}_{i}= \sum_{i}y_{i}$ and $h^{0}\ /\ \sum_{i}h^{0}_{i}=1$.

\subsubsection{* Itérations on ``$h$" - Non-multiplicative form.}
We have at our disposal $x^{k}$, $h^{k}$; here, $q=X^{k}h^{k}$.\\
For a divergence invariant with respect to ``$q$" denoted $DI\left(p\|q\right)$, the iterative algorithm on ``$h$" is written :
\begin{equation}
	h^{k+1}_{l}=h^{k}_{l}+\alpha^{k}_{l}h^{k}_{l}\left[-\frac{\partial DI\left(p\|q\right)}{\partial h^{k}_{l}}\right]
	\label{eq.iterhDCV}
\end{equation}
The gradient with respect to ``$h$" is given by:
\begin{equation}
	\frac{\partial DI\left(p\|q\right)}{\partial h_{l}}=\sum_{i}\frac{\partial DI\left(p_{i}\|q_{i}\right)}{\partial q_{i}}\frac{\partial q_{i}}{\partial h_{l}}
\end{equation}
With 
\begin{equation}
	\frac{\partial q_{i}}{\partial h_{l}}=X_{il}
\end{equation}
We have:
\begin{equation}
	\frac{\partial DI\left(p\|q\right)}{\partial h_{l}}=\sum_{i}\frac{\partial DI\left(p_{i}\|q_{i}\right)}{\partial q_{i}}X_{il}=\sum_{i}X^{T}_{li}\frac{\partial DI\left(p_{i}\|q_{i}\right)}{\partial q_{i}}
\end{equation}
And finally:
\begin{equation}
	\frac{\partial DI\left(p\|q\right)}{\partial h_{l}}=\left[X^{T}\frac{\partial DI\left(p\|q\right)}{\partial q}\right]_{l}
\end{equation}
The iterative algorithm on ``$h$" is then written with $x\equiv x^{k}\Leftrightarrow X\equiv X^{k}$ and $q=X^{k}h^{k}$:
\begin{equation}
	h^{k+1}_{l}=h^{k}_{l}+\alpha^{k}_{l}h^{k}_{l}\left[-X^{T}\frac{\partial DI\left(p\|q\right)}{\partial q}\right]_{l}
\end{equation}
As in all analogous algorithms, as indicated in \textbf{Chapter 10}, we must in a first step calculate the step size $\left[\alpha^{k}_{l}\right]_{Max}$ ensuring the non-negativity of all the components of $h^{k+1}$, then calculate by a one-dimensional search method over the interval $\left\{0,\left[\alpha^{k}_{l}\right]_{Max}\right\}$, the descent step size $\alpha^{k}$ ensuring the convergence of the algorithm:
\begin{equation}
	h^{k+1}_{l}=h^{k}_{l}+\alpha^{k}h^{k}_{l}\left[-X^{T}\frac{\partial DI\left(p\|q\right)}{\partial q}\right]_{l}
	\label{eq.iterhDCVbis}
\end{equation} 
For the whole ``$h$" vector, we can write in a more explicit way:
\begin{equation}
	h^{k+1}=h^{k}+\alpha^{k}h^{k}\odot \left[-X^{T}\frac{\partial DI\left(p\|q\right)}{\partial q}\right]
\end{equation}
In this writing the symbol $\odot$ is the pointwise product of two vectors.\\
With such an algorithm, considering the property (\ref{eq.Pfond}), we obtain $\sum_{l}h^{k+1}_{l}=\sum_{l}h^{k}_{l}$, whatever the descent step size $\alpha^{k}$, as shown by the following simple calculation.

\subsubsection{* Démonstration of $\sum_{l}h^{k+1}_{l}=\sum_{l}h^{k}_{l}$}
From (\ref{eq.iterhDCV}), we have:
\begin{equation}
	\sum_{l}h^{k+1}_{l}=\sum_{l}h^{k}_{l}+\alpha^{k}\sum_{l}\left\{h^{k}_{l}\left[-X^{T}\frac{\partial DI\left(p\|q\right)}{\partial q}\right]_{l}\right\}
\end{equation} 
By explaining:
\begin{equation}
	\left[X^{T}\frac{\partial DI\left(p\|q\right)}{\partial q}\right]_{l}=\sum_{m}X^{T}_{lm}\left[\frac{\partial DI\left(p\|q\right)}{\partial q}\right]_{m}
\end{equation}
Putting in the previous equation, it comes:
\begin{equation}
\sum_{l}h^{k+1}_{l}=\sum_{l}h^{k}_{l}+\alpha^{k}\sum_{l}\left\{h^{k}_{l}\left[-\sum_{m}X^{T}_{lm}\left[\frac{\partial DI\left(p\|q\right)}{\partial q}\right]_{m}\right]\right\}
\end{equation}
By swapping the order of the summations, it comes:
\begin{equation}
\sum_{l}h^{k+1}_{l}=\sum_{l}h^{k}_{l}+\alpha^{k}\sum_{m}\left[-\frac{\partial DI\left(p\|q\right)}{\partial q}\right]_{m}\sum_{l}h^{k}_{l}X^{T}_{lm}	
\end{equation}
That is:
\begin{equation}
\sum_{l}h^{k+1}_{l}=\sum_{l}h^{k}_{l}+\alpha^{k}\sum_{m}\left[-\frac{\partial DI\left(p\|q\right)}{\partial q}\right]_{m}\sum_{l}X_{ml}	h^{k}_{l}
\end{equation}
Given the property (\ref{eq.Pfond}) the second term of the second member of the previous equation is zero, hence the announced result.

\subsubsection{* Itérations on ``$h$" - Multiplicative form.}
In order to obtain a multiplicative form of the algorithm, from the relations (\ref{eq.iterhDCV}) (\ref{eq.iterhDCVbis}), we strictly apply the SGM method and we carry out in a first step the decomposition:
\begin{equation}
\left[-\frac{\partial DI\left(p\|q\right)}{\partial h^{k}}\right]_{l}=\left[-X^{T}\frac{\partial DI\left(p\|q\right)}{\partial q}\right]_{l}=U^{k}_{l}-V^{k}_{l}	
\end{equation}
Avec $U^{k}_{l}>0$, $V^{k}_{l}>0$.\\
We then write in the next step:
\begin{equation}
	\widetilde{h}^{k+1}_{l}=h^{k}_{l}+\alpha^{k}h^{k}_{l}\left[\frac{U^{k}_{l}}{V^{k}_{l}}-1\right]_{l}
\end{equation}
Then, by choosing a descent step size equal to 1 for any ``$k$", we get:
\begin{equation}
	\widetilde{h}^{k+1}_{l}=h^{k}_{l}\left[\frac{U^{k}_{l}}{V^{k}_{l}}\right]_{l}
\end{equation}
At this point, we don't have $\sum_{l}\widetilde{h}^{k+1}_{l}=\sum_{l}h^{k}_{l}=1$.\\
However, the invariance properties of the divergences allow to carry out, in a last step, a normalization which does not modify the value of the divergence, and we finally have:
\begin{equation}
	h^{k+1}_{l}=\frac{\widetilde{h}^{k+1}_{l}}{\sum_{l}\widetilde{h}^{k+1}_{l}}
\end{equation}

\subsubsection{* Itérations on ``$x$" - Non multiplicative form.}
At this point we have $x^{k}$, $h^{k+1}$; here, $q=H^{k+1}x^{k}$.\\
The iterative algorithm on ``$x$" is written as follows:
\begin{equation}
	x^{k+1}_{l}=x^{k}_{l}+\beta^{k}_{l}x^{k}_{l}\left[-\frac{\partial DI\left(p\|q\right)}{\partial x^{k}_{l}}\right]
\end{equation}
The gradient with respect to ``$x$" is written:
\begin{equation}
	\frac{\partial DI\left(p\|q\right)}{\partial x_{l}}=\sum_{i}\frac{\partial DI\left(p_{i}\|q_{i}\right)}{\partial q_{i}}\frac{\partial q_{i}}{\partial x_{l}}
\end{equation}
With: 
\begin{equation}
	\frac{\partial q_{i}}{\partial x_{l}}=H_{il}
\end{equation}
It comes:
\begin{equation}
	\frac{\partial DI\left(p\|q\right)}{\partial x_{l}}=\sum_{i}\frac{\partial DI\left(p_{i}\|q_{i}\right)}{\partial q_{i}}H_{il}=\sum_{i}H^{T}_{li}\frac{\partial DI\left(p_{i}\|q_{i}\right)}{\partial q_{i}}
\end{equation}
And finally:
\begin{equation}
	\frac{\partial DI\left(p\|q\right)}{\partial x_{l}}=\left[H^{T}\frac{\partial DI\left(p\|q\right)}{\partial q}\right]_{l}
\end{equation}
The iterative algorithm on ``$x$" is then written with $h\equiv h^{k+1}\Leftrightarrow H\equiv H^{k+1}$ and $q=H^{k+1}x^{k}$:
\begin{equation}
	x^{k+1}_{l}=x^{k}_{l}+\beta^{k}_{l}x^{k}_{l}\left[-H^{T}\frac{\partial DI\left(p\|q\right)}{\partial q}\right]_{l}
\end{equation}
As in all analogous algorithms, as indicated in \textbf{Chapter 10}, we must in a first step calculate the step size $\left[\beta^{k}_{l}\right]_{Max}$ which ensures the non-negativity of all the components of $x^{k+1}$, then calculate by a one-dimensional search method over the interval $\left\{0,\left[\beta^{k}_{l}\right]_{Max}\right\}$, the descent step size $\beta^{k}$ ensuring the convergence of the algorithm:
\begin{equation}
	x^{k+1}_{l}=x^{k}_{l}+\beta^{k}x^{k}_{l}\left[-H^{T}\frac{\partial DI\left(p\|q\right)}{\partial q}\right]_{l}
	\label{eq.iterx}
\end{equation} 
A demonstration similar to the one developed for ``$h$" makes it easy to establish that with this type of algorithm one obtains $\sum_{l}x^{k+1}_{l}=\sum_{l}x^{k}_{l}$, regardless of the descent step size $\beta^{k}$.\\
An expression of this algorithm for the ``$x$" vector as a whole, is:
\begin{equation}
	x^{k+1}=x^{k}+\beta^{k}x^{k}\odot\left[-H^{T}\frac{\partial DI\left(p\|q\right)}{\partial q}\right]
\end{equation}
In this expression, the symbol $odot $ represents the pointwise product.\\ 
At this point, we have obtained $h^{k+1}$, $x^{k+1}$.

\subsubsection{* Iterations on ``$x$" - Multiplicative form.}
To obtain this form, we strictly follow what has been described for ``$h$" ( the decomposition here is about $\left[-H^{T}\frac{\partial DI\left(p\|q\right)}{\partial q}\right]_{l}$), until we obtain $\widetilde{x}^{k+1}_{l}$, after which, the normalization step is written:
\begin{equation}
	x^{k+1}_{l}=\frac{\widetilde{x}^{k+1}_{l}}{\sum_{l}\widetilde{x}^{k+1}_{l}}\sum_{l}y_{l}
\end{equation}
Such normalization will not change the value of the divergence.

\section{Analysis of 2 particular cases.}
\subsection{Invariant ``Beta'' divergence.}
In \textbf{Chapter 5}, we developed the invariance aspects concerning the Beta divergence.\\
The initial divergence is given by  (\ref{eq.BC}).\\
The expression of the invariance factor $K$ is obtained explicitly, it is therefore the nominal value which is written as follows (\ref{eq.K0BI}):
\begin{equation}
	K_{0}=\frac{\sum_{i}p_{i}q_{i}^{\lambda-1}}{\sum_{i}q^{\lambda}_{i}}	
\end{equation}
And the corresponding divergence invariant by scale change on ``$q$" is written (\ref{eq.BCIpq2}):
\begin{equation}	
BI\left(p\|q\right)=\frac{1}{\lambda\left(\lambda-1\right)}\left[\sum_{i}p_{i}^{\lambda}-\sum_{i}K^{\lambda}_{0}q_{i}^{\lambda}\right]
\end{equation}
Its gradient with respect to ``$q$" has been previously written in the form (\ref{eq.GBCIpqs}):
\begin{equation}
		\frac{\partial BI\left(p\|q\right)}{\partial q_{j}}=K^{\lambda}_{0}\left[q^{\lambda-1}_{j}-\left(K_{0}\right)^{-1}p_{j}q^{\lambda-2}_{j}\right]	
\end{equation}
The algorithm can then be written:\\
Initialization:
\begin{equation}
	h^{0}_{i} \ \ / \ \ \sum_{i}h^{0}_{i}=1\ \ ;\ \ \ \ 
	x^{0}_{i} \ \ / \ \ \sum_{i}x^{0}_{i}=\sum_{i}y_{i}
\end{equation}

\textbf{* Iteration on ``$h$" according to (\ref{eq.iterhDCV}):}\\

We know $h^{k}$, $x^{k}\rightarrow X^{k}\equiv X$, $p=y$, $q=X^{k}h^{k}\left(=H^{k}x^{k}\right)$.\\
The invariance factor is iteration dependent and is written as:
\begin{equation}
	K^{k}_{0}=\frac{\sum_{i}p_{i}q_{i}^{\lambda-1}}{\sum_{i}q^{\lambda}_{i}}
		\label{eq.KoBkactu}
\end{equation}
\begin{equation}
h^{k+1}_{l}=h^{k}_{l}+\alpha^{k}\left(K^{k}_{0}\right)^{\lambda}h^{k}_{l}\left\{\sum_{j}X^{T}_{lj}\left[\left(K^{k}_{0}\right)^{-1}p^{\lambda}_{j}q^{\lambda-2}_{j}-q^{\lambda-1}_{j}\right]\right\}	
\end{equation}
At this point, we have obtained $h^{k+1}\rightarrow H\equiv H^{k+1}$, we know $p=y$, $x^{k}$, $q=H^{k+1}x^{k}$, then, updating $K^{k}_{0}$ according to (\ref{eq.KoBkactu}), we iterate on ``$x$''.\\

\textbf{* Iteration on ``$x$" according to (\ref{eq.iterx}):}\\
\begin{equation}
x^{k+1}_{l}=x^{k}_{l}+\alpha^{k}\left(K^{k}_{0}\right)^{\lambda}x^{k}_{l}\left\{\sum_{j}H^{T}_{lj}\left[\left(K^{k}_{0}\right)^{-1}p^{\lambda}_{j}q^{\lambda-2}_{j}-q^{\lambda-1}_{j}\right]\right\}		
\end{equation}

\subsection{Invariant mean square deviation.}
For the root mean square deviation which corresponds to $\lambda=2$ for the Beta divergence, the nominal invariance factor $K_{0}$ can be explicitly computed and is written as follows:
\begin{equation}
K_{0}=\frac{\sum_{l}p_{l}q_{l}}{\sum_{l}q^{2}_{l}}	
\end{equation}
The corresponding invariant divergence is given by:
\begin{equation}
	MCI\left(p\|q\right)=\frac{1}{2}\sum_{i}\left[p^{2}_{i}-\left(\frac{\sum_{l}p_{l}q_{l}}{\sum_{l}q^{2}_{l}}\right)^{2}q^{2}_{i}\right]
\end{equation}
Its gradient with respect to ``$q$" is written as follows:
\begin{equation}
	\frac{\partial MCI\left(p\|q\right)}{\partial q_{j}}=\left(\frac{\sum_{l}p_{l}q_{l}}{\sum_{l}q^{2}_{l}}\right)^{2}q_{j}-\left(\frac{\sum_{l}p_{l}q_{l}}{\sum_{l}q^{2}_{l}}\right)p_{j}=K^{2}_{0}q_{j}-K_{0}p_{j}
\end{equation}
According to the relation (\ref{eq.iterhDCV}), after discussion about the descent step size and taking into account the variarion of $K_{0}$ during iterations, it comes for the entire vector ``$h$":
\begin{equation}
	h^{k+1}_{l}=h^{k}_{l}+\alpha^{k}h^{k}_{l}\left\{X^{T}\left[K^{k}_{0}p-\left(K^{k}_{0}\right)^{2}q\right]\right\}_{l}	
\end{equation}
In this expression, $x\equiv x^{k}\Leftrightarrow X\equiv X^{k}$, $q=H^{k}x^{k}=X^{k}h^{k}$, and $K^{k}_{0}$ is expressed as:
\begin{equation}
K^{k}_{0}=\frac{\sum_{l}y_{l}\left(H^{k}x^{k}\right)_{l}}{\sum_{l}\left(H^{k}x^{k}\right)^{2}_{l}}	
\end{equation}

We then iterate on ``$x$" according to the relationship (\ref{eq.iterx}); after discussion about the descent step size and  with $h\equiv h^{k+1}\Leftrightarrow H\equiv H^{k+1}$ et $q=H^{k+1}x^{k}$, it comes for the entire vector ``$x$":
\begin{equation}
	x^{k+1}_{l}=x^{k}_{l}+\beta^{k}x^{k}_{l}\left\{H^{T}\left[K^{k}_{0}p-\left(K^{k}_{0}\right)^{2}q\right]\right\}_{l}
\end{equation} 

The invariance factor $K^{k}_{0}$ variable during iterations, will be expressed by:
\begin{equation}
K^{k}_{0}=\frac{\sum_{l}y_{l}\left(H^{k+1}x^{k}\right)_{l}}{\sum_{l}\left(H^{k+1}x^{k}\right)^{2}_{l}}	
\end{equation}

\subsection{Invariant Alpha divergence.}
In \textbf{Chapter 5}, we discussed the invariance aspects concerning the Alpha divergence.\\
The initial divergence is given by  (\ref{eq.AC}).\\
The expression of the invariance factor $K$ is obtained explicitly, so it is the nominal value which is written as (\ref{eq.K0AI}):
\begin{equation}
	K_{0}=\left(\frac{\sum_{i}p_{i}^{\lambda}q_{i}^{1-\lambda}}{\sum_{i}q_{i}}\right)^{\frac{1}{\lambda}}	
\end{equation}
And the corresponding divergence invariant by scale change on ``$q$" is written (\ref{eq.ACIpq2}):
\begin{equation}	
AI\left(p\|q\right)=\frac{1}{\lambda-1}\left[\sum_{i}K_{0}q_{i}-\sum_{i}p_{i}\right]
\end{equation}
Its gradient with respect to ``$q$" has been written in the form (\ref{eq.GACIpqs}):
\begin{equation}
		\frac{\partial AI\left(p\|q\right)}{\partial q_{j}}=\frac{1}{\lambda}\left[K_{0}-K^{1-\lambda}_{0}p^{\lambda}_{j}q^{-\lambda}_{j}\right]	
\end{equation}
With the initialization:
\begin{equation}
	h^{0}_{i} \ \ / \ \ \sum_{i}h^{0}_{i}=1\ \ ;\ \ \ \ 
	x^{0}_{i} \ \ / \ \ \sum_{i}x^{0}_{i}=\sum_{i}y_{i}
\end{equation}
The algorithm can then be described as follows:\\

\textbf{* Iteration on ``$h$" according to (\ref{eq.iterhDCV}):}\\

We know $h^{k}$, $x^{k}\rightarrow X^{k}\equiv X$, $p=y$, $q=X^{k}h^{k}$.\\
The iteration-dependent invariance factor is written as:
\begin{equation}
	K^{k}_{0}=\left(\frac{\sum_{i}p_{i}^{\lambda}q_{i}^{1-\lambda}}{\sum_{i}q_{i}}\right)^{\frac{1}{\lambda}}
		\label{eq.Kokactu}
\end{equation}
The iteration is written for ``$h$" as a whole:
\begin{equation}
h^{k+1}=h^{k}+\frac{\alpha^{k}K^{k}_{0}}{\lambda}h^{k}\odot\left\{X^{T}\left[\left(K^{k}_{0}\right)^{-\lambda}p^{\lambda}q^{-\lambda}-1\right]\right\}	
\end{equation}
The symbol ``$\odot$" corresponds to the pointwise multiplication of the two vectors.\\
At this point, we obtained $h^{k+1}$.\\

\textbf{* Iteration on ``$x$" according to (\ref{eq.iterx}):}\\

We have $h^{k+1}\rightarrow H^{k+1}\equiv H$, $x^{k}$, $p=y$, $q=H^{k+1}x^{k}$.\\
The iteration-dependent invariance factor is updated according to (\ref{eq.Kokactu}), and it comes:
\begin{equation}
x^{k+1}=x^{k}+\frac{\alpha^{k}K^{k}_{0}}{\lambda}x^{k}\odot\left\{H^{T}\left[\left(K^{k}_{0}\right)^{-\lambda}p^{\lambda}q^{-\lambda}-1\right]\right\}		
\end{equation}
The symbol ``$\odot$" corresponds to the pointwise multiplication of the two vectors.\\
At this point, we obtained $h^{k+1}$, $x^{k+1}$.

\subsection{Kullback-Leibler invariant divergence.}
This divergence corresponds to the case $\lambda=1$ of the Alpha divergence as well as the Beta divergence. This is the point in common between these two types of divergence.\\
The nominal invariance factor $K_{0}$ can be explicitly derived and is written:
\begin{equation}
K_{0}=\frac{\sum_{l}p_{l}}{\sum_{l}q_{l}}	
\end{equation}
The corresponding invariant divergence is given by:
\begin{equation}
	KLI\left(p\|q\right)=\sum_{l}p_{l}\sum_{i}\frac{p_{i}}{\sum_{l}p_{l}}\log\left(\frac{\sum_{l}q_{l}}{\sum_{l}p_{l}}\right)\frac{p_{i}}{q_{i}}
\end{equation}
Its gradient with respect to ``$q$" is written as follows:
\begin{equation}
	\frac{\partial KLI\left(p\|q\right)}{\partial q_{j}}=\frac{\sum_{l}p_{l}}{\sum_{l}q_{l}}-\frac{p_{j}}{q_{j}}=K_{0}-\frac{p_{j}}{q_{j}}
\end{equation}
According to the relationship (\ref{eq.iterhDCV}) and after discussing the step size, he comes for ``$h$" as a whole:
\begin{equation}
	h^{k+1}=h^{k}+\alpha^{k}h^{k}\odot\left[X^{T}\left(\frac{p}{q}-K^{k}_{0}\right)\right]	
\end{equation}
In this expression, $p\equiv y$, $x\equiv x^{k}\Leftrightarrow X\equiv X^{k}$ and $q=H^{k}x^{k}=X^{k}h^{k}$.\\
The invariance factor $K^{k}_{0}$ will be expressed by:
\begin{equation}
K^{k}_{0}=\frac{\sum_{l}y_{l}}{\sum_{l}\left(X^{k}h^{k}\right)_{l}}	
\end{equation}

We then iterate on ``$x$" according to the relation (\ref{eq.iterx}); after discussion on the descent step size, it comes for ``$x$" as a whole:
\begin{equation}
	x^{k+1}=x^{k}+\beta^{k}x^{k}\odot\left[H^{T}\left(\frac{p}{q}-K^{k}_{0}\right)\right]
\end{equation} 
with $p\equiv y$, $h\equiv h^{k+1}\Leftrightarrow H\equiv H^{k+1}$,  $q=H^{k+1}x^{k}$ and $K^{k}_{0}$ given by:
\begin{equation}
K^{k}_{0}=\frac{\sum_{l}y_{l}}{\sum_{l}\left(H^{k+1}x^{k}\right)_{l}}	
\end{equation}

\section{Regularization - Scale invariant divergences.}
In the context of deconvolution or blind deconvolution, we classically impose a property of smoothness of the solution; this is the regularization in the sense of Tikhonov \cite{tikhonov1974methods}.\\
The 2 classical methods used to introduce such a property consist in using as a penalty term either the Euclidean norm of the deviation between the solution and its mean value, or the Euclidean norm of the Laplacian of the solution.\\
In a first step, we develop the expressions of the penalty terms having the invariance properties, which will be associated with the ``\textit{data attachment}" term.

\subsection{Euclidean norm of the solution.}
The non-invariant penalty term for the variable ``$x$" is written as:
\begin{equation}
	DRQ_{x}=\sum_{i}\left(c_{i}-x_{i}\right)^{2}
	\label{eq.MCninv}
\end{equation}
Where $c_{i}=1/N^{2}\ \forall i$ and where $x_{i}$ is the image of dimension $N*N$  lexicographically ordered.\\
First we need to obtain a form that is invariant with respect to ``$x$".\\
To do that, we derive the invariance factor, which is written:
\begin{equation}
	K_{0}=\frac{\sum_{j}c_{j}x_{j}}{\sum_{j}x^{2}_{j}}
\end{equation}
The corresponding invariant divergence is written as follows:
\begin{equation}
	DRQI_{x}=\sum_{i}\left(c_{i}-K_{0}x_{i}\right)^{2}
	\label{eq.MCreg}
\end{equation}
After a few calculations, the gradient with respect to ``$x$" is expressed as follows:
\begin{equation}
	\frac{\partial DRQI_{x}}{\partial x_{l}}=K^{2}_{0}x_{l}-K_{0}c_{l}
\end{equation}
We can easily verify that we have:
\begin{equation}
	\sum_{l}x_{l}\frac{\partial DRQI_{x}}{\partial x_{l}}=0
	\label{eq.propinvar}
\end{equation}

\subsection{Use of the Laplacian norm of the solution.}
The classical non-invariant penalty term for the variable ``$x$" is written as:
\begin{equation}
	DRL_{x}=\left\|Lx\right\|^{2}=\left\|x-Tx\right\|^{2}
	\label{eq.Lapninv}
\end{equation}
The term ``$Lx$" corresponds to the convolution of the table $\textbf{X}$ by the Laplacian mask: 
$
\begin{pmatrix}0&-1/4&0\\ 
-1/4&1&-1/4\\ 
0&-1/4&0\\ 
\end{pmatrix}.  
$\\
The resulting table is written in lexicographical order.\\
We can thus write the vector ``$Lx$" in the form $\left(x-Tx\right)$ where the term ``$Tx$" is a vector which, in lexicographic order, corresponds to the convolution of $\textbf{X}$ by the mask:
$
\begin{pmatrix}0&1/4&0\\ 
1/4&0&1/4\\ 
0&1/4&0\\ 
\end{pmatrix}.
$\\
The ``$x$" vector corresponds to the lexicographical writing of the table $\textbf{X}$.\\
We will note that under these conditions, the matrix $T$ is symmetrical, it is normalized to 1 in columns, i.e. $\sum_{i}T_{ij}=1\ \forall j$, thus $\sum_{i}\left(Tx\right)_{i}=\sum_{i}x_{i}$.\\
However, the penalty term so described is not invariant with respect to ``$x$".\\
This aspect of the problem was discussed in  \textbf{Chapter 9}.\\
In order to have the property of invariance, one can use, to express the discrepancy between ``$x$" and ``$Tx$", the divergences $LAI$ or $LBI$, which are invariant to both arguments and which are given respectively by (\ref{eq.LAIlap}) and (\ref{eq.LBIlap}).\\
The corresponding gradients with respect to ``$x$" are given respectively by (\ref{eq.gradLAIlap}) and (\ref{eq.gradLBIlap}).\\
We use here their expressions; for $LAI_{a}\left(x\|Tx\right)$, we have (\ref{eq.gradLAIlap}): 
\begin{align}
\frac{\partial LAI_{a}\left(x\|Tx\right)}{\partial x_{l}}=
&\frac{1}{a}\left[\frac{\sum_{i}T_{il}}{\sum_{i}\left(Tx\right)_{i}}-\frac{\sum_{i}T_{il}x_{i}^{a}\left(Tx\right)_{i}^{-a}}{\sum_{i}x_{i}^{a}\left(Tx\right)_{i}^{1-a}}\right]+\nonumber \\
&\frac{1}{a-1}\left[\frac{x_{l}^{a-1}\left(Tx\right)_{l}^{1-a}}{\sum_{i}x_{i}^{a}\left(Tx\right)_{i}^{1-a}}-\frac{1}{\sum_{i}x_{i}}\right]
\end{align}
The case $a=1$ which corresponds to the basis of a Kullback-Leibler divergence imply the computation of a limit which leads after some calculations to the expression:
\begin{equation}
\frac{\partial LAI_{1}\left(x\|Tx\right)}{\partial x_{l}}=\frac{1}{\sum_{i}x_{i}}\left[\sum_{i}T_{il}-\sum_{i}T_{il}\frac{x_{i}}{\left(Tx\right)_{i}}\right]
\label{eq.gradLAI1lap}	
\end{equation}
Or else by simplifying, taking into account the properties of the matrix $T$:
\begin{equation}
\frac{\partial LAI_{1}\left(x\|Tx\right)}{\partial x_{l}}=\frac{1}{\sum_{i}x_{i}}\left[1-\sum_{i}T_{il}\frac{x_{i}}{\left(Tx\right)_{i}}\right]
\label{eq.gradLAI1laps}	
\end{equation}
One can observe that, as always with this type of divergence,\\
$\sum_{l}x_{l}\frac{\partial }{\partial x_{l}}=0$.\\

Similarly, as regards $ LBI_{b}\left(x\|Tx\right)$, we have (\ref{eq.gradLBIlap}):
\begin{align}
\frac{\partial LBI_{b}\left(x\|Tx\right)}{\partial x_{l}}=
&\frac{1}{b-1}\left[\frac{x_{l}^{b-1}}{\sum_{i}x_{i}^{b}}-\frac{\left(Tx\right)_{l}^{b-1}}{\sum_{i}x_{i}\left(Tx\right)_{i}^{b-1}}\right]+\nonumber\\
&\left[\frac{\sum_{i}T_{il}\left(Tx\right)_{i}^{b-1}}{\sum_{i}\left(Tx\right)_{i}^{b}}-\frac{\sum_{i}T_{il}x_{i}\left(Tx\right)_{i}^{b-2}}{\sum_{i}x_{i}\left(Tx\right)_{i}^{b-1}}\right]
\end{align}
The special case $b=1$ leads to the same result as the divergence $LAI_{1}$, whereas the case $b=2$ which corresponds to a basis divergence which is the mean square deviation, does not imply a passage to the limit and leads to:
\begin{equation}
\frac{\partial LBI_{2}\left(x\|Tx\right)}{\partial x_{l}}=\frac{x_{l}}{\sum_{i}x^{2}_{i}}-\frac{\sum_{i}T_{il}\left(Tx\right)_{i}}{\sum_{i}\left(Tx\right)^{2}_{i}}
\label{eq.gradLBI2lap}	
\end{equation}
Or else, given the properties of the matrix $T$ ($T=T^{T}$):
\begin{equation}
\frac{\partial LBI_{2}\left(x\|Tx\right)}{\partial x_{l}}=\frac{x_{l}}{\sum_{i}x^{2}_{i}}-\frac{\left[T\left(Tx\right)\right]_{l}}{\sum_{i}\left(Tx\right)^{2}_{i}}
\label{eq.gradLBI2laps}	
\end{equation}
We have, again, as always with this type of divergence:\\
$\sum_{l}x_{l}\frac{\partial }{\partial x_{l}}=0$.

\section{Regularized algorithm.}
The problem here is to minimize under constraints, with respect to the two unknowns ``$x$" and ``$h$", an invariant composite divergence of the form:
\begin{equation}
	DCI=DI\left(y\|\left[h,x\right]\right)+\gamma DRI\left(h\right)+\mu DRI\left(x\right)
\end{equation}
The coefficients ``$\gamma$" and ``$\mu$" are the positive regularization factors.\\
For the $DI\left(y\|\left[h,x\right]\right)$ divergence representing the attachment to the data, we can use any divergence invariant with respect to $q\equiv Hx\equiv Xh$.\\
In contrast, as far as the terms of regularization are concerned, there are two possible cases.\\
* either this term translates the discrepancy to a constant, which is the case for example of the mean square deviation proposed in the previous section, then an invariant form of this term i.e. (\ref{eq.MCreg}) is usable.\\
This is valid for any invariant divergence between a constant and the variable considered.\\
* or the variable is involved in the two arguments of the divergence expressing the penalty, this is the case of a Laplacian type regularization, then one must use invariant divergences with respect to its both arguments, as it is the case for divergences of type $LAI$ or $LBI$ developed in the preceding section.\\
\textbf{These considerations being taken into account, the algorithmic method developed in \textbf{section (12.2)} applies to the composite divergence $DCI$ without any change}.

\section{Blind deconvolution - Non-invariant divergences.}
In the previous section, it was shown that the use of invariant divergences allows us to obtain multiplicative algorithms while maintaining sum constraints on the unknowns, subject to a normalization step.\\
In order to propose such algorithms without resorting to this procedure, we develop here Blind Deconvolution algorithms using non-invariant divergences.
In order to take into account the sum constraints on ``$h$" and on ``$x$", we rely on the variable change method proposed in the \textbf{section 10.3.2}, and more specifically on the algorithm (\ref{eq.algocontraint}).\\
The general structure of the algorithm remains as proposed in the previous section, which consists of iterating alternately on the two unknowns.\\
Considering a non-invariant divergence:
\begin{equation}
	D\left(p\|q\right)=\sum_{i}d\left(p_{i}\|q_{i}\right)
\end{equation}
With $p_{i}\equiv y_{i}$ et $q_{i}\equiv \left(Hx\right)_{i}\equiv \left(Xh\right)_{i}$.\\

\textbf{* Iterations on ``$h$" - Non-multiplicative form.}\\

We have $x^{k}\ \Leftrightarrow\ X^{k}$ and $h^{k}$, we are searching for $h^{k+1}$.\\
After performing the variable change $h_{j}=\frac{t_{j}}{\sum_{n}t_{n}}$, and according to (\ref{eq.algocontraint}), we obtain the algorithm:
\begin{equation}
	h^{k+1}_{l}=h^{k}_{l}+\delta^{k}h^{k}_{l}\left\{\left[-\frac{\partial D\left(y\|X^{k}h^{k}\right)}{\partial h^{k}_{l}}\right]-\sum_{m}h^{k}_{m}\left[-\frac{\partial D\left(y\|X^{k}h^{k}\right)}{\partial h^{k}_{m}}\right]\right\}
	\label{eq.iterationh}
\end{equation}
We can easily verify that we have $\ \sum_{l}h^{k+1}_{l}=\sum_{l}h^{k}_{l}$; as a consequence, considering the initialization $\sum_{j}h^{0}_{j}=1$, we will have $\sum_{l}h^{k+1}_{l}=1$.\\

\textbf{* Iterations on ``$h$" - Multiplicative form.}\\

Unlike invariant divergences, the use of non-invariant divergences allows multiplicative forms of the algorithms to appear, without the need to use normalization.\\
This is done by shifting the components of the opposite of the gradient involved in (\ref{eq.iterationh}), in order to make them all positive.
\begin{equation}
	\left[-\frac{\partial D\left(y\|X^{k}h^{k}\right)}{\partial h^{k}}\right]_{d}=\left[-\frac{\partial D\left(y\|X^{k}h^{k}\right)}{\partial h^{k}}\right]-\underbrace{\min}_{l} \left[-\frac{\partial D\left(y\|X^{k}h^{k}\right)}{\partial h^{k}_{l}}\right]+\epsilon
\end{equation}
Obviously, given the properties of ``$h$'', such an offset will not change the algorithm.\\
We can then write an algorithm that will allow us to obtain a purely multiplicative form:
\begin{equation}
	h^{k+1}_{l}=h^{k}_{l}+\delta^{k}h^{k}_{l}\left\{\frac{\left[-\frac{\partial D\left(y\|X^{k}h^{k}\right)}{\partial h^{k}_{l}}\right]_{d}}{\sum_{m}h^{k}_{m}\left[-\frac{\partial D\left(y\|X^{k}h^{k}\right)}{\partial h^{k}_{m}}\right]_{d}}-1\right\}
	\label{eq.iterationhmult}
\end{equation}
The purely multiplicative form is obtained, as always, by taking a descent step size $\delta^{k}=1\ \forall k$, however, there is no guarantee of convergence of such an algorithm over which we no longer have any control.\\

\textbf{* Iterations on ``$x$"- Non-multiplicative form.}\\

We know $h^{k+1}\ \Leftrightarrow\ H^{k+1}$ and $x^{k}$, we want to obtain $x^{k+1}$.\\
After having carried out the change of variables $x_{j}=\frac{u_{j}}{\sum_{n}u_{n}}\sum_{i}y_{i}$, referring to (\ref{eq.algocontraint}), one obtains the algorithm:
\begin{equation}
	x^{k+1}_{l}=x^{k}_{l}+\beta^{k}x^{k}_{l}\left\{\sum_{i}y_{i}\left[-\frac{\partial D\left(y\|H^{k+1}x^{k}\right)}{\partial x^{k}}\right]_{l}-\sum_{m}x^{k}_{m}\left[-\frac{\partial D\left(y\|H^{k+1}x^{k}\right)}{\partial x^{k}}\right]_{m}\right\}
	\label{eq.iterationx}
\end{equation}

Taking into account that $\sum_{j}x_{j}=\sum_{j}y_{j}$, we will have $\sum_{l}x^{k+1}_{l}=\sum_{l}x^{k}_{l}=\sum_{j}y_{j}$.\\

\textbf{* Iterations on ``$x$" - Multiplicative form.}\\

As proposed for the variable ``$h$", a multiplicative algorithm can be obtained;
to do so, we shift the components of the opposite gradient in  involved in (\ref{eq.iterationx}), in order to make them all positive.
\begin{equation}
	\left[-\frac{\partial D\left(y\|H^{k+1}x^{k}\right)}{\partial x^{k}}\right]_{d}=\left[-\frac{\partial D\left(y\|H^{k+1}x^{k}\right)}{\partial x^{k}}\right]-\underbrace{\min}_{l} \left[-\frac{\partial D\left(y\|H^{k+1}x^{k}\right)}{\partial x^{k}_{l}}\right]+\epsilon
\end{equation}
Obviously, an offset like that will not change the algorithm.\\
We can then write an algorithm that will allow us to obtain a purely multiplicative form:
\begin{equation}
	x^{k+1}_{l}=x^{k}_{l}+\beta^{k}x^{k}_{l}\left\{\frac{\sum_{i}y_{i}\left[-\frac{\partial D\left(y\|H^{k+1}x^{k}\right)}{\partial x^{k}_{l}}\right]_{d}}{\sum_{m}x^{k}_{m}\left[-\frac{\partial D\left(y\|H^{k+1}x^{k}\right)}{\partial x^{k}_{m}}\right]_{d}}-1\right\}
	\label{eq.iterationxmult}
\end{equation}
The purely multiplicative form is obtained, as always, by taking a descent step size $\beta^{k}=1\ \forall k$, but there's no convergence guarantee for such an algorithm over which we no longer have any control.\\

\subsection{Regularization.}
In this case, with $q\equiv Hx\equiv Xh$, we have to minimize a composite divergence:
\begin{equation}
	DC=D\left(y\|q\right)+\gamma DR_{x}\left(x\right)+\mu DR_{h}\left(h\right)
\end{equation}
Decomposing the problem as described in the previous paragraph leads us to:\\
1 - For known "$x$", minimize with respect to "$h$", under constraints, a divergence of the form:
\begin{equation}
		DC_{h}=D\left(y\|Xh\right)+\gamma DR_{h}\left(h\right)
\end{equation}
2 - For known "$h$", minimize with respect to "$x$", under constraints, a divergence of the form:
\begin{equation}
	DC_{x}=D\left(y\|Hx\right)+\mu DR_{x}\left(x\right)
\end{equation}

\subsubsection{Expressions of regularization terms.}
For both "$h$" and "$x$", we use here a smoothness constraint regularization of the solution which is written:\\
 either in the form (\ref{eq.MCninv}):
\begin{equation}
	DRQ_{x}=\frac{1}{2}\sum_{i}\left(c-x_{i}\right)^{2}
\end{equation}
Where $c=\frac{\sum_{i}y_{i}}{N}$; "$N$" is the number of components of "$y$" and/or "$x$", i.e. the number of pixels in the measured and/or reconstructed image.\\

with:
\begin{equation}
	\frac{\partial DRQ_{x}}{\partial x_{l}}=x_{l}-c
\end{equation}
, or in the form (\ref{eq.Lapninv}):
\begin{equation}
	DRL_{x}=\frac{1}{2}\left\|Lx\right\|^{2}=\left\|x-Tx\right\|^{2}
\end{equation}
With, taking into account the symmetry of the matrix "$T$":
\begin{equation}
	\frac{\partial DRL_{x}}{\partial x_{l}}=x^{2}_{l}-2\left(Tx\right)_{l}+\left[T\left(Tx\right)\right]_{l}
\end{equation}

\subsubsection{Algorithms.}
The regularized algorithms are based on the expressions (\ref{eq.iterationh}) and (\ref{eq.iterationx}), replacing the unregularized divergence $D$ by the composite divergences respectively $DC_{h}$ and $DC_{x}$.\\
The expressions of the involved gradients are then written as follows:
\begin{equation}
\frac{\partial DC_{h}\left(y\|X^{k}h^{k}\right)}{\partial h^{k}_{l}}=\frac{\partial D\left(y\|X^{k}h^{k}\right)}{\partial h^{k}_{l}}+\gamma \frac{\partial DR_{h}\left(h^{k}\right)}{\partial h^{k}_{l}}	
\end{equation}
and
\begin{equation}
\frac{\partial DC_{x}\left(y\|H^{k+1}x^{k}\right)}{\partial x^{k}_{l}}=\frac{\partial D\left(y\|H^{k+1}x^{k}\right)}{\partial x^{k}_{l}}+\mu \frac{\partial DR_{x}\left(x^{k}\right)}{\partial x^{k}_{l}}	
\end{equation}
In these expressions, the regularization terms $DR$ are either $DRQ$ or $DRL$.\\
So the non-multiplicative regularized algorithms are written in the form:
\begin{equation}
	h^{k+1}_{l}=h^{k}_{l}+\delta^{k}h^{k}_{l}\left\{\left[-\frac{\partial DC_{h}}{\partial h^{k}_{l}}\right]-\sum_{m}h^{k}_{m}\left[-\frac{\partial DC_{h}}{\partial h^{k}_{m}}\right]\right\}
	\label{eq.iterationhconvreg}
\end{equation}
\begin{equation}
	x^{k+1}_{l}=x^{k}_{l}+\beta^{k}x^{k}_{l}\left\{\sum_{i}y_{i}\left[-\frac{\partial DC_{x}}{\partial x^{k}_{l}}\right]-\sum_{m}x^{k}_{m}\left[-\frac{\partial DC_{x}}{\partial x^{k}_{m}}\right]\right\}
	\label{eq.iterationxconvreg}
\end{equation}

\subsubsection{* Multiplicative forms.}
They are obtained by following the procedure indicated in the previous section; the components of the opposite of the gradients $\left[-\frac{\partial DC_{h}}{\partial h^{k}_{l}}\right]$ and $\left[-\frac{\partial DC_{x}}{\partial x^{k}_{l}}\right]$ are shifted in order to obtain components that are all positive, and we build the algorithms:
\begin{equation}
	h^{k+1}_{l}=h^{k}_{l}+\delta^{k}h^{k}_{l}\left\{\frac{\left[-\frac{\partial DC_{h}}{\partial h^{k}_{l}}\right]_{d}}{\sum_{m}h^{k}_{m}\left[-\frac{\partial DC_{h}}{\partial h^{k}_{m}}\right]_{d}}-1\right\}
\end{equation}
\begin{equation}
	x^{k+1}_{l}=x^{k}_{l}+\beta^{k}x^{k}_{l}\left\{\frac{\sum_{i}y_{i}\left[-\frac{\partial DC_{x}}{\partial x^{k}_{l}}\right]_{d}}{\sum_{m}x^{k}_{m}\left[-\frac{\partial DC_{x}}{\partial x^{k}_{m}}\right]_{d}}-1\right\}
\end{equation}
Purely multiplicative forms are obtained, as always, by taking the descent step size $\beta^{k}=\delta^{k}=1\ \forall k$, but there's no guarantee that such algorithms, over which we have no control, will converge.\\

\setcounter{table}{0}  \setcounter{equation}{0}  \setcounter{figure}{0} \setcounter{chapter}{0} \setcounter{section}{0}
 
\appendix

\chapter{Appendix 1 - Generalized Logarithm.}
 \textbf{A reminder on the Generalized Logarithm function.}\\ 

The function ``Generalized Logarithm" or ``Deformed Logarithm" (d-logarithm, Box-Cox transformation) is written as follows:
\begin{equation}
	\log_{d}\left(x\right)=\frac{x^{1-d}-1}{1-d}
\end{equation}
The limit cases are the function $\log x$ when $d\rightarrow 1$ and the linear function $\left(x-1\right)$ when $d\rightarrow 0$.\\
The inverse function is the ``Generalized exponential" or ``Deformed exponential" function:
\begin{equation}
	\exp_{d}\left(x\right)=\left[1+\left(1-d\right)x\right]^{\frac{1}{1-d}}
\end{equation}
This last function is obviously only defined if the quantity in square brackets is positive or zero.\\

\chapter{Appendix 2 - Means.}
 \textbf{Recalls on means.}\\ 

For a set of real numbers $a_{j}\geq 0 \ \ j=1...N$ and a set of weights $w_{j}\geq 0$ we define:\\

\textbf{A. Weighted arithmetic mean.}
\begin{equation}
	M_{a}=\frac{\sum_{j}w_{j}a_{j}}{\sum_{i=1}^{N}w_{i}}=\sum_{j}c_{j}a_{j}\ \ ;\ \ c_{j}=\frac{w_{j}}{\sum_{i=1}^{N}w_{i}}\ \ ;\ \sum_{j}c_{j}=1
\end{equation}
In the particular case where all the weights are equal, we have $c_{j}=\frac{1}{N}$ and the unweighted arithmetic mean is written as follows:
\begin{equation}
	M_{a}=\frac{\sum_{j}a_{j}}{N}
\end{equation}

\textbf{B. Weighted geometric mean.}
\begin{equation}
	M_{g}=\left[\prod_{j}a_{j}^{w_{j}}\right]^{\frac{1}{\sum_{i=1}^{N}w_{i}}}=\exp\left[\frac{1}{\sum_{i=1}^{N}w_{i}}\sum_{j}w_{j}\log a_{j}\right]
\end{equation}
Which can also be written:
\begin{equation}
	M_{g}=\prod_{j}a_{j}^{c_{j}}=\exp\left[\sum_{j}c_{j}\log a_{j}\right]\ \ ;\ \ c_{j}=\frac{w_{j}}{\sum_{i=1}^{N}w_{i}}\ \ ;\ \ \sum_{j}c_{j}=1 
\end{equation}
In the particular case where all the weights are equal, we have $c_{j}=\frac{1}{N}$ and the unweighted geometric mean is written as follows:
\begin{equation}
	M_{g}=\prod_{j}a_{j}^{\frac{1}{N}}=\left[\prod_{j}a_{j}\right]^{\frac{1}{N}}
\end{equation}

\textbf{C. Weighted harmonic mean.}\\

\begin{equation}
	M_{h}=\frac{\sum_{i=1}^{N}w_{i}}{\sum_{j}\frac{w_{j}}{a_{j}}}=\frac{1}{\sum_{j}\frac{c_{j}}{a_{j}}}\ \ ;\ \ c_{j}=\frac{w_{j}}{\sum_{i=1}^{N}w_{i}}\ \ ;\ \ \sum_{j}c_{j}=1
\end{equation}
In the particular case where all the weights are equal, we have $c_{j}=\frac{1}{N}$ and the unweighted harmonic mean is written as follows:
\begin{equation}
	M_{h}=\frac{N}{\sum_{j}\frac{1}{a_{j}}}
\end{equation}

\textbf{D. Generalized mean with exponent ``t" - Power mean - Hôlder mean.}\\

For non-zero real ``t", we define:
\begin{equation}
	M_{t}=\left[\frac{1}{N}\sum_{j}a_{j}^{t}\right]^{\frac{1}{t}}
\end{equation}
- $ t=1\ \ \Rightarrow$ Arithmetic mean\\ 
- $ t=2\ \ \Rightarrow$ Quadratic mean $M_{q}$ (Square root mean)\\ 
- $ t=-1\ \ \Rightarrow$ Harmonic mean\\ 
- $ t\rightarrow 0\ \ \Rightarrow$ one tends towards the geometric mean ( caution, it is a passage to the limit.\\
- $ t\rightarrow+\infty \ \ \Rightarrow$ one tends towards : $\max\left[a_{j}\right]$\\
- $ t\rightarrow-\infty \ \ \Rightarrow$ one tends towards : $\min\left[a_{j}\right]$\\

The weighted versions are obtained by the relation:
\begin{equation}
	M_{t}=\left[\frac{1}{\sum_{i=1}^{N}w_{i}}\sum_{j}w_{j}a_{j}^{t}\right]^{\frac{1}{t}}
\end{equation}
or also:
\begin{equation}
	M_{t}=\left[\sum_{j}c_{j}a_{j}^{t}\right]^{\frac{1}{t}}\ \ ;\ \ c_{j}=\frac{w_{j}}{\sum_{i=1}^{N}w_{i}}\ \ ;\ \ \sum_{j}c_{j}=1
\end{equation}
An interesting inequality: if $t1<t2$ we have $M_{t1}\leq M_{t2}$.\\

This allows us to recover the inequalities between the different means:\\
\begin{equation}
	\min\leq M_{h}\leq M_{g}\leq M_{a}\leq M_{q}\leq \max
\end{equation}

This notion can be extended to the generalized mean in the meaning of a function ``$f$".\\

\textbf{E. Generalized ``$f$" mean.}\\

Let ``$f$" be an injective and continuous function, the weighted average in the meaning of the function ``$f$" is defined by:
\begin{equation}
	M_{f}=f^{-1}\left[\frac{\sum_{j}w_{j}f\left(a_{j}\right)}{\sum_{i=1}^{N}w_{i}}\right]
\end{equation}
or also:
\begin{equation}
	M_{f}=f^{-1}\left[\sum_{j}c_{j}f\left(a_{j}\right)\right]\ \ ;\ \ c_{j}=\frac{w_{j}}{\sum_{i=1}^{N}w_{i}}\ \ ;\ \ \sum_{j}c_{j}=1	
\end{equation}
If all the weights $w_{j}$ are equal, we obtain the unweighted version:
\begin{equation}
	M_{f}=f^{-1}\left[\frac{\sum_{j}f\left(a_{j}\right)}{N}\right]	
\end{equation}

\textbf{Remarks:}\\
- the function ``$f$" is assumed to be injective to ensure that ``$f^{-1}$" exists.\\
- the function ``$f$" is assumed to be continuous to ensure that the quantities $\left[\frac{\sum_{j}f\left(a_{j}\right)}{N}\right]$ or  $\left[\frac{\sum_{j}f\left(a_{j}\right)}{\sum_{i=1}^{N}w_{i}}\right]$ belong to the domain of ``$f^{-1}$".\\
- if ``$f$" is injective and continuous, it is strictly monotonous, which results in:
\begin{equation}
	\min\left[a_{j}\right]\leq M_{f}\leq \max\left[a_{j}\right]
\end{equation}
The various classical means can be considered as special cases:\\
- if $f\left(x\right)=x$ we recover the arithmetic mean.\\
- if $f\left(x\right)=\log x$ we recover the geometric mean.\\
- if $f\left(x\right)=\frac{1}{x}$ we recover the harmonic mean.\\
- if $f\left(x\right)=x^{t}$ we recover the Hölder mean  (generalized mean with exponent "t").\\

\textbf{Application to our problem:}\\
In our applications, we have 2 data fields ``$p$" and ``$q$" whose homologous points are noted $p_{i}$ and $q_{i}$; the means considered are made always between 2 homologous points, and are then extended to the whole field by simply summing the means between homologous points.\\
Therefore, in all relations concerning means, we always have N=2.\\
In the weighted means, we will use weights in the form of $\left(\alpha\right)$ and $\left(1-\alpha\right)$, $0<\alpha<1$.\\
And so we have:\\
* Weighted arithmetic mean:
\begin{equation}
	M_{a}\left(p,q\right)=\sum_{i}\left[\alpha p_{i}+\left(1-\alpha\right)q_{i}\right]
\end{equation}
* Weighted geometric mean:
\begin{equation}
	M_{g}\left(p,q\right)=\sum_{i}p_{i}^{\alpha}q_{i}^{1-\alpha}
\end{equation}
* Weighted harmonic mean:
\begin{equation}
	M_{h}\left(p,q\right)=\sum_{i}\frac{p_{i}q_{i}}{\left(1-\alpha\right)p_{i}+\alpha q_{i}}
\end{equation}
* Weighted quadratic mean:
\begin{equation}
	M_{q}\left(p,q\right)=\sqrt{\sum_{i}\alpha p_{i}^{2}+\left(1-\alpha\right)q_{i}^{2}}
\end{equation}

\chapter{Appendix 3 - Generalization of GHOSH.}
 \textbf{Generalization of GHOSH et al. \cite{ghosh2013}}\\

 With the change of denomination of the parameters, and considering that it is necessary to introduce the multiplicative factor already mentioned, the divergence of Ghosh et al. is written as that of BHHJ, by replacing in this one: $1\Leftrightarrow A$ and $\beta\Leftrightarrow B$, or equivalently, in our divergences: $\lambda\Leftrightarrow A+B$ and $\lambda-1 \Leftrightarrow B$.\\
This leads to the generalized divergence:
\begin{equation}	
GH(p\|q)=\frac{1}{A\left(A+B\right)}\left[\sum_{i}q_{i}^{A+B}+\frac{A}{B}\sum_{i}p_{i}^{A+B}-\frac{A+B}{B}\sum_{i}p_{i}^{A}q_{i}^{B}\right]
\end{equation}
The gradient with respect to $``q"$ will be written as follows:
\begin{equation}
\frac{\partial GH\left(p\|q\right)}{\partial q_{i}}=\frac{1}{A}q_{i}^{B-1}\left(q_{i}^{A}-p_{i}^{A}\right)
\end{equation}
The dual divergence will be written:
\begin{equation}	GH(q\|p)=\frac{1}{A\left(A+B\right)}\left[\sum_{i}p_{i}^{A+B}+\frac{A}{B}\sum_{i}q_{i}^{A+B}-\frac{A+B}{B}\sum_{i}q_{i}^{A}p_{i}^{B}\right]
\end{equation}
The gradient with respect to $``q"$ will be expressed as:
\begin{equation}
		\frac{\partial GH\left(q\|p\right)}{\partial q_{i}}=\frac{1}{B}q_{i}^{A-1}\left(q_{i}^{B}-p_{i}^{B}\right)
\end{equation}
If, from the divergence of Ghosh et al. we want to obtain the invariant form of this divergence by scale change on ``$q$", by applying the appropriate procedure, we obtain the factor $K$:
\begin{equation}
	K=\left(\frac{\sum_{i}p_{i}^{A}q_{i}^{B}}{\sum_{i}q_{i}^{A+B}}\right)^{\frac{1}{A}}     
\end{equation}
This leads after simplification to the invariant divergence:
\begin{equation}		
GHI(p\|q)=\frac{1}{B}\left[\sum_{i}p_{i}^{A+B}-\left(\frac{\sum_{i}p_{i}^{A}q_{i}^{B}}{\sum_{i}q_{i}^{A+B}}\right)^{\frac{A+B}{A}}\sum_{i}q_{i}^{A+B}\right]
\end{equation}
We can easily verify that if we replace ``$q$" by ``$Kq$", nothing changes, so we have well obtained the required result.\\
This expression is the difference of 2 positive terms (disregarding the constant multiplicative factor $\frac{1}{B}$), one can thus apply on each term of this difference the same increasing function, the ``Generalized Logarithm" for example, then the extreme form, the ``Log" which will lead to:
\begin{equation}
LGHI(p\|q)=\frac{1}{B}\log\left(\sum_{i}p_{i}^{A+B}\right)-\frac{A+B}{AB}\log\left(\sum_{i}p_{i}^{A}q_{i}^{B}\right)+\frac{1}{A}\log\left(\sum_{i}q_{i}^{A+B}\right)
\end{equation}
We can see, which is true in all cases, that this divergence is not only invariant by a change of scale on $``q"$, but also by a change of scale on $``p"$.\\
At this point, we can re-examine the dual Ghosh's divergence which is written:
\begin{equation}	
GH(q\|p)=\frac{1}{A\left(A+B\right)}\left[\sum_{i}p_{i}^{A+B}+\frac{A}{B}\sum_{i}q_{i}^{A+B}-\frac{A+B}{B}\sum_{i}q_{i}^{A}p_{i}^{B}\right]
\end{equation}
It is rendered invariant by scale change with respect to $``q"$ by introducing the invariance factor:
\begin{equation}
	K=\left(\frac{\sum_{i}q_{i}^{A}p_{i}^{B}}{\sum_{i}q_{i}^{A+B}}\right)^{\frac{1}{B}}
\end{equation}
We so obtain:
\begin{equation}	
GHI(q\|p)=\frac{1}{A\left(A+B\right)}\left[\sum_{i}p_{i}^{A+B}-\left(\frac{\sum_{i}q_{i}^{A}p_{i}^{B}}{\sum_{i}q_{i}^{A+B}}\right)^{\frac{A+B}{B}}\sum_{i}q_{i}^{A+B}\right]
\end{equation}
Hence, the Logarithmic form:
\begin{equation}	
LGHI(q\|p)=\frac{1}{A\left(A+B\right)}\left[\log\sum_{i}p_{i}^{A+B}-\frac{A+B}{B}\log\sum_{i}q_{i}^{A}p_{i}^{B}+\frac{A}{B}\log\sum_{i}q_{i}^{A+B}\right]
\end{equation}

\chapter{Appendix 4 - Gradient of invariant divergences.}
\textbf{Computation of the gradient of an invariant divergence; influence of the invariance factor.}\\

Considering a divergence of $D\left(p,q\right)$, which we make invariant with respect to ``$q$" by introducing the invariance factor $K\left(p,q\right)$ , we attempt here to answer the following question:\\
What happens to the expression of the gradient with respect to ``$q$" depending on whether we use the \textit{``nominal"} invariance factor $K_{0}$, or another acceptable expression of $K$, in this case $K_{1}=\frac{\sum_{i}p_{i}}{\sum_{i}q_{i}}$ ?\\
The two cases to be distinguished are therefore:\\

* Either $K=K_{0}$, where $K_{0}$ is the ``\textit{nominal}" invariance factor , obtained by solving the equation:
\begin{equation}
	\sum_{i}\frac{\partial D\left(p_{i}\|Kq_{i}\right)}{\partial K}=0
\end{equation}
Then, the gradient $\frac{\partial D\left(p\|Kq\right)}{\partial q_{j}}$ is given by (\ref{eq.GradDpKQ3S}).\\ 
\begin{equation}
\frac{\partial D\left(p\|K_{0}q\right)}{\partial q_{j}}=K_{0}	\left[\frac{\partial D\left(p\|Kq\right)}{\partial \left(Kq_{j}\right)}\right]_{K=K_{0}}
\label{eq.GradDpKQ3S2}
\end{equation}

* Or, $K=K_{1}$, in this case, the gradient $\frac{\partial D\left(p\|Kq\right)}{\partial q_{j}}$ is given by (\ref{eq.GradDpKQ3C})
\begin{equation}
\frac{\partial D\left(p\|Kq\right)}{\partial q_{j}}=\frac{\partial K_{1}}{\partial q_{j}}\left[\sum_{i}\frac{\partial D\left(p_{i}\|Kq_{i}\right)}{\partial K}\right]_{K=K_{1}}+K_{1}\left[\frac{\partial D\left(p\|Kq\right)}{\partial\left(Kq_{j}\right)}\right]_{K=K_{1}}
\end{equation}
Note that this expression can also be written as follows:
\begin{equation}
\frac{\partial D\left(p\|Kq\right)}{\partial q_{j}}=\frac{\partial K_{1}}{\partial q_{j}}\left[\sum_{i}q_{i}\frac{\partial D\left(p_{i}\|Kq_{i}\right)}{\partial \left(Kq_{i}\right)}\right]_{K=K_{1}}+K_{1}\left[\frac{\partial D\left(p\|Kq\right)}{\partial\left(Kq_{j}\right)}\right]_{K=K_{1}}
\end{equation}
We observe that when $K\neq K_{0}$, the first term of the second member of this expression is no longer zero.

In both cases, we still have:
\begin{equation}
	\sum_{j}q_{j}\frac{\partial D\left(p\|Kq\right)}{\partial q_{j}}=0
\end{equation}
If we consider the case of a linear model, that is to say $q_{j}=\left(Hx\right)_{j}$, then:
\begin{equation}
	\frac{\partial q_{j}}{\partial x_{l}}=h_{j l}\ \ \Rightarrow\ \ \frac{\partial D\left(p\|Kq\right)}{\partial x_{l}}=\sum_{j}\frac{\partial D\left(p\|Kq\right)}{\partial q_{j}} \frac{\partial q_{j}}{\partial x_{l}}
\end{equation}
and:
\begin{equation}
\frac{\partial D\left(p\|Kq\right)}{\partial x_{l}}=\sum_{j}h_{j l}\frac{\partial K}{\partial q_{l}}\left[\sum_{i}q_{i}\frac{\partial D\left(p_{i}\|Kq_{i}\right)}{\partial \left(Kq_{i}\right)}\right]+K\sum_{j}h_{j l}\left[\frac{\partial D\left(p\|Kq\right)}{\partial\left(Kq_{l}\right)}\right]
\end{equation}
If $K=K_{0}$, the first term will be zero.\\
 
\textbf{Applications to some divergences.}\\
If we use:
\begin{equation}
	K=K1=\frac{\sum_{j}p_{j}}{\sum_{j}q_{j}}\ \ \ \Rightarrow\ \ \ \frac{\partial K}{\partial q_{l}}=\frac{\partial K_{1}}{\partial q_{l}}=-\frac{\sum_{j}p_{j}}{\left(\sum_{j}q_{j}\right)^{2}}
\end{equation}
In a first step, it comes:
\begin{equation}
\frac{\partial D\left(p\|K_{1}q\right)}{\partial x_{l}}=\frac{\sum_{j}p_{j}}{\sum_{j}q_{j}}\left\{-\frac{1}{\sum_{j}q_{j}}\left[\sum_{i}q_{i}\frac{\partial D\left(p\|K_{1}q\right)}{\partial \left(K_{1}q_{i}\right)}\right]\left[H^{T}1\right]_{l}+\left[H^{T}\frac{\partial D\left(p\|K_{1}q\right)}{\partial\left(K_{1}q\right)}\right]_{l}\right\}
\label {eq.Basegrad}	
\end{equation}
We can now specify the divergence under consideration:\\

\textbf{1 - Mean Square Deviation.}\\
We have here:
\begin{equation}
 D\left(p\|Kq\right)=MC\left(p\|Kq\right)=\frac{1}{2}\sum_{i}\left(p_{i}-Kq_{i}\right)^{2}	
\end{equation}
Using the nominal invariance factor: $K=K_{0}=\frac{\sum_{i}p_{i}q_{i}}{\sum_{i}q^{2}_{i}}$, we obtain:
\begin{equation}
\frac{\partial MC\left(p\|K_{0}q\right)}{\partial q_{j}}=-K_{0}	\left(p_{j}-K_{0}\:q_{j}\right)=-\frac{\sum_{i}p_{i}q_{i}}{\sum_{i}q^{2}_{i}}	\left(p_{j}-\frac{\sum_{i}p_{i}q_{i}}{\sum_{i}q^{2}_{i}}\:q_{j}\right)
\end{equation}
and it comes:
\begin{equation}
\frac{\partial MC\left(p\|K_{0}q\right)}{\partial x_{l}}=-\frac{\sum_{i}p_{i}q_{i}}{\sum_{i}q^{2}_{i}}\left[	H^{T}\left(p-\frac{\sum_{i}p_{i}q_{i}}{\sum_{i}q^{2}_{i}}\:q\right)\right]_{l}	
\end{equation}
With another invariance factor $K_{1}$, we have:
\begin{equation}
	MC\left(p\|K_{1}q\right)=\frac{1}{2}\sum_{j}\left(p_{j}-K_{1}q_{j}\right)^{2}\ \ \Rightarrow\ \ \frac{\partial MC\left(p\|K_{1}q\right)}{\partial\left(K_{1}q_{i}\right)}=-\left(p_{i}-K_{1}q_{i}\right)
\end{equation}
from where, in a first step, with $K_{1}=\frac{\sum_{j}p_{j}}{\sum_{j}q_{j}}$, we have :
\begin{equation}
\frac{\partial MC\left(p\|K_{1}q\right)}{\partial x_{l}}=\frac{\sum_{j}p_{j}}{\sum_{j}q_{j}}\left\{\frac{1}{\sum_{j}q_{j}}\left[\sum_{i}q_{i}\left(p_{i}-K_{1}q_{i}\right)\right]\left[H^{T}1\right]_{l}-\left[H^{T}\left(p-K_{1}q\right)\right]_{l}\right\}	
\end{equation}
so, all simplifications made:
\begin{equation}
\frac{\partial MC\left(p\|K_{1}q\right)}{\partial x_{l}}=\frac{\left(\sum_{j}p_{j}\right)^{2}}{\sum_{j}q_{j}}\left\{\left[\sum_{i}\bar{q_{i}}\left(\bar{p_{i}}-\bar{q_{i}}\right)\right]\left[H^{T}1\right]_{l}-\left[H^{T}\left(\bar{p}-\bar{q}\right)\right]_{l}\right\}	
\end{equation}\\

\textbf{2 - Kullback-Leibler divergence.}\\
In this case, $K_{0}$ et $K_{1}$ are the same, so:
\begin{equation}	KL\left(p\|K_{1}q\right)=\sum_{i}p_{i}\log\frac{p_{i}}{K_{1}q_{i}}+K_{1}q_{i}-p_{i}=\sum_{i}d\left(p_{i}\|K_{1}q_{i}\right)
\end{equation}
\begin{equation}
	\frac{d\left(p_{i}\|K_{1}q_{i}\right)}{\partial\left(K_{1}q_{i}\right)}=-\frac{p_{i}}{K_{1}q_{i}}+1
\end{equation}
From the equation (\ref {eq.Basegrad}), we have:
\begin{equation}
\frac{\partial KL\left(p\|K_{1}q\right)}{\partial x_{l}}=\frac{\sum_{j}p_{j}}{\sum_{j}q_{j}}\left\{-\frac{1}{\sum_{j}q_{j}}\left[\sum_{i}q_{i}\left(-\frac{p_{i}}{K_{1}q_{i}}+1\right)\right]\left[H^{T}1\right]_{l}+\left[H^{T}\left(-\frac{p}{K_{1}q}+1\right)\right]_{l}\right\}	
\end{equation}
With $K_{1}=\frac{\sum_{j}p_{j}}{\sum_{j}q_{j}}=K_{0}$, we have, all simplifications made:
\begin{equation}
\frac{\partial KL\left(p\|K_{1}q\right)}{\partial x_{l}}=\frac{\sum_{j}p_{j}}{\sum_{j}q_{j}}\left\{-\underbrace{\left[\sum_{i}\bar{q_{i}}\left(1-\frac{\bar{p_{i}}}{\bar{q_{i}}}\right)\right]}_{=0}\left[H^{T}1\right]_{l}+\left[H^{T}\left(1-\frac{\bar{p}}{\bar{q}}\right)\right]_{l}\right\}	
\end{equation}
We can see that the first term of the expression between braces is zero; this is related to the fact that the invariance factor used here corresponds to the one that would be computed explicitly for this divergence, so it is in fact the \textit{``nominal"} invariance factor $K_{0}$ for the Kullback-Leibler divergence and the gradient has a simplified expression (\ref{eq.GradDpKQ3S}) (\ref{eq.GradDpKQ3S2}), then, finally:
\begin{equation}
\frac{\partial KL\left(p\|K_{1}q\right)}{\partial x_{l}}=\frac{\sum_{j}p_{j}}{\sum_{j}q_{j}}\left\{\left[H^{T}\left(1-\frac{\bar{p}}{\bar{q}}\right)\right]_{l}\right\}	
\end{equation}\\
	
\textbf{3 - Neyman Chi2.}\\
We have:
\begin{equation}	
\chi_{N}^{2}\left(p\|Kq\right)=\sum_{i}\frac{\left(p_{i}-Kq_{i}\right)^{2}}{Kq_{i}}=\sum_{i}d\left(p_{i}\|Kq_{i}\right)
\end{equation}
\begin{equation}
	\frac{\partial \chi_{N}^{2}\left(p\|Kq\right)}{\partial Kq_{i}}=1-\frac{p^{2}_{i}}{\left(Kq_{i}\right)^{2}}
	\label{eq.gradchi2}
\end{equation}
 From equation (\ref {eq.Basegrad}), it comes:
\begin{equation}
\frac{\partial \chi_{N}^{2}\left(p\|K_{1}q\right)}{\partial x_{l}}=\frac{\sum_{j}p_{j}}{\sum_{j}q_{j}}\left\{\left[H^{T}\left(1-\frac{p^{2}}{\left(K_{1}q\right)^{2}}\right)\right]_{l}-\frac{1}{\sum_{j}q_{j}}\left[\sum_{i}q_{i}\left(1-\frac{p^{2}_{i}}{\left(K_{1}q_{i}\right)^{2}}\right)\right]\left[H^{T}1\right]_{l}\right\}	
\end{equation}
With $K_{1}=\frac{\sum_{j}p_{j}}{\sum_{j}q_{j}}$, we have, all simplifications made:
\begin{equation}
\frac{\partial \chi_{N}^{2}\left(p\|K_{1}q\right)}{\partial x_{l}}=\frac{\sum_{j}p_{j}}{\sum_{j}q_{j}}\left\{\left[H^{T}\left(1-\frac{\bar{p^{2}}}{\bar{q^{2}}}\right)\right]_{l}-\left[\sum_{i}\bar{q_{i}}\left(1-\frac{\bar{p_{i}}^{2}}{\bar{q_{i}}^{2}}\right)\right]\left[H^{T}1\right]_{l}\right\}	
\end{equation}
With $K=K_{0}=\sqrt{\frac{\sum_{j}p^{2}_{j}q^{-1}_{j}}{\sum_{j}q_{j}}}$, we have, from equation (\ref{eq.gradchi2}) 
\begin{equation}
\frac{\partial \chi_{N}^{2}\left(p\|K_{0}q\right)}{\partial K_{0}q_{i}}=1-\left(\frac{p_{i}}{K_{0}q_{i}}\right)^{2}	
\end{equation}
Then, with (\ref{eq.GradDpKQ3S2}):
\begin{equation}
\frac{\partial \chi_{N}^{2}\left(p\|K_{0}q\right)}{\partial q_{j}}=	K_{0}\left[1-\left(\frac{p_{i}}{K_{0}q_{i}}\right)^{2}\right]
\end{equation}
And finally:
\begin{equation}
\frac{\partial \chi_{N}^{2}\left(p\|K_{0}q\right)}{\partial x_{l}}=	K_{0}\left[H^{T}.1\right]_{l}-\frac{1}{K_{0}}\left[H^{T}\left(\frac{p^{2}}{q^{2}}\right)\right]_{l}
\end{equation}

\chapter{Appendix 5 - Expressions of the invariance factor.}
\textbf{Here we give the expressions of the \textit{``nominal"} invariance factor $K_{0}\left(p,q\right)$ for some particular divergences.}\\

\textbf{General form:}
\begin{equation}
K_{0}\left(p,q\right)=\left(\frac{\sum_{j}p^{\alpha}_{j}q^{\beta}_{j}}{\sum_{j}p^{\delta}_{j}q^{\gamma}_{j}}\right)^{\mu} \ \ ;\ \  \alpha+\beta=\delta+\gamma \ \ ;\ \  \mu=\frac{1}{\gamma-\beta}	
\end{equation}
In a more synthetic form, one can write, taking into account the relations between the parameters:
\begin{equation}
K_{0}\left(p,q\right)=\left[\sum_{i}\left(\frac{p_{i}}{q_{i}}\right)^{\alpha-\delta}\frac{p^{\delta}_{i}q^{\gamma}_{i}}{\sum_{j}p^{\delta}_{j}q^{\gamma}_{j}}\right]^{\frac{1}{\alpha-\delta}}
\label{eq.Kgene3}	
\end{equation}
It is, for quantities of the form $\left(p_{i}/q_{i}\right)$, a weighted generalized mean of the order ``$t$" with the exponent $t=\alpha-\delta$, and weighting factors $w_{i}$ such that $\sum_{i}w_{i}=1$, which are written as follows:
\begin{equation}
	w_{i}=\frac{p^{\delta}_{i}q^{\gamma}_{i}}{\sum_{j}p^{\delta}_{j}q^{\gamma}_{j}}
	\label{eq.fpond}
\end{equation}\\
In the case of the dual Kullback-Leibler divergence, the invariance factor appears as a generalized mean based on a ``$\psi$'' function of the form:
\begin{equation}
	K_{0}=\psi^{-1}\left[\sum_{i}w_{i}\psi\left(\frac{p_{i}}{q_{i}}\right)\right]
\end{equation}
With $\psi\left(x\right)=\log x$.\\
The general form of the weighting factors (\ref{eq.fpond}) remains unchanged.\\
 
\textbf{1 - Kullback-Leibler divergence:}
\begin{equation}
K_{0}\left(p,q\right)=\frac{\sum _{j}p_{j}}{\sum _{j}q_{j}}	
\end{equation}
Here, we have: $\alpha=1, \beta=0, \gamma=1, \delta=0$,\\
then $\alpha-\delta=1, \gamma-\delta=1$.\\
\begin{equation}
	w_{i}=\frac{q_{i}}{\sum_{j}q_{j}}	
\end{equation}\\

\textbf{2 - Dual Kullback-Leibler divergence:}
\begin{equation}
	K_{0}=\psi^{-1}\left[\sum_{i}w_{i}\psi\left(\frac{p_{i}}{q_{i}}\right)\right]
\end{equation}
With $\psi\left(x\right)=\log x$.
\begin{equation}
	w_{i}=\frac{q_{i}}{\sum_{j}q_{j}}	
\end{equation}\\

\textbf{3 - Mean square deviation:}
\begin{equation}
K_{0}\left(p,q\right)=\frac{\sum _{j}p_{j}q_{j}}{\sum _{j}q^{2}_{j}}	
\end{equation}
We have here: $\alpha=1, \beta=1, \gamma=2, \delta=0$,\\
so, $\alpha-\delta=1, \gamma-\delta=2$.\\
\begin{equation}
	w_{i}=\frac{q^{2}_{i}}{\sum_{j}q^{2}_{j}}	
\end{equation}\\

\textbf{4 - Neyman's Chi2:}
\begin{equation}
K_{0}\left(p,q\right)=\sqrt{\frac{\sum _{j}p^{2}_{j}q^{-1}_{j}}{\sum _{j}q_{j}}}	
\end{equation}
Here: $\alpha=2, \beta=-1, \gamma=1, \delta=0$,\\
then $\alpha-\delta=2, \gamma-\delta=1$.\\
\begin{equation}
	w_{i}=\frac{q_{i}}{\sum_{j}q_{j}}	
\end{equation}\\

\textbf{5 - Pearson's Chi2:}
\begin{equation}
K_{0}\left(p,q\right)=\frac{\sum _{j}q_{j}}{\sum _{j}p^{-1}_{j}q^{2}_{j}}	
\end{equation}
We have here: $\alpha=0, \beta=1, \gamma=2, \delta=-1$,\\
then $\alpha-\delta=1, \gamma-\delta=3$.\\
\begin{equation}
	w_{i}=\frac{p^{-1}_{i}q^{2}_{i}}{\sum_{j}p^{-1}_{j}q^{2}_{j}}	
\end{equation}\\

\textbf{6 - Arithmetic-Geometric mean}
\begin{equation}
K_{0}\left(p,q\right)=\left(\frac{\sum_{j}\sqrt{p_{j}q_{j}}}{\sum_{j}q_{j}}\right)^{2}	
\end{equation}
We have: $\alpha=1/2, \beta=1/2, \gamma=1, \delta=0$,\\
then $\alpha-\delta=1/2, \gamma-\delta=1$.\\
\begin{equation}
	w_{i}=\frac{q^{3}_{i}}{\sum_{j}q^{3}_{j}}	
\end{equation}\\

\textbf{7 - Alpha divergence:}
\begin{equation}
K_{0}\left(p,q\right)=\left(\frac{\sum_{j}p^{\lambda_{a}}_{j}q^{1-\lambda_{a}}_{j}}{\sum_{j}q_{j}}\right)^{\frac{1}{\lambda_{a}}}	
\end{equation}
If $\lambda_{a}=1$ we recover $K_{0}$ corresponding to the K.L. divergence.\\
We have here: $\alpha=\lambda_{a}, \beta=1-\lambda_{a}, \gamma=1, \delta=0$,\\
so: $\alpha-\delta=\lambda_{a}, \gamma-\delta=1$.\\
\begin{equation}
	w_{i}=\frac{q_{i}}{\sum_{j}q_{j}}	
\end{equation}\\

\textbf{8 - Beta divergence:}
\begin{equation}
K_{0}\left(p,q\right)=\frac{\sum_{j}p_{j}q^{\lambda_{b}-1}_{j}}{\sum_{j}q^{\lambda_{b}}_{j}}	
\end{equation}
If $\lambda_{b}=1$ we recover $K_{0}$ corresponding to the K.L. divergence.\\
If $\lambda_{b}=2$ we recover $K_{0}$ corresponding to the mean square deviation.\\
Here: $\alpha=1, \beta=\lambda_{b}-1, \gamma=\lambda_{b}, \delta=0$,\\
then: $\alpha-\delta=1-\lambda_{b}, \gamma-\delta=\lambda_{b}$.\\
\begin{equation}
	w_{i}=\frac{q^{\lambda_{b}}_{i}}{\sum_{j}q^{\lambda_{b}}_{j}}	
\end{equation}\\

\chapter{Appendix 6 -  Effects of the invariance factor.}
\textbf{Here we show the effects of using an invariance factorsuch as $K=\frac{\sum_{j}p_{j}}{\sum_{j}q_{j}}$ on some classical divergences.}\\

These effects can be summarized by:\\

\textbf{Rule:\\
 Except for a multiplicative coefficient, which depends only on $\sum_{j}p_{j}$, the invariant divergence can be obtained by replacing in the initial divergence, the variables $p_{i}$ and $q_{i}$ by the normalized variables $\bar{p}_{i}=\frac{p_{i}}{\sum_{j}p_{j}}$ and $\bar{q}_{i}=\frac{q_{i}}{\sum_{j}q_{j}}$.}\\

\textbf{1 - Mean square deviation:}\\ 
We have the base divergence:
\begin{equation}
	EQM\left(p\|q\right)=\sum_{i}\left(p_{i}-q_{i}\right)^{2}
\end{equation}
With $K=\frac{\sum_{j}p_{j}}{\sum_{j}q_{j}}$, we have:
\begin{equation}
EQM\left(p\|Kq\right)=\sum_{i}\left(p_{i}-\frac{\sum_{j}p_{j}}{\sum_{j}q_{j}}q_{i}\right)^{2}	
\end{equation}
then:
\begin{equation}
EQMI\left(p\|q\right)=\left(\sum_{j}p_{j}\right)^{2}\sum_{i}\left(\bar{p}_{i}-\bar{q}_{i}\right)^{2}	
\end{equation}
which can be reduced to an invariant divergence with respect to both arguments:
\begin{equation}
EQMI\left(p\|q\right)=\sum_{i}\left(\bar{p}_{i}-\bar{q}_{i}\right)^{2}
\label{eq.EQMInorm}	
\end{equation}

\textbf{2 - Kullback-Leibler divergence:}\\
The initial divergence is written:
\begin{equation}
	KL\left(p\|q\right)=\sum_{i}p_{i}\log\frac{p_{i}}{q_{i}}+q_{i}-p_{i}
\end{equation}
With $K=\frac{\sum_{j}p_{j}}{\sum_{j}q_{j}}$, it comes:
\begin{equation}	KL\left(p\|Kq\right)=\sum_{i}p_{i}\log\frac{p_{i}}{\frac{\sum_{j}p_{j}}{\sum_{j}q_{j}}q_{i}}+\frac{\sum_{j}p_{j}}{\sum_{j}q_{j}}q_{i}-p_{i}
\end{equation}
That is:
\begin{equation}	KL\left(p\|Kq\right)=KLI\left(p\|q\right)=\sum_{j}p_{j}\left[\sum_{i}\bar{p}_{i}\log\frac{\bar{p}_{i}}{\bar{q}_{i}}+\bar{q}_{i}-\bar{p}_{i}\right]
\end{equation}
which can be reduced to an invariant divergence with respect to both arguments:
\begin{equation}	
KLI\left(p\|q\right)=\left[\sum_{i}\bar{p}_{i}\log\frac{\bar{p}_{i}}{\bar{q}_{i}}+\bar{q}_{i}-\bar{p}_{i}\right]=\left[\sum_{i}\bar{p}_{i}\log\frac{\bar{p}_{i}}{\bar{q}_{i}}\right]
\label{eq.KLInorm}
\end{equation}
It is a Kullback-Leibler between normalized variables.\\

\textbf{3 - Neyman's Chi2:}\\
The initial divergence is written:
\begin{equation}
	\chi^{2}_{N}\left(p\|q\right)=\sum_{i}\frac{\left(p_{i}-q_{i}\right)^{2}}{q_{i}}
\end{equation}
With $K=\frac{\sum_{j}p_{j}}{\sum_{j}q_{j}}$, we have:
\begin{equation}	\chi^{2}_{N}\left(p\|Kq\right)=\sum_{i}\frac{\left(p_{i}-\frac{\sum_{j}p_{j}}{\sum_{j}q_{j}}q_{i}\right)^{2}}{\frac{\sum_{j}p_{j}}{\sum_{j}q_{j}}q_{i}}
\end{equation}
That is:
\begin{equation}
	\chi^{2}_{N}I\left(p\|q\right)=\sum_{j}p_{j}\sum_{i}\frac{\left(\bar{p}_{i}-\bar{q}_{i}\right)^{2}}{\bar{q}_{i}}
\end{equation}
which can be reduced to an invariant divergence with respect to both arguments:
\begin{equation}
	\chi^{2}_{N}I\left(p\|q\right)=\sum_{i}\frac{\left(\bar{p}_{i}-\bar{q}_{i}\right)^{2}}{\bar{q}_{i}}
	\label{eq.chi2Nnorm}
\end{equation}

\textbf{4 - Pearson's Chi2:}\\
The initial divergence is written:
\begin{equation}
	\chi^{2}_{P}\left(p\|q\right)=\sum_{i}\frac{\left(q_{i}-p_{i}\right)^{2}}{p_{i}}
\end{equation}
With $K=\frac{\sum_{j}p_{j}}{\sum_{j}q_{j}}$, we obtain:
\begin{equation}
	\chi^{2}_{P}I\left(p\|q\right)=\sum_{i}\frac{\left(\frac{\sum_{j}p_{j}}{\sum_{j}q_{j}}q_{i}-p_{i}\right)^{2}}{p_{i}}
\end{equation}
That is:
\begin{equation}
	\chi^{2}_{P}I\left(p\|q\right)=\sum_{j}p_{j}\sum_{i}\frac{\left(\bar{q}_{i}-\bar{p}_{i}\right)^{2}}{\bar{p}_{i}}
\end{equation}
which can be reduced to an invariant divergence with respect to both arguments:
\begin{equation}
	\chi^{2}_{P}I\left(p\|q\right)=\sum_{i}\frac{\left(\bar{q}_{i}-\bar{p}_{i}\right)^{2}}{\bar{p}_{i}}
	\label{chi2Pnorm}
\end{equation}

\chapter{Appendix 7 - Inequalities between invariant divergences.}

\textbf{We show here, for some classic divergences denoted $D\left(p,q\right)$, that if $K_{0}\left(p,q\right)$ is the \textit{``nominal"} invariance factor, and if $K_{1}\left(p,q\right)\neq K_{0}\left(p,q\right)$ is another possible invariance factor, we have:
\begin{equation}
	D\left(p\|K_{0}q\right)\leq D\left(p\|K_{1}q\right)
\end{equation}
This is an obvious consequence of the definition of the \textit{``nominal"} invariance factor.\\
On the other hand, we have the relationship:
\begin{equation}
	D\left(p\|K_{1}q\right)-D\left(p\|K_{0}q\right)=D\left(K_{0}\|K_{1}\right)
\end{equation}
This is the basic divergence taken between the scalar quantities $K_{0}$ and $K_{1}$.}\\ 

\textbf{1 - Mean square deviation:}\\
We have here:
\begin{equation}
MC\left(p\|K_{0}q\right)=\sum_{i}\left(p_{i}-K_{0}q_{i}\right)^{2},\ \ \ MC\left(p\|K_{1}q\right)=\sum_{i}\left(p_{i}-K_{1}q_{i}\right)^{2}	
\end{equation}
We can write:
\begin{equation}
MC_{K_{1}}-MC_{K_{0}}=\left(K_{1}-K_{0}\right)\left[\left(K_{0}+K_{1}\right)\left(\sum_{i}q^{2}_{i}\right)-2\sum_{i}p_{i}q_{i}\right]	
\end{equation}
For this divergence the nominal invariance factor $K_{0}$ is computed analytically and is written as follows:
\begin{equation}
K_{0}=\frac{\sum_{j}p_{j}q_{j}}{\sum_{j}q^{2}_{j}}	
\end{equation}
This immediately leads to:
\begin{equation}
MC_{K_{1}}-MC_{K_{0}}=\sum_{i}q^{2}_{i}\left(K_{1}-K_{0}\right)^{2}\geq 0	
\end{equation}
Consequently:
\begin{equation}
	MC\left(p\|K_{0}q\right)\leq MC\left(p\|K_{1}q\right)
\end{equation}
With an equality (at $0$ moreover) for $p_{i}=q_{i}\ \forall i$, with then: $K_{0}=K_{1}=1$.\\
We also note, moreover, that:
\begin{equation}
	MC\left(p\|K_{1}q\right)-MC\left(p\|K_{0}q\right)\approx MC\left(K_{0}\|K_{1}\right)
\end{equation}

\textbf{2 - Kullback-Leibler divergence:}\\
We have here:
\begin{equation}
KL\left(p\|K_{0}q\right)=\sum_{i}\left(p_{i}\log\frac{p_{i}}{K_{0}q_{i}}+K_{0}q_{i}-p_{i}\right) 	
\end{equation}
\begin{equation}
	KL\left(p\|K_{1}q\right)=\sum_{i}\left(p_{i}\log\frac{p_{i}}{K_{1}q_{i}}+K_{1}q_{i}-p_{i}\right)
\end{equation}
We can write:
\begin{equation}
KL_{K_{1}}-KL_{K_{0}}=\sum_{i}p_{i}\left[\log\frac{K_{0}}{K_{1}}+\left(K_{1}-K_{0}\right)\frac{\sum_{i}q_{i}}{\sum_{i}p_{i}}\right]	
\end{equation}
For this divergence the nominal invariance factor $K_{0}$ is derived analytically and is written as follows:
\begin{equation}
K_{0}=\frac{\sum_{j}p_{j}}{\sum_{j}q_{j}}	
\end{equation}
This allows to write:
\begin{equation}
KL_{K_{1}}-KL_{K_{0}}=\sum_{i}q_{i}\left[K_{0}\log\frac{K_{0}}{K_{1}}+K_{1}-K_{0}\right]	
\end{equation}
So $\left[KL_{KL_{1}}-KL_{K_{0}}\right]$ is, with the exception of the multiplicative factor $\sum_{i}q_{i}$, equal to the Kullback-Leibler divergence between the scalars factors $K_{0}$ et $K_{1}$ denoted $KL\left(K_{0}\|K_{1}\right)\geq 0$.\\ 

Consequently:
\begin{equation}
KL_{K_{0}}\leq KL_{K_{1}}	
\end{equation}
With equality (to zero) if $p_{i}=q_{i}\ \forall i$, i.e. equivalently $K_{0}=K_{1}=1$.\\

\textbf{3 - Pearson's  Chi2:}\\
In this case, we have:
\begin{equation}
	\chi^{2}_{P}\left(p\|K_{0}q\right)=\sum_{i}\frac{\left(p_{i}-K_{0}q_{i}\right)^{2}}{p_{i}},\ \ \   \chi^{2}_P{}\left(p\|K_{1}q\right)=\sum_{i}\frac{\left(p_{i}-K_{1}q_{i}\right)^{2}}{p_{i}} 
\end{equation}
We can write:
\begin{equation}
\chi^{2}_{P,K_{1}}-\chi^{2}_{P,K_{0}}=\left(K_{1}-K_{0}\right)\sum_{i}\frac{\left(K_{0}+K_{1}\right)q^{2}_{i}-2p_{i}q_{i}}{p_{i}}
\end{equation}
For this divergence the nominal invariance factor $K_{0}$ is derived analytically and is written as follows:
\begin{equation}
	K_{0}=\frac{\sum_{j}q_{j}}{\sum_{j}q^{2}_{j}p^{-1}_{j}}
\end{equation}
This leads to:
\begin{equation}
	\chi^{2}_{P,K_{1}}-\chi^{2}_{P,K_{0}}=\sum_{i}q_{i}\frac{\left(K_{0}-K_{1}\right)^{2}}{K_{0}}
\end{equation}
This expression is obviously positive, and therefore:
\begin{equation}
	\chi^{2}_{P,K_{0}}\leq\chi^{2}_{P,K_{1}}
\end{equation}
With equality (to zero) if $p_{i}=q_{i}\ \forall i$, then $K_{0}=K_{1}=1$.\\
It can be seen that with the exception of a multiplicative factor, we have:
\begin{equation}
\chi^{2}_{P,K_{1}}-\chi^{2}_{P,K_{0}}\approx \chi^{2}_{K_{0},K_{1}}	
\end{equation}

\textbf{4 - Neyman's Chi2:}\\
In this case, we have:
\begin{equation}
	\chi^{2}_{N}\left(p\|K_{0}q\right)=\sum_{i}\frac{\left(p_{i}-K_{0}q_{i}\right)^{2}}{K_{0}q_{i}},\ \ \   \chi^{2}_N{}\left(p\|K_{1}q\right)=\sum_{i}\frac{\left(p_{i}-K_{1}q_{i}\right)^{2}}{K_{1}q_{i}} 
\end{equation}
We can write:
\begin{equation}
\chi^{2}_{N,K_{1}}-\chi^{2}_{N,K_{0}}=\frac{\left(K_{0}-K_{1}\right)}{K_{0}K_{1}}\left(\sum_{i}p^{2}_{i}q^{-1}_{i}-K_{0}K_{1}\sum_{i}q_{i}\right)
\end{equation}
For this divergence the nominal invariance factor $K_{0}$ is derived analytically and is written as follows:
\begin{equation}
	K_{0}=\left(\frac{\sum_{j}p^{2}_{j}q^{-1}_{j}}{\sum_{j}q_{j}}\right)^{\frac{1}{2}}
\end{equation}
We obtain immediately:
\begin{equation}
\chi^{2}_{N,K_{1}}-\chi^{2}_{N,K_{0}}=\sum_{j}q_{j}\frac{\left(K_{0}-K_{1}\right)}{K_{0}K_{1}}\left(K^{2}_{0}-K_{0}K_{1}\right)=\sum_{j}q_{j}\frac{\left(K_{0}-K_{1}\right)^{2}}{K_{1}}	
\end{equation}
This expression is obviously positive, and therefore:
\begin{equation}
	\chi^{2}_{N,K_{0}}\leq\chi^{2}_{N,K_{1}}
\end{equation}
With equality (to zero) if $p_{i}=q_{i}\ \forall i$, then $K_{0}=K_{1}=1$.\\
As in the previous examples, we have:
\begin{equation}
	\chi^{2}_{N}\left(p\|K_{1}q\right)-\chi^{2}_{N}\left(p\|K_{0}q\right)\approx \chi^{2}_{N}\left(K_{0}\|K_{1}\right)
\end{equation}

\chapter{Appendix 8 - Order relations between divergences.}
\textbf{We show on a few examples that for a divergence $D\left(p\|q\right)$ rendered invariant by using its nominal factor $K_{0}$, that is $D\left(p\|K_{0} q\right)$, we have:
\begin{equation}
	 D\left(p\|K_{0}q\right)\leq D\left(p\|q\right)
\end{equation}
with equality to "$0$" if $p_{i}=q_{i}\ \ \forall i$.}\\
\textbf{This is a consequence of the definition of the nominal invariance factor, because $D\left(p\|q\right)$ corresponds to $K=1$.}\\

\textbf{1 - Mean square deviation.}\\
Here, we have:
\begin{equation}
MC\left(p\|q\right)=\sum_{i}\left(p_{i}-q_{i}\right)^{2},\ \ \ MC\left(p\|K_{0}q\right)=\sum_{i}\left(p_{i}-K_{0}q_{i}\right)^{2}	
\end{equation}
We can write:
\begin{equation}
MC-MC_{K_{0}}=\left(1-K_{0}\right)\left[\left(1+K_{0}\right)\left(\sum_{i}q^{2}_{i}\right)-2\sum_{i}p_{i}q_{i}\right]	
\end{equation}
With: 
\begin{equation}
	K_{0}=\frac{\sum_{i} p_{i}q_{i}}{\sum_{i}q^{2}_{i}}
\end{equation}
it comes:
\begin{equation}
MC-MC_{K_{0}}=\sum_{i}q^{2}_{i}\left(1-K_{0}\right)^{2}\geq 0	
\end{equation}

\textbf{2 - Kullback-Leibler divergence.}\\
We have:
\begin{equation}
KL\left(p\|q\right)=\sum_{i}\left(p_{i}\log\frac{p_{i}}{q_{i}}+q_{i}-p_{i}\right),\ \ \ KL\left(p\|K_{0}q\right)=\sum_{i}\left(p_{i}\log\frac{p_{i}}{K_{0}q_{i}}+K_{0}q_{i}-p_{i}\right)	
\end{equation}
We can write:
\begin{equation}
KL-KL_{K_{0}}=\sum_{i}q_{i}\left[\frac{\sum_{i}p_{i}}{\sum_{i}q_{i}}\log K_{0}+\left(1-K_{0}\right)\right]	
\end{equation}
For this divergence the nominal invariance factor $K_{0}$ is derived analytically and is written as follows:
\begin{equation}
K_{0}=\frac{\sum_{j}p_{j}}{\sum_{j}q_{j}}>0	
\end{equation}
Then:
\begin{equation}
KL-KL_{K_{0}}=\sum_{i}q_{i}\left(K_{0}\log K_{0}+1-K_{0}\right)\geq 0	
\end{equation}

\textbf{3 - Pearson's Chi2.}\\
In this case, we have:
\begin{equation}
	\chi^{2}_{P}\left(p\|q\right)=\sum_{i}\frac{\left(p_{i}-q_{i}\right)^{2}}{p_{i}},\ \ \   \chi^{2}_{P}\left(p\|K_{0}q\right)=\sum_{i}\frac{\left(p_{i}-K_{0}q_{i}\right)^{2}}{p_{i}} 
\end{equation}
We can write:
\begin{equation}
\chi^{2}_{P}-\chi^{2}_{P,K_{0}}=\left(1-K_{0}\right)\sum_{i}\frac{\left(1+K_{0}\right)q^{2}_{i}-2p_{i}q_{i}}{p_{i}}
\end{equation}
For this divergence the nominal invariance factor $K_{0}$ is derived analytically and is written as follows:
\begin{equation}
	K_{0}=\frac{\sum_{j}q_{j}}{\sum_{j}q^{2}_{j}p^{-1}_{j}}>0
\end{equation}
Then:
\begin{equation}
\chi^{2}_{P}-\chi^{2}_{P,K_{0}}=\left(\sum_{i}q^{2}_{i}p^{-1}_{i}\right)\left(1-K_{0}\right)\left[\left(1+K_{0}\right)-2K_{0}\right]
\end{equation}
Or, also:
\begin{equation}
\chi^{2}_{P}-\chi^{2}_{P,K_{0}}=\left(\sum_{i}q^{2}_{i}p^{-1}_{i}\right)\left(1-K_{0}\right)^{2}
\end{equation}
And:
\begin{equation}
\chi^{2}_{P}-\chi^{2}_{P,K_{0}}=\sum_{i}q_{i}\frac{1}{K_{0}}\left(1-K_{0}\right)^{2}\geq 0
\end{equation}

\textbf{4 - Neyman's Chi2.}\\
We have, in this case:
\begin{equation}
	\chi^{2}_{N}\left(p\|q\right)=\sum_{i}\frac{\left(p_{i}-q_{i}\right)^{2}}{q_{i}},\ \ \   \chi^{2}_{N}\left(p\|K_{0}q\right)=\sum_{i}\frac{\left(p_{i}-K_{0}q_{i}\right)^{2}}{K_{0}q_{i}} 
\end{equation}
We can write:
\begin{equation}
\chi^{2}_{N}-\chi^{2}_{N,K_{0}}=\left(1-\frac{1}{K_{0}}\right)\sum_{i}p^{2}_{i}q^{-1}_{i}+\left(1-K_{0}\right)\sum_{i}q_{i}
\end{equation}
Then:
\begin{equation}
\chi^{2}_{N}-\chi^{2}_{N,K_{0}}=\sum_{i}q_{i}\left[\frac{\sum_{i}p^{2}_{i}q^{-1}_{i}}{\sum_{i}q_{i}}\left(1-\frac{1}{K_{0}}\right)+\left(1-K_{0}\right)\right]
\end{equation}
For this divergence the nominal invariance factor $K_{0}$ is derived analytically and is written as follows:
\begin{equation}
	K_{0}=\left(\frac{\sum_{j}p^{2}_{j}q^{-1}_{j}}{\sum_{j}q_{j}}\right)^{\frac{1}{2}}>0
\end{equation}
Consequently, we have:
\begin{equation}
\chi^{2}_{N}-\chi^{2}_{N,K_{0}}=\sum_{i}q_{i}\left[K^{2}_{0}\left(1-\frac{1}{K_{0}}\right)+\left(1-K_{0}\right)\right]=\sum_{i}q_{i}\left(1-K_{0}\right)^{2}\geq 0
\end{equation}

\textbf{4 - Alpha Divergence.}\\
Using the expression Alpha divergence (\ref{eq.AC}), it comes:
\begin{equation}	AC-AC_{K_{0}}=\frac{1}{\lambda_{a}\left(\lambda_{a}-1\right)}\left[\sum_{i}p^{\lambda_{a}}_{i}q^{1-\lambda_{a}}_{i}\left(1-K^{1-\lambda_{a}}_{0}\right)-\left(1-K_{0}\right)\right]
\end{equation}
For this divergence the nominal invariance factor $K_{0}$ is derived analytically and is written as follows:
\begin{equation}
	K_{0}=\left(\frac{\sum_{j}p^{\lambda_{a}}_{j}q^{1-\lambda_{a}}_{j}}{\sum_{j}q_{j}}\right)^{\frac{1}{\lambda_{a}}}>0
\end{equation}
Finally, we have:
\begin{equation}	AC-AC_{K_{0}}=\frac{\sum_{i}q_{i}}{\lambda_{a}\left(\lambda_{a}-1\right)}\left[K^{\lambda_{a}}_{0}-\lambda_{a}K_{0}-1+\lambda_{a}\right]
\end{equation}
That expression is always positive if $K_{0}\geq 0$.\\
The specific case of ${\lambda_{a}\rightarrow 0}$ leads to:
\begin{equation}
K_{0}-\log K_{0}-1\geq 0	
\end{equation}
In a similar way, $\lambda_{a}\rightarrow 1$ leads to:
\begin{equation}
K_{0}\log K_{0}-K_{0}+1\geq 0	
\end{equation}

\textbf{4 - Beta Divergence.}\\
Using the expression of the Beta divergence (\ref{eq.BC}), the derivation of the nominal invariance factor leads to:
\begin{equation}
	K_{0}=\frac{\sum_{i}p_{i}q_{i}^{\lambda_{b}-1}}{\sum_{i}q_{i}^{\lambda_{b}}}
\end{equation} 
and we obtain after a few calculations:
\begin{equation}	BC-BC_{K_{0}}=\frac{\sum_{i}q^{\lambda_{b}}_{i}}{\lambda_{b}\left(\lambda_{b}-1\right)}\left[K^{\lambda_{b}}_{0}-\lambda_{b}K_{0}-1+\lambda_{b}\right]\geq 0
\end{equation}
This expression is analogous to that obtained for the Alpha Divergence and the particular cases can be derived from it immediately.
 
\chapter{Appendix 9 - Regularization problems for the NMF.}
In this appendix, we indicate the difficulties that arise when constructing purely multiplicative regularized algorithms for NMF, as they have been proposed in the literature \cite{cichocki2009} \cite{pauca2006nonnegative} \cite{berry2007algorithms}.\\
Non-regularized multiplicative algorithms as they are proposed in the literature dedicated to NMF \cite{cichocki2009} \cite{lee2001algorithms} are immediately deduced from the SGM method described in Chapter 10.\\
All these algorithms only take into account the non-negativity constraint.\\
The method proposed in the literature can be summarized as follows:\\
For a divergence $D\left(Y\|HX\right)$  representing the data attachment term, the algorithmic method always consists, as previously indicated, of alternately performing a $H$ descent step, followed by a $X$ descent step.\\
Considering that we have obtained $H^{k+1}$ from $H^{k},\ X^{k}$, the basic algorithm can be written, with respect to the variable $X$ for example, with $X^{0}> 0$, in the form:
\begin{equation}
	X^{k+1}_{nm}=X^{k}_{nm}+\alpha^{k}X^{k}_{nm}\left[-\frac{\partial D}{\partial X^{k}}\right]_{nm}
\end{equation}
With $U^{k}_{nm}>0$ and $V^{k}_{nm}>0$, it is always possible to write:
\begin{equation}
\left[-\frac{\partial D}{\partial X^{k}}\right]_{nm}=U^{k}_{nm}-V^{k}_{nm}	
\end{equation}
This leads to:
\begin{equation}
	X^{k+1}_{nm}=X^{k}_{nm}+\alpha^{k}X^{k}_{nm}\left[U^{k}_{nm}-V^{k}_{nm}\right]
\end{equation}
From this, we deduce an algorithm which is at the base of purely multiplicative algorithms:
\begin{equation}
	X^{k+1}_{nm}=X^{k}_{nm}+\alpha^{k}X^{k}_{nm}\left[\frac{U^{k}_{nm}}{V^{k}_{nm}}-1\right]
\end{equation}
The initial direction of displacement $X^{k}\left[-\frac{\partial D}{\partial X^{k}}\right]$ has of course been modified, and becomes $\frac{X^{k}}{V_{k}}\left[-\frac{\partial D}{\partial X^{k}}\right]$, but remains a descent direction because $V^{k}$ is positive.
After which, by choosing a descent step size $\alpha^{k}=1\ \forall k$, we obtain the purely multiplicative form:
\begin{equation}
	X^{k+1}_{nm}=X^{k}_{nm}\left[\frac{U^{k}_{nm}}{V^{k}_{nm}}\right]
	\label{eq.mult}
\end{equation}
In the case where the term of attachment to the data is given by the mean square (RMS) deviation:
\begin{equation}
	D\left(Y\|HX\right)=\sum_{ij}\left[Y_{ij}-\left(HX\right)_{ij}\right]^{2}
\end{equation}
The multiplicative algorithm is written from (\ref{eq.mult}), with $H\equiv H^{k+1}$:
\begin{equation}
	X^{k+1}_{nm}=X^{k}_{nm}\left[\frac{\left(H^{T}Y\right)}{\left(H^{T}HX^{k}\right)}\right]_{nm}
\end{equation}
This algorithm was proposed by Lee and Seung \cite{lee2001algorithms}. It is the direct transposition of the ISRA algorithm \cite{daube1986}, well known in the field of deconvolution.\\
Similarly, for a data attachment term, which is the Kullback-Leibler divergence, we obtain a purely multiplicative algorithm, which is the immediate transposition of the Richardson-Lucy (RL) algorithm \cite{richardson1972} \cite{lucy1974} well known in the field of deconvolution, which is written, with $H\equiv H^{k+1}$:
\begin{equation}
	X^{k+1}_{nm}=\frac{X^{k}_{nm}}{\sum_{l}H^{T}_{nl}}\left[H^{T}\frac{Y}{HX^{k}}\right]_{nm}
\end{equation}
As in all descent algorithms, the convergence is ensured by a computation of the descent step size, so that with a fixed descent step size for all iterations, and equal to 1, nothing is guaranteed for convergence, even if non-negativity is ensured subject to the choice of a positive initial iterate.\\
In the case of deconvolution, the convergence of such algorithms has been demonstrated for the two divergences mentioned above.\\
In the case of NMF, multiplicative algorithms have been proposed for many of the divergences listed in the book by Cichocki et al. \cite{cichocki2009}.\\

\textbf{* Problems related to regularization.}\\

When we consider the regularized algorithms, things get a little worse.\\
For example, again with the $X$ variable, it is now a matter of minimizing under a non-negativity constraint, a composite divergence:
\begin{equation}
	J\left(Y,HX\right)=D\left(Y\|HX\right)+\mu\; DR\left(X\right)
\end{equation}
where $DR\left(X\right)$ is the penalty term, and ``$\mu$" the positive regularization factor.
Basically, knowing $X^{k}$ and $H^{k+1}$, the algorithm for the variable $X$ is written as follows:
\begin{equation}
	X^{k+1}_{nm}=X^{k}_{nm}+\alpha^{k}X^{k}_{nm}\left[-\frac{\partial D}{\partial X^{k}}+\mu\left(-\frac{\partial DR}{\partial X^{k}}\right)\right]_{nm}
\end{equation}
The method generally proposed in the literature consists of decomposing $-\frac{\partial D}{\partial X^{k}}$ in a difference of 2 positive terms $U^{k}$ and $V^{k}$: $-\frac{\partial D}{\partial X^{k}}=U^{k}-V^{k}$ as previously proposed, but on the other hand such a decomposition is not applied to the regularization term $-\frac{\partial DR}{\partial X^{k}}$, which means to implicitly consider that this term is always positive.\\

\textbf{This is false. It is only true for certain penalty functions (such as $DR\left(X\right)=\left\|X\right\|^{2}$ for example).}\\

From there, one thing leads to another; the proposed algorithm is first written:
\begin{equation}
	X^{k+1}_{nm}=X^{k}_{nm}+\alpha^{k}X^{k}_{nm}\left[U^{k}_{nm}-V^{k}_{nm}+\mu \left(-\frac{\partial DR}{\partial X^{k}}\right)\right]
\end{equation}
Then:
\begin{equation}
	X^{k+1}_{nm}=X^{k}_{nm}+\alpha^{k}X^{k}_{nm}\left[\frac{U^{k}_{nm}}{V^{k}_{nm}+\mu \left(\frac{\partial DR}{\partial X^{k}}\right)_{nm}}-1\right]
\end{equation}
And finally, by taking a descent step size equal to 1, we obtain a multiplicative algorithm which is written as follows:
\begin{equation}
	X^{k+1}_{nm}=X^{k}_{nm}+\alpha^{k}X^{k}_{nm}\left[\frac{U^{k}_{nm}}{V^{k}_{nm}+\mu \left(\frac{\partial DR}{\partial X^{k}}\right)_{nm}}\right]
\end{equation}
\textbf{It is absolutely clear (in principle) that all this is only valid if $\left[V^{k}_{nm}+\mu \left(\frac{\partial DR}{\partial X^{k}}\right)_{nm}\right]$ is strictly positive; otherwise we have transformed a descent algorithm into an ascent algorithm.....}\\

Authors who apply this type of algorithm prudently point out that if the regularization factor ``$\mu$" is small, everything is fine; it is just as clear that everything is fine regardless of ``$\mu$" if the penalty term is $DR\left(X\right)=\left\|X\right\|^{2}$, because in this case the gradient of the penalty term is positive.\\
It should be mentioned that this method is the transposition to the NMF of a rather old algorithm designated in the literature as One Step Late (OSL) by P.J.Green \cite{green1990bayesian}, in which the author clearly mentions the limits of use of his algorithm.\\

\textbf{In fact, in order for this to work in all cases, it is necessary to proceed as described in \textbf{Chapter 10} regarding the SGM method}.\\

\textbf{$\ast$ The opposite of the gradient of the penalty term must be decomposed into the form:}

\begin{equation}
	\left(-\frac{\partial DR}{\partial X^{k}}\right)_{nm}=UR^{k}_{nm}-VR^{k}_{nm}
\end{equation}
With $UR^{k}_{nm}\geq 0$ and $VR^{k}_{nm}\geq 0$.\\
 Then, we first write the algorithm in the form:
\begin{equation}
	X^{k+1}_{nm}=X^{k}_{nm}+\alpha^{k}X^{k}_{nm}\left[U^{k}_{nm}-V^{k}_{nm}+\mu\; UR^{k}_{nm}-\mu\; VR^{k}_{nm}+\epsilon-\epsilon\right]
\end{equation}
We have: $U^{k}_{nm}+\mu\; UR^{k}_{nm}+\epsilon>0$ and $V^{k}_{nm}+\mu\; VR^{k}_{nm}+\epsilon>0$.\\

From this, we deduce the algorithm that will be the basis of the purely multiplicative algorithms:
\begin{equation}
	X^{k+1}_{nm}=X^{k}_{nm}+\alpha^{k}X^{k}_{nm}\left[\frac{U^{k}_{nm}+\mu\; UR^{k}_{nm}+\epsilon}{V^{k}_{nm}+\mu\; VR^{k}_{nm}+\epsilon}-1\right]
\end{equation}
Which does not cause any problems with regard to the descent properties since $V^{k}_{nm}+\mu\; VR^{k}_{nm}+\epsilon>0$.\\
After that, a purely multiplicative algorithm is obtained by taking a descent step size equal to 1 for all the iterations:
\begin{equation}
	X^{k+1}_{nm}=X^{k}_{nm}\left[\frac{U^{k}_{nm}+\mu\; UR^{k}_{nm}+\epsilon}{V^{k}_{nm}+\mu\; VR^{k}_{nm}+\epsilon}\right]
\end{equation}
Obviously, the same reasoning applies for iterations on $H$, leading to $H^{k+1}$ knowing $H^{k}$ and $X^{k}$.\\ 
Nevertheless, the convergence of this type of algorithm is not demonstrated in all generality because of the particular choice of the descent step size.

\bibliographystyle{plain}

\begin{thebibliography}{}

\end{thebibliography}


\begin{thebibliography}{10}

\bibitem{ali1966}
S.M. Ali and S.D. Silvey.
\newblock A general class of coefficients of divergence of one distribution
  from another.
\newblock {\em Journal of the Royal Statistical Society. Series B
  (Methodological)}, pages 131--142, 1966.

\bibitem{amari2009}
S.I. Amari.
\newblock Alpha divergence is unique, belonging to both $f$-divergence and
  {B}regman divergence classes.
\newblock {\em Information Theory, IEEE Transactions on}, 55(11):4925--4931,
  2009.

\bibitem{andrews1977digital}
H.C. Andrews and B.~R. Hunt.
\newblock Digital image restoration.
\newblock {\em Prentice-Hall Signal Processing Series, Englewood Cliffs:
  Prentice-Hall, 1977}, 1977.

\bibitem{arimoto1971}
S.~Arimoto.
\newblock Information-theoretical considerations on estimation problems.
\newblock {\em Information and Control}, 19(3):181--194, 1971.

\bibitem{armijo1966}
L.~Armijo.
\newblock Minimization of functions having {L}ipschitz continuous first partial
  derivatives.
\newblock {\em Pacific Journal of mathematics}, 16(1):1--3, 1966.

\bibitem{arndt2001}
C.~Arndt.
\newblock {\em Information measures; information and its description in science
  and engineering}.
\newblock Springer, 2001.

\bibitem{ayers1988iterative}
G.R. Ayers and J.C. Dainty.
\newblock Iterative blind deconvolution method and its applications.
\newblock {\em Optics letters}, 13(7):547--549, 1988.

\bibitem{basseville1989}
M.~Basseville.
\newblock Distance measures for signal processing and pattern recognition.
\newblock {\em Signal processing}, 18(4):349--369, 1989.

\bibitem{basseville2013}
M.~Basseville.
\newblock Divergence measures for statistical data processing - an annotated
  bibliography.
\newblock {\em Signal Processing}, 93(4):621--633, 2013.

\bibitem{basseville1996}
M.~Basseville.
\newblock Information: entropies, divergences et moyennes.
\newblock Technical report, IRISA 1996.

\bibitem{basu1998}
A.~Basu, I.R. Harris, N.L. Hjort, and M.C. Jones.
\newblock Robust and efficient estimation by minimising a density power
  divergence.
\newblock {\em Biometrika}, 85(3):549--559, 1998.

\bibitem{ben1989}
A.~Ben-Tal, A.~Charnes, and M.~Teboulle.
\newblock Entropic means.
\newblock {\em Journal of Mathematical Analysis and Applications},
  139(2):537--551, 1989.

\bibitem{beran1977}
R.~Beran.
\newblock Minimum {H}ellinger distance estimates for parametric models.
\newblock {\em The Annals of Statistics}, pages 445--463, 1977.

\bibitem{bercher2008}
J-F. Bercher.
\newblock On some entropy functionals derived from {R}{\'e}nyi information
  divergence.
\newblock {\em Information Sciences}, 178(12):2489--2506, 2008.

\bibitem{berry2007algorithms}
M.~W. Berry, M.~Browne, A.~N. Langville, V.~P. Pauca, and R.~J. Plemmons.
\newblock Algorithms and applications for approximate nonnegative matrix
  factorization.
\newblock {\em Computational statistics \& data analysis}, 52(1):155--173,
  2007.

\bibitem{bertero1998}
M.~Bertero and P.~Boccacci.
\newblock {\em Introduction to inverse problems in imaging}.
\newblock CRC press, 1998.

\bibitem{bertsekas1995}
D.P. Bertsekas.
\newblock {\em Non linear programming}.
\newblock Athena scientific, Belmont, Mass., 1995.

\bibitem{bhatia2013}
P.K. Bhatia and S.~Singh.
\newblock On a new csisz\"ar\'{}s f-divergence measure.
\newblock {\em Cybernetics and information technologies}, 13(2):43--57, 2013.

\bibitem{boyd2004}
S.~Boyd and L.~Vandenberghe.
\newblock {\em Convex optimization}.
\newblock Cambridge university press, 2004.

\bibitem{bregman1967}
L.M. Bregman.
\newblock The relaxation method of finding the common point of convex sets and
  its application to the solution of problems in convex programming.
\newblock {\em USSR computational mathematics and mathematical physics},
  7(3):200--217, 1967.

\bibitem{burbea1982}
J.~Burbea and C.R. Rao.
\newblock On the convexity of higher order {J}ensen differences based on
  entropy functions (corresp.).
\newblock {\em Information Theory, IEEE Transactions on}, 28(6):961--963, 1982.

\bibitem{censor1997}
Y.~Censor and S.A. Zenios.
\newblock {\em Parallel optimization: Theory, algorithms, and applications}.
\newblock Oxford University Press on Demand, 1997.

\bibitem{cichocki2010}
A.~Cichocki and S.I. Amari.
\newblock Families of $\alpha$-$\beta$-and $\gamma$-divergences: Flexible and
  robust measures of similarities.
\newblock {\em Entropy}, 12(6):1532--1568, 2010.

\bibitem{cichocki2011}
A.~Cichocki, S.~Cruces, and S.I. Amari.
\newblock Generalized alpha-beta divergences and their application to robust
  nonnegative matrix factorization.
\newblock {\em Entropy}, 13(1):134--170, 2011.

\bibitem{cichocki2009}
A.~Cichocki, R.~Zdunek, A.H. Phan, and S.I. Amari.
\newblock {\em Nonnegative matrix and tensor factorizations: applications to
  exploratory multi-way data analysis and blind source separation}.
\newblock John Wiley \& Sons, 2009.

\bibitem{csiszar1974}
I.~Csisz\"ar.
\newblock Information measures: A critical survey.
\newblock In {\em Transactions of the Seventh Prague Conference on Information
  Theory, Statistical Decision Functions, Random Processes}, pages 73--86,
  1974.

\bibitem{csiszar1967}
I.~Csisz\"ar et~al.
\newblock Information-type measures of difference of probability distributions
  and indirect observations.
\newblock {\em Studia Sci. Math. Hungar.}, 2:299--318, 1967.

\bibitem{culioli2012introduction}
J.C. Culioli.
\newblock {\em Introduction {\`a} l'optimisation}.
\newblock Ed.Ellipses, 2012.

\bibitem{daube1986}
M.E. Daube-Witherspoon and G.~Muehllehner.
\newblock An iterative image space reconstruction algorthm suitable for volume
  {ECT}.
\newblock {\em Medical Imaging, IEEE Transactions on}, 5(2):61--66, 1986.

\bibitem{de1987}
A.R. De~Pierro.
\newblock On the convergence of the iterative image space reconstruction
  algorithm for volume {ECT}.
\newblock {\em IEEE TRANS. MED. IMAG.}, 6(2):174--175, 1987.

\bibitem{demoment1985deconvolution}
G.~Demoment.
\newblock D{\'e}convolution des signaux.
\newblock Technical report, Cours Supelec 1985.

\bibitem{dempster1977}
A.P. Dempster, N.M. Laird, and D.B. Rubin.
\newblock Maximum likelihood from incomplete data via the {EM} algorithm.
\newblock {\em Journal of the royal statistical society. Series B
  (methodological)}, pages 1--38, 1977.

\bibitem{dragomir2001}
S.S. Dragomir, V.~Gluscevic, and C.E.M. Pearce.
\newblock Approximations for the {C}sisz\"ar's f divergence via mid point
  inequalities.
\newblock {\em Inequality Theory and Applications}, 1:139--154, 2001.

\bibitem{eguchi2010}
S.~Eguchi and S.~Kato.
\newblock Entropy and divergence associated with power function and the
  statistical application.
\newblock {\em Entropy}, 12(2):262--274, 2010.

\bibitem{fevotte2011}
C.~F{\'e}votte and J.~Idier.
\newblock Algorithms for nonnegative matrix factorization with the
  $\beta$-divergence.
\newblock {\em Neural computation}, 23(9):2421--2456, 2011.

\bibitem{fujisawa2008}
H.~Fujisawa and S.~Eguchi.
\newblock Robust parameter estimation with a small bias against heavy
  contamination.
\newblock {\em Journal of Multivariate Analysis}, 99(9):2053--2081, 2008.

\bibitem{furuichi2012}
S.~Furuichi and F.C. Mitroi.
\newblock Mathematical inequalities for some divergences.
\newblock {\em Physica A: Statistical Mechanics and its Applications},
  391(1):388--400, 2012.

\bibitem{ghosh2013}
A.~Ghosh, I.R. Harris, A.~Maji, A.~Basu, and L.~Pardo.
\newblock A generalized divergence for statistical inference.
\newblock Technical report, BIRU/2013/3, Bayesian and Interdisciplinary
  Research Unit, Indian Statistical Institute, Kolkata, India, 2013.

\bibitem{green1990bayesian}
P.~J. Green.
\newblock Bayesian reconstructions from emission tomography data using a
  modified em algorithm.
\newblock {\em IEEE transactions on medical imaging}, 9(1):84--93, 1990.

\bibitem{grendar2010}
M.~Grend{\'a}r and R.K. Niven.
\newblock The {P}\'olya information divergence.
\newblock {\em Information Sciences}, 180(21):4189--4194, 2010.

\bibitem{hardy1952}
G.H. Hardy, J.E. Littlewood, and G.~P{\'o}lya.
\newblock {\em Inequalities}.
\newblock Cambridge university press, 1952.

\bibitem{havrda1967}
J.~Havrda and F.~Charv{\'a}t.
\newblock Quantification method of classification processes. concept of
  structural $ \alpha $-entropy.
\newblock {\em Kybernetika}, 3(1):30--35, 1967.

\bibitem{hellinger1909}
E.~Hellinger.
\newblock Neue begr{\"u}ndung der theorie quadratischer formen von
  unendlichvielen ver{\"a}nderlichen.
\newblock {\em Journal f{\"u}r die reine und angewandte Mathematik},
  136:210--271, 1909.

\bibitem{henze1999}
N.~Henze and M.D. Penrose.
\newblock On the multivariate runs test.
\newblock {\em Annals of statistics}, pages 290--298, 1999.

\bibitem{hildebrand1987}
F.B. Hildebrand.
\newblock {\em Introduction to numerical analysis}.
\newblock Courier Corporation, 1987.

\bibitem{hiriart2012}
J-B. Hiriart-Urruty and C.~Lemar{\'e}chal.
\newblock {\em Fundamentals of convex analysis}.
\newblock Springer Science \& Business Media, 2012.

\bibitem{hoyer2004}
P.O. Hoyer.
\newblock Non-negative matrix factorization with sparseness constraints.
\newblock {\em Journal of machine learning research}, 5(Nov):1457--1469, 2004.

\bibitem{Idier01a}
J.~Idier, editor.
\newblock {\em Approche bay\'esienne pour les probl\`emes inverses}.
\newblock Trait\'e IC2, S\'erie traitement du signal et de l'image, Herm\`es,
  Paris, nov. 2001.

\bibitem{itakura1968}
F.~Itakura and S.~Saito.
\newblock Analysis synthesis telephony based on the maximum likelihood method.
\newblock In {\em Proceedings of the 6th International Congress on Acoustics},
  volume~17, pages C17--C20. pp. C17--C20, 1968.

\bibitem{jansson1984deconvolution}
P.~A. Jansson.
\newblock {\em Deconvolution. With applications in spectroscopy}.
\newblock New York: Academic Press, 1984, edited by Jansson, Peter A., 1984.

\bibitem{jansson2014deconvolution2}
P.~A. Jansson.
\newblock {\em Deconvolution of images and spectra}.
\newblock Courier Corporation, 2014.

\bibitem{jeffreys1946}
H.~Jeffreys.
\newblock An invariant form for the prior probability in estimation problems.
\newblock In {\em Proceedings of the Royal Society of London A: Mathematical,
  Physical and Engineering Sciences}, volume 186, pages 453--461. The Royal
  Society, 1946.

\bibitem{jensen1906}
J.L. Jensen.
\newblock Sur les fonctions convexes et les in{\'e}galit{\'e}s entre les
  valeurs moyennes.
\newblock {\em Acta Mathematica}, 30(1):175--193, 1906.

\bibitem{kapur1967}
J.N. Kapur.
\newblock Generalized entropy of order $\alpha$ and type $\beta$.
\newblock In {\em The Math. Seminar}, volume~4, pages 78--82, 1967.

\bibitem{kapur1983}
J.N. Kapur.
\newblock Nonadditive measures of entropy and distributions of statistical
  mechanics.
\newblock {\em Indian J. Pure Appl. Math}, 14:1372--1387, 1983.

\bibitem{knockaert1993}
L.~Knockaert.
\newblock A class of statistical and spectral distance measures based on
  {B}ose-{E}instein statistics.
\newblock {\em IEEE Trans. on Signal Processing}, 11:3171--3174, 1993.

\bibitem{kullback1951}
S.~Kullback and R.A. Leibler.
\newblock On information and sufficiency.
\newblock {\em The annals of mathematical statistics}, 22(1):79--86, 1951.

\bibitem{lange1984}
K.~Lange, R.~Carson, et~al.
\newblock {EM} reconstruction algorithms for emission and transmission
  tomography.
\newblock {\em J Comput Assist Tomogr}, 8(2):306--16, 1984.

\bibitem{lanteri1994blind}
H.~Lanteri, C.~Aime, H.~Beaumont, and P.~Gaucherel.
\newblock Blind deconvolution using the richardson-lucy algorithm.
\newblock In {\em Optics in Atmospheric Propagation and Random Phenomena},
  volume 2312, pages 182--192. International Society for Optics and Photonics,
  1994.

\bibitem{lanteri1995comparison}
H.~Lanteri, M.~Barilli, H.~Beaumont, C.~Aime, P.~Gaucherel, and H.~Touma.
\newblock Comparison of several algorithms for blind deconvolution: analysis of
  noise effects.
\newblock In {\em Optics in Atmospheric Propagation and Adaptive Systems},
  volume 2580, pages 275--287. International Society for Optics and Photonics,
  1995.

\bibitem{lanteri2002}
H.~Lanteri, M.~Roche, and C.~Aime.
\newblock Penalized maximum likelihood image restoration with positivity
  constraints: multiplicative algorithms.
\newblock {\em Inverse problems}, 18(5):1397, 2002.

\bibitem{lanteri2001}
H.~Lant{\'e}ri, M.~Roche, O.~Cuevas, and C.~Aime.
\newblock A general method to devise maximum-likelihood signal restoration
  multiplicative algorithms with non-negativity constraints.
\newblock {\em Signal Processing}, 81(5):945--974, 2001.

\bibitem{lanteri2013}
H.~Lant{\'e}ri and C.~Theys.
\newblock La {G}amma divergence de {F}ujisawa et {E}guchi. {U}ne alternative
  pour imposer la contrainte de somme. {A}pplication au d{\'e}m{\'e}lange
  lin{\'e}aire.
\newblock In {\em Actes du 24e Colloque GRETSI sur le Traitement du Signal et
  des Images}, 2013.

\bibitem{lanteri2015x}
H.~Lant{\'e}ri, C.~Theys, and C.~Aime.
\newblock Divergences invariantes par changement d'{\'e}chelle pour la
  reconstruction de signaux ou d'images . {A}pplication \'a
  l'interf{\'e}rom\'etrie optique.
\newblock In {\em Actes du 25e Colloque GRETSI sur le Traitement du Signal et
  des Images}, 2015.

\bibitem{lanteri2015}
H.~Lant{\'e}ri, C.~Theys, and C.~Aime.
\newblock Scale invariant divergences for signal and image reconstruction.
\newblock In {\em Signal Processing Conference (EUSIPCO), 2015 23rd European},
  pages 2162--2166. IEEE, 2015.

\bibitem{lanteri2011}
H.~Lant{\'e}ri, C.~Theys, and C.~Richard.
\newblock Minimisation d'une forme g{\'e}n{\'e}rale de divergence sous
  contrainte de non-n{\'e}gativit{\'e} application \`a la factorisation en
  matrices non-n{\'e}gatives.
\newblock In {\em Actes du 23e Colloque GRETSI sur le Traitement du Signal et
  des Images}, 2011.

\bibitem{lee2001algorithms}
D.~D. Lee and H.~S. Seung.
\newblock Algorithms for non-negative matrix factorization.
\newblock In {\em Advances in neural information processing systems}, pages
  556--562, 2001.

\bibitem{lin1991}
J.~Lin.
\newblock Divergence measures based on the {S}hannon entropy.
\newblock {\em Information Theory, IEEE Transactions on}, 37(1):145--151, 1991.

\bibitem{lucy1974}
L.B. Lucy.
\newblock An iterative technique for the rectification of observed
  distributions.
\newblock {\em The astronomical journal}, 79:745, 1974.

\bibitem{michel1994}
O.J.J. Michel, R.G. Baraniuk, and P.~Flandrin.
\newblock Time-frequency based distance and divergence measures.
\newblock In {\em Time-Frequency and Time-Scale Analysis, 1994., Proceedings of
  the IEEE-SP International Symposium on}, pages 64--67. IEEE, 1994.

\bibitem{mihoko2002}
M.~Mihoko and S.~Eguchi.
\newblock Robust blind source separation by $\beta$ divergence.
\newblock {\em Neural computation}, 14(8):1859--1886, 2002.

\bibitem{neemuchwala2007}
H.~Neemuchwala, A.~Hero, S.~Zabuawala, and P.~Carson.
\newblock Image registration methods in high-dimensional space.
\newblock {\em International Journal of Imaging Systems and Technology},
  16(5):130--145, 2007.

\bibitem{neyman1949}
J.~Neyman.
\newblock Contribution to the theory of the x2 test.
\newblock In {\em Proceedings of the Berkeley symposium on mathematical
  statistics and probability}, volume~1, pages 239--273. University of
  California Press, 1949.

\bibitem{osterreicher1996}
F.~{\"O}sterreicher.
\newblock On a class of perimeter-type distances of probability distributions.
\newblock {\em Kybernetika}, 32(4):389--393, 1996.

\bibitem{osterreicher2002}
F.~{\"O}sterreicher.
\newblock Csisz\"ar\'{}s f-divergences-basic properties.
\newblock {\em RGMIA Res. Rep. Coll}, 2002.

\bibitem{osterreicher2013}
F.~{\"O}sterreicher.
\newblock Distances based on the perimeter of the risk set of a testing
  problem.
\newblock {\em Austrian Journal of Statistics}, 42(1):3--19, 2013.

\bibitem{pardo1999}
J.A. Pardo and M.C. Pardo.
\newblock Small-sample comparisons for the {R}ukhin goodness-of-fit-statistics.
\newblock {\em Statistical Papers}, 40(2):159--174, 1999.

\bibitem{pardo2005}
L.~Pardo.
\newblock {\em Statistical inference based on divergence measures}.
\newblock CRC Press, 2005.

\bibitem{pauca2006nonnegative}
V.~P. Pauca, J.~Piper, and R.~J. Plemmons.
\newblock Nonnegative matrix factorization for spectral data analysis.
\newblock {\em Linear algebra and its applications}, 416(1):29--47, 2006.

\bibitem{Polyanin:2001ul}
A.D. Polyanin, V.F. Zaitsev, and A.~Moussiaux.
\newblock {\em {Handbook of First-Order Partial Differential Equations}}.
\newblock CRC Press, November 2001.

\bibitem{renyi1955}
A.~R{\'e}nyi.
\newblock On a new axiomatic theory of probability.
\newblock {\em Acta Mathematica Hungarica}, 6(3-4):285--335, 1955.

\bibitem{renyi1961}
A.~R\'enyi.
\newblock Fourth berkeley symp.
\newblock {\em Math. Stat. Probability}, 1:547, 1961.

\bibitem{richardson1972}
W.H. Richardson.
\newblock Bayesian-based iterative method of image restoration.
\newblock {\em JOSA}, 62(1):55--59, 1972.

\bibitem{rockafellar2015}
R.T. Rockafellar.
\newblock {\em Convex analysis}.
\newblock Princeton university press, 2015.

\bibitem{roddier1971distributions}
F.~Roddier.
\newblock {\em Distributions et transformations de Fourier {\`a} l'usage des
  physiciens et des ing{\'e}nieurs}.
\newblock Ediscience, 1971.

\bibitem{rukhin1994}
A.L. Rukhin.
\newblock Optimal estimator for the mixture parameter by the method of moments
  and information affinity.
\newblock In {\em Transactions of the 12th Prague Conference on Information
  Theory}, pages 214--219, 1994.

\bibitem{salicru1994}
M.~Salicr{\'u}.
\newblock Measures of information associated with {C}sisz\"ar's divergences.
\newblock {\em Kybernetika}, 30(5):563--573, 1994.

\bibitem{shannon1948}
C.~E. Shannon.
\newblock A mathematical theory of communication.
\newblock {\em Bell system technical journal}, 27(3):379--423, 1948.

\bibitem{sharma1975}
B.D. Sharma and D.P. Mittal.
\newblock New nonadditive measures of entropy for discrete probability
  distributions.
\newblock {\em J. Math. Sci}, 10:28--40, 1975.

\bibitem{shepp1982}
L.A. Shepp and Y.~Vardi.
\newblock Maximum likelihood reconstruction for emission tomography.
\newblock {\em Medical Imaging, IEEE Transactions on}, 1(2):113--122, 1982.

\bibitem{taneja2001}
I.J. Taneja.
\newblock Generalized information measures and their applications. on-line
  book, 2001.
\newblock {\em URL www. mtm. ufsc. br/taneja/book/book. html}.

\bibitem{taneja1989}
I.J. Taneja.
\newblock On generalized information measures and their applications.
\newblock {\em Advances in Electronics and Electron Physics}, 76:327--413,
  1989.

\bibitem{tikhonov1974methods}
A.N. Tikhonov and V.Y. Arsenin.
\newblock {\em Methods of solving incorrect problems}.
\newblock Nauka, Moscow, 1974.

\bibitem{titterington1985}
D.M. Titterington.
\newblock General structure of regularization procedures in image
  reconstruction.
\newblock {\em Astronomy and Astrophysics}, 144:381, 1985.

\bibitem{tsallis1988}
C.~Tsallis.
\newblock Possible generalization of {B}oltzmann-{G}ibbs statistics.
\newblock {\em Journal of statistical physics}, 52(1-2):479--487, 1988.

\bibitem{uchida2005}
M.~Uchida and H.~Shioya.
\newblock An extended formula for divergence measures using invariance.
\newblock {\em Electronics and Communications in Japan (Part III: Fundamental
  Electronic Science)}, 88(4):35--42, 2005.

\bibitem{zaccheo1996}
T.S. Zaccheo and R.A. Gonsalves.
\newblock Iterative maximum-likelihood estimators for positively constrained
  objects.
\newblock {\em JOSA A}, 13(2):236--242, 1996.

\bibitem{zhang2004}
J.~Zhang.
\newblock Divergence function, duality, and convex analysis.
\newblock {\em Neural Computation}, 16(1):159--195, 2004.

\end{thebibliography}

\end{document}